\documentclass[11pt,a4paper]{article}

\usepackage{fullpage,wrapfig,relsize}
\input{manuscript/preamble}

\title{Nonparametric Maximum Likelihood Estimation for Mixtures of Logistic Regressions}
\author{
  Muheng Li\textsuperscript{1}
  \quad Peize Zhang\textsuperscript{2}
  \quad Xin Bing\textsuperscript{1,3}\\[0.5em]
  \begin{tabular}{@{}l@{}}
    \textsuperscript{1}\,Department of Statistical Sciences,
    University of Toronto\\
    \textsuperscript{2}\,Department of Computer Science,
    University of Toronto\\
    \textsuperscript{3}\,Department of Statistics,
    University of Illinois Urbana-Champaign
  \end{tabular}
}
\date{}

\begin{document}
\maketitle

\begin{abstract}
    Mixtures of logistic regressions provide a flexible framework for modeling binary responses, extending standard logistic regression to accommodate unobserved heterogeneity while retaining an interpretable regression structure. In this paper, we study such models with an unspecified mixing distribution.  We analyze a sieve nonparametric maximum likelihood estimator (NPMLE) and derive nonasymptotic bounds on the estimation error of the resulting mixture regression function under both fixed and random designs. Our results reveal an explicit tradeoff between statistical estimation error and the approximation error induced by the sieve radius, from which we derive convergence rates under several tail conditions on the mixing distribution.

   For computation, we propose a finite-particle algorithm that approximates
    the sieve NPMLE by optimizing the locations and weights of a finitely
    supported mixing distribution. The method offers numerical advantages
    over existing algorithms and comes with rigorous convergence guarantees.
    For any fixed dataset, we establish global convergence of the likelihood
    to the restricted NPMLE value when the number of particles exceeds a
    finite threshold. When the true mixing distribution is finitely supported,
    we further establish uniform global convergence guarantees and derive
    error bounds that separate algorithmic and statistical errors.
    Our analysis also yields a denoising bound for the mixing distribution
    and an identifiability result that complements existing results
    in the literature.

\end{abstract}

\newpage
\tableofcontents
\section{Introduction}

Logistic regression is a ubiquitous approach to binary classification, with applications across a broad range of scientific and decision-making problems. Let $(x_1,Y_1),\ldots,(x_n,Y_n)$ be $n$ observations, where $x_i\in\cX\subseteq\RR^d$ is a covariate vector and $Y_i\in\{0,1\}$ is a binary response for each $i\in[n]:=\{1,\ldots,n\}$. Logistic regression models
\begin{equation}\label{model_logis}
    \PP(Y_i=1\mid x_i)=\sigma(x_i^\T\beta),
\end{equation}
independently for $i\in[n]$, where $\sigma(t)=1/(1+e^{-t})$ is the sigmoid function mapping $t\in\RR$ to $(0,1)$ and $\beta\in \RR^d$ is an unknown parameter vector. The primary goals are to estimate the regression function $x\mapsto\sigma(x^\T\beta)$ and to classify future observations. Logistic regression is widely used in applications such as medical diagnosis, spam filtering, credit risk assessment, and click-through rate prediction. More recently, it has also played an important role in reward modeling for reinforcement learning from human feedback (RLHF). In a typical preference-learning setting, for each prompt $s_i$, a human rater is asked to indicate a preference between two candidate responses $a_i^+$ and $a_i^-$, yielding the binary outcome $Y_i= 1\{a_i^+ \text{ is preferred to } a_i^-\}$. The Bradley–Terry model provides a standard formulation in which \eqref{model_logis} holds with $x_i=\phi(s_i,a_i^+)-\phi(s_i,a_i^-)$, where $\phi$ is a pretrained feature mapping. The fitted reward model can subsequently be used to fine-tune a large language model toward better alignment with human preferences \citep{wang2023aligning}.

In modern data applications, however, a single logistic regression model as in \eqref{model_logis} can be too restrictive to capture heterogeneity across observations or individuals. When the relationship between the covariates and the binary response varies across individuals or subpopulations, imposing a common coefficient vector $\beta$ may lead to model misspecification and fail to capture such heterogeneity \citep{ichimura1998maximum}.  This issue is particularly salient in the RLHF application discussed above. Empirical studies of human preference datasets have found substantial and systematic disagreement among annotators, reflecting heterogeneity in individual preferences rather than merely random variation in their responses \citep{zhang2025diverging}. Recent theoretical work further shows that, in the presence of such heterogeneity, fitting a single reward model can fail to accurately represent population-level preferences, with implications for the subsequent policy-optimization stage of RLHF \citep{aouad2026sign}.

A natural extension of \eqref{model_logis} is therefore to allow the regression coefficient to vary across observations. Specifically, suppose that each observation $(x_i,Y_i)$ is associated with an unobserved coefficient $\beta_i$, where $\beta_1,\ldots,\beta_n$ are independently drawn from an unknown distribution $G\in\cP(\RR^d)$. Marginalizing over $\beta_i$ gives the regression function
\begin{equation}\label{def_p_G}
    p_G(x)
    :=\PP_G(Y=1\mid x)
    =\int \sigma(x^\T\beta)\,\rd G(\beta),
    \qquad x\in\cX.
\end{equation}
The mixing distribution $G$ provides a flexible representation of latent heterogeneity in the population, with different values of $\beta$ corresponding to different response patterns or preferences. The standard logistic regression model \eqref{model_logis} is  the special case in which $G=\delta_\beta$ is a Dirac measure.

Models of the form \eqref{def_p_G} have a long history, with different approaches placing different structural restrictions on the mixing distribution $G$. A common parametric approach assumes that the random coefficients, or a subset such as random intercepts and slopes, follow a Gaussian distribution, leading to logistic random-effects models  \citep{molenberghs2005models}. Although computationally convenient, the Gaussian assumption can be restrictive when the underlying population is asymmetric, heavy-tailed, or multimodal; misspecification of the random-effects distribution can in turn lead to biased estimation and inaccurate statistical inference \citep{marquart2019misspecification}. Another widely used approach assumes that $G$ is a finite discrete distribution,
$
G=\sum_{k=1}^K w_k\delta_{\beta_k},
$
where $w_k\ge 0$ with $\sum_{k=1}^K w_k = 1$ are the mixing weights and $\beta_1,\ldots,\beta_K\in\RR^d$ are the support points, leading to a finite mixture of logistic regressions \citep{wedel1995mixture}. Such models provide greater flexibility than a single logistic regression and are commonly fitted using likelihood-based procedures such as the Expectation-Maximization algorithm \citep{grun2007fitting}. Nevertheless, they require specification or selection of the number of components $K$. Estimating $K$ entails a nontrivial model-selection problem, while an incorrectly specified $K$ imposes an additional source of model misspecification. More fundamentally, restricting $G$ to a discrete distribution with a prespecified, and often relatively small, number of support points may be undesirable when little prior information is available about the underlying structure and complexity of population heterogeneity.

An alternative approach is to leave the mixing distribution $G$ in \eqref{def_p_G} unspecified and estimate the resulting regression function $p_G$ by maximizing the (nonparametric) likelihood. Specifically,  for any $G\in\cP(\RR^d)$, define the conditional log-likelihood
\begin{equation}\label{def_llh}
    \ell_n(G)
    :=\frac{1}{n}\sum_{i=1}^n
    \left\{
        Y_i\log p_G(x_i)
        +(1-Y_i)\log\big(1-p_G(x_i)\big)
    \right\}.
\end{equation}
The nonparametric maximum likelihood estimator (NPMLE) is obtained by maximizing $\ell_n(G)$ over $G\in\cP(\RR^d)$.
Nonparametric maximum likelihood estimation of mixing distributions has a long history dating back to the seminal work of \citet{kiefer1956consistency}, with subsequent developments establishing fundamental structural and computational properties of mixture likelihoods \citep[e.g.,][]{lindsay1983geometry}. Under the logistic mixture models, the existing literature of NPMLE has largely focused on more structured settings. For example, \citet{follmann1989generalizing} considered nonparametric mixing over a random intercept, while \citet{lesperance2014nonparametric} developed a computational method for repeated binary observations sharing a random intercept and a random slope within each subject.

More recently, substantial progress has been made in understanding both the statistical and computational properties of the NPMLE for Gaussian mixture models (GMMs). In particular, finite-sample convergence rates for the mixture density estimated by the NPMLE have been established under a range of mixing distributions \citep[e.g.,][]{zhang2009generalized,saha2020nonparametric,soloff2025multivariate,polyanskiy2020self}, while computational procedures for approximating the NPMLE have also been developed \citep[e.g.,][]{koenker2014convex,feng2018approximate,zhang2024efficient}. Related to our regression setting, \citet{jiang2025nonparametric} study the NPMLE for mixtures of linear regressions, establishing finite-sample guarantees for conditional density estimation under fixed and random designs and proposing a computational procedure for approximating the NPMLE.

By contrast, both finite-sample statistical theory and provable computational guarantees remain much less developed for logistic mixture regression, particularly when the entire coefficient vector $\beta\in\RR^d$ is allowed to vary according to an unrestricted mixing distribution. The goal of this paper is to develop finite-sample statistical guarantees for a sieve NPMLE under the logistic mixture model \eqref{def_p_G}, for both fixed and random designs, and to provide a computationally tractable procedure with provable guarantees for computing the sieve NPMLE.

\subsection{Contributions}\label{sec:contrib}

We summarize our main contributions in this section.

\paragraph{Non-asymptotic error bounds for the sieve NPMLE}  Our first contribution is to establish finite-sample convergence rates in \cref{sec_theory} for estimating the mixture regression function $p^*:=p_{G^*}$, with $G^*\in \cP(\RR^d)$, via a sieve NPMLE under both fixed and random designs.

The unrestricted NPMLE over $G\in\cP(\RR^d)$ need not exist,
as approaching the likelihood supremum may require support points
to diverge to infinity. For $r>0$, we therefore restrict candidate
mixing distributions to $\cB_r:=\{\beta\in\RR^d:\|\beta\|_2\le r\}$ and define the sieve NPMLE
\begin{equation}\label{def_hat_p_n_r}
    \wh G_r  ~ \in ~  \argmax_{G\in\cP(\cB_r)}\, \ell_n(G),
    \qquad \wh p_r ~ := ~ p_{\wh G_r}.
\end{equation}
The radius $r$ controls an estimation--approximation tradeoff:
increasing $r$ reduces approximation error but increases the
complexity of the sieve class and hence the estimation error bound.
Our nonasymptotic analysis allows $r=r_n$ to grow with $n$,
so that the sieve can approximate mixing distributions with
unbounded support.

For the statistical analysis, define the unrestricted mixture regression class $\cF:=\{p_G:G\in\cP(\RR^d)\}$ and its restricted counterpart
\begin{equation}\label{def_F_r}
    \cF_r := \{p_G:G\in\cP(\cB_r)\}.
\end{equation}
Since $G^*$ remains
unrestricted, $p^*$ need not belong to $\cF_r$.  For fixed design, \cref{thm:mse_rate_fixed} establishes the following high-probability bound on the in-sample mean squared error:
\begin{equation}\label{bd_sieve_intro}
    P_n(\wh p_r-p^*)^2   ~ \lesssim ~  \epsilon_{\est}(r,R_{x,n}) + \epsilon_{\app}(r,P_n),
\end{equation}
where $P_n$ denotes the empirical average over the design points and  $R_{x,n}:=\max_{i\in[n]}\|x_i\|_2$. The term $\epsilon_{\est}(r,R_{x,n})$ captures the estimation error within the sieve $\cF_r$ and increases with the sieve radius $r$, the dimension $d$, and the covariate radius $R_{x,n}$. When these quantities are fixed, it scales as $1/n$ up to logarithmic factors in $n$. In contrast, $\epsilon_{\app}(r,P_n)$ captures the approximation error incurred by restricting the unrestricted class $\cF$ to $\cF_r$ and vanishes as $r\to\infty$. We control this approximation error through the minimal empirical $\chi^2$-divergence between the Bernoulli distributions induced by $p^*$ and functions $p\in\cF_r$, and further relate it to weak moments of $G^*$ in \cref{lem_approx_error}. Balancing $\epsilon_{\est}(r,R_{x,n})$ and $\epsilon_{\app}(r,P_n)$ yields explicit rates under three standard conditions on $G^*$—bounded support, exponential tail decay, and a finite weak $p$-th moment—as summarized in \cref{cor:fixed_tail_conditions}. With suitable choices of $r=r_n$, the bounded-support and exponential-tail rates have the same polynomial dependence on $n$, up to logarithmic factors, as corresponding NPMLE rates for Gaussian mixture models \citep{zhang2009generalized,saha2020nonparametric}; see \cref{rem_rates_GMM} for a detailed comparison. Moreover, the results of \cref{thm:mse_rate_fixed} accommodate any {\em approximate maximizers} of the sieve likelihood in \eqref{def_hat_p_n_r}, which is particularly important  for our computational results.

A key ingredient in both analyses is a metric entropy bound for $\cF_r$ under the sup-norm that depends only polylogarithmically on the inverse covering radius; see \cref{thm:metric_entropy} in \cref{sec_entropy_logistic_mixtures}. We establish this bound by using Chebyshev polynomials to
approximate the sigmoid function on bounded intervals.
Integrating these approximations against the mixing distribution
places the resulting functions in a common finite-dimensional
polynomial space, whose covering number can be controlled.
\cref{sec_theory_proof} explains the argument and its relation to existing
entropy bounds for Gaussian location mixtures.

\paragraph{Finite-particle algorithm for approximating the sieve NPMLE}
Our second contribution is to develop the reweighted particle-gradient algorithm in \cref{sec:particle_prediction_formulation} for computing the sieve NPMLE. Computing $\wh G_r$ in \eqref{def_hat_p_n_r}  requires optimization over the infinite-dimensional space $\cP(\cB_r)$ of mixing distributions. A key structural property of mixture likelihoods is that an NPMLE can be chosen to have finite support \citep{lindsay1983geometry}; see \eqref{eq:potential_directional_derivative} for details.  Motivated by this fact, we represent the mixing distribution by an $M$-particle probability measure
$$
    G=\sum_{j=1}^M w_j\delta_{\beta_j},
    \qquad
    w\in\Delta^{M-1} := \Bigl\{v \in [0,1]^M:  \sum_{j=1}^M v_j=1\Bigr\},\quad
    \beta_j\in\cB_r,
$$
where $M$ is a potentially large number of particles that controls the flexibility of the finite-particle approximation. Although the log-likelihood is concave in the mixing distribution $G$ over the convex set $\cP(\cB_r)$, it is generally not concave in the particle locations $\beta_1,\ldots,\beta_M$, and the set of measures supported on at most $M$ points is not convex. Our algorithm alternates between updating the particle locations $\beta_j$ in an inner loop and optimizing the mixing weights $w_j$ in an outer loop. Specifically, at each outer iteration $t\ge 1$, we hold the weights $w^{(t)}$ fixed and perform $H\ge 1$ projected gradient-ascent steps on the particle locations, producing the inner-loop iterates $\beta_1^{(t,h)},\ldots,\beta_M^{(t,h)}$ for $h\in[H]$. We then set $\beta_j^{(t+1)}:=\beta_j^{(t,H)}$ and, holding these locations fixed, maximize the likelihood over the weights in $\Delta^{M-1}$ to obtain $w^{(t+1)}$. The resulting particle measure is
\begin{equation}\label{def_G_t}
    G_{t+1}=\sum_{j=1}^M w_j^{(t+1)}\delta_{\beta_j^{(t+1)}}.
\end{equation}
To initialize the algorithm, we randomly sample the particle locations from $\cB_r$ and then optimize the initial weights. After $t$ outer iterations, the algorithm has performed $Ht$ particle-location updates and $t+1$ weight optimizations, including the initial weight optimization.

Existing computational approaches to mixture NPMLEs can be broadly grouped according to how the support points of the mixing distribution are constructed and updated. The first approach sequentially adds new support points, guided by likelihood-gradient or related optimality criteria, and re-optimizes the mixing weights over the enlarged support set \citep{lesperance2007fast,jagabathula2020conditional}. The second approach prespecifies a finite set of candidate support points and optimizes the mixing weights over this set \citep{fox2011simple,train2008em,kim2020fast,zhang2024efficient}. The number of candidate points may depend on the sample size \citep{fox2016simple}.  Alternatively, the candidate support may be constructed from the data; for example, \citet{jiang2025nonparametric} generate candidate support points for mixtures of linear regressions by fitting coefficients to randomly selected subsets of the observations. The third approach allows the support points to move during optimization, possibly together with their weights, as in the Wasserstein–Fisher–Rao gradient-flow method of \citet{yan2024learning}, developed in the context of Gaussian mixture models, and the neural-network formulation of mixed logit studied by \citet{wang2024neural}.

Our method is most closely related to this third approach: it alternates between exact optimization of the mixing weights and projected gradient-ascent updates of all particle locations. Unlike the second approach, the particle locations are updated during optimization; unlike the first approach, all particle locations are updated jointly rather than added one at a time. As detailed below, an important distinction from the third approach is that our global convergence guarantees apply directly to the proposed algorithm with a sufficiently large but {\em finite} number of particles, whereas the convergence analyses in \citet{yan2024learning,wang2024neural} rely on measure-space or mean-field limiting dynamics as $M\to\infty$.

\paragraph{Convergence guarantees for the finite-particle algorithm}
Our third contribution is to establish global convergence guarantees for the proposed algorithm, ranging from deterministic fixed-data results to population and finite-sample guarantees that are uniform over the data distribution.
For any fixed dataset, \cref{thm:reweighted_particle_global} of \cref{sec:particle_fixed_data} establishes global algorithmic convergence for the proposed algorithm. In particular, once the number of particles $M$ exceeds a threshold $M_0$, the iterate $G_t$ in \eqref{def_G_t} of the proposed algorithm satisfies, for all $t\ge 1$,
$$
    \sup_{G\in\cP(\cB_r)}\ell_n(G)-\ell_n(G_t)
    ~ \le ~
    \frac{C}{t+c},
$$
for suitable $C,c>0$ that do not depend on $M$. Thus, although the likelihood is nonconcave in the finite-particle parameters, a sufficiently large $M$ guarantees convergence of the likelihood values to the global optimum of the sieve NPMLE problem
at rate $\cO(t^{-1})$. Importantly, this guarantee holds for any fixed dataset without assumptions on the true mixing distribution $G^*$. The price of this generality is that the required threshold $M_0$ and constants $C,c$ depend on the data, so the result alone
does not yield an explicit bound accounting for both optimization and statistical errors.

We obtain more quantitative guarantees in \cref{sec:particle_population} when the true mixing distribution $G^*$ is finitely supported and both $G^*$ and $P_X$ satisfy suitable identifiability and regularity conditions. At the population level, \cref{thm:particle_population} shows that, for a sufficiently large but finite number of particles $M$, the induced regression function $p_{G_t}$ satisfies, for $t\ge 1$,
\begin{equation}\label{rate_popu}
    P_X(p_{G_t}-p^*)^2 \le  C \left\{ (1+\eta Ht)^{-2}+e^{-ct} \right\},
\end{equation}
where $c,C>0$ are some constants and $\eta$ is the step size used in the inner-loop gradient-ascent updates of the particle locations. The two terms on the right hand side of \eqref{rate_popu} decay at rates $\cO((Ht)^{-2})$ and $\cO(e^{-ct})$, respectively, reflecting the distinct roles of the $H$ inner location updates and the outer weight-optimization iterations: increasing $H$ accelerates the decay of the first term through the product $Ht$, but does not affect the exponential term $e^{-ct}$. For fixed $H$, the overall rate in \eqref{rate_popu} is $\cO(t^{-2})$, and \cref{prop:particle_inverse_square_sharp} further shows that this rate in general cannot be improved for the proposed algorithm.

We extend the population-level guarantee to finite samples in \cref{thm:particle_uniform_empirical}. In particular, with high probability, the iterates satisfy, for $t\ge 1$,
\begin{equation}\label{rate_emp}
    P_X(p_{G_t}-p^*)^2
    \le
    C
    \left\{
       (1+\eta Ht)^{-2}+e^{-ct}
    \right\}
    +
    \epsilon_{\est}(r,R_X),
\end{equation}
where $R_X$ is the support radius of the covariate distribution $P_X$.
The first term represents the algorithmic error due to optimization, while the second term is the statistical error that we have seen in \eqref{bd_sieve_intro}.  In particular, once $Ht\gtrsim 1/\sqrt{\epsilon_{\est}(r,R_X)}$ and $t\gtrsim\log\{1/\epsilon_{\est}(r,R_X)\}$, the algorithmic error is dominated by the statistical error and the overall upper bound in \eqref{rate_emp} is of order $\epsilon_{\est}(r,R_X)$.

The analysis of the proposed algorithm requires several ingredients to overcome the nonconvexity arising from the joint optimization of the particle weights and locations. First, we establish a quantitative denoising inequality that lower bounds the estimation error $P_X(p_G-p_{G^*})^2$ of the regression function by a suitable discrepancy between $G$ and $G^*$. Second, we use this inequality to construct ascent directions for the particle locations and to show that, once the likelihood gap of an iterate relative to the maximum is below a suitable threshold, the projected location updates align with these directions and achieve the $\cO(t^{-2})$ rate described above. Third, exact optimization of the weights at each outer iteration ensures that the likelihood is nondecreasing, allowing us to control the effect of the location updates separately from that of the weight updates. Finally, we show that, with a sufficiently large number of randomly initialized particles, the algorithm reaches the required likelihood region with high probability. Together, these ingredients yield the global convergence guarantees described above. See \cref{sec:particle_proof_techniques} for a detailed discussion of the proof strategy.

As a byproduct, the denoising inequality also implies identifiability of finite mixtures of logistic regressions under the random-design setting. Our result strengthens and complements existing identifiability results for finite logistic mixtures \citep{auder2021least} and uniqueness results for finite sums of sigmoid functions \citep{sussmann1992uniqueness}; see \cref{sec_ident} for a detailed comparison.\\

Extensive simulation studies in \cref{sec_sim} evaluate the finite-sample performance of the proposed algorithm and corroborate our theoretical findings. We examine in \cref{sec_sim_reweight} how the frequency of weight optimization in the outer loop as well as the number of inner-loop steps affects estimation accuracy and computational cost of our proposed algorithm. In \cref{sec_sim_cv} we investigate the effects of the sieve radius $r$ and the number of particles $M$ on the estimation error and assess their selection via cross-validation. The proposed algorithm is further compared with existing NPMLE-based computational methods in \cref{sec_sim_npmle},  an EM-based method and several machine learning approaches in \cref{sec_sim_benchmarks}.  Overall, the numerical results demonstrate superior performance of the proposed algorithm across a range of settings. All the proofs  are deferred to the Appendix.

\section{NPMLE-based regression function estimation}\label{sec_theory}

In this section, we analyze the estimation error of the $r$-sieve NPMLE $\wh p_{r}$ in \eqref{def_hat_p_n_r} under both fixed and random designs.
A central step of both analyses is to establish suitable entropy bounds for the restricted logistic-mixture class $\cF_r$ given in \eqref{def_F_r}.
This allows us to quantify the complexity of the model class and to control the associated empirical process for analyzing the estimation errors. We first state the entropy bound in \cref{sec_entropy_logistic_mixtures}.

\subsection{Entropy bounds for logistic mixtures}
\label{sec_entropy_logistic_mixtures}

For any $\epsilon>0$,  let $\cN(\epsilon, \cF, \|\cdot\|)$ denote the $\epsilon$-covering number of a function class $\cF$ under any pseudometric $\|\cdot\|$, that is, the smallest integer $N$ such that there exist $f_1,\ldots, f_N\in \cF$ satisfying $\sup_{f\in \cF}\min_{j\in [N]} \|f-f_j\| \le \epsilon$. Let $\|\cdot\|_{\i,\cX}$ be a pseudonorm on $\cF$ such that for any $f\in \cF$,
$$
    \|f\|_{\i,\cX} := \sup_{x\in \cX}|f(x)|.
$$
The entropy studied in the following is $\log \cN(\epsilon, \cF_r, \|\cdot\|_{\i,\cX})$ with $\cF_r$ given in \eqref{def_F_r}.

A generic convex-hull argument gives an empirical $\ell_2$-norm entropy
bound with polynomial dependence on $1/\epsilon$
\citep[Lemma~2.6.11]{vanderVaartWellner1996}.
Such bounds exploit only the boundedness and Lipschitz continuity
of the sigmoid function $\sigma(\cdot)$ and are too coarse for the near-parametric
risk bounds sought here.

For Gaussian location mixtures, sharper entropy bounds are
obtained through Taylor expansion and moment matching to approximate mixture densities by finite mixtures \citep[Lemma~1]{zhang2009generalized}; see also \citet[Lemmas~D.2--D.3]{saha2020nonparametric}.  For logistic mixtures, we instead approximate the sigmoid
function $\sigma(\cdot)$ on bounded intervals using Chebyshev polynomials.
Integrating these approximations against the mixing distribution $G\in\cP(\cB_r)$
yields a common finite-dimensional polynomial space that
approximates $\mathcal F_r$ uniformly, leading to the entropy
bound below.

Let $R_x := \sup_{x\in\cX}\|x\|_2$. The following theorem gives our entropy bound for
$\log\mathcal N(\epsilon,\mathcal F_r,\|\cdot\|_{\infty,\cX})$.
Its proof appears in \cref{app:proof_metric_entropy}, while
\cref{sec_theory_proof} provides a sketch and explains how
the argument differs from existing analyses of Gaussian
location mixtures.

\begin{theorem}
\label{thm:metric_entropy}
For any $r>0$, set
$
    \bar R_r:= 1\vee (r R_x).
$
Then there exists a universal constant $C>0$ such that for every
$\epsilon\in(0,1]$,
$$
\log \mathcal N(\epsilon,\mathcal F_r,\|\cdot\|_{\i,\cX})
~\le~
\bigg[
C\bar R_r\log\!\Big(\frac{C\bar R_r}{\epsilon}\Big)
\bigg]^d
\log\!\Big(\frac{C}{\epsilon}\Big).
$$
\end{theorem}

The polylogarithmic dependence on $1/\epsilon$ is the key ingredient for deriving near-$n^{-1}$ estimation-error bounds for the NPMLE in the next two subsections. An extension to mixing distributions with unbounded support and finite exponential moments is given in \cref{app:entropy_exponential_moments}.

General sup-norm entropy bounds for functions that are analytic and uniformly
bounded on a complex neighborhood of the unit cube are also available in
\citet[Theorem~2.7.16]{vanderVaartWellner2023}. Applying this result to
$\cF_r$ requires rescaling the domain and bounding the analytic extensions
on a common complex neighborhood. The argument above instead approximates
the univariate function $\sigma(\cdot)$ before integrating against $G$, giving
the explicit dependence on $r$, $R_x$, and $d$ in
\cref{thm:metric_entropy}.

\subsection{Estimation error under fixed design}
\label{sec_theory_fix}

Throughout this subsection, we treat $x_1,\ldots,x_n$ as fixed design points, so the labels $Y_i$ are independent Bernoulli random variables with success probabilities $p^*(x_i)=p_{G^*}(x_i)$ for $i\in[n]$.
For any $\wh p\in \cF$,   its in-sample mean squared estimation error (MSE) is
\begin{equation}\label{def_pred_emp}
   P_n (\wh p - p^*)^2   := {1\over n}\sum_{i=1}^n\bigl[\wh p(x_i) - p^*(x_i)\bigr]^2,
\end{equation}
where, for any function $f:\cX\to\RR$, we write $P_n f := \frac{1}{n}\sum_{i=1}^n f(x_i).$

For any $r>0$, the estimation error of the sieve NPMLE depends on both the complexity of $\cF_r$ and how well it approximates $p^*$. The entropy bound in \cref{thm:metric_entropy}, applied to $\cX=\{x_1,\ldots,x_n\}$, gives the estimation error term
\begin{equation}\label{def_epsilon_1}
    \epsilon_{\est}(r, R_{x,n}) :=  \left(1 + rR_{x,n}\right)^d\left\{\log[n(1 + rR_{x,n})]\right\}^{d+1} \frac{1}{n},
\end{equation}
where $R_{x,n}:=\max_{i\in [n]}\|x_i\|_2$.
To quantify the error from restricting the mixing distribution to $\cB_r$, define the approximation error by the smallest empirical $\chi^2$-divergence between the Bernoulli distributions induced by $p^*$ and functions in $\cF_r$:
\begin{equation}\label{def_chisquare_emp}
\epsilon_{\app}(r, P_n)
:=
\inf_{q\in \cF_r}
\chi^2_n(p^*,q)
:=
\inf_{q\in \cF_r}
P_n \frac{(p^*-q)^2}{q(1-q)}.
\end{equation}
A larger $r$ increases $\epsilon_{\est}(r,R_{x,n})$ but decreases $\epsilon_{\app}(r,P_n)$. We explain in \cref{sec_theory_proof} how the entropy bound and the likelihood comparison give these two error terms.

The following theorem gives an upper bound on the MSE of any approximate solution to the sieve NPMLE in \eqref{def_hat_p_n_r} in terms of both estimation and approximation errors defined above.  The full proof is given in \cref{app:proof_mse_rate_fixed}, and we offer a sketch in \cref{sec_theory_proof}. %

\begin{theorem}
\label{thm:mse_rate_fixed}
Assume $p^*\in\mathcal F$. Fix any $r>0$ and any $\tau \ge 0$. Let
$\wt p_r\in\mathcal F_r$ be any measurable estimator, and choose
$\wt G_r\in\cP(\cB_r)$ such that $p_{\wt G_r}=\wt p_r$. Suppose that
\begin{equation}
\label{eq:fixed_design_oracle_lr}
 \ell_n(\wt G_r)
~ \ge~  \sup_{G\in\cP(\cB_r)}\ell_n(G)
 -\tau.
\end{equation}
Then there exists a constant $C>0$, depending only on $d$, such that, for every $u\ge1$,
\begin{equation}
\label{eq:fixed_design_lr_bound}
\PP\left\{
P_n (\wt p_r-p^*)^2
>
C\left(
\epsilon_{\est}(r, R_{x,n}) +
\epsilon_{\app}(r, P_n) + \tau +
\frac{u}{n}
\right)
\right\}
\le e^{-u}.
\end{equation}
\end{theorem}

\cref{thm:mse_rate_fixed} separates the statistical error from the optimization error $\tau$ in \eqref{eq:fixed_design_oracle_lr}. In particular, an approximate maximizer has the same estimation error rate as the exact sieve NPMLE whenever $\tau$ is no larger than the sum of the estimation and approximation terms. For ease of presentation, we state the corollaries below for the exact sieve NPMLEs, for which $\tau=0$.

The estimation error $\epsilon_{\est}(r,R_{x,n})$ increases with the sieve radius $r$ and the covariate norm $R_{x,n}$. For fixed $r$ and bounded $R_{x,n}$, it is of nearly parametric order $n^{-1}$, up to logarithmic factors. In contrast, the approximation error $\epsilon_{\app}(r,P_n)$ decreases as $r$ increases. Its dependence on $G^*$ arises through $p^*=p_{G^*}$. In particular, if $G^*(\cB_r)=1$, then $\epsilon_{\app}(r,P_n)=0$. For $G^*$ with unbounded support, $\epsilon_{\app}(r,P_n)$ can be bounded in terms of the tail behavior of $G^*$.

Specifically, for $\beta\sim G^*$, define the weak moments
\citep[Section~2.3]{zhang2009generalized} of $\|\beta\|_2$ by
\begin{equation}\label{def_weak_moments}
    \mu^w_p(G^*)
    =
    \left\{
    \sup_{r>0}\ r^pG^*(\cB_r^c)
    \right\}^{1/p}
    =
    \left\{
    \sup_{r>0}\ r^p\int_{\|\beta\|_2>r}\,\rd G^*(\beta)
    \right\}^{1/p}
    \quad \forall~p>0.
\end{equation}
The term ``weak'' reflects that, by Markov's inequality, a finite $p$-th moment implies a finite weak $p$-th moment, while the converse need not hold. The following lemma, proved in \cref{app:proof_lem_approx_error}, bounds $\epsilon_{\app}(r,P_n)$ in terms of the weak moments of $G^*$. Let $\cB_r^c:=\RR^d\setminus\cB_r$.
\begin{lemma}\label{lem_approx_error}
For any $r>0$, with the convention $a/0:=+\infty$ for all $a\ge0$, one has
$$
    \epsilon_{\app}(r,P_n)
    ~ \le~
    \min\left\{
        G^*(\cB_r^c),\,
        \inf_{0<s\le r}
        \frac{(1+e^{s R_{x,n}})}{2G^*(\cB_s)} \{G^*(\cB_r^c)\}^2
    \right\},
$$
Moreover,
$
    G^*(\cB_r^c)
    ~ \le~  \inf_{p>0}\left({\mu^w_p(G^*)/ r}\right)^p.
$
\end{lemma}

\cref{lem_approx_error} bounds the approximation error in terms of the mass of $G^*$ outside $\cB_r$. Combining this result with \cref{thm:mse_rate_fixed} yields the following rates under three tail conditions on $G^*$. In each case, we choose $r$
to balance the estimation and approximation error bounds. The proof appears in \cref{app:proof_fixed_tail_conditions}. To simplify the presentation, we assume $\log R_{x,n}\lesssim\log n$ and $rR_{x,n}\gtrsim1$ for the choices of $r$ in this corollary.

\begin{corollary}
\label{cor:fixed_tail_conditions}
Assume $p^*\in\cF$. Then for sufficiently large $n$, the following holds for the sieve NPMLE $\wh p_{r}$ in \eqref{def_hat_p_n_r} with probability at least $1-n^{-1}$.
\begin{enumerate}[label=(\roman*),noitemsep]
\item \emph{(Bounded support).}
If $G^*(\cB_R)=1$ for some positive $R<\i$, then
$$
    P_n(\wh p_R-p^*)^2 ~ \lesssim~  (R R_{x,n})^d  \frac{ (\log n)^{d+1}}{n} .
$$

\item \emph{(Exponential decay).}
Suppose that
\begin{equation}\label{def_exp_G}
    G^* \in \cG(a,q,b) := \left\{
        G\in \cP(\RR^d): \int_{\RR^d}\exp(a\|\beta\|_2^q)\,\rd G(\beta)\le b
    \right\}
\end{equation}
for some $a>0$, $q>0$, and $1\le b <\i$. Then by setting
$
    r_n  = \{\log(n)/a\}^{1/q},
$
one has
$$
P_n(\wh p_{r_n}-p^*)^2
~ \lesssim~  \left({R_{x,n}  \over a^{1/q}}\right)^d \frac{(\log n)^{d/q + d+1} }{n} .
$$

\item \emph{(Finite weak $p$-th moment).}
If there exists some $p>0$ such that $\mu^w_p(G^*)=\mu \in (0,\i)$, by setting
$
    r_n = \mu^{2p/(2p+d)}n^{1/(2p+d)},
$
one has
$$
P_n(\wh p_{r_n}-p^*)^2
~ \lesssim ~
\left[
R_{x,n}^d
(\log n)^{d+1}
+1+e^{2^{1/p}\mu R_{x,n}}
\right]\frac{\mu^{2pd/(2p+d)}}{n^{2p/(2p+d)}}.
$$
\end{enumerate}
\end{corollary}

The exponential decay family in $(ii)$ allows sub-Weibull tails of any order $q>0$, including sub-Gaussian tails for $q=2$ and sub-exponential tails for $q=1$.

Treating all tail constants as fixed, and assuming $R_{x,n}=\cO(1)$,
 the following table summarizes the simplified rates in \cref{cor:fixed_tail_conditions} for all three tail conditions on $G^*$.
\begin{table}[H]
\centering
\caption{Simplified error rates under fixed design for different tail conditions on $G^*$}
\label{tab_rates_fixed}
\begin{tabular}{l|ccc}
\toprule
Tail condition on $G^*$
& Bounded support
& Exponential tail
& Finite weak $p$-th moment \\
\midrule
Order of $r_n$
& $1$
& $(\log n)^{1/q}$
& $n^{1/(2p+d)}$ \\[1ex]
Rate for $P_n(\wh p_{r_n}-p^*)^2$
& $\displaystyle \frac{(\log n)^{d+1}}{n}$
& $\displaystyle \frac{(\log n)^{d/q+d+1}}{n}$
& $\displaystyle \frac{(\log n)^{d+1}}{n^{2p/(2p+d)}}$ \\
\bottomrule
\end{tabular}
\end{table}

\begin{remark}[Comparison with Gaussian location mixtures]\label{rem_rates_GMM}
    It is informative to compare the results in \cref{tab_rates_fixed} with the corresponding results for Gaussian location mixtures, first established in the univariate case by \citet{zhang2009generalized} and subsequently extended to the multivariate setting by \citet{saha2020nonparametric}. Under bounded support or exponential tail decay of $G^*$, our bounds exhibit the same polynomial rates in $n$, with differences in the logarithmic factors. In particular, dimension-dependent powers of $\log n$ also appear in NPMLE rates for estimating Gaussian location mixture densities.  Under bounded support, both results yield the rate $(\log n)^{d+1}/n$. Under exponential-type tails of order $0<q\le2$, \citet{saha2020nonparametric} obtain $(\log n)^{d/q+d/2+1}/n$, compared with $(\log n)^{d/q+d+1}/n$ here. Under moment conditions, their squared Hellinger error bound has the polynomial rate $n^{-p/(p+d)}$, whereas our squared estimation error bound for the regression function has the rate $n^{-2p/(2p+d)}$, up to logarithmic factors.

    A more substantive difference concerns adaptation.  For Gaussian location mixtures, the unrestricted NPMLE, which optimizes over all mixing distributions, adapts to the tail behavior of $G^*$. In contrast, our guarantees for the sieve NPMLE over $\cF_r$ require choosing the radius $r$ according to that tail behavior. This distinction reflects an important structural difference between the two models. In addition to the global analyticity and quadratic-exponential decay of the Gaussian kernel discussed earlier, a Gaussian location mixture assumes
    $$
        Z=\beta+W,\qquad \beta\sim G^*,\qquad W\sim N(0,I_d),
    $$
    where $\beta$ and $W$ are independent, and only $Z\in\RR^d$ is observed. Thus, the observations are noisy versions of latent draws from $G^*$ and consequently carry direct information about the location and spread of its mass. In the logistic mixture model considered here, by contrast, the design points  $x_1,\ldots,x_n$ are fixed and therefore carry no information about $G^*$; information about $G^*$ is available only indirectly through the binary responses $Y_i$ at these design points. Consequently, the data provide no analogous direct information about where the mass of $G^*$ is concentrated or how widely it is spread, motivating the use of the sieve radius $r$ to control the range of candidate mixing distributions. Our simulation studies in \cref{sec_sim_cv} suggest that this choice is less of a concern in practice, as $r$ can be reliably selected via cross-validation.
\end{remark}

\subsection{Estimation error under random design}
\label{sec_theory_rand}

We now treat the design points $x_1,\ldots,x_n$ as realizations of random vectors  $X_1,\ldots,X_n$ that are i.i.d. from $P_X\in \cP(\cX)$. For any $\wh p\in\cF$, its $L_2(P_X)$ mean-squared error is defined as
\begin{equation}\label{def_pred_L2}
    P_X(\wh p  - p^*)^2 := \int \bigl[\wh p(x)-p^*(x)\bigr]^2 \,\rd P_X(x).
\end{equation}
We impose the following sub-Gaussian assumption on $P_X$.
\begin{assumption}
\label{assum:subgaussian_x}
For $X\sim P_X$, there exist constants $m_x\ge 0$ and $\sigma>0$ such that $\|\EE X\|_2\le m_x,$
and
$
\EE\exp(u^\T(X-\EE X))
\le
\exp({\sigma^2\|u\|_2^2/2})
$
for all $u\in\RR^d.$
\end{assumption}
Under \cref{assum:subgaussian_x}, a tail bound of $\|X-\EE X\|_2$ and a union bound over the observations give the following high-probability bound on $R_{x,n}$. For any $u>0$, by letting
\begin{equation}\label{def_Tn}
 T_u := m_x+ C_{\sigma} \sqrt{d+u+\log n}
\end{equation}
for some sufficiently large constant $C_{\sigma}>0$ depending only on $\sigma$, one has $\PP(R_{x,n}\ge T_u)\le e^{-u}$. The radius $T_u$ determines the compact region on which we apply the entropy bound in \cref{thm:metric_entropy}, leading to the estimation term $\epsilon_{\est}(r,T_u)$.

The population analogue of the approximation error in \eqref{def_chisquare_emp} is
\begin{equation}\label{def_chisquare_pop}
   \epsilon_{\app}(r, P_X)
   :=
   \inf_{q\in \cF_r}\chi^2_{P_X}(p^*,q)
   :=
   \inf_{q\in \cF_r}
    P_X{(p^*-q)^2\over q(1-q)}.
\end{equation}
The following theorem gives the $L_2(P_X)$ MSE bound for approximate sieve NPMLEs. Its full proof is given in \cref{app:proof_mse_rate_random}, and a sketch appears in \cref{sec_theory_proof}.

\begin{theorem}
\label{thm:mse_rate_random}
Assume $p^*\in\cF$ and \cref{assum:subgaussian_x} holds. Fix any $r>0$ and any
$\tau\ge0$. Let $\wt p_r\in\cF_r$ be any measurable estimator satisfying
\eqref{eq:fixed_design_oracle_lr}. There exists a constant $C>0$, depending only on $d$, such that, for every $u\ge1$,
\begin{equation}
\label{eq:random_design_l2_bound}
\PP\left\{
P_X(\wt p_r-p^*)^2
>
C\left(
\epsilon_{\est}(r,T_u)
+
\epsilon_{\app}(r, P_X)
+
\tau
+
{u\over n}
\right)
\right\}
\le e^{-u}.
\end{equation}
\end{theorem}

Compared with \cref{thm:mse_rate_fixed}, the bound uses $T_u$ in place of $R_{x,n}$ in the estimation error, while the approximation error is taken with respect to $P_X$. The analysis requires additional arguments beyond those used for the fixed design, as discussed in \cref{sec_theory_proof}.

Similar to the fixed-design case, the approximation error
$\epsilon_{\app}(r,P_X)$ can be bounded in terms of the weak moments of $G^*$ given in \eqref{def_weak_moments}.
In particular, the same proof as that of \cref{lem_approx_error} gives that, for any $r>0$,
\begin{equation}
\label{eq:pop_chisquare_tail_bound}
    \epsilon_{\app}(r, P_X)
    \le  \min\left\{
        G^*(\cB_r^c),\,
        \inf_{0<s\le r}
        \frac{1+\EE e^{s \|X\|_2}}{2G^*(\cB_s)} \{G^*(\cB_r^c)\}^2
    \right\}.
\end{equation}
The expectation is finite under \cref{assum:subgaussian_x}.

The following corollary states simplified estimation error bounds for the random-design setting under the same three tail conditions on $G^*$ as in \cref{cor:fixed_tail_conditions}. For simplicity, we only state the results for exact restricted NPMLEs in \eqref{def_hat_p_n_r}. The proof appears in \cref{app:proof_random_tail_conditions}.

\Needspace{8\baselineskip}
\begin{corollary}
\label{cor:random_tail_conditions}
Assume $p^*\in\cF$ and \cref{assum:subgaussian_x} holds. For sufficiently large $n$, the following results hold with probability at least $1-n^{-1}$.

\begin{enumerate}[label=(\roman*),noitemsep]
\item \emph{(Bounded support).}
If there exists some $0<R<\i$ such that $G^*(\cB_R)=1$, then
$$
P_X(\wh p_R-p^*)^2
~ \lesssim ~
R^d (d+\log n)^{d/2}
\frac{(\log n)^{d+1}}{n}.
$$

\item \emph{(Exponential decay).}
Suppose $G^*\in \cG(a,q,b)$ for some $a>0$, $q>0$, and $b\ge1$.
By setting
$
    r_n  =   \{ \log(n)/a \}^{1/q},
$
one has
$$
    P_X(\wh p_{r_n}-p^*)^2 ~ \lesssim ~
    \frac{
    (d+\log n)^{d/2}
    }{a^{d/q}}
    \frac{(\log n)^{d/q+d+1}}{n}.
$$

\item \emph{(Finite weak $p$-th moment).}
If there exists some $p>0$ such that $\mu^w_p(G^*)=\mu\in(0,\i)$, by setting
$
    r_{n} =     \mu^{2p/(2p+d)} n^{1/(2p+d)}(d+\log n)^{-d/(4p+2d)},
$
one has
$$
P_X(\wh p_{r_n}-p^*)^2  ~ \lesssim ~
\mu^{2pd\over 2p+d}
(d+\log n)^{pd\over 2p+d}
\frac{(\log n)^{d+1}}{n^{2p/(2p+d)} }.
$$
\end{enumerate}
The implicit constants may depend on $d,m_x,\sigma$ in all three parts,
and additionally on $a,q,b$ in part (ii) and on $p,\mu$ in part (iii),
but not on $n$.
\end{corollary}

Treating all tail constants and $d$ as fixed, the rates in \cref{cor:random_tail_conditions} simplify as follows.
\begin{table}[H]
\centering
\caption{Simplified estimation error rates under random design for different tail conditions on $G^*$}
\label{tab_rates_random}
\begin{tabular}{l|ccc}
\toprule
Tail condition on $G^*$
& Bounded support
& Exponential tail
& Finite weak $p$-th moment \\
\midrule
Order of $r_n$
& $1$
& $(\log n)^{1/q}$
& $\displaystyle n^{1/(2p+d)}(\log n)^{-d/(4p+2d)}$ \\[1ex]
Rate for $P_X(\wh p_{r_n}-p^*)^2$
& $\displaystyle \frac{(\log n)^{1+3d/2}}{n}$
& $\displaystyle \frac{(\log n)^{1 + d/q+3d/2}}{n}$
& $\displaystyle
\frac{(\log n)^{1+d+dp/(2p+d)}}{n^{2p/(2p+d)}}$ \\
\bottomrule
\end{tabular}
\end{table}
Compared with \cref{tab_rates_fixed}, the additional logarithmic factors in the rates arise from truncating the sub-Gaussian $P_X$. If $P_X(\cB_T)=1$ for some fixed $T$, then no truncation is needed and the argument in \cref{sec_theory_proof}, applied on $\cB_T$, yields the same rates as in \cref{tab_rates_fixed}.

\subsection{Discussion of the proof techniques for  \cref{thm:metric_entropy,thm:mse_rate_fixed,thm:mse_rate_random}}\label{sec_theory_proof}

In this section, we discuss the main proof techniques underlying the entropy bound as well as the estimation-error bounds under both fixed and random designs. %

\paragraph{Proof techniques for the entropy bound in \cref{thm:metric_entropy}}

In \cref{rem_approx_GMM} of \cref{app:proof_metric_entropy},
we review the approximation argument for Gaussian location
mixtures \citep[e.g.,][]{zhang2009generalized,saha2020nonparametric}.
This argument relies on two properties of the Gaussian density:
its two-sided exponential tail decay, which controls contributions
from distant mixing locations, and its globally convergent Taylor
expansion, which, together with moment matching, controls the
remaining contribution. However, neither property holds for the sigmoid function: it tends to one in its positive tail and, due to its complex poles, does not admit a globally convergent Taylor expansion. Nevertheless, it is analytic on every compact real interval and extends holomorphically to a suitable complex neighborhood. We exploit this local analytic structure through Chebyshev polynomial approximation, as described below.

\begin{proof}[Proof sketch of \cref{thm:metric_entropy}]
To approximate $x\mapsto \int \sigma(x^\T \beta)\,\rd G(\beta)$ for $G \in  \cP(\cB_r)$, we apply Chebyshev polynomial approximation to the sigmoid function $\sigma(t)$ uniformly over $t\in [-rR_x,rR_x]$ as $|x^\T \beta|\le rR_x$. Since the extended $\sigma(z)$ has complex poles at $z = (2k+1)i\pi$ for $k\in \ZZ$ as shown in \cref{lem:logistic_meromorphic}, the key step is to identify a so-called Bernstein ellipse in the complex domain and establish the analyticity of $\sigma(z)$ within this ellipse. This is the result of \cref{lem:logistic_chebyshev} and we give a detailed review of Chebyshev approximation in \cref{app:complex_chebyshev_prelim}. The semi-minor axis of such an ellipse determines both the rate of convergence of the resulting approximation, and the polynomial degree required to achieve a given accuracy. In \cref{cor:degree_choice_logistic} we show that for any $\epsilon\in (0,1/2]$, there is a polynomial
$P_k$ with degree $k\le C\bar R_r\log(C\bar R_r/\eps)$ such that
$$
    \sup_{|t|\le rR_x}|\sigma(t)-P_k(t)|\le \eps/4.
$$
Consequently,  we obtain a class of functions
$
\wt\cF_{r,k} := \{x\mapsto \int P_k(x^\T\beta) \,\rd G(\beta): G\in  \cP(\cB_r)\}
$
that approximates $\cF_r$ with error at most $\eps/4$. Since $P_k(x^\T\beta)$ is a
$d$-variate polynomial in $x$ of total degree at most $k$, the
approximating function class $\wt\cF_{r,k}$ lies in a polynomial space on $\cX$
of dimension at most $\binom{d+k}{d}$. An $\eps/2$-cover of this class,
with each polynomial center replaced by its corresponding logistic mixture,
gives an $\eps$-cover of $\cF_r$.
By \citet[Example~5.8]{wainwright2019high},
$$
\log \cN(\eps,\cF_r,\|\cdot\|_{\i,\cX})
\le
\binom{d+k}{d}\log(13/\eps) \le  (k+1)^d \log(13/\eps)
$$
which leads to the entropy bound in \cref{thm:metric_entropy}.
\end{proof}

\paragraph{Proof techniques for the error bounds in \cref{thm:mse_rate_fixed,thm:mse_rate_random}}

We first explain the proof techniques for the fixed-design setting and then discuss the additional arguments needed to extend the analysis to the random-design setting.

\begin{proof}[Proof sketch of \cref{thm:mse_rate_fixed}]
For $G,H\in\cP(\RR^d)$, define the squared Hellinger distance between the Bernoulli distributions with success probabilities $p_G(x)$ and $p_H(x)$ by
\begin{equation}\label{def_h_sq_main}
    h^2(p_G,p_H)(x)
    =1-\sqrt{p_G(x)p_H(x)}-\sqrt{(1-p_G(x))(1-p_H(x))},
    \qquad x\in\cX.
\end{equation}
Its empirical version is $h_n^2(p_G,p_H):=P_n h^2(p_G, p_H)$. For simplicity, we focus on the case $\tau=0$ in \cref{thm:mse_rate_fixed}, where $\wh p_r = p_{\wh G_r}$ is the exact sieve NPMLE defined in \eqref{def_hat_p_n_r}. Since \cref{lem:bernoulli_hellinger_x} in \cref{app:proof_mse_rate_fixed} gives $P_n(\wh p_r-p^*)^2\le4h_n^2(\wh p_r,p^*)$ with $p^*=p_{G^*}\in \cF$, it suffices to bound $h_n^2(\wh p_r,p^*)$.

Let $\gamma>0$ be  any deterministic sequence to be specified later. Start with
\begin{align}\nonumber
    & \PP\left\{
        h_n^2(\wh p_r,p^*) \ge \gamma
    \right\}   =   \inf_{H\in\cP(\cB_r)} \PP\left\{
        h_n^2(\wh p_r,p^*) \ge \gamma,  \ell_n(\wh G_r)
    \ge  \ell_n(H)
    \right\}\\
    & \qquad \le    \PP\left\{
        h_n^2(\wh p_r,p^*) \ge \gamma,  \ell_n(\wh G_r) -\ell_n(G^*)
    \ge  -\delta
    \right\} +  \inf_{H\in\cP(\cB_r)}\PP\left\{ \ell_n(G^*) - \ell_n(H)   \ge \delta\right\}\label{decomp_h_proof}
\end{align}
for some $\delta>0$ to be chosen. The second term in \eqref{decomp_h_proof} arises as  $G^*$ need not belong to $\cP(\cB_r)$.  By the independence of $Y_1,\ldots, Y_n$ and the definition of $\ell_n$ in \eqref{def_llh}, one has
\begin{align*}
    \EE\exp\{n[\ell_n(G^*)-\ell_n(H)]\}
    &=\prod_{i=1}^n\left\{1+
      {\{p^*(x_i)-p_H(x_i)\}^2\over p_H(x_i)\{1-p_H(x_i)\}}\right\}\\
    &\le\exp\left\{nP_n{(p^*-p_H)^2\over p_H(1-p_H)}\right\}.
\end{align*}
This explains why the approximation error is measured by $\chi^2$-divergence rather than the squared error alone.
Pick any $u>0$. By Markov's inequality, taking the infimum over $H\in\cP(\cB_r)$ and choosing $\delta = u/n + \epsilon_{\app}(r,P_n)$,  the second term in \eqref{decomp_h_proof} equals
\begin{equation}\label{bd_approx_proof}
    \inf_{H\in\cP(\cB_r)} \PP\left\{ \ell_n(G^*)-\ell_n(H)  \ge {u\over n} + \epsilon_{\app}(r,P_n)  \right\} \le e^{-u}.
\end{equation}
For the first term in \eqref{decomp_h_proof}, for any fixed $G\in\cP(\cB_r)$ satisfying $h_n^2(p_G,p^*) \ge \gamma$, we have
\begin{align}\label{bd_moment_proof}\nonumber
    \EE\exp\left\{{n\over2}[\ell_n(G)-\ell_n(G^*)]\right\}
    &=\prod_{i=1}^n\left[
      \sqrt{p_G(x_i)p^*(x_i)}+
      \sqrt{\{1-p_G(x_i)\}\{1-p^*(x_i)\}}\right]\\
    &\le\exp\left\{-n h_n^2(p_G,p^*)\right\} ~\le \exp\{-n  \gamma\}.
\end{align}
In other words, a function $p_G\in\cF_r$ with large $h_n^2(p_G,p^*)$ is unlikely to have a log-likelihood close to that of $p^*$. By Markov's inequality, \eqref{bd_moment_proof} and choosing $\gamma = u/n +   \delta / 2$,
\begin{align*}
    \PP\left\{
        h_n^2(p_G,p^*) \ge \gamma,\,  \ell_n(G) -\ell_n(G^*)
    \ge     -\delta
    \right\} \le e^{-u}.
\end{align*}
What remains is to make the above bound uniform over $p_G\in \cF_r$ to accommodate $\wh p_r$. To this end, take an $\epsilon$-cover of $\cF_r$ under $\|\cdot\|_{\i,\cX}$ with $\epsilon=n^{-2}$ and $\cX = \{x_1,\ldots,x_n\}$. On the event $\{h_n^2(\wh p_r,p^*)\ge\gamma\}$, the estimator $\wh p_r$ lies within $2\epsilon$ of some $p_j\in\cF_r$ with $h_n^2(p_j,p^*)\ge\gamma$, where each $p_j$ is fixed independently of $Y_1,\ldots,Y_n$ and belongs to a collection of at most $\cN(\epsilon,\cF_r,\|\cdot\|_{\i,\cX})$ functions. Consequently, the likelihood of $\wh p_r$ is at most that obtained by replacing $p_j(x_i)$ and $1-p_j(x_i)$ with $p_j(x_i)+2\epsilon$ and $1-p_j(x_i)+2\epsilon$, and the calculation in \eqref{bd_moment_proof}, applied with these enlarged values, gives the bound $\exp\{-n\gamma+2n\sqrt{2\epsilon}\}$. A union-bound argument over this collection, together with the entropy bound in \cref{thm:metric_entropy} with $R_{x,n}$ in lieu of $R_x$, gives
\begin{equation}\label{bd_est_proof}
    \PP\left\{
        h_n^2(\wh p_r,p^*) \ge \gamma,\,  \ell_n(\wh G_r) -\ell_n(G^*)
    \ge    -\delta
    \right\} \le e^{-u}
\end{equation}
for $\gamma \asymp u/n +  \epsilon_{\app}(r,P_n) + \epsilon_{\est}(r,R_{x,n})$, with implicit constants depending only on $d$.

Combining  \eqref{bd_approx_proof} and \eqref{bd_est_proof} gives $\PP\{h_n^2(\wh p_r,p^*) \lesssim  \gamma \} \ge  1-2e^{-u}$. Replacing $u$ by $u+\log 2$ gives the stated probability bound, after enlarging the constant since $u\ge1$. We defer the  full proof of \cref{thm:mse_rate_fixed} to \cref{app:proof_mse_rate_fixed}.
\end{proof}

\begin{proof}[Proof sketch of \cref{thm:mse_rate_random}]
The preceding likelihood comparison for fixed design remains the starting point, with additional steps to handle the randomness and unboundedness of $X$.  First suppose $P_X(\cB_T)=1$ for some finite $T\ge 0$.
For $p,q\in\cF$, define
$
    h_X^2(p,q):=P_Xh^2(p,q),
$
the squared Hellinger distance between the joint distributions of $(Y,X)$ corresponding to $p$ and $q$ under $P_X$. Then
$
    P_X(p-p^*)^2\le4h_X^2(p,p^*).
$
To bound $h_X^2(\wh p_r, p^*)$ from above, we follow the arguments for the fixed-design setting to obtain the decomposition \eqref{decomp_h_proof}. The argument used to prove \eqref{bd_approx_proof} controls the approximation error, which in this case is $\epsilon_{\app}(r,P_X)$. To obtain an analogue of \eqref{bd_est_proof}, we invoke the uniform likelihood-ratio inequality of \citet[Theorem~1]{wong1995probability}. This requires controlling the bracketing number of the joint densities induced by $\cF_r$ on $\cB_T$, which, as shown in \cref{lem:hellinger_bracketing_transfer_random}, can be achieved using our entropy bound in \cref{thm:metric_entropy}.

When $P_X$ is sub-Gaussian, take $T=T_u$ from \eqref{def_Tn}, so that $\{R_{x,n}\le T_u\}$ holds with probability at least $1-e^{-u}$. Conditional on this event, the observations $(X_1,Y_1),\ldots,(X_n,Y_n)$ remain i.i.d., with covariate distribution $P_X^{(T)}:=P_X(\cdot\mid\cB_T)$ and the same regression function $p^*$. The estimator $\wh p_r$ still satisfies \eqref{eq:fixed_design_oracle_lr}, so the preceding result for the bounded design applies to $P_X^{(T)}(\wh p_r-p^*)^2$. Moreover, since $\epsilon_{\app}(r,P_X^{(T)})\le2\epsilon_{\app}(r,P_X)$ using  $P_X(\cB_T)\ge1/2$ and
$$
    P_X(\wh p_r-p^*)^2
    \le P_X^{(T)}(\wh p_r-p^*)^2+P_X(\cB_T^c),
$$
the result of \cref{thm:mse_rate_random} follows. The full proof is given in  \cref{app:proof_mse_rate_random}.
\end{proof}

\section{A finite-particle approximation to the NPMLE}
\label{sec:computation_npmle}
Our proposed reweighted particle-gradient algorithm (RPG) for approximating the
restricted NPMLE is stated in  \cref{sec:particle_prediction_formulation}. Its global convergence for any fixed dataset is established in \cref{sec:particle_fixed_data}, while the results in \cref{sec:particle_population} account for the randomness of the data and provide  high probability result of global convergence of RPG when the true mixing distribution is finite.

\subsection{Particle representation and the proposed algorithm}
\label{sec:particle_prediction_formulation}

For any $G\in\cP(\cB_r)$, recall that $p_G:\cX\to(0,1)$ denotes the regression function defined in \eqref{def_p_G}. Let $\bm{p}_G=(p_G(x_1),\ldots,p_G(x_n))^\T \in (0,1)^n$ be its empirical vector at $n$ observed covariates. In particular, the Dirac measure $\delta_\beta$ at any $\beta\in\cB_r$ gives
$
    \bm{p}_{\delta_\beta}
    =(\sigma(x_1^\T\beta),\ldots,\sigma(x_n^\T\beta))^\T,
$
and
\begin{equation*}
    \bm{p}_G=\int_{\cB_r}\bm{p}_{\delta_\beta}\,\rd G(\beta).
\end{equation*}
Consequently,  the vectors $\bm{p}_G$, for all $G\in \cP(\cB_r)$, form the convex hull of $\{\bm{p}_{\delta_\beta}:\beta\in\cB_r\}$. Since the
log-likelihood $\ell_n(G)$ in \eqref{def_llh} depends on $G$ only through
$\bm{p}_G$ and is strictly concave as a function of $\bm{p}_G$, its maximizing
prediction vector $\bm{p}_{\wh G_r}$ is unique, although the mixing distribution $\wh G_r$ representing it
need not be.
Nevertheless, Carath\'eodory's theorem \citep{lindsay1983geometry} ensures that $\wh G_r$ can be chosen in the form
\begin{equation}\label{eq:caratheodory_npmle_prediction}
    \wh G_r=\sum_{\ell=1}^{L^*}a_\ell^*\delta_{\beta_\ell^*},
\end{equation}
where $L^*\le n+1$, $a^*\in\Delta^{L^*-1}$, and $\beta_\ell^*\in\cB_r$ for $\ell\in[L^*]$. Thus, even when the true mixing distribution $G^*$ is continuous, a restricted NPMLE with at most $n+1$ atoms exists. This is the motivation of maximizing $\ell_n(G)$ over finite-particle mixing distributions with a sufficiently large number $M$ of particles:
\begin{equation}
    G=\sum_{j=1}^M w_j\delta_{\beta_j},
    \qquad
    w=(w_1,\ldots,w_M)^\T\in\Delta^{M-1},
    \qquad
    \beta_1,\ldots,\beta_M\in\cB_r.
    \label{eq:particle_measure}
\end{equation}

\paragraph{Location updates} We start by  describing how $\ell_n(G)$ changes when an infinitesimal amount of mass is moved in $G$, to motivate the proposed updates of the particle locations $\beta_j$ and the weights in $w$. For any $p_G$, define its empirical residual vector by $\bm{r}(p_G)=(r_1(p_G),\ldots,r_n(p_G))^\T\in \RR^n$, where, for each $i\in[n]$,
\begin{equation}
    r_i(p_G)
    ~ :=~ n\frac{\partial\ell_n(G)}{\partial p_G(x_i)}
    ~ =~ \frac{Y_i}{p_G(x_i)}-\frac{1-Y_i}{1-p_G(x_i)}.
    \label{eq:particle_residual_potential}
\end{equation}
For any atom $\beta\in\cB_r$, replacing $G$ by $(1-\epsilon)G+\epsilon\delta_\beta$ changes its prediction vector from $\bm{p}_G$ to $(1-\epsilon)\bm{p}_G+\epsilon\bm{p}_{\delta_\beta}$. Therefore, by writing  $\langle u,v\rangle_n:=n^{-1}\sum_{i=1}^n u_iv_i$ for any $u,v\in \RR^n$,
\begin{align}\label{eq:potential_directional_derivative}
    \left.\frac{d}{d\epsilon}
    \ell_n\bigl((1-\epsilon)G+\epsilon\delta_\beta\bigr)
    \right|_{\epsilon=0}
    &~ =~ -\langle\bm{r}(p_G),\bm{p}_G\rangle_n
      +\langle\bm{r}(p_G),\bm{p}_{\delta_\beta}\rangle_n.
\end{align}
By further writing
\begin{equation} \label{eq:particle_potential}
g_{n,G}(\beta)
    :=\langle\bm{r}(p_G),\bm{p}_{\delta_\beta}\rangle_n,
    \qquad
    \bar g_{n,G}
    :=  \langle\bm{r}(p_G),\bm{p}_G\rangle_n = \int_{\cB_r}g_{n,G}(\beta)\,\rd G(\beta),
\end{equation}
we call $g_{n,G}(\beta)$ the empirical potential at $\beta$, and $\bar g_{n,G}$ its average under the current mixing distribution $G$. Their difference, $g_{n,G}(\beta)-\bar g_{n,G}$, gives the first-order likelihood gain from moving an infinitesimal amount of mass from $G$ to $\beta$. Since $\bar g_{n,G}$ does not depend on $\beta$, maximizing this gain amounts to maximizing $g_{n,G}(\beta)$ over $\beta\in \cB_r$.

The commonly used Frank--Wolfe method for computing NPMLE builds on this observation by searching for a global maximizer of $g_{n,G}$ over $\cB_r$ and adding the resulting atom \citep{frank1956algorithm,holloway1974extension}. However, maximizing $g_{n,G}$ is itself a challenging problem due to its non-concavity.  We refer to \cref{rem:particle_fixed_data_comparison} for discussion on related algorithms.

Our proposed method instead moves {\em all} $M$ particles in
\eqref{eq:particle_measure} along the gradient of the empirical potential,
\begin{equation}\label{eq:grad_potential}
    \nabla g_{n,G}(\beta)
    =\frac{1}{n}\sum_{i=1}^n
      r_i(p_G)\sigma'(x_i^\T\beta)x_i,
\end{equation}
evaluated at each particle's current location $\beta=\beta_j$, $j\in[M]$. This is   detailed in \eqref{eq:reweighted_particle_update} below.
When differentiating $g_{n,G}(\beta)$ with respect to $\beta$,
we hold $G$ fixed, and primes on $\sigma$ denote derivatives
with respect to its scalar argument.

\paragraph{Weights update} After updating the particles, we hold $\bbeta=(\beta_1,\ldots,\beta_M)\in\RR^{d\times M}$ fixed and update the weights  by maximizing $\ell_n$. The set of optimal weights, relative to $\bbeta$, is given by
\begin{equation}\label{def_Mn}
    \cM_n(\bbeta)~ :=~ \argmax_{w\in\Delta^{M-1}}
    ~ \ell_n\left(\sum_{j=1}^M w_j\delta_{\beta_j}\right).
\end{equation}
Since the objective is concave in $w$ and $\Delta^{M-1}$ is convex, this is a convex optimization problem. We call any $w\in\cM_n(\bbeta)$ {\em exactly reweighted} at $\bbeta$ and such weights need not be unique.
The exact reweighting allows us to distinguish the gain from updating the particle locations from that obtained by updating the weights.  Indeed, for any $w\in\cM_n(\bbeta)$ with $G=\sum_{j=1}^M w_j\delta_{\beta_j}$, the KKT conditions for \eqref{def_Mn} (see, also, \cref{lem:app_particle_weight_kkt}) assert
\begin{equation}
    g_{n,G}(\beta_j)\le\bar g_{n,G},
    \qquad
    w_j\bigl(g_{n,G}(\beta_j)-\bar g_{n,G}\bigr)=0,
    \qquad j\in[M].
    \label{eq:exact_reweight_t_kkt}
\end{equation}
Thus, every active particle, with $w_j>0$, satisfies $g_{n,G}(\beta_j)=\bar g_{n,G}$, while every inactive particle with $w_\ell = 0$ satisfies $g_{n,G}(\beta_\ell)\le\bar g_{n,G}$. Consequently, no reweighting over the current locations can increase the likelihood, so any further increase requires moving the particles.\\

Combining both updates, our proposed method is an iterative procedure of outer and inner iterations. Each outer iteration consists of $H\ge1$ inner steps that update the particle locations with the weights held fixed, followed by one exact reweighting. Specifically, at outer iteration $t\ge0$, let $w^t$  and $\bbeta^t=(\beta_1^t,\ldots,\beta_M^t)\in\RR^{d\times M}$ denote the current weights and  particle locations. Set $\beta_j^{t,0}=\beta_j^t$ for $j\in[M]$. For each inner step $h=0,\ldots,H-1$, define
$$
    G^{t,h}=\sum_{j=1}^M w_j^t\delta_{\beta_j^{t,h}}
$$
with the weight $w^t$ held fixed,
and update the particle locations by
\begin{equation}
    \beta_j^{t,h+1}
    =\Pi_{\cB_r}\left\{
      \beta_j^{t,h}+\eta\nabla g_{n,G^{t,h}}(\beta_j^{t,h})
    \right\},
    \qquad j\in[M].
    \label{eq:reweighted_particle_update}
\end{equation}
Here $\Pi_{\cB_r}$ denotes Euclidean projection onto $\cB_r$,
and $\eta>0$ is the step size. The potential gradient
$\nabla g_{n,G^{t,h}}$ is recomputed at each inner step.

For fixed weights, the ordinary location gradient satisfies
$\nabla_{\beta_j}\ell_n(G^{t,h})=w_j^t\,\nabla g_{n,G^{t,h}}(\beta_j^{t,h})$.
The update in \eqref{eq:reweighted_particle_update} omits the
multiplicative factor $w_j^t$, allowing particles with small
weights to move without this attenuation and zero-weight
particles to potentially become active after moving and
reweighting. When all weights in $w^t$ are positive and held fixed,   \eqref{eq:reweighted_particle_update} can also be interpreted as a projected Wasserstein gradient-ascent step for the likelihood function \citep{yan2024learning}.

After $H$ inner steps, we set the particle locations as $\bbeta^{t+1}=(\beta_1^{t,H}, \ldots,\beta_M^{t,H})$ and update the weights to conclude the $t$-th outer iteration by
\begin{equation}\label{eq:exact_reweight}
    w^{t+1}\in\cM_n(\bbeta^{t+1}),
\end{equation}
giving the next iterate
$$
    G_{t+1}=\sum_{j=1}^M w_j^{t+1}\delta_{\beta_j^{t+1}}.
$$
Iterations in \eqref{eq:reweighted_particle_update}--\eqref{eq:exact_reweight} can be run for a specified number of outer iterations or terminated when the change in the parameters or the log-likelihood between successive outer iterations is sufficiently small. See \cref{sec_sim} for the stopping rule used in our simulations.

For initialization, our theory in \cref{sec:particle_fixed_data,sec:particle_population}  suggests drawing the locations $\beta_1^0,\ldots,\beta_M^0$ independently and uniformly within $\cB_r$ and computing $w^0\in\cM_n(\bbeta^0)$. We summarize the whole procedure in \cref{alg:fully_corrective_particle_gradient}.

\begin{algorithm}[ht]
\caption{Reweighted particle-gradient algorithm (RPG)}
\label{alg:fully_corrective_particle_gradient}
\begin{algorithmic}[1]
\Require Data $\{(x_i,Y_i)\}_{i=1}^n$, radius $r$, particle count $M$,
number of blocks $T$, block length $H$, and step size $\eta$.
\State Initialize $\bbeta^0=(\beta_1^0,\ldots,\beta_M^0)\in\cB_r^M$.
\State Compute $w^0\in\cM_n(\bbeta^0)$.
\For{$t=0,\ldots,T-1$}
    \State Set $\bbeta^{t,0}=\bbeta^t$.
    \For{$h=0,\ldots,H-1$}
        \State Compute $\bm{p}_{G^{t,h}}$, $\bm{r}(\bm{p}_{G^{t,h}})$, and
        $\nabla g_{n,G^{t,h}}(\beta_j^{t,h})$.
        \State Update all locations by \eqref{eq:reweighted_particle_update}.
    \EndFor
    \State Set $\bbeta^{t+1}=\bbeta^{t,H}$.
    \State Compute $w^{t+1}\in\cM_n(\bbeta^{t+1})$.
\EndFor
\State \Return $G_T=\sum_{j=1}^M w_j^T\delta_{\beta_j^T}$ and $p_{G_T}$.
\end{algorithmic}
\end{algorithm}

\subsection{Convergence to the restricted NPMLE for fixed data}
\label{sec:particle_fixed_data}

Throughout this subsection, the observations $\{Y_i,x_i\}_{i\in[n]}$ are fixed. We establish the global convergence of the RPG iterates $\{G_t\}_{t\ge1}$ by showing that its induced regression function $p_{G_t}$ converges to the restricted NPMLE $\wh p_r = p_{\wh G_r}$ as $t\to\infty$.

We first establish a key property that $\ell_n(G_t)$ is nondecreasing in $t$. The next lemma bounds from below the likelihood gain of each outer iteration  in terms of the following {\em energy} of its first location update,
\begin{equation}
    E_t:=\frac{1}{\eta^2}\sum_{j=1}^M
    w_j^t\|\beta_j^{t,1}-\beta_j^{t,0}\|_2^2,
    \qquad t\ge0.
    \label{eq:particle_gap_energy}
\end{equation}
Note that $\eta^2 E_t$ is the weighted sum of squared particle displacements in the first inner step of outer iteration $t$.  Recall that $R_{x,n}:=\max_{i\in[n]}\|x_i\|_2$ and set $\rho_{n,r}:=\sigma(-rR_{x,n})$.

\begin{samepage}
\begin{lemma}
\label{lem:particle_one_block_ascent}
If $R_{x,n}>0$ and the step size  in \eqref{eq:reweighted_particle_update}  satisfies $0<\eta\le  8\rho_{n,r}^2 / R_{x,n}^2$,
then the iterates $\{G_t\}_{t\ge 0}$ in
\cref{alg:fully_corrective_particle_gradient} satisfy
\begin{equation}
    \ell_n(G_{t+1})-\ell_n(G_t) ~ \ge ~  {\eta\over2}E_t,
    \qquad t\ge 0.
    \label{eq:first_step_progress}
\end{equation}
\end{lemma}
\end{samepage}

\cref{lem:particle_one_block_ascent} shows that $\ell_n(G_t)$ is nondecreasing along the proposed iterates $\{G_t\}_{t\ge1}$. However, this does not guarantee global convergence, as the algorithm may converge to a suboptimal stationary point. Indeed, although $\ell_n(G)$ is concave in both the prediction vector $\bm p_G$ and the mixing measure $G$, it is generally {\em not} concave in the particle locations. Moreover, the set of mixing measures with at most $M$ atoms is generally {\em nonconvex}. Consequently, a stationary point of the finite-particle \cref{alg:fully_corrective_particle_gradient} need not maximize $\ell_n$. The main difficulty in establishing global convergence to the restricted NPMLE is therefore to show that the energy $E_t$ remains bounded away from zero whenever $\ell_n(G_t)$ is away from the optimum.

For fixed data, we prove, in \cref{lem:app_particle_fixed_order_energy} of \cref{app:particle_fixed_data_auxiliary}, the following key inequality
\begin{equation}\label{lb_energy}
    E_t ~ \ge ~ c\left\{\ell_n(\wh G_r)-\ell_n(G_t)\right\}^2,
\end{equation}
whenever the likelihood gap $\ell_n(\wh G_r)-\ell_n(G_t)$ is small. Together with \cref{lem:particle_one_block_ascent}, this gives
$$
    \ell_n(G_{t+1})-\ell_n(G_t)
    ~ \ge~  \frac{c\,\eta }{2}\left\{\ell_n(\wh G_r)-\ell_n(G_t)\right\}^2.
$$
A recursive argument readily yields an $\cO(t^{-1})$ convergence rate for $\ell_n(\wh G_r)-\ell_n(G_t)$, provided that the initial gap $\ell_n(\wh G_r)-\ell_n(G_0)$ is sufficiently small. The latter holds for sufficiently large $M$ under the proposed initialization, as shown in \cref{lem:app_particle_fixed_initialization} of \cref{app:sec_initial_fixed}. The following theorem summarizes these results.

\begin{theorem}
\label{thm:reweighted_particle_global}
Fix $r>0$ and $H\ge1$. Suppose that $R_{x,n}>0$ and $0<\eta\le  8\rho_{n,r}^2 / R_{x,n}^2$. For every $\delta\in(0,1)$, there exists a finite $M_0$ such that, whenever $M\ge M_0$, the following holds with probability at least $1-\delta$ over the random initialization:
\begin{equation}
    \ell_n(\wh G_r)-\ell_n(G_t) ~ \le ~  {C\over t+c},
    \qquad t\ge0,
    \label{eq:particle_global_rate}
\end{equation}
for some $C,c>0$ that do not depend on $M$. In particular, as $t\to\infty$,
$$
    \ell_n(G_t)\to\ell_n(\wh G_r),
    \qquad
    \bm{p}_{G_t}\to\bm{p}_{\wh G_r}.
$$
\end{theorem}
\begin{proof}[Proof sketch]
    To prove the key inequality \eqref{lb_energy}, we first establish a positive gap between $\ell_n(\wh G_r)$ and the likelihood values attained at suboptimal stationary points of the finite-particle algorithm. This follows from the subanalytic structure of the stationary set, which yields only finitely many stationary likelihood values (see \cref{lem:app_particle_stationary_values}). If the initial $\ell_n(G_0)$ lies within this gap, monotonicity in \cref{lem:particle_one_block_ascent} keeps the iterates $\ell_n(G_t)$ in this region and excludes suboptimal stationary limits, thereby ensuring global convergence. Obtaining the quantitative lower bound on $E_t$ in \eqref{lb_energy}, however, requires further analysis based on the analytic structure of the likelihood (see \cref{lem:app_particle_projected_curve,lem:app_particle_fixed_order_energy}).

    To ensure that $\ell_n(G_0)$ lies in the gap, we choose $M\ge M_0$ so that, with high probability, each atom of a finite representation of the restricted NPMLE has a nearby initial particle. Exact reweighting then makes the initial likelihood gap small enough for the bound to apply; see \cref{lem:app_particle_fixed_initialization}. This initialization requirement determines $M_0$, which depends on $r$, the data, and $\log(1/\delta)$. We defer the full proof to \cref{app:particle_fixed_data_analysis}.
\end{proof}

\cref{thm:reweighted_particle_global} establishes that $\ell_n(G_t)$ converges to the global maximum over $\cP(\cB_r)$ despite the nonconvexity of the finite-particle formulation. It holds for any fixed dataset and imposes no restriction on $G^*$. In particular, $G^*$ need not have finite support.
Choosing a sufficiently large number of particles $M$ is crucial for global convergence. More importantly, the convergence rate is independent of $M$ once it exceeds a certain threshold.

The constants $C$ and $c$ in \eqref{eq:particle_global_rate} depend on the step size $\eta$ and the data, in particular, the sample size $n$, through the quantity $c$ in the energy lower bound \eqref{lb_energy}.
The preceding proof of \eqref{lb_energy} relies heavily on the observed data being fixed, making it difficult to extend the result uniformly over data realizations. Establishing such a uniform bound is substantially more challenging and requires a different approach. In the next section, we obtain such control under further conditions on $G^*$ and the covariates.

\begin{remark}[Existing algorithms for computing the NPMLE]
\label{rem:particle_fixed_data_comparison}
As mentioned before \eqref{eq:grad_potential}, conditional-gradient methods, also known as Frank--Wolfe methods, successively add a support point by maximizing the potential $g_{n,G}(\beta)$ over $\beta\in\cB_r$. Its fully corrective variant  further  reoptimizes the weights over the expanded support \citep{jagabathula2020conditional,jiang2025nonparametric}. When the potential maximization is solved globally, or approximately with  suitable decay errors \citep[Theorem~2]{jiang2025nonparametric}, these methods attain an $\cO(1/t)$ convergence rate for the likelihood gap \citep[Theorem~1]{jaggi2013revisiting}, matching our \cref{thm:reweighted_particle_global}. However, since $g_{n,G}(\beta)$ is nonconcave in $\beta$, its global maximization, or even sufficiently accurate approximation, is a challenging problem on its own. In contrast, \cref{thm:reweighted_particle_global} shows that RPG achieves global convergence of the likelihood with a sufficiently large but finite number of particles, by simultaneously updating their locations via
projected gradient steps followed by exact reweighting. In \cref{sec_sim_npmle}, we also demonstrate numerical advantages of RPG over these conditional-gradient methods, including computational gains in high dimensions.

Another related approach optimizes $G$ directly over probability measures, with particle implementations that move support points as our location updates do. For the Wasserstein--Fisher--Rao method of \citet{yan2024learning}, full-support initialization and suitable step-size conditions ensure that, if the measure-valued iterates converge weakly, their limit is globally optimal. However, this result provides neither a convergence rate nor a finite-particle guarantee. Quantitative convergence guarantees are available for another measure-space approach, mean-field Langevin dynamics, which optimizes an entropy-regularized objective and converges exponentially in continuous time under a uniform log-Sobolev inequality \citep{nitanda2022convex}. Its finite-particle convergence bounds nevertheless include an additional approximation term of order $1/M$ \citep{suzuki2023convergence,nitanda2024improved}. In contrast, \eqref{eq:particle_global_rate} applies directly to finite-particle iterates, yielding a convergence bound independent of
$M$ for $M\ge M_0$, without any particle-approximation error term.
\end{remark}

\subsection{Global convergence and statistical guarantees}
\label{sec:particle_population}

In this section we establish global convergence and statistical guarantees for the proposed RPG algorithm.  We consider the random-design setting with $x_1,\ldots,x_n$ being i.i.d. realizations from the covariate distribution $P_X$. We impose the following conditions on $P_X$ and the true mixing distribution $G^*$.

\begin{assumption}
\label{assum:particle_population}
For a fixed dimension $d$, the following conditions hold.
\begin{itemize}[leftmargin=7mm,itemsep=0mm]
\item[(i)] \emph{(Mixing distribution.)} For a fixed positive integer $K_0$ and fixed constants $w_0,\allowbreak s_-,\allowbreak s_+,\allowbreak r_\beta>0$, the true mixing distribution has the form $G^*=\sum_{a=1}^{K^*}w_a^*\delta_{\beta_a^*}$, where $K^*\le K_0$, $w_a^*\ge w_0$ and $\|\beta_a^*\|_2\le r_\beta$ for all $a\in[K^*]$. The atoms further satisfy
\begin{equation}
    \min_{a\ne b}\|\beta_a^*-\beta_b^*\|_2\ge s_-,
    \qquad
    \min_{a,b}\|\beta_a^*+\beta_b^*\|_2\ge s_+.
    \label{eq:population_atom_separation}
\end{equation}

\item[(ii)] \emph{(Covariate distribution.)} For fixed constants $0<r_X\le R_X$ and $m_X>0$, $P_X$ is supported on $\cB_{R_X}$ and has a density
$f_X$ satisfying $f_X(x)\ge m_X$ for all $x\in\cB_{r_X}$.
\end{itemize}
\end{assumption}
\begin{remark}
\label{rem:particle_separation}
Part $(i)$ of \cref{assum:particle_population} assumes a finitely supported $G^*$, which is  commonly used to model heterogeneous binary responses \citep{wedel1995mixture,grun2007fitting}. Our RPG algorithm, however, does not require knowledge of either $K^*$ or $K_0$.  Weights bounded away from zero and separated atoms are also common assumptions in finite mixture models. The second condition in \eqref{eq:population_atom_separation}, referred to as the {\em reflection condition},
is specific to logistic mixture models. By excluding the origin and reflected pairs from the true support, it rules out an ambiguity arising from sigmoid symmetry: since $\sigma(t)+\sigma(-t)=1$, one has $p_{(\delta_\beta+\delta_{-\beta})/2}=p_{\delta_0}\equiv1/2$ for any $\beta$. Similar restrictions on coefficient directions are imposed in random-coefficient binary choice models \citep{gautier2013nonparametric,gautier2018adaptive};
see \cref{sec_ident} for a detailed discussion.

Part $(ii)$ of \cref{assum:particle_population} imposes standard conditions on $P_X$: bounded covariates are assumed in analyses of the NPMLE for mixtures of linear regressions \citep[Theorem~4]{jiang2025nonparametric}, and density lower bounds are common in classification \citep[Definition~2.2]{audibert2007fast}. In our analysis, bounded covariates and the sieve radius keep $p_G(X)$ uniformly away from zero and one, yielding uniform bounds on the log-likelihood function and its derivatives. The density lower bound ensures that the estimation error of the regression function under $P_X$ controls the integrated estimation error over $\cB_{r_X}$ with respect to Lebesgue measure.

Finally, fixing the constants in \cref{assum:particle_population} allows us to establish uniform guarantees over the resulting class. The constants $w_0$, $s_-$, $s_+$, $r_X$, and $m_X$ may be arbitrarily small, although the bounds given below may deteriorate as they approach zero.
\end{remark}

Under \cref{assum:particle_population}, we are able to derive a denoising result in \cref{thm:app_population_inverse} of \cref{app:particle_denoising} that controls a discrepancy between $G$ and $G^*$ by the squared error $P_X(p_G-p_{G^*})^2$, uniformly over $G\in\cP(\cB_r)$. This is stronger than identifiability and is detailed in \cref{sec_ident}. The denoising bound is crucial to establishing an energy lower bound of the form \eqref{lb_energy}, with the constant $c$ depending only on the quantities specified in \cref{assum:particle_population}. Specifically, at each iteration $t$, we use the discrepancy between $G_t$ and $G^*$ to construct a feasible reference mixing distribution that reduces the population likelihood gap $\cL(G^*)-\cL(G_t)$. We then compare the next iterate $G_{t+1}$ with this reference to derive the lower bound on $E_t$. See \cref{sec:particle_proof_techniques} for details.

Throughout this subsection, for a fixed $\tau_0>0$ and with $r_\beta$ given in   \cref{assum:particle_population}$(i)$, we choose the sieve radius as
\begin{equation}
    r=r_\beta+\tau_0.
    \label{eq:particle_radius_choice}
\end{equation}
The sample-level convergence for the RPG algorithm in \cref{thm:particle_uniform_empirical} is our main result. Its proof builds on the population-level result in \cref{thm:particle_population}, which we present first.

\subsubsection{Population-level convergence guarantees}
\leavevmode\par

For any $G\in\cP(\cB_r)$, define the population-level log-likelihood by
\begin{equation}
    \cL(G) :=\EE_X\bigl[p^*(X)\log p_G(X)
       +\{1-p^*(X)\}\log\{1-p_G(X)\}\bigr].
    \label{eq:population_likelihood_gap}
\end{equation}
For any $\bbeta = (\beta_1,\ldots,\beta_M)$, its corresponding set of optimal weights (with respect to $\cL$) is
$$
    \cM(\bbeta)~ :=~ \argmax_{w\in\Delta^{M-1}}
    ~ \cL\biggl(\sum_{j=1}^M w_j\delta_{\beta_j}\biggr).
$$
The population-level counterpart of \cref{alg:fully_corrective_particle_gradient} replaces the empirical average in the location updates \eqref{eq:reweighted_particle_update} by an expectation under $P_X$ and uses the above population exact reweighting step $w^t\in\cM(\bbeta^t)$ in place of \eqref{eq:exact_reweight}. The mixing distribution $G_t$ at the end of the $t$-th outer iteration is obtained after $Ht$ location updates and $t+1$ exact reweightings: one at initialization and one at the end of each completed outer iteration. The following theorem establishes its global convergence to the true mixing distribution $G^*$.

\begin{theorem}
\label{thm:particle_population}
Under \cref{assum:particle_population}, with $r$ chosen as in
\eqref{eq:particle_radius_choice}, fix an integer $H\ge1$, a step size
$0<\eta\le 27\sigma(-rR_X)/(4R_X^2)$, and
$\delta\in(0,1)$. There are constants $C,c>0$ and a finite $M_0$
such that, for every $M\ge M_0$, with probability at least $1-\delta$ over the random initialization,
\begin{equation}
    2P_X(p_{G_t}-p^*)^2
    \le\cL(G^*)-\cL(G_t)
    \le C\left\{{1\over(1+\eta Ht)^2}+e^{-ct}\right\},
    \quad t\ge0.
    \label{eq:uniform_population_rate}
\end{equation}
Here $C,c$ depend only on $d,\tau_0$ and the constants in
\cref{assum:particle_population}; $M_0$ may additionally depend on
$\log(1/\delta)$. None of these quantities depends on $M,H$, or $\eta$.
\end{theorem}
\begin{proof}
    The proof is given in \cref{app:particle_uniform_analysis}, where the dependence of $c$, $C$, and $M_0$ on the relevant quantities can be tracked. We provide a proof sketch in \cref{sec:particle_proof_techniques}.
\end{proof}

The effect of the step size in \cref{thm:particle_population} is classical: $\eta$ must be sufficiently small to ensure convergence, but smaller values lead to slower convergence. The bound in \eqref{eq:uniform_population_rate} also reveals interesting separate effects of the number of inner steps $H$ and the number of outer iterations $t$, as detailed in the following remark.

\begin{remark}[Reweighting frequency and convergence]
\label{rem:particle_converged_locations}
Exact reweighting optimizes all weights over $\Delta^{M-1}$ and can be computationally expensive for large $M$, whereas projected gradient-ascent updates of particle locations are typically less costly. To balance these costs, for fixed $\eta$, \cref{thm:particle_population} suggests choosing $t=\lceil C\log H\rceil$ for a sufficiently large constant $C>0$, making the exponential term in \eqref{eq:uniform_population_rate} negligible and yielding a likelihood gap of order $(Ht)^{-2}$. This requires only $O(\log H)$ exact reweightings over $Ht$ location updates, suggesting that longer inner loops can reduce computational overhead without sacrificing the theoretical convergence rate. Our simulations in \cref{sec_sim_reweight} further illustrate this trade-off.

A special case is to run the location updates until convergence before each reweighting. In this case, we show in \cref{app:particle_reweighting_remark_proof}  that there exists some $\kappa\in(0,1)$ such that
the likelihood gap decreases geometrically in the number of outer iterations:
$$
    \cL(G^*)-\cL(G_t)
    \le(1-\kappa)^t\{\cL(G^*)-\cL(G_0)\}.
$$
\end{remark}

For fixed $H$ and $\eta$, the population-level RPG algorithm achieves an $\cO(t^{-2})$ convergence rate in \eqref{eq:uniform_population_rate}. The next proposition gives a matching lower bound on its convergence rate, even in a one-dimensional model with a single true atom.

\begin{proposition}
\label{prop:particle_inverse_square_sharp}
There is a fixed one-dimensional model with a single true atom satisfying
\cref{assum:particle_population} and \eqref{eq:particle_radius_choice}
for which the following holds.
Fix an integer $H\ge1$ and a sufficiently small step size $\eta>0$.
There are constants $c>0$ and $\delta_0\in(0,1)$, independent of $M$,
and a deterministic sequence of iteration indices $T_M\to\infty$ as $M\to\infty$,
such that, for every sufficiently large $M$, the population-level RPG algorithm satisfies
$$
    \PP\left\{\cL(G^*)-\cL(G_{T_M})
       \ge {c\over(1+T_M)^2}\right\}\ge \delta_0.
$$
The same bound holds with the likelihood difference replaced by
$P_X(p_{G_{T_M}}-p^*)^2$.
\end{proposition}

The lower bound holds with probability at least $\delta_0>0$ over the random
initialization, where $\delta_0$ does not decrease as $M$ increases.
For fixed $H$, $\eta$, and $\delta<\delta_0$, \cref{prop:particle_inverse_square_sharp} rules out an $o(t^{-2})$ convergence rate for RPG that is independent of $M$ and holds with probability at least $1-\delta$ for every sufficiently large $M$. In this sense, the $\cO(1/t^2)$ rate in \cref{thm:particle_population} cannot be improved.

\begin{remark}[Comparison with overparameterization under finite Gaussian mixtures]
\label{rem:particle_overparam_em}
Global convergence has also been studied for EM and gradient EM in
overparameterized Gaussian location mixtures, where the number of fitted
components exceeds the true number. \citet{dwivedi2020singularity}
analyze EM when a symmetric two-component mixture with known mixing
weights and isotropic covariance is fitted to a single true Gaussian.
Extending beyond two fitted components, \citet{xu2024toward} establish
global convergence of population gradient EM when an arbitrary number
$M>1$ of components is fitted to a single true Gaussian, with an
$\cO(t^{-1/2})$ population likelihood gap. For a true $K^*$-component
mixture, \citet{zhou2025global} obtain an $\cO(t^{-2})$ rate when fitting
$M\gg K^*$ components using gradient EM with the exact reweighting steps, under
separation conditions that depend on $M$ and a step size that decreases
polynomially in $M$.

These results concern Gaussian location mixtures, where the mixing
distribution is identifiable from the mixture density. Our setting
instead concerns the restricted NPMLE for logistic mixtures,
where sigmoid symmetry introduces an additional identifiability issue
(see \cref{rem:particle_separation}). Our analysis relies on a novel
quantitative denoising bound established in
\cref{thm:app_population_inverse} of \cref{app:particle_denoising}.
Moreover, \cref{thm:particle_population} establishes the same
$\cO(t^{-2})$ rate, but with the step-size restriction, separation
conditions on $G^*$, and convergence constants all {\em independent} of $M$
once $M\ge M_0$. Beyond such population guarantee,
\cref{thm:reweighted_particle_global} establishes convergence to the
restricted NPMLE for any fixed dataset without requiring $G^*$ to be
finitely supported.
\end{remark}

\subsubsection{Sample-level convergence guarantees}
\label{sec:particle_uniform_empirical}

Building on the population  convergence of the RPG algorithm in \cref{thm:particle_population}, we now establish its sample-level counterpart.
For any $\delta\in(0,1)$, set the statistical error as
\begin{equation}
    \varepsilon_n(\delta)
    ~ :=~ 1\wedge {A\over n}
    \left[\{\log(Bn)\}^{d+1}+\log{12\over\delta}\right],
    \label{eq:particle_statistical_floor}
\end{equation}
where $A\ge1$ and $B\ge e$ depend only on $d,\tau_0$ and the constants in
\cref{assum:particle_population}. The following theorem establishes the global convergence of the RPG iterates $\{G_t\}_{t\ge1}$ in \cref{alg:fully_corrective_particle_gradient}. Its proof is given in \cref{app:proof_thm:particle_uniform_empirical} with a proof sketch offered in \cref{sec:particle_proof_techniques}.

\begin{theorem}
\label{thm:particle_uniform_empirical}
Grant \cref{assum:particle_population} and choose $r$ as in \eqref{eq:particle_radius_choice}. Fix an integer $H\ge1$, a step size $0<\eta\le27\sigma(-rR_X)/(4R_X^2)$, and $\delta\in(0,1)$. There exist finite $n_0$ and $M_0$ such that, for every $n\ge n_0$ and $M\ge M_0$, with probability at least $1-\delta$ over the data and the random initialization, simultaneously for all $t\ge0$,
\begin{equation*}
    \max\left\{2 P_X(p_{G_t}-p^*)^2,~ \ell_n(\wh G_r)-\ell_n(G_t)\right\}
    ~ \le~  C_{\mathrm{opt}}
        \left\{\frac{1}{(1+\eta Ht)^2}+e^{-ct}\right\}
        +C_{\mathrm{stat}}\, \varepsilon_n(\delta).
\end{equation*}
Here $C_{\mathrm{opt}},C_{\mathrm{stat}},c,n_0,M_0$ depend only on $d,\tau_0$ and the constants in \cref{assum:particle_population}, except that $n_0,M_0$ may additionally depend on $\log(1/\delta)$.
\end{theorem}

\cref{thm:particle_uniform_empirical} bounds both the empirical likelihood gap and the $L_2(P_X)$ MSE of the RPG iterates by an optimization term and a statistical term. The former vanishes as $t\to\infty$, reflecting the population convergence guarantee in \cref{thm:particle_population}; the latter arises from empirical deviations in the likelihood and the quantities governing the location updates. Once $Ht\gtrsim 1/\sqrt{\varepsilon_n(\delta)}$ and $t\gtrsim \log(1/\varepsilon_n(\delta))$, the statistical term dominates, and further iterations do not improve the guaranteed error order $\varepsilon_n(\delta)$.

\subsection{Identifiability of finite mixtures of logistic regressions}
\label{sec_ident}

A by-product of our convergence analysis is an identifiability
result that is of its own interest. Specifically, the regression function $p_{G^*}$ uniquely determines the mixing distribution $G^*$ under a subset of conditions in \cref{assum:particle_population}. This is stated   in the following proposition with its proof given in \cref{app:particle_identifiability}.

\begin{proposition}[Identifiability]
\label{prop:identifiability}
Let $G^*=\sum_{a=1}^{K^*}w_a^*\delta_{\beta_a^*}$, where $K^*<\infty$, $w_a^*>0$, $\beta_a^*\in\cB_r$, and $\beta_a^*+\beta_b^*\ne0$ for all $a,b\in[K^*]$. Suppose that $P_X$ satisfies \cref{assum:particle_population}(ii). Then for any $G\in\cP(\cB_r)$ satisfying $p_G=p_{G^*}$, $P_X$-almost surely, we have $G = G^*$.
\end{proposition}

The competing distribution $G$ in \cref{prop:identifiability} can be any probability measure on $\cB_r$; in particular, it need not be finitely supported. For identifiability of a finitely supported $G^*$, condition  \eqref{eq:population_atom_separation} in \cref{assum:particle_population}
is relaxed in \cref{prop:identifiability} to exclude only the zero atom and exact reflected pairs. The uniform denoising bound in
\cref{thm:app_population_inverse} of \cref{app:particle_denoising}, however,
requires the stronger condition
$\min_{a,b}\|\beta_a^*+\beta_b^*\|_2\ge s_+$
in \eqref{eq:population_atom_separation}. It does not follow from \cref{prop:identifiability} and requires a separate argument, given in \cref{app:particle_denoising}.

\begin{remark}[Existing identifiability results]
 \citet[Section~6]{auder2021least}
establish identifiability of finite logistic mixtures with $\operatorname{supp}(P_X)=\RR^d$ when the weights $w^*$ are all positive and the coefficient matrix $\bbeta^* = (\beta_1^*\ \cdots\ \beta_{K^*}^*)\in \RR^{d\times K^*}$ has rank $K^*.$ The latter implies  $K^*\le d$.
In contrast, \cref{prop:identifiability} does not require such linear
independence and allows $K^*>d$.

Uniqueness results for finite sums of hyperbolic tangent functions also provide a related identifiability result for our model. Indeed, the identity $2p_{G^*}(x)-1=\sum_{a=1}^{K^*}w_a^*\tanh(x^\T\beta_a^*/2)$ connects our model to the results of \citet{sussmann1992uniqueness,albertini1993uniqueness} and \citet{vlacic2021affine}. Under the reflection condition in \cref{prop:identifiability}, these results imply that $G^*$ is identifiable among finitely supported probability measures. Specifically, if $G$ is finitely supported and $p_G=p_{G^*}$ on a nonempty open set, real analyticity extends this equality to all of $\mathbb R^d$, and the cited uniqueness results
then imply $G=G^*$. In contrast, \cref{prop:identifiability} allows arbitrary competing
$G\in\cP(\cB_r)$, including those with infinite support.

For univariate two-component logistic mixtures, \citet{do2025strong} establish first-order identifiability when $\beta_1^*\ne\beta_2^*$, $\beta_1^*+\beta_2^*\ne0$, and $\operatorname{supp}(P_X)$ contains a nonempty open interval. Their result rules out nontrivial perturbations of the weights and locations that leave $p_{G^*}$ unchanged to first order. This local property supports parameter estimation but does not establish uniqueness of $G^*$ among arbitrary $G\in\cP(\cB_r)$, as in \cref{prop:identifiability}.

Finally, restrictions similar in spirit to our reflection condition also appear in the identification of random-coefficient binary-choice models \citep{gautier2013nonparametric,gautier2018adaptive}. We describe these distinct models and compare their restrictions with ours in \cref{app:particle_identifiability}.
\end{remark}

\subsection{Proof techniques of \cref{thm:particle_population,thm:particle_uniform_empirical}}
\label{sec:particle_proof_techniques}

\begin{proof}[Proof sketch of \cref{thm:particle_population}]
Since an analogue of \cref{lem:particle_one_block_ascent} holds for the population likelihood, $\cL(G_t)$ is nondecreasing. The heart of the proof is therefore to derive a quantitative lower bound analogous to \eqref{lb_energy} for the energy at each outer iteration $t$.

To this end, if the iterate $G_t=\sum_{j=1}^M w_j^t\delta_{\beta_j^t}$ is already close to $G^*$ under a suitable discrepancy measure, we can construct a set of candidate particles $\bbeta^{t,\dagger}=(\beta_1^{t,\dagger},\ldots,\beta_M^{t,\dagger})$, called the {\em reference configuration}, such that its corresponding mixing measure $G_t^\dagger=\sum_{j=1}^M w_j^t\delta_{\beta_j^{t,\dagger}}$ has a smaller likelihood gap than $G_t$:
\begin{equation}\label{lh_gap_ref}
    \cL(G^*)-\cL(G_t^\dagger)
    ~ \le ~ (1-\alpha)\left\{\cL(G^*)-\cL(G_t)\right\}
\end{equation}
for some $\alpha\in(0,1)$. Establishing this inequality is one difficulty in our analysis. It also ensures that the reference configuration yields a nonvanishing likelihood improvement whenever the current likelihood gap is positive.

The second challenge is to show that the reference configuration can be approached through the particle gradient-ascent updates in \eqref{eq:reweighted_particle_update}. This is generally not true unless $G_t$ is already in a neighborhood of $G^*$. In this case, we prove that, for each inner step $h=0,1,\ldots,H$,
\begin{equation}\label{local_ref}
    \sum_{j=1}^M w_j^t
    \|\beta_j^{t,h}-\beta_j^{t,\dagger}\|_2^2
     ~ \lesssim ~ \sqrt{\cL(G^*)-\cL(G_t)}.
\end{equation}
Summing over the $H$ inner steps and combining the resulting bound with \eqref{lh_gap_ref}, we deduce
$$
    \sum_{h=1}^H E_{t,h}
    ~ \gtrsim ~ \alpha^2H
    \left\{\cL(G^*)-\cL(G_t)\right\}^{3/2},
$$
where $E_{t,h}$ is the energy of the $h$-th inner step at outer iteration $t$, defined analogously to \eqref{eq:particle_gap_energy}. This further implies that for some constant $\kappa\in(0,1)$,
\begin{equation}\label{lh_one_step}
    \cL(G_{t+1})-\cL(G_t)
    ~ \ge ~ \kappa\eta H
    \{\cL(G^*)-\cL(G_t)\}^{3/2}.
\end{equation}
A recursive argument subsequently leads to  the term of order $\cO\{(\eta Ht)^{-2}\}$ in \eqref{eq:uniform_population_rate}. The other term is easier to obtain from the nondecreasing property of $\cL(G_t)$.

The last difficulty is to show that $G_t$ remains close to $G^*$ for \eqref{lh_gap_ref} and \eqref{local_ref} to hold. Choosing a sufficiently large $M$ ensures that the randomly initialized $\bbeta^0$, followed by exact reweighting for $w^0$, makes the initial likelihood gap $\cL(G^*)-\cL(G_0)$ small (see \cref{lem:app_population_initialization}). Monotonicity then ensures that $\cL(G^*)-\cL(G_t)$ remains small for all $t\ge1$. To translate this likelihood gap into control of the discrepancy between $G_t$ and $G^*$, we invoke the denoising bound in \cref{thm:app_population_inverse} of \cref{app:particle_denoising}, which holds uniformly over $G\in\cP(\cB_r)$ and hence applies to every iterate $G_t$. The discrepancy measure  is introduced in \eqref{eq:app_population_gauge} to accommodate  differing numbers of fitted and true components.
\end{proof}

\begin{proof}[Proof sketch of \cref{thm:particle_uniform_empirical}]
The sample-level proof extends the population-level argument to the empirical likelihood through a series of uniform concentration inequalities established in \cref{lem:app_particle_empirical_comparison} of \cref{app:proof_thm:particle_uniform_empirical}. These inequalities control the deviations between the empirical and population likelihoods, as well as the errors in evaluating the location updates. When the population likelihood gap exceeds $C\varepsilon_n(\delta)$ for some large constant $C>0$, they preserve both the improvement toward the reference configuration in \eqref{lh_gap_ref} and the bound on the weighted squared particle distances in \eqref{local_ref}. Following the argument leading to \eqref{lh_one_step} and using the concentration bounds relating $\cL$ and $\ell_n$, we obtain
$$
    \ell_n(G_{t+1})-\ell_n(G_t)\ge\kappa\{\ell_n(\wh G_r)-\ell_n(G_t)\}\min\left\{\eta H\sqrt{\ell_n(\wh G_r)-\ell_n(G_t)},1\right\}
$$
whenever $\ell_n(\wh G_r)-\ell_n(G_t)>C\varepsilon_n(\delta)$. The resulting recursion gives the algorithmic error in  \cref{thm:particle_uniform_empirical}. Once the empirical likelihood gap reaches the statistical scale, monotonicity keeps it below that scale, and the likelihood comparison gives the estimation error bound. \end{proof}

\section{Simulation studies}\label{sec_sim}

We conduct simulation studies to evaluate the finite-sample performance of the proposed RPG algorithm. In \cref{sec_sim_reweight}, we examine how allocating a fixed total number of iterations between the number of outer-loop iterations and inner-loop updates affects estimation error and computational cost.  In \cref{sec_sim_cv} we examine how the sieve radius $r$ and the number of particles $M$ affect the estimation error, and assess their selection via cross-validation. In \cref{sec_sim_npmle} we compare RPG with alternative algorithms for computing the NPMLE while in \cref{sec_sim_benchmarks} we benchmark RPG against  the EM algorithm and several standard machine-learning methods. %

\paragraph{Data generation process} We consider the finite mixing distribution $G^*=\sum_{a=1}^{K^*}w_a^*\delta_{\beta_a^*}$, where the mixing weights $w^*$ are drawn from a symmetric Dirichlet distribution with concentration parameter $1$, and the location parameters $\{\beta_a^*\}_{a=1}^{K^*}$ are sampled independently and uniformly from $\cB_{r^*}$. The covariates $x_1,\ldots,x_n$ are drawn independently from $N(0,\sigma_X^2 I_d)$ and the responses $Y_1,\ldots, Y_n$ are drawn independently according to model \eqref{def_p_G} using the regression function $p^*:= p_{G^*}$.

\paragraph{Evaluation metrics} To evaluate the estimation error of an estimator $\wh p$ of $p^*$, we use the mean squared error (MSE), $P_{n_{\mathrm{test}}}(\wh p-p^*)^2$, evaluated on a test set of size $n_{\mathrm{test}}$. We write $n_{\mathrm{tot}}=n+n_{\mathrm{test}}$ for the total sample size, where $n$ is the training sample size. Each simulation setting is evaluated over $50$ independent Monte Carlo replications. The baseline configuration is $n_{\mathrm{tot}}=20{,}000$, $n_\mathrm{test}=2{,}000$, $K^*=4$, $d=3$, $r=r^*=3$, $\sigma_X=5$, and $M=64$. For comparisons of computational cost, we also report wall-clock training time.

\paragraph{Specifications of the proposed RPG in \cref{alg:fully_corrective_particle_gradient}} We set the step size of the location updates in \eqref{eq:reweighted_particle_update} to $\eta=1$. We fix the total number of location updates at $B=TH=10{,}000$, with $T$ being the number of outer-loop iterations in \eqref{eq:exact_reweight} and $H$ the number of inner steps of location updates  in \eqref{eq:reweighted_particle_update}. In \cref{sec_sim_reweight}, we examine how allocating between $T$ and $H$ affects estimation error and computational cost. Based on the results therein, we set $H = 1000$ and $T = 10$ for the remaining experiments. %
We use SLSQP \citep{kraft1988software} to solve the reweighting step in \eqref{eq:exact_reweight}, with a tolerance of $10^{-12}$ and at most $1000$ iterations. %
In addition to the choice of $B$, we employ an early stopping criterion for the inner loop: the log-likelihood is evaluated every $20$ location updates, and the inner loop is terminated if the relative change falls below $10^{-8}$ for five consecutive evaluations or if $H$ location updates are reached. At the outer-loop, we terminate the algorithm if the relative change in the log-likelihood stays below  $10^{-7}$ for five consecutive outer iterations or when $T$ outer iterations are reached.

\subsection{Effect of reweighting frequency}\label{sec_sim_reweight}

We examine the effect of allocating inner- and outer-loop iterations
on estimation accuracy and computational cost, fixing the total number
of location updates at $B=TH=10{,}000$. We vary $H$ and set $T=B/H$,
disabling both stopping criteria to ensure exactly $B$ updates.
Within each replication, all choices of $H$ share the same data
and initial particle locations.

\begin{figure}[ht]
    \centering
    \includegraphics[width=1\linewidth]{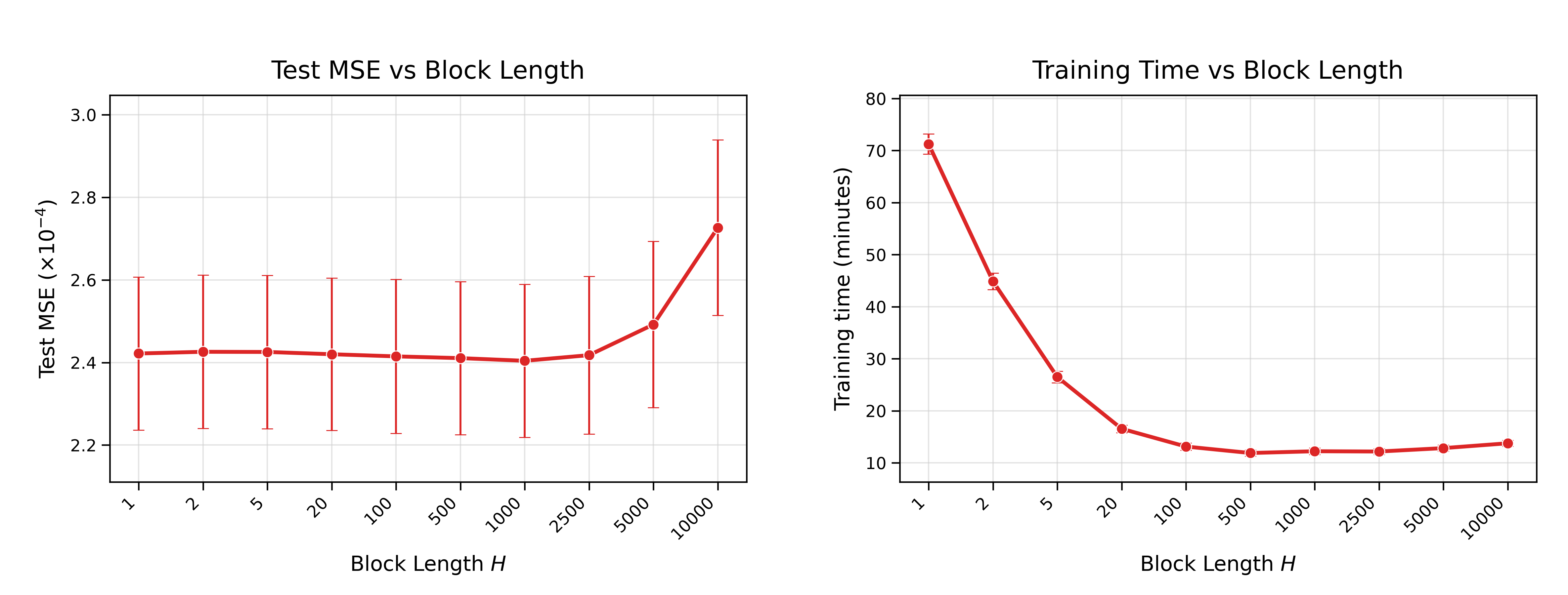}
    \caption{Test MSE and computational cost as functions of the block length $H$.}
    \label{fig:block_length}
\end{figure}

Figure~\ref{fig:block_length} depicts the average (test) MSE and computation time for each $H$ over $50$ replications, with error bars indicating one standard error. The MSE remains nearly unchanged for $H\leq 2500$, but increases for $H=5000$ and $H=10{,}000$. This behavior is qualitatively consistent with \cref{thm:particle_population,thm:particle_uniform_empirical}, which allow multiple particle location updates $H\geq 1$ within each inner loop, provided that the number $T$ of outer-loop iterations is sufficiently large. When $H$ becomes very large, however, the corresponding $T$ is only $2$ or $1$, resulting in too few outer-loop iterations and, consequently, too few weight optimizations. On the computational side, with $B$ fixed, increasing $H$ from $1$ to $1000$ reduces the mean computation time from approximately $71.3$ to $12.2$ minutes, while changing the MSE by less than $1\%$. This substantial reduction in computation time is primarily due to the decrease in the number of weight optimizations from $10{,}001$ to $11$, since solving the weight optimization problem in \eqref{eq:exact_reweight} is relatively expensive for moderate or large $M$. Thus, choosing a relatively large $H$ can substantially reduce computational cost without sacrificing estimation accuracy, provided that the corresponding number $T$ of outer-loop iterations is not too small.

\begin{figure}[ht]
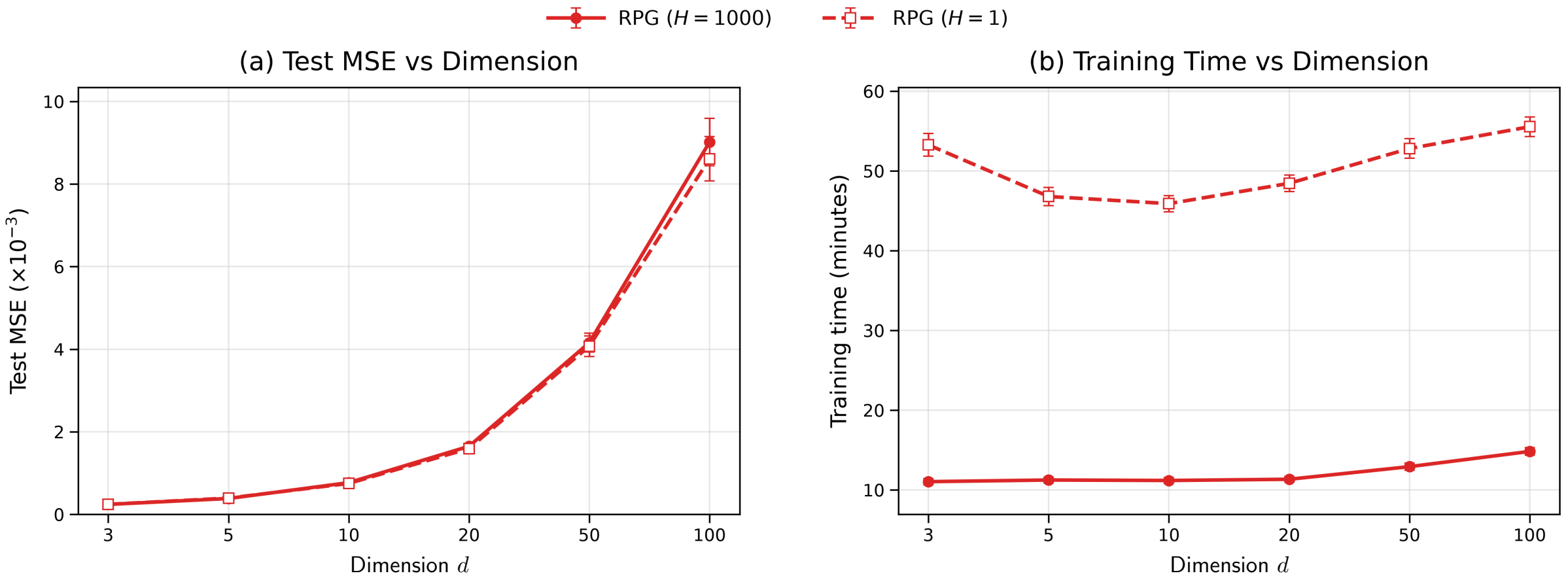

    \centering
    \simulationfigure{two_H}
    \caption{Estimation error and computational cost as functions of the dimension $d$.  We compare $(H,T)=(1000,10)$ and $(H,T)=(1,10000)$.}
    \label{fig:block_length_dimension}
\end{figure}

Figure~\ref{fig:block_length_dimension} further compares $H=1$ and $H=1000$ as the dimension $d$ increases, again with $B$ fixed. Across all dimensions considered, $H=1000$ substantially reduces the computation time while yielding similar test MSEs, with relative differences below $5\%$. This suggests that the computational advantage of less frequent weight optimization persists in higher dimensions, with little difference in estimation accuracy. Based on these results, we use $H=1000$ and $T=10$ for the remaining experiments, with both stopping criteria enabled and $B=10{,}000$.

\subsection{Hyperparameter sensitivity and selection via cross-validation}\label{sec_sim_cv}

We examine the sensitivity of RPG to the number of particles $M$ and the sieve radius $r$, and evaluate their selection via cross-validation (CV). We set $K^*=8$ and consider $d\in\{3,50\}$. In each replication, $90\%$ of the observations are used for training and the remaining $10\%$ are for testing. We consider $M\in\{2,4,8,16,32,64\}$ and $r\in\{1,2,3,5,8,12\}$, and use a 5-fold CV to select  $(r,M)$ based on the log-likelihood of held-out data. For each selected pair, RPG is then refitted using the full training set and evaluated on the test set.

\begin{figure}[ht]
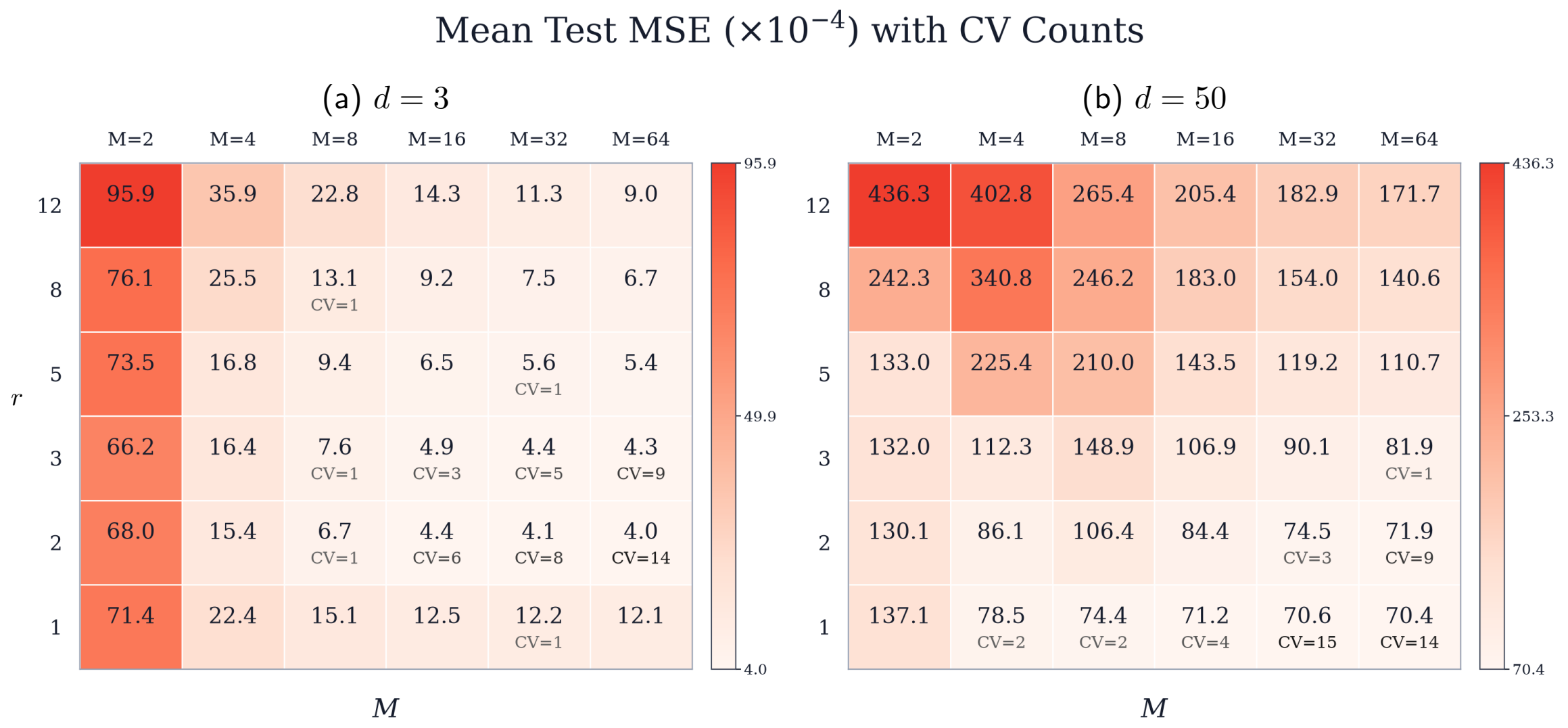

    \centering
    \simulationfigure{cv}
    \caption{Mean test MSE and CV selection counts over the grid of $(r,M)$ for (a) $d=3$ and (b) $d=50$, with $K^*=8$ and $n_{\mathrm{tot}}=20000$, averaged over 50 simulation replications. The MSE values are reported in units of $10^{-4}$, with darker shading indicating larger values. A label $\mathrm{CV}=t$ indicates that the corresponding parameter pair was selected by five-fold CV in $t$ of the 50 replications; cells without a CV label were never selected.}
    \label{fig:cv}
\end{figure}

Figure~\ref{fig:cv} shows the average test MSE for each pair $(M,r)$ and the frequency with which each pair is selected by CV across $50$ replications. The concentration of CV selections in regions of the grid with relatively low test MSE suggests CV as a reasonable approach for selecting both the sieve radius and number of particles in practice.

When $d=3$, $47$ of the $50$ selected pairs have $r\in\{2,3\}$, close to the generating radius $r^*=3$. When $d=50$, CV selects $r=1$ in $37$ replications, $r=2$ in $12$, and $r=3$ in one, indicating a preference for smaller sieve radii in higher dimensions. This preference may reflect a regularization effect and is consistent with the tradeoff in \cref{thm:mse_rate_random}: increasing $r$ reduces the approximation error but increases the estimation error, whose dependence on $r$ becomes stronger as $d$ increases. Thus, although the true radius is $r^*=3$, a smaller sieve radius can yield a lower test MSE in higher dimensions. %

Turning to the number of particles, increasing $M$ generally reduces the test MSE, although the improvement becomes smaller as $M$ keeps increasing. Notably, CV selects $M>K^*$ in $47$ of the $50$ replications when $d=3$ and in $46$ when $d=50$, suggesting that overspecifying the number of particles can be beneficial in both settings, in line with our theory in \cref{sec:computation_npmle}.
\subsection{Comparison of other methods for computing the NPMLE}\label{sec_sim_npmle}

Computing the NPMLE becomes increasingly challenging as the dimension $d$ grows. Common approaches include fixed-grid methods and conditional-gradient methods. The former discretize the parameter space, but the number of grid points required for a sufficiently fine discretization grows rapidly with $d$, making these methods computationally expensive in moderate and high dimensions. We therefore focus on the latter which optimize support locations directly. We compare RPG with two such alternatives as the dimension $d$ and covariate scale $\sigma_X$ vary.

The first comparator is an equally weighted particle-gradient method (EWPG), which only updates the particle locations with all mixing weights fixed at $1/M$. Since EWPG otherwise follows the same particle-location updates as RPG, this comparison isolates the effect of the weight optimization step in \eqref{eq:exact_reweight}. The second comparison method is the fully corrective Frank--Wolfe (FCFW) algorithm \citep{jaggi2013revisiting,lesperance2007fast,jagabathula2020conditional}, which searches for a new support point at each iteration and then reoptimizes the mixing weights over all accumulated support points. Finding each new support point requires solving a linear maximization oracle (LMO) subproblem, which itself  is  computationally challenging  due to nonconvexity. We approximately solve the LMO subproblem using projected gradient ascent with five random restarts, a step size of $1$, and at most $200$ iterations per restart, retaining the solution with the largest objective value. RPG and EWPG each use $M=64$ particles, while FCFW is allowed up to $64$ support points for a comparable specification across methods.

\begin{figure}[ht]
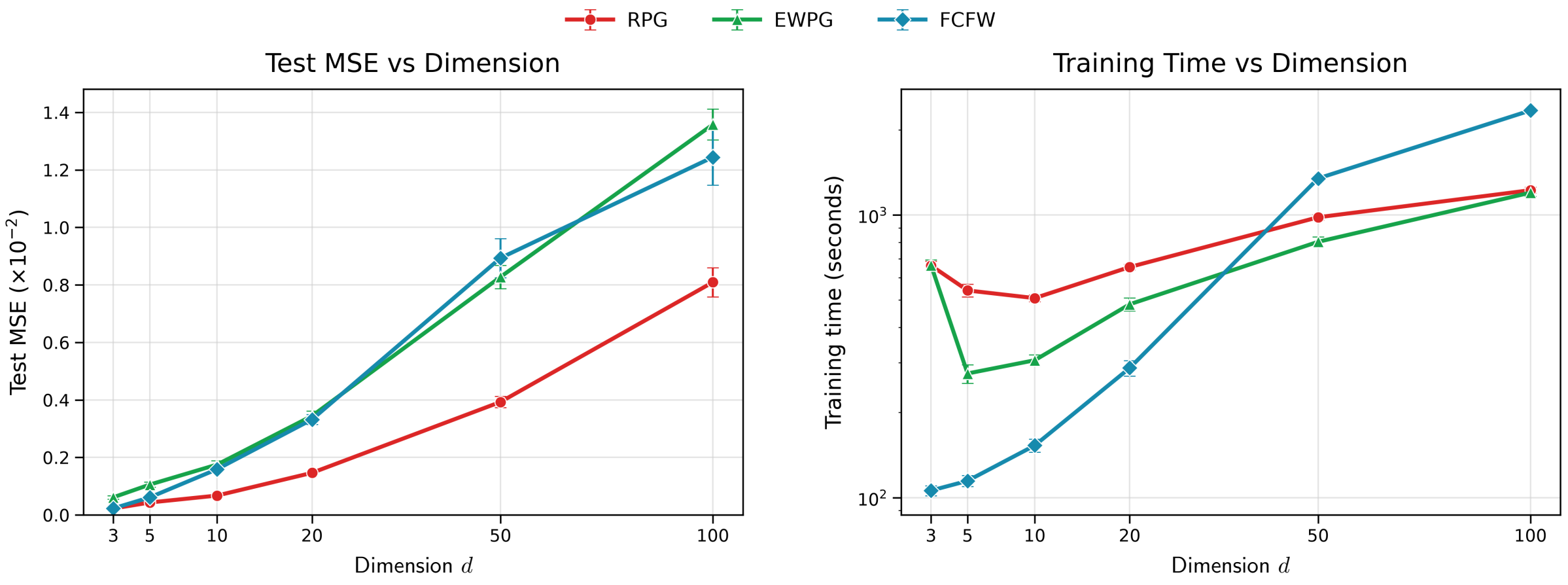

    \centering
    \simulationfigure{fcfw}
    \caption{MSE and computational cost as functions of the dimension $d$.}
    \label{fig:fcfw}
\end{figure}

Figure~\ref{fig:fcfw} compares the test MSE and training time of the three methods as $d$ varies, with both metrics averaged over $50$ replications. RPG achieves lower mean test MSE than both EWPG and FCFW in higher dimensions. Its advantage over EWPG reflects the benefit of additionally optimizing the mixing weights. At the largest dimension considered, RPG and EWPG have similar training times despite the additional weight optimization in RPG. One possible explanation is that weight optimization directly adjusts the contribution of each particle, reducing the number of location updates needed to satisfy the stopping criterion.
FCFW is faster in lower dimensions, but its training time
increases substantially with $d$. This may be due in part to
the nonconvex LMO subproblem, which requires searching over the
$d$-dimensional parameter space at each iteration. Moreover,
once a support point is added, FCFW adjusts its weight but
leaves its location fixed. In contrast, RPG maintains a fixed
collection of particles whose locations and weights are updated
throughout the optimization. These differences may contribute
to RPG's lower estimation error and shorter training time in
higher dimensions.

\begin{figure}[ht]
    \centering
    \includegraphics[width=1\linewidth]{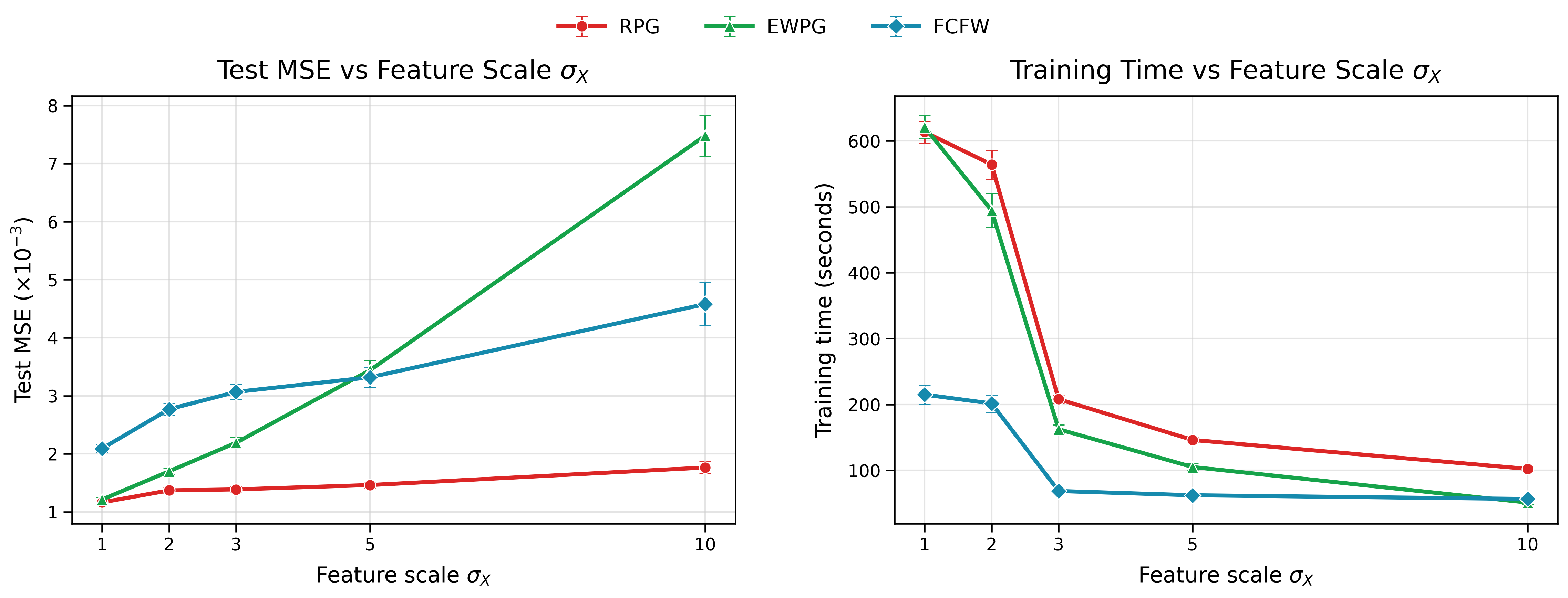}
    \caption{MSE and computational cost as functions of covariate scale $\sigma_X$.}
    \label{fig:vary_sigma}
\end{figure}

Figure~\ref{fig:vary_sigma} compares the three methods as the covariate scale $\sigma_X$ varies, with $d=20$ and the remaining parameters fixed at their baseline values. The mean test MSE of RPG increases only modestly with $\sigma_X$, whereas those of EWPG and FCFW increase more substantially. For EWPG, this increasing gap may be due to its fixed particle weights, which limit its ability to adjust the fitted mixture through location updates alone. For FCFW, larger $\sigma_X$ leads to more pronounced saturation of the logistic function, which may make the gradient-based LMO less effective. Overall, RPG is less sensitive to increasing covariate scale than the two alternative methods. On the computational side, the mean training time of RPG decreases from approximately $614$ seconds at $\sigma_X=1$ to $102$ seconds at $\sigma_X=10$. One possible explanation is that stronger saturation leads to smaller changes in the likelihood during particle updates, causing the inner-loop stopping criterion to be reached earlier. FCFW remains faster throughout, but has consistently higher mean test MSE.

Overall, these comparisons show that combining particle location updates with weight optimization improves estimation accuracy while remaining computationally competitive, particularly in higher dimensions and at larger covariate scales.

\subsection{Comparison with broader benchmarks}\label{sec_sim_benchmarks}
We now compare RPG with the expectation-maximization (EM) algorithm for fitting finite logistic mixtures and three general-purpose machine learning methods. The former examines the effect of specifying the number of mixture components, while the latter assesses the benefit of explicitly modeling the mixture structure when estimating the regression function.

\paragraph{Comparison with the EM algorithm}
We compare RPG with the EM algorithm for fitting finite mixtures of logistic regressions \citep[e.g.,][]{grun2007fitting}. At each EM iteration, the E-step computes the posterior component probabilities, while the M-step updates the mixing weights and fits a weighted logistic regression for each component. To enforce the same parameter constraint as RPG, the component parameters are restricted to $\cB_r$ and updated using projected Newton steps with backtracking line search. Each component update uses at most $50$ Newton steps and terminates early if the absolute change in its weighted average log-likelihood falls below $10^{-8}$ for five consecutive steps. Each EM run uses at most $25{,}000$ iterations and terminates if the relative change in the average training log-likelihood falls below $10^{-7}$ for ten consecutive iterations. Any iteration that decreases the average training log-likelihood by more than $10^{-10}$ is rejected, at which point the run is terminated.
Since EM is sensitive to initialization, we use five random initializations for each specification, with component parameters drawn independently and uniformly from $\cB_r$ and the mixing weights initialized uniformly. We retain the solution with the highest training log-likelihood, while including all five runs in the reported computation time.
We consider two benchmarks, Oracle-$K^*$ EM that uses the true  $K^*$, and Fixed-$K$ EM that uses $K=6$ throughout.

\begin{figure}[ht]
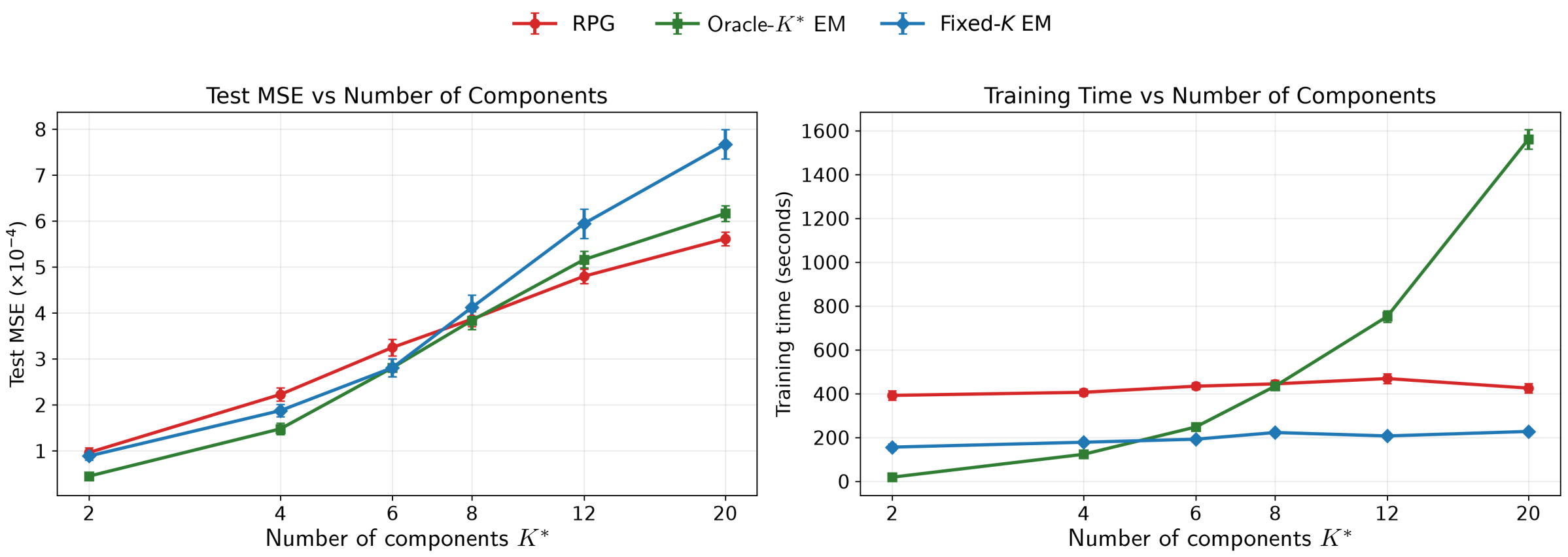

    \centering
    \simulationfigure{em}
    \caption{MSE and computational cost as functions of  $K^*$.  }
    \label{fig:em}
\end{figure}

We vary $K^*\in\{2,4,6,8,12,20\}$ and report the average test MSEs and computation times over $50$ replications. Figure~\ref{fig:em} shows that RPG performs comparably to Oracle-$K^*$ EM despite not knowing $K^*$. Oracle-$K^*$ EM has slightly lower mean test MSE for smaller $K^*$, while RPG performs   slightly better as $K^*$ increases. When $K^*\leq6$, Fixed-$K$ EM performs comparably to Oracle-$K^*$ EM, suggesting that overspecifying the number of components has little adverse effect on MSE. In contrast, when $K^*>6$, Fixed-$K$ EM underspecifies the number of mixture components and incurs increasingly larger MSE.

The right panel shows that the training time of Oracle-$K^*$ EM increases sharply with $K^*$, while that of RPG remains relatively stable. Together, the MSE and timing results highlight a tradeoff in selecting the number of components for EM: using more components can protect against underfitting but increases computational cost, whereas using fewer components reduces computation at the risk of larger estimation errors. The sharp increase in the training time of Oracle-$K^*$ EM further suggests that choosing a conservatively large number of components to accommodate an unknown $K^*$ could be even more computationally demanding. By contrast, RPG does not require knowledge of $K^*$ and, for the larger mixtures considered here, achieves lower MSE in less training time than Oracle-$K^*$ EM.

\paragraph{Comparison with other machine learning approaches} We also compare RPG with random forests (RF), support vector machines with a radial basis function kernel (RBF-SVM), and extreme gradient boosting (XGBoost). Unlike RPG, these methods estimate the regression function $p^*$  without modeling the mixing distribution $G^*$. RF uses $500$ trees, and XGBoost uses $400$ boosting rounds with a maximum tree depth of $2$ and a learning rate of $0.1$. RBF-SVM uses standardized covariates, with regularization parameter $C=0.1$ and kernel parameter $\gamma=1/d$. %
We vary $d$, $K^*$, and $n_{\mathrm{tot}}$ separately, with the remaining parameters fixed at their baseline values.

\begin{figure}[ht]
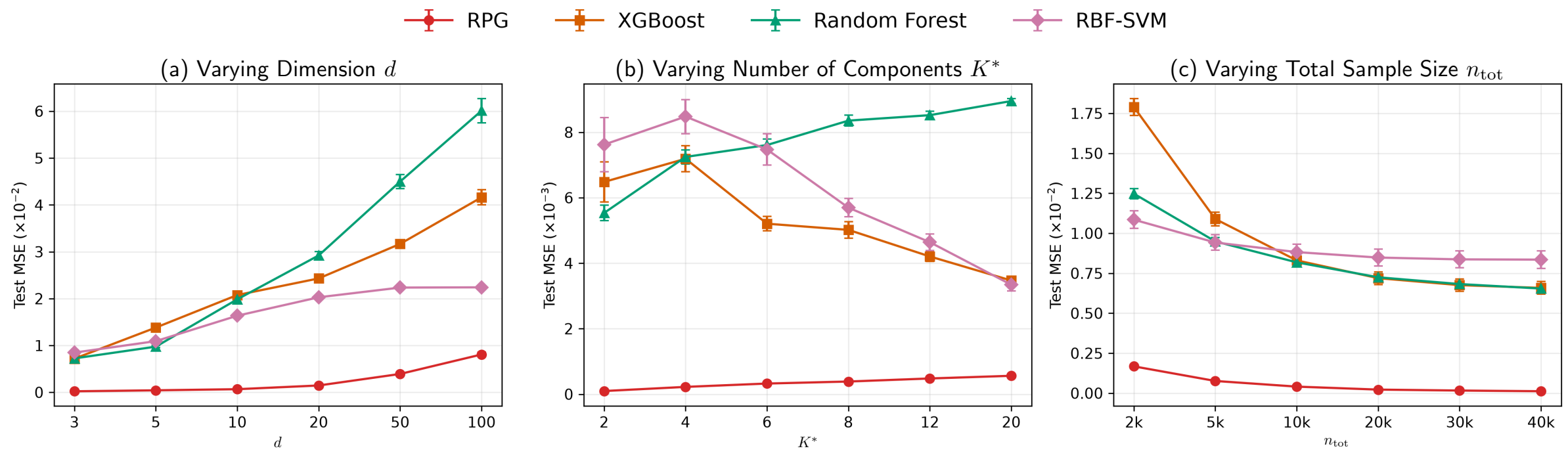

    \centering
    \simulationfigure{ml}
    \caption{Mean test MSE of RPG, XGBoost, RF, and RBF-SVM. Panel (a) uses $n_{\mathrm{tot}}=20000$ and $K^*=4$. Panel (b) uses $n_{\mathrm{tot}}=20000$ and $d=3$. Panel (c) uses   $d=3$ and $K^*=4$.}
    \label{fig:ml}
\end{figure}

Figure~\ref{fig:ml} depicts the average test MSEs over $50$ replications. RPG achieves the lowest mean test MSE across all settings considered. In panel~(a), the performance gap widens as $d$ increases, particularly relative to RF and XGBoost. The increase in RPG's MSE with $d$ is qualitatively consistent with  \cref{thm:particle_uniform_empirical}. In panel~(b), RPG's MSE increases only modestly with $K^*$ and remains substantially lower than those of the three benchmark methods. In panel~(c), all methods improve as $n_{\mathrm{tot}}$ increases, while RPG consistently achieves the lowest MSE, consistent with the decay of the estimation error in \cref{thm:particle_uniform_empirical}. Overall, these results highlight the benefit of explicitly modeling the underlying mixture structure under the specified mixture model.

\setlength{\bibsep}{5pt}
\bibliographystyle{plainnat}
\bibliography{refs}

\newpage
\appendix

The proof of the entropy bound in \cref{thm:metric_entropy} is provided in \cref{app:proof_metric_entropy}. The proofs of the estimation error bounds for fixed and random designs are given in \cref{app:proof_mse_rate_fixed} and \cref{app:proof_mse_rate_random}, respectively. \cref{app:particle_fixed_data_analysis} establishes the convergence of the RPG algorithm for a fixed dataset. Identifiability and the denoising bound are presented in \cref{app_sec_ident}. The convergence and statistical error bounds for RPG in finite logistic mixtures are established in \cref{app:particle_uniform_analysis} and \cref{app:proof_thm:particle_uniform_empirical} for the population- and sample-level results, respectively. Finally, the proof of the convergence lower bound in \cref{prop:particle_inverse_square_sharp} is provided in \cref{app:particle_lower_bound_construction}.

\section{Proofs for the entropy bounds}
\label{app:proof_metric_entropy}

We prove Theorem~\ref{thm:metric_entropy} by first approximating the mixture function class $\cF$ by a polynomial, and then covering the resulting
truncated mixture class inside a finite-dimensional polynomial space. The former established in \cref{app:sec_approx_F} builds on analytic approximation of the sigmoid function by Chebyshev
polynomials in \cref{app:general_aux_analytic}.  The latter resorts to a standard volumetric covering bound for bounded
subsets of finite-dimensional normed spaces. Before we proceed, we first recall
some standard notions from complex analysis and Chebyshev approximation theory.

\subsection{Preliminaries on complex analysis and Chebyshev approximation}
\label{app:complex_chebyshev_prelim}

In this subsection we briefly recall several standard notions from complex
analysis and Chebyshev approximation that will be used in the sequel.

A complex number $z\in\mathbb C$ is written as
$$
z=x+iy,
\qquad x,y\in\mathbb R,
$$
where $x = \Re z$ is the real part of $z$, $y=\Im z$ is the imaginary part of $z$, and
$i=\sqrt{-1}$.

A complex-valued function is said to be \emph{holomorphic} on an
open set if it is complex differentiable at every point of that set.  It is
said to be \emph{meromorphic} if it is holomorphic except at isolated poles.   For $a>0$, denote the horizontal strip of half-width $a$ by
$$
S_a:=\{x+iy\in\mathbb C:\ |y|<a\}.
$$

Our goal is to approximate a function $g$ on $[-1,1]$ by polynomials.
When $g$ extends holomorphically to a complex neighborhood of the interval,
the rate of approximation can be controlled by the size of this neighborhood.
We describe the relevant neighborhoods by Bernstein ellipses, defined using
the standard Joukowski map \citep[Chapter~8]{trefethen2019approximation}
$$
J(w):=\frac12\left(w+w^{-1}\right),\qquad\forall ~  w\in \CC\setminus\{0\}.
$$
For a parameter $\rho>1$, the image of the circle $\{w\in\mathbb C:|w|=\rho\}$
under $J$ is the \emph{Bernstein ellipse} with parameter $\rho$, which we
denote by
$$
E_\rho:=J(\{w\in\mathbb C:|w|=\rho\}).
$$
Its closed interior will be denoted by $\overline E_\rho$. This ellipse has
foci at $\pm1$ and semi-major and semi-minor axes
$$
a_\rho=\frac{\rho+\rho^{-1}}{2},
\qquad
b_\rho=\frac{\rho-\rho^{-1}}{2}.
$$
We use Chebyshev polynomials to construct the approximation.
For $j\ge 0$, let $T_j$ denote the $j$th Chebyshev polynomial of the first kind,
defined by
$$
T_j(x):=\cos(j\arccos x),
\qquad x\in[-1,1].
$$
We approximate $g$ by truncating its expansion in these polynomials.
The following standard result gives this expansion and bounds its coefficients
when $g$ is holomorphic on a neighborhood of $\overline E_\rho$; see
\citet[Theorems~8.1 and~8.2]{trefethen2019approximation}.

\begin{lemma}
\label{lem:chebyshev_projection_analytic}
Let $\rho>1$, and let $g$ be holomorphic on an open neighborhood of
$\overline E_\rho$. Then $g$ admits a Chebyshev expansion
$$
g(x)=\frac{\gamma_0}{2}+\sum_{j=1}^\infty \gamma_j T_j(x),
\qquad x\in[-1,1],
$$
where the coefficients satisfy
$$
|\gamma_j|\le 2M_\rho\,\rho^{-j},
\qquad j\ge 1
$$
with $M_\rho:=\sup_{z\in\overline E_\rho}|g(z)|$.
\end{lemma}
Consequently, for any integer $k\ge 0$, the truncated Chebyshev expansion
\begin{equation}\label{def_g_k}
g_k(x):=\frac{\gamma_0}{2}+\sum_{j=1}^k \gamma_j T_j(x)
\end{equation}
satisfies
$$
\sup_{x\in[-1,1]} |g(x)-g_k(x)|
\le
\sum_{j=k+1}^\infty |\gamma_j|
\le
2M_\rho\sum_{j=k+1}^\infty \rho^{-j}
=
\frac{2M_\rho}{\rho-1}\rho^{-k}.
$$

\subsection{Polynomial approximation to the sigmoid function over  compact intervals}
\label{app:general_aux_analytic}

Fix $R>0$. We aim to apply the Chebyshev approximation result in \cref{lem:chebyshev_projection_analytic} to the sigmoid function $\sigma(z) = 1/(1+e^{-z})$ for $z\in [-R,R]$. To this end, we rescale the
sigmoid function by setting
\begin{equation}\label{def_scaled_sigmoid}
    g(t):=\sigma(Rt),\qquad \forall ~ t\in [-1,1].
\end{equation}
To find an open neighborhood of a Bernstein ellipse on which $g$ is holomorphic, we first identify the strip of $\CC$ on which the sigmoid function $\sigma$ is
holomorphic.

\begin{lemma}
\label{lem:logistic_meromorphic}
The sigmoid function $\sigma:\CC\to \CC$
is meromorphic on $\mathbb C$, with poles at
\begin{equation}\label{zeroes}
z=(2m+1)\pi i,\qquad m\in\mathbb Z.
\end{equation}
In particular, $\sigma$ is holomorphic on the horizontal strip $S_\pi$.
\end{lemma}

\begin{proof}
The reciprocal of the sigmoid function is
$1+e^{-z}$,
which is entire. Therefore the poles of $\sigma(\cdot)$ are exactly the zeros of
$1+e^{-z}$, which  satisfy
$$
e^{-z}=-1=e^{(2m+1)\pi i},
\qquad m\in\mathbb Z.
$$
Solving the above gives the set of zeros as in \eqref{zeroes}.
Since none of these poles lies in $S_\pi$, the function $\sigma$ is
holomorphic on $S_\pi$.
\end{proof}

   We proceed to identify a Bernstein
ellipse such that $g$ in \eqref{def_scaled_sigmoid} is holomorphic on its open neighborhood, and to control
the magnitude of $g$ on that ellipse.

\begin{lemma}
\label{lem:logistic_chebyshev}
Let $R>0$, and define
$$
\rho_R
:=
\frac{\pi}{2R}+\sqrt{1+\frac{\pi^2}{4R^2}}.
$$
For any integer $k\ge 0$, there exists a polynomial $P_k$ of degree at most $k$ such that
$$
\sup_{z\in[-R,R]} |\sigma(z)-P_k(z)|
\le
\frac{4R}{\pi}\rho_R^{-k}.
$$
\end{lemma}

\begin{proof}
For any $t\in[-1,1]$, recall from \eqref{def_scaled_sigmoid} that
$
g(t)=\sigma(Rt).
$
With $g_k$ as in \eqref{def_g_k}, we only need to prove
\begin{equation}\label{appro_unit}
\sup_{t\in[-1,1]} |g(t)-g_k(t)|
\le
\frac{4R}{\pi}\rho_R^{-k}
\end{equation}
as the claim follows by setting
$$
P_k(z):=g_k(z/R),\qquad z\in[-R,R]
$$
and the change of variables $t=z/R$.

To prove \eqref{appro_unit}, we first verify that $g$ is holomorphic on an open neighborhood of
$\overline E_{\rho_R}$. By the definition of $\rho_R$,
$$
\frac{\rho_R-\rho_R^{-1}}{2}=\frac{\pi}{2R},
$$
so the semi-minor axis of $\overline E_{\rho_R}$ is $\pi/(2R)$. Hence every
point $t=x+iy\in\overline E_{\rho_R}$ satisfies
$
|y|\le  \pi / (2R),
$
and therefore
$
|\Im(Rt)|\le \pi / 2<\pi.
$
Thus
$$
R\overline E_{\rho_R} := \left\{R t: t\in \overline E_{\rho_R}\right\}\subset S_\pi.
$$
Since $R\overline E_{\rho_R}$ is compact and $S_\pi$ is open, there exists
an open neighborhood $U$ of $R\overline E_{\rho_R}$ such that
$U\subset S_\pi$.
By Lemma~\ref{lem:logistic_meromorphic}, $\sigma$ is holomorphic on $S_\pi$,
hence in particular on $U$. Therefore $g(t)=\sigma(Rt)$ is holomorphic on
an open neighborhood of $\overline E_{\rho_R}$.

We next show that $|g(t)|\le 1$ for all $t\in\overline E_{\rho_R}$.
Let $t=x+iy\in\overline E_{\rho_R}$. As above, $|Ry|\le \pi/2$ implies
$
\cos(Ry)\ge 0,
$
so that
\begin{align*}
|1+e^{-Rt}|^2
&=
\bigl(1+e^{-R(x+iy)}\bigr)\bigl(1+e^{-R(x-iy)}\bigr) \\
&=
1+2e^{-Rx}\cos(Ry)+e^{-2Rx}\ge 1.
\end{align*}
Therefore
$$
|g(t)|=\frac{1}{|1+e^{-Rt}|}\le 1.
$$

Applying Lemma~\ref{lem:chebyshev_projection_analytic} with $\rho=\rho_R$
yields
$$
\sup_{t\in[-1,1]} |g(t)-g_k(t)|
\le
\frac{2}{\rho_R-1}\rho_R^{-k}.
$$
Finally, from
$$
\rho_R=\frac{\pi}{2R}+\sqrt{1+\frac{\pi^2}{4R^2}}
>
1+\frac{\pi}{2R},
$$
we obtain
$$
\frac{2}{\rho_R-1}<\frac{4R}{\pi}.
$$
This proves \eqref{appro_unit} and hence completes the proof.
\end{proof}

The preceding lemma gives the approximation error as a function of the degree
$k$. For a
prescribed approximation accuracy, the following corollary states an explicit degree choice.

\begin{corollary}
\label{cor:degree_choice_logistic}
Let $R>0$ and $\epsilon\in(0,1]$. Define
$$
k_*(R,\epsilon)
:=
\begin{cases}
0, & \text{if } \dfrac{8R}{\pi\epsilon}\le 1, \\[1.2ex]
\left\lceil \dfrac{2R+\pi}{\pi}\log\!\left(\dfrac{8R}{\pi\epsilon}\right)\right\rceil,
& \text{if } \dfrac{8R}{\pi\epsilon}>1.
\end{cases}
$$
Then there exists a polynomial $P_k$ of degree at most
$k_*(R,\epsilon)$ such that
$$
\sup_{z\in[-R,R]} |\sigma(z)-P_k(z)|\le \frac{\epsilon}{2}.
$$
\end{corollary}

\begin{proof}
If $8R/(\pi\epsilon)\le 1$, then Lemma~\ref{lem:logistic_chebyshev} with
$k=0$ gives
$$
\sup_{z\in[-R,R]} |\sigma(z)-P_0(z)|
\le
\frac{4R}{\pi}
\le
\frac{\epsilon}{2}.
$$
Now suppose $8R/(\pi\epsilon)>1$. By Lemma~\ref{lem:logistic_chebyshev}, it is
enough to require
$$
\frac{4R}{\pi}\rho_R^{-k}\le \frac{\epsilon}{2},
$$
or equivalently
$$
k\ge \frac{1}{\log\rho_R}\log\!\left(\frac{8R}{\pi\epsilon}\right).
$$
Since
$
\rho_R>1+ \pi/(2R),
$
the elementary inequality $\log(1+x)\ge x/(1+x)$ for $x>0$ yields
$$
\log\rho_R
>
\log\!\left(1+\frac{\pi}{2R}\right)
\ge
\frac{\pi}{2R+\pi}.
$$
Hence
$$
\frac{1}{\log\rho_R}\le \frac{2R+\pi}{\pi},
$$
so the displayed value of $k_*(R,\epsilon)$ is sufficient.
\end{proof}

To control the dimension of the polynomial space in the proof of
\cref{thm:metric_entropy}, we use the degree bound
\begin{equation}\label{degree_choice_upper_bound}
k_*(R,\epsilon)
\le
C(1\vee R)\log\!\left(\frac{C(1\vee R)}{\epsilon}\right)
\qquad
\text{for all }R>0,\ \epsilon\in(0,1],
\end{equation}
where $C>0$ is a universal constant.

\subsection{Polynomial approximation to the logistic mixture class}\label{app:sec_approx_F}

We use the polynomial approximation in \cref{cor:degree_choice_logistic}
to approximate every function in $\cF_r$ within a common finite-dimensional
space.

Fix $r>0$ and a nonempty bounded set $\cX\subset\RR^d$.
Recall that $R_x=\sup_{x\in\cX}\|x\|_2$, and set $R=rR_x$.
Let $\mathbb N_0:=\{0,1,2,\ldots\}$.
For a multi-index $\alpha=(\alpha_1,\ldots,\alpha_d)\in\mathbb N_0^d$, write
$$
|\alpha|:=\alpha_1+\cdots+\alpha_d,\qquad
\alpha!:=\alpha_1!\cdots\alpha_d!,\qquad
x^\alpha:=x_1^{\alpha_1}\cdots x_d^{\alpha_d}.
$$
For an integer $k\ge0$, let $V_k$ be the space of functions on $\cX$ of the form
$$
x\mapsto\sum_{|\alpha|\le k}c_\alpha x^\alpha,\qquad c_\alpha\in\RR.
$$
Polynomials that agree on $\cX$ define the same element of $V_k$, so
$\|\cdot\|_{\infty,\cX}$ is a norm on this space. Since the monomials
$x^\alpha$ with $|\alpha|\le k$ span $V_k$,
$$
\dim(V_k)\le
\#\{\alpha\in\mathbb N_0^d:|\alpha|\le k\}
=\binom{d+k}{d}.
$$
The following proposition gives the polynomial approximation to $\cF_r$
and the covering bound for the resulting class.

\begin{proposition}
\label{prop:finite_dim_approx_class}
Fix $\epsilon\in(0,1]$. Suppose $P_k$ is a polynomial of degree at most $k$
such that
$$
\sup_{z\in[-R,R]}|\sigma(z)-P_k(z)|\le\frac{\epsilon}{2}.
$$
For each $G\in\cP(\cB_r)$, define
$$
\wt p_{G,k}(x):=\int_{\cB_r}P_k(x^\T\beta)\,\rd G(\beta),
\qquad x\in\cX,
$$
and let $\wt\cF_{r,k}:=\{\wt p_{G,k}:G\in\cP(\cB_r)\}$.
Then $\wt\cF_{r,k}\subset V_k$, and
$$
\sup_{G\in\cP(\cB_r)}\|p_G-\wt p_{G,k}\|_{\infty,\cX}
\le\frac{\epsilon}{2},\qquad
\sup_{f\in\wt\cF_{r,k}}\|f\|_{\infty,\cX}
\le1+\frac{\epsilon}{2}.
$$
Moreover, for every $u>0$,
$$
\log\cN(u,\wt\cF_{r,k},\|\cdot\|_{\infty,\cX})
\le
\binom{d+k}{d}\log\!\left(1+\frac{4(1+\epsilon/2)}{u}\right).
$$
\end{proposition}

\begin{proof}
Write $P_k(z)=\sum_{j=0}^k a_jz^j$ with $a_j\in\RR$.
By the multinomial theorem, for every $\beta\in\cB_r$,
$$
P_k(x^\T\beta)
=\sum_{j=0}^k a_j(x^\T\beta)^j
=\sum_{j=0}^k\sum_{|\alpha|=j}
a_j\frac{j!}{\alpha!}\beta^\alpha x^\alpha.
$$
Since $|\beta^\alpha|\le r^{|\alpha|}$ for $\beta\in\cB_r$,
integrating the finite sum against $G$ gives
$$
\wt p_{G,k}(x)
=\sum_{|\alpha|\le k}
a_{|\alpha|}\frac{|\alpha|!}{\alpha!}
\left(\int_{\cB_r}\beta^\alpha\,\rd G(\beta)\right)x^\alpha.
$$
Thus $\wt p_{G,k}\in V_k$ for every $G\in\cP(\cB_r)$.

For any $x\in\cX$ and $\beta\in\cB_r$,
$|x^\T\beta|\le\|x\|_2\|\beta\|_2\le rR_x=R$.
By the choice of $P_k$,
\begin{align*}
|p_G(x)-\wt p_{G,k}(x)|
&=\left|\int_{\cB_r}\{\sigma(x^\T\beta)-P_k(x^\T\beta)\}\,\rd G(\beta)\right|\\
&\le\int_{\cB_r}|\sigma(x^\T\beta)-P_k(x^\T\beta)|\,\rd G(\beta)
\le\frac{\epsilon}{2}.
\end{align*}
Taking the supremum over $x\in\cX$ and $G\in\cP(\cB_r)$ proves the
approximation bound. Since $0\le p_G(x)\le1$, we also have
$$
\|\wt p_{G,k}\|_{\infty,\cX}
\le\|p_G\|_{\infty,\cX}+\|p_G-\wt p_{G,k}\|_{\infty,\cX}
\le1+\frac{\epsilon}{2}.
$$
Fix $u>0$. By \citet[Example~5.8]{wainwright2019high}, there exist
$q_1,\ldots,q_M\in V_k$ such that
$$
M\le\left\{1+\frac{4(1+\epsilon/2)}{u}\right\}^{\dim(V_k)},
\qquad
\sup_{f\in\wt\cF_{r,k}}\min_{j\in[M]}
\|f-q_j\|_{\infty,\cX}\le\frac u2.
$$
Our definition of the covering number requires the centers to belong to
$\wt\cF_{r,k}$. For each $j$, choose
$f_j\in\wt\cF_{r,k}$ with $\|f_j-q_j\|_{\infty,\cX}\le u/2$ if such a
function exists, and discard $q_j$ otherwise. For any
$f\in\wt\cF_{r,k}$, there exists $j\in[M]$ such that
$\|f-q_j\|_{\infty,\cX}\le u/2$. For this $j$, $f_j$ is defined, and
$$
\|f-f_j\|_{\infty,\cX}
\le\|f-q_j\|_{\infty,\cX}+\|q_j-f_j\|_{\infty,\cX}
\le\frac u2+\frac u2=u.
$$
Thus the chosen functions $f_j$ form a $u$-cover of $\wt\cF_{r,k}$
with at most $M$ centers, and hence
$$
\begin{aligned}
\log\cN(u,\wt\cF_{r,k},\|\cdot\|_{\infty,\cX})
&\le\dim(V_k)\log\!\left(1+\frac{4(1+\epsilon/2)}{u}\right)\\
&\le\binom{d+k}{d}\log\!\left(1+\frac{4(1+\epsilon/2)}{u}\right).
\end{aligned}
$$
This completes the proof.
\end{proof}

\subsection{Proof of \cref{thm:metric_entropy}}

\begin{proof}
Fix $\epsilon\in(0,1]$. If $R_x=0$, all functions in $\cF_r$ equal $1/2$
on $\cX$, so the covering number is one and the claim follows.
Thus assume $R_x>0$ and write $R=rR_x$.
By \cref{cor:degree_choice_logistic}, taking $k=k_*(R,\epsilon/2)$ gives
a polynomial $P_k$ such that
$$
\sup_{z\in[-R,R]}|\sigma(z)-P_k(z)|\le\frac{\epsilon}{4}.
$$
We use the notation in \cref{prop:finite_dim_approx_class}.
Applying this proposition with $\epsilon/2$ in place
of $\epsilon$, we can choose an $(\epsilon/2)$-cover
$\{\wt p_{G_j,k}:j\in[M]\}$ of $\wt\cF_{r,k}$ under
$\|\cdot\|_{\infty,\cX}$, with $G_j\in\cP(\cB_r)$, such that
$$
\log M
\le\binom{d+k}{d}\log\!\left(1+\frac{8(1+\epsilon/4)}{\epsilon}\right)
\le\binom{d+k}{d}\log\!\left(\frac{13}{\epsilon}\right).
$$

We next show that $\{p_{G_j}:j\in[M]\}\subset\cF_r$ is an
$\epsilon$-cover of $\cF_r$. For any $G\in\cP(\cB_r)$, there exists
$j\in[M]$ such that
$\|\wt p_{G,k}-\wt p_{G_j,k}\|_{\infty,\cX}\le\epsilon/2$.
The approximation bound in \cref{prop:finite_dim_approx_class} then gives
\begin{align*}
\|p_G-p_{G_j}\|_{\infty,\cX}
&\le\|p_G-\wt p_{G,k}\|_{\infty,\cX}
  +\|\wt p_{G,k}-\wt p_{G_j,k}\|_{\infty,\cX}
  +\|\wt p_{G_j,k}-p_{G_j}\|_{\infty,\cX}\\
&\le\frac{\epsilon}{4}+\frac{\epsilon}{2}+\frac{\epsilon}{4}
=\epsilon.
\end{align*}
Hence $\log\cN(\epsilon,\cF_r,\|\cdot\|_{\infty,\cX})\le\log M$.

Finally, by \eqref{degree_choice_upper_bound}, there exists a universal
constant $C>0$ such that
$$
k
\le
C(1\vee R)\log\!\left(\frac{C(1\vee R)}{\epsilon}\right).
$$
Increasing $C$ if necessary gives
$$
\binom{d+k}{d}
\le (k+1)^d \le
\bigg[
C(1\vee R)\log\!\Bigl(\frac{C(1\vee R)}{\epsilon}\Bigr)
\bigg]^d,
$$
and therefore
$$
\log\cN(\epsilon,\cF_r,\|\cdot\|_{\infty,\cX})
\le
\bigg[
C(1\vee R)\log\!\Bigl(\frac{C(1\vee R)}{\epsilon}\Bigr)
\bigg]^d
\log\!\Bigl(\frac{C}{\epsilon}\Bigr).
$$
This proves the theorem.
\end{proof}

\begin{remark}[Existing approximation arguments under GMMs]\label{rem_approx_GMM}
    We recall the approximation argument for Gaussian location mixtures
\citep[e.g.][]{zhang2009generalized,saha2020nonparametric}, where the goal is to approximate the following mixture density uniformly over a fixed compact set $S\subset\RR^d$:
$$
f_G(x) =
\int \phi(x-\theta)\,\rd G(\theta) =
\int_{\theta:\|x-\theta\|_2\le a}\phi(x-\theta)\,\rd G(\theta) + \int_{\theta:\|x-\theta\|_2>a}\phi(x-\theta)\,\rd G(\theta),
$$
where $\phi(z)=(2\pi)^{-d/2}e^{-\|z\|_2^2/2}$ and $a\ge 1$ is arbitrary. The second integral on the far right is bounded by $(2\pi)^{-d/2}e^{-a^2/2}$, which is of order $\epsilon$ when $a$ is chosen appropriately with $a\asymp\sqrt{\log(1/\epsilon)}$ as $\epsilon\to0$. For the first integral, Taylor expansion of $e^{-u}$ around zero, with $u=\|z\|_2^2/2$, gives $\phi(z)=H_{2k}(z)+R_{2k}(z)$ for any $k\in\NN$, where
$$
    H_{2k}(z)=\frac{1}{(2\pi)^{d/2}}
    \sum_{j=0}^k(-1)^j\frac{\|z\|_2^{2j}}{j!\,2^j},
    \qquad
    \sup_{\|z\|_2\le3a}|R_{2k}(z)|
    \le
    \frac{1}{(2\pi)^{d/2}}
    \frac{(3a)^{2(k+1)}}{(k+1)!\,2^{k+1}}.
$$
Choosing $k\asymp\log(1/\epsilon)$ with a sufficiently large proportionality constant makes the remainder $|R_{2k}(z)|$ uniformly of order $\epsilon$. Combining moment matching on partition sets with Gaussian tail decay, the argument in \citet[Lemmas~D.2--D.3]{saha2020nonparametric} yields a discrete probability measure $G'=\sum_{\ell=1}^L w_\ell\delta_{\theta_\ell}$ such that
$$
    \sup_{x\in S}|f_G(x)-f_{G'}(x)|
    \lesssim
    \epsilon+
    \sup_{x\in S}
    \left|
    \int_{\theta:\|x-\theta\|_2\le a}
    H_{2k}(x-\theta)\,\rd (G-G')(\theta)
    \right|
    \lesssim\epsilon.
$$
Here, $L$ is controlled by the partition size and polynomial degree. Discretizing the weights and support locations then gives the desired sup-norm covering bound on $S$.

Clearly, the above argument crucially relies on the two-sided exponential decay of the Gaussian density as well as its globally convergent Taylor expansion.
\end{remark}

\subsection{Entropy bounds under exponential moment conditions}
\label{app:entropy_exponential_moments}

We extend \cref{thm:metric_entropy} to mixing distributions in
$\cG(a,q,b)$, defined in \eqref{def_exp_G}, by truncating their tails.
The following corollary shows that the entropy bound for the corresponding
regression function class remains polylogarithmic in $1/\eps$.

\begin{corollary}
\label{cor:metric_entropy_light_tail_beta}
For $a>0$, $q>0$, and
$b\ge1$, let
$
\cF_{a,q,b} := \{ p_G: G \in  \cG(a,q,b)  \}.
$
For any $\eps\in(0,1]$, set $\bar R_\eps:=1\vee(r_\eps R_x) $ where $r_\eps
:=
1\vee
\left\{
\frac1a\log\left(\frac{4b}{\eps}\right)
\right\}^{1/q}$. Then there exists a universal constant $C>0$ such that
$$
\log\cN
\left(
\eps,\cF_{a,q,b},\|\cdot\|_{\i,\cX}
\right)
\le
\left[
C\bar R_\eps
\log\left(\frac{C\bar R_\eps}{\eps}\right)
\right]^d
\log\left(\frac{C}{\eps}\right).
$$
\end{corollary}

\begin{proof}
For any $G$ such that $\int_{\RR^d}\exp(a\|\beta\|_2^q)\,\rd G(\beta)\le b$, Markov's inequality gives
$$
G(\cB_{r_\eps}^c)
\le
b\exp(-a r_\eps^q)
\le
\eps/4.
$$
Move the mass outside $\cB_{r_\eps}$ to the origin:
$$
G^{(\eps)}
:=
G|_{\cB_{r_\eps}}
+
G(\cB_{r_\eps}^c)\delta_0
$$
where $G|_{\cB{r_\eps}}$ denotes the restriction of $G$ to $\cB_{r_\eps}$.
Then $G^{(\eps)}\in\cP(\cB_{r_\eps})$, so
$p_{G^{(\eps)}}\in\cF_{r_\eps}$. Since $0\le\sigma(\cdot)\le1$,
$$
\|p_G-p_{G^{(\eps)}}\|_{\i,\cX}
\le
G(\cB_{r_\eps}^c)
\le
\eps/4.
$$
Thus an $\eps/2$-cover of $\cF_{r_\eps}$ is an $\eps$-cover of
$\cF_{a,q,b}$. The result follows by applying
\cref{thm:metric_entropy} to $\cF_{r_\eps}$.
\end{proof}

\section{Proofs for bounding the estimation error under fixed design}
\label{app:proof_mse_rate_fixed}

\subsection{Bernoulli likelihood and Hellinger distance}
\label{app:fixed_design_preliminaries}

In the fixed-design setting, recall that $Y_1,\ldots,Y_n$ are independent,
with $Y_i\mid x_i\sim\Bern(p^*(x_i))$ for $i\in[n]$. For any
$x\in\cX$ and $y\in \{0,1\}$, write the likelihood of $Y=y \mid x$ at any $p\in \cF$ as
\begin{equation}
\label{def_Lf}
    L(p; y, x) := p(x)^y(1-p(x))^{1-y}.
\end{equation}
For any $a,b\in(0,1)$, define the squared Hellinger distance between
$\Bern(a)$ and $\Bern(b)$ by
\begin{equation}
\label{eq_hsquare_xi}
    \begin{split}
        h^2(a,b)
        &:=
        {1\over 2}
        \left[
        \left(\sqrt a-\sqrt b\right)^2
        +
        \left(\sqrt{1-a}-\sqrt{1-b}\right)^2
        \right] \\
        &~ =
        1
        -
        \sqrt{ab}
        -
        \sqrt{(1-a)(1-b)}.
    \end{split}
\end{equation}
For any $p,p'\in\cF$, write
$h^2(p,p')(x):=h^2(p(x),p'(x))$ for $x\in\cX$, and define their
empirical squared Hellinger distance by
\begin{equation}
\label{eq_hsquare_emp}
    h_n^2(p,  p')
    := P_n h^2(p,p') =
    {1\over n}
    \sum_{i=1}^n
    h^2(p(x_i), p'(x_i)) .
\end{equation}

The following lemma bounds the squared difference between two success
probabilities in terms of their squared Hellinger distance.

\begin{lemma}
\label{lem:bernoulli_hellinger_x}
For any $p,q\in (0,1)$,
$$
    {1\over 2}\left(p-q\right)^2
    \le
    2\,h^2(p,q) .
$$
\end{lemma}

\begin{proof}
By \eqref{eq_hsquare_xi} and the basic inequality
$$
(\sqrt a-\sqrt b)^2
=
\frac{(a-b)^2}{(\sqrt a+\sqrt b)^2}
\ge
\frac{(a-b)^2}{4},
\qquad a,b\in(0,1),
$$
one has
\begin{align}
\label{basic_ineq_hellinger_L2}
2h^2(p,q)
&=
(\sqrt{p}-\sqrt{q})^2
+
(\sqrt{1-p}-\sqrt{1-q})^2  \ge
{1\over 2}(p-q)^2,
\end{align}
which is the stated inequality.
\end{proof}

\subsection{Proof of \cref{thm:mse_rate_fixed}}
\label{app:fixed_design_theorem_proof}

\begin{proof}
Fix $u\ge1$, and set
\begin{equation}
\label{eq:fixed_design_lr_radius}
\gamma_{n,u,\tau}^2
:=
A\left(\epsilon_{\est}(r, R_{x,n}) + \epsilon_{\app}(r, P_n)+\tau+\frac{u}{n}\right),
\end{equation}
where $A\ge 1$ is a sufficiently large constant depending only on $d$, to be chosen below.
If $\gamma_{n,u,\tau}^2\ge1$, then the claim is trivial after enlarging the final
constant, since
$
P_n (\wt p_r-p^*)^2\le1.
$
Thus we assume $\gamma_{n,u,\tau}^2<1$.

Let
$
\mathcal X_n:=\{x_1,\ldots,x_n\}.
$
Take $\epsilon:=n^{-2}$, and let
$
f_1,\ldots,f_{N_\epsilon}
$
be an $\epsilon$-cover of $\mathcal F_r$ with respect to
$\|\cdot\|_{\infty,\mathcal X_n}$.  By
\cref{thm:metric_entropy}, with $\bar R_r$ replaced by
$$
    \bar R_{r,n}:=1\vee (rR_{x,n}),
$$
and since
$n\ge2$ and $\bar R_{r,n}\ge1$,
\begin{equation}
\label{eq:fixed_design_cover_size}
\log N_\epsilon
\le
C_1\bar R_{r,n}^d\{\log(n\bar R_{r,n})\}^{d+1}
\le
C_1 n \epsilon_{\est}(r,R_{x,n}) .
\end{equation}
Here $C_1>0$ depends only on $d$.
Let $J\subseteq[N_\epsilon]$ be the set of indices $j$ for which there exists
one $p_j\in\mathcal F_r$ satisfying
\begin{equation}
\label{eq:fixed_design_representatives}
\|f_j-p_j\|_{\infty,\mathcal X_n}\le\epsilon,
\qquad
h_n^2(p_j,p^*)
\ge {1\over 4}\gamma_{n,u,\tau}^2 .
\end{equation}
For every $j\in J$, fix one such $p_j$, independently of
$Y_1,\ldots,Y_n$. The second inequality in
\eqref{eq:fixed_design_representatives} is needed for the likelihood-ratio
bound below and need not hold for the original covering centers $f_j$.
If $J$ is empty, every $p\in\cF_r$ satisfies
$h_n^2(p,p^*)<\gamma_{n,u,\tau}^2/4$, and the claim follows from
\cref{lem:bernoulli_hellinger_x}. Thus we assume $J\ne\varnothing$.

On the event
$$
\left\{
h_n^2(\wt p_r,p^*)
\ge
    {1\over 4}\gamma_{n,u,\tau}^2
\right\},
$$
there must exist  $j\in[N_\epsilon]$ such that
$$
\|\wt p_r-f_j\|_{\infty,\mathcal X_n}\le\epsilon .
$$
By \eqref{eq:fixed_design_representatives}, we also have  $j\in J$, and by the triangle inequality,
$$
\|\wt p_r-p_j\|_{\infty,\mathcal X_n}\le2\epsilon .
$$
Therefore, for every $i=1,\ldots,n$,
$$
\wt p_r(x_i)\le p_j(x_i)+2\epsilon,
\qquad
1-\wt p_r(x_i)\le1-p_j(x_i)+2\epsilon.
$$
Fix any $q\in\cF_r$, and choose $H\in\cP(\cB_r)$ with $p_H=q$.
On this event, the same index $j$ therefore satisfies
\begin{align*}
\prod_{i=1}^n
\frac{(p_j(x_i)+2\epsilon)^{Y_i}(1-p_j(x_i)+2\epsilon)^{1-Y_i}}
{L(q;Y_i,x_i)}
&\ge \prod_{i=1}^n\frac{L(\wt p_r;Y_i,x_i)}{L(q;Y_i,x_i)}\\
&=\exp\{n[\ell_n(\wt G_r)-\ell_n(H)]\}
\ge\exp(-n\tau),
\end{align*}
where the last inequality follows from \eqref{eq:fixed_design_oracle_lr}.
Consequently, for every fixed $q\in\cF_r$,
\begin{equation}
\label{eq:fixed_design_cover_event}
\begin{aligned}
&\PP\left\{h_n^2(\wt p_r,p^*)\ge\gamma_{n,u,\tau}^2/4\right\}\\
&\quad\le\PP\left\{
\max_{j\in J}\prod_{i=1}^n
\frac{(p_j(x_i)+2\epsilon)^{Y_i}(1-p_j(x_i)+2\epsilon)^{1-Y_i}}
{L(q;Y_i,x_i)}
\ge\exp(-n\tau)
\right\}.
\end{aligned}
\end{equation}
To bound the probability on the right hand side of
\eqref{eq:fixed_design_cover_event}, for each $j\in J$, write
\begin{align*}
&\prod_{i=1}^n
\frac{(p_j(x_i)+2\epsilon)^{Y_i}(1-p_j(x_i)+2\epsilon)^{1-Y_i}}
{L(q;Y_i,x_i)}\\
&\quad=
\left[\prod_{i=1}^n
\frac{(p_j(x_i)+2\epsilon)^{Y_i}(1-p_j(x_i)+2\epsilon)^{1-Y_i}}
{L(p^*;Y_i,x_i)}\right]
\left[\prod_{i=1}^n\frac{L(p^*;Y_i,x_i)}{L(q;Y_i,x_i)}\right].
\end{align*}
The product can be at least $\exp(-n\tau)$ only if the first factor is at
least $\exp(-n\tau-n\gamma_{n,u,\tau}^2/4)$ or the second factor is at
least $\exp(n\gamma_{n,u,\tau}^2/4)$. Thus, by
\eqref{eq:fixed_design_cover_event} and a union bound, for every fixed
$q\in\cF_r$,
\begin{align*}
&\PP\left\{h_n^2(\wt p_r,p^*)\ge\gamma_{n,u,\tau}^2/4\right\}\\
&\quad\le\underbrace{\PP\left\{
\max_{j\in J}\prod_{i=1}^n
\frac{
(p_j(x_i)+2\epsilon)^{Y_i}
(1-p_j(x_i)+2\epsilon)^{1-Y_i}
}
{L(p^*; Y_i, x_i)}
\ge
\exp\left(- n\tau -{1\over 4}n\gamma_{n,u,\tau}^2 \right)
\right\}}_{=:~\rI}\\
&\qquad+\PP\left\{
 \prod_{i=1}^n
\frac{L(p^*; Y_i, x_i)
}
{L(q; Y_i, x_i)}
\ge
\exp\left({1\over 4}n\gamma_{n,u,\tau}^2\right)
\right\}.
\end{align*}
Since $\rI$ does not depend on $q$, taking the infimum over
$q\in\cF_r$ gives
\begin{align*}
&\PP\left\{h_n^2(\wt p_r,p^*)\ge\gamma_{n,u,\tau}^2/4\right\}\\
&\quad\le\inf_{q\in\cF_r}\left[
\rI+\PP\left\{
\prod_{i=1}^n\frac{L(p^*;Y_i,x_i)}{L(q;Y_i,x_i)}
\ge\exp\left({1\over4}n\gamma_{n,u,\tau}^2\right)
\right\}\right]\\
&\quad=\rI+\inf_{q\in\cF_r}\PP\left\{
\prod_{i=1}^n\frac{L(p^*;Y_i,x_i)}{L(q;Y_i,x_i)}
\ge\exp\left({1\over4}n\gamma_{n,u,\tau}^2\right)
\right\}\\
&\quad=:\rI+\rII.
\end{align*}
We proceed to bound each term separately.

\paragraph{Bounding $\rI$}
By the union
bound, $|J| \le N_\epsilon$ and Markov's inequality applied to the square root of the likelihood
ratio,
\begin{align}
\rI \le N_\epsilon
\exp\left({n\tau \over 2}+  {n\gamma_{n,u,\tau}^2 \over 8}\right)
\max_{j\in J}
\EE
\prod_{i=1}^n
\left[
\frac{
(p_j(x_i)+2\epsilon)^{Y_i}
(1-p_j(x_i)+2\epsilon)^{1-Y_i}
}
{L(p^*; Y_i, x_i)}
\right]^{1/2}.
\label{eq:fixed_design_markov}
\end{align}
To bound the expectation in \eqref{eq:fixed_design_markov}, by independence
across $i\in[n]$, for every fixed $j\in J$,
\begin{align*}
&\EE
\prod_{i=1}^n
\left[
\frac{
(p_j(x_i)+2\epsilon)^{Y_i}
(1-p_j(x_i)+2\epsilon)^{1-Y_i}
}
{L(p^*; Y_i, x_i)}
\right]^{1/2}  \\
&\quad =
\prod_{i=1}^n
\left\{
\sqrt{(p_j(x_i)+2\epsilon)p^*(x_i)}
+
\sqrt{(1-p_j(x_i)+2\epsilon)(1-p^*(x_i))}
\right\}.
\end{align*}
Using $\sqrt{a+b}\le\sqrt a+\sqrt b$ and the fact that $p^*(x_i) \in (0,1)$, the $i$-th factor is bounded by
\begin{align*}
    &\sqrt{p_j(x_i)p^*(x_i)}
+
\sqrt{(1-p_j(x_i))(1-p^*(x_i))}
\\&\quad+
\sqrt{2\epsilon}\left(\sqrt{p^*(x_i)}+\sqrt{1-p^*(x_i)}\right) \\
\le&\sqrt{p_j(x_i)p^*(x_i)}
+
\sqrt{(1-p_j(x_i))(1-p^*(x_i))}
+
2\sqrt{2\epsilon} \\
=& 1-h^2(p_j(x_i),p^*(x_i)) + 2\sqrt{2\epsilon} &&\text{by \eqref{eq_hsquare_xi}}\\
\le & \exp \left\{
-h^2(p_j(x_i),p^*(x_i)) + 2\sqrt{2\epsilon}
\right\}.
\end{align*}
The last step uses $1+x\le e^x$. By \eqref{eq_hsquare_emp} and
\eqref{eq:fixed_design_representatives}, we obtain
\begin{align*}
\EE
\prod_{i=1}^n
\left[
\frac{
(p_j(x_i)+2\epsilon)^{Y_i}
(1-p_j(x_i)+2\epsilon)^{1-Y_i}
}
{L(p^*; Y_i, x_i)}
\right]^{1/2}
&\le
\exp\left\{
-nh_n^2(p_j,p^*)
+2n\sqrt{2\epsilon}
\right\} \\
&\le
\exp\left\{
-\frac{1}{4}n \gamma_{n,u,\tau}^2
+2n\sqrt{2\epsilon}
\right\}.
\end{align*}
Since $\epsilon=n^{-2}$, substituting this bound into \eqref{eq:fixed_design_markov} gives
\begin{align}
\rI
&\le
\exp\left\{
\log N_\epsilon
+{n\tau \over 2}
-\frac{n \gamma_{n,u,\tau}^2}{8}
+2n\sqrt{2\epsilon}
\right\}\nonumber\\
&\le \exp\left\{
C_1n\,\epsilon_{\est}(r,R_{x,n})+\frac{n\tau}{2}
-\frac{n \gamma_{n,u,\tau}^2}{8}+2\sqrt{2}\right\} &&\text{by \eqref{eq:fixed_design_cover_size}}\nonumber\\
&\le
{1\over 2}\exp(-u) &&\text{by \eqref{eq:fixed_design_lr_radius} and choosing $A$ large}. \label{eq:fixed_design_exp_bound_I}
\end{align}

\paragraph{Bounding $\rII$}

By Markov's inequality,
\begin{align*}
   \rII &\le \inf_{q\in \cF_r} \exp\left(-n\gamma_{n,u,\tau}^2/4\right)  \EE
 \prod_{i=1}^n
\frac{L(p^*; Y_i, x_i)
}
{L(q; Y_i, x_i)} \\
& = \inf_{q\in \cF_r}\exp\left(-n\gamma_{n,u,\tau}^2/4\right)
\prod_{i=1}^n \left\{
    {[p^*(x_i)]^2 \over q(x_i)} +  {[1-p^*(x_i)]^2 \over 1-q(x_i)}
\right\}\\
&= \inf_{q\in \cF_r}\exp\left(-n\gamma_{n,u,\tau}^2/4\right)
\prod_{i=1}^n \left\{ 1+ {(p^*(x_i) - q(x_i))^2 \over q(x_i) (1-q(x_i))}
\right\}.
\end{align*}
Using $1+x \le e^x$ and \eqref{def_chisquare_emp} gives
\begin{equation}\label{eq:fixed_design_exp_bound_II}
    \rII \le
   \exp \left\{ n \epsilon_{\app}(r, P_n)- {n\over 4}\gamma_{n,u,\tau}^2
\right\}\le {1\over 2}\exp(-u).
\end{equation}
Combining \eqref{eq:fixed_design_exp_bound_I} and \eqref{eq:fixed_design_exp_bound_II} gives
\begin{equation}\label{eq:fixed_design_hellinger_bound}
    \PP\left\{
    h_n^2(\wt p_r,p^*)
    \ge
    \frac{1}{4}\gamma_{n,u,\tau}^2
    \right\}
    \le e^{-u}.
\end{equation}
Finally, since Lemma~\ref{lem:bernoulli_hellinger_x} gives
$$
P_n (\wt p_r-p^*)^2
~ \le ~
4h_n^2(\wt p_r,p^*),
$$
the proof is complete.
\end{proof}

\subsection{Proof of \cref{lem_approx_error}}\label{app:proof_lem_approx_error}

\begin{proof}
Let $\Pi_n$ denote the orthogonal projection onto $\operatorname{span}\{x_1,\ldots,x_n\}\subseteq\RR^d$.
Let $T_r:\RR^d\to\cB_r$ be defined by
$$
    T_r(\beta)
    :=
    \begin{cases}
    \Pi_n\beta, & \|\Pi_n\beta\|_2\le r,\\[0.5ex]
    r{\Pi_n\beta\over \|\Pi_n\beta\|_2}, & \|\Pi_n\beta\|_2>r.
    \end{cases}
$$
Let $G_r:=T_r\#G^*$ be the pushforward of $G^*$, and set
$$
    q_r(x):=\int \sigma(x^\T\beta)\,\rd G_r(\beta)
    =
    \int \sigma(x^\T T_r(\beta))\,\rd G^*(\beta).
$$
Since $G_r\in\cP(\cB_r)$, we have $q_r\in\cF_r$. Therefore, by the
definition of $\epsilon_{\app}(r, P_n)$,
$$
    \epsilon_{\app}(r, P_n)
    \le
    {1\over n}\sum_{i=1}^n
    \chi^2
    \left(
    p^*(x_i),q_r(x_i)
    \right).
$$
By the joint convexity of the $\chi^2$-divergence,
$$
    \chi^2\!\left(
    p^*(x_i),q_r(x_i)
    \right)
    \le
    \int
    \chi^2\!\left(
    \sigma(x_i^\T\beta),
    \sigma(x_i^\T T_r(\beta))
    \right)
    \,\rd G^*(\beta).
$$
If $\|\Pi_n\beta\|_2\le r$, then
$$
    x_i^\T T_r(\beta)
    =
    x_i^\T\Pi_n\beta
    =
    x_i^\T\beta,
$$
so the integrand is zero. If $\|\Pi_n\beta\|_2>r$, write
$$
    t=x_i^\T\beta=x_i^\T\Pi_n\beta,
    \qquad
    s=x_i^\T T_r(\beta)
    =
    {r\over \|\Pi_n\beta\|_2}t .
$$
Then $s$ and $t$ have the same sign and $|s|\le |t|$. If $t\ge s\ge0$,
then
$$
    0\le \sigma(t)-\sigma(s)\le 1-\sigma(s),
    \qquad
    \sigma(s)\ge {1\over2},
$$
and hence
$$
    \chi^2(\sigma(t),\sigma(s))
    =
    {(\sigma(t)-\sigma(s))^2\over \sigma(s)(1-\sigma(s))}
    \le 1.
$$
The case $t\le s\le0$ is symmetric. Consequently,
$$
    \epsilon_{\app}(r, P_n)
    \le
    G^*(\|\Pi_n\beta\|_2>r).
$$
Finally, since $\|\Pi_n\beta\|_2\le \|\beta\|_2$,
$$
    G^*(\|\Pi_n\beta\|_2>r)
    \le
    G^*(\|\beta\|_2>r)
    =
    G^*(\cB_r^c).
$$
The weak-moment bound follows immediately from the definition of
$\mu^w_p(G^*)$:
$$
    G^*(\cB_r^c)
    \le
    \left({\mu^w_p(G^*)\over r}\right)^p .
$$

To prove the quadratic bound, fix any $0<s\le r$ such that
$G^*(\cB_s)>0$. For each $i\in[n]$, we have
$$
\begin{aligned}
    |p^*(x_i)-q_r(x_i)|
    &\le \int_{\|\Pi_n\beta\|_2>r}
    \left|\sigma(x_i^\T\beta)-\sigma(x_i^\T T_r(\beta))\right|
    \,\rd G^*(\beta)\\
    &\le {1\over2}G^*(\|\Pi_n\beta\|_2>r)
    \le {1\over2}G^*(\cB_r^c),
\end{aligned}
$$
since the two sigmoid values are on the same side of $1/2$ and their
difference is zero when $\|\Pi_n\beta\|_2\le r$.
For $\beta\in\cB_s$, one has $T_r(\beta)=\Pi_n\beta$ and
$\|T_r(\beta)\|_2\le s$, so that
$$
    |x_i^\T T_r(\beta)|
    \le \|x_i\|_2\|T_r(\beta)\|_2
    \le s\|x_i\|_2.
$$
Since $\sigma$ is increasing and $1-\sigma(u)=\sigma(-u)$, we have
$$
\begin{aligned}
    q_r(x_i)
    &\ge \int_{\cB_s}\sigma(x_i^\T T_r(\beta))\,\rd G^*(\beta) \ge \sigma(-s\|x_i\|_2)G^*(\cB_s)
    = {G^*(\cB_s)\over1+e^{s\|x_i\|_2}},\\
    1-q_r(x_i)
    &\ge \int_{\cB_s}\sigma(-x_i^\T T_r(\beta))\,\rd G^*(\beta) \ge \sigma(-s\|x_i\|_2)G^*(\cB_s)
    = {G^*(\cB_s)\over1+e^{s\|x_i\|_2}}.
\end{aligned}
$$
It follows that
$$
\begin{aligned}
    q_r(x_i)(1-q_r(x_i))
    &\ge {1\over2}\min\{q_r(x_i),1-q_r(x_i)\} \ge {G^*(\cB_s)\over2(1+e^{s\|x_i\|_2})},
\end{aligned}
$$
and thus
$$
    {(p^*(x_i)-q_r(x_i))^2\over q_r(x_i)(1-q_r(x_i))}
    \le {1+e^{s\|x_i\|_2}\over2G^*(\cB_s)}
    \left\{G^*(\cB_r^c)\right\}^2.
$$
Averaging over $i\in[n]$ and using $\|x_i\|_2\le R_{x,n}$ gives the
quadratic bound. Combining it with
$\epsilon_{\app}(r,P_n)\le G^*(\cB_r^c)$ proves the minimum bound.
\end{proof}

\subsection{Proof of \cref{cor:fixed_tail_conditions}}
\label{app:proof_fixed_tail_conditions}

\begin{proof}
Throughout the proof, we apply \cref{thm:mse_rate_fixed} with $\tau=0$
and $u=\log n$.
The bounded-support case follows from \cref{lem_approx_error}, since
$G^*(\cB_R)=1$ implies $\epsilon_{\app}(r, P_n)=0$ for every $r\ge R$.
Applying \cref{thm:mse_rate_fixed} gives the first claim.

For the exponential-decay case, Markov's inequality gives, for every $r>0$,
$$
    G^*(\cB_r^c)
    =
    G^*(\|\beta\|_2>r)
    \le
    \exp(-ar^q)
    \int_{\RR^d}\exp(a\|\beta\|_2^q)\,\rd G^*(\beta)
    \le
    b\exp(-ar^q).
$$
Thus $\epsilon_{\app}(r, P_n)\le b\exp(-ar^q)$. With
$$
    r_n
    =
    \left\{
    {1\over a}\log n
    \right\}^{1/q},
$$
we have $b\exp(-ar_n^q)=b/n$. Applying
\cref{thm:mse_rate_fixed} with $\tau=0$ proves the second claim.

For the finite weak $p$-th moment case, the definition of $\mu^w_p(G^*)=\mu$
gives
$$
    G^*(\cB_r^c)\le \left({\mu\over r}\right)^p,
    \qquad r>0.
$$
Taking $K=2^{1/p}\mu$ gives $G^*(\cB_K)\ge1/2$. By the second bound in
\cref{lem_approx_error}, for every $r\ge K$,
$$
    \epsilon_{\app}(r,P_n)
    \le \left(1+e^{KR_{x,n}}\right)\left({\mu\over r}\right)^{2p}.
$$
With
$$
    r_n:=\mu^{\frac{2p}{2p+d}}n^{\frac{1}{2p+d}},
$$
we have $r_n\ge K$ for sufficiently large $n$, and
$$
    {r_n^d\over n}
    =\left({\mu\over r_n}\right)^{2p}
    =\mu^{\frac{2pd}{2p+d}}n^{-\frac{2p}{2p+d}}.
$$
Applying \cref{thm:mse_rate_fixed} with $\tau=0$ proves the third claim.
\end{proof}

\section{Proofs for bounding the estimation error under random design}
\label{app:proof_mse_rate_random}

We prove \cref{thm:mse_rate_random} by first establishing a bounded-design
version under $P_X(\cB_T) = 1$ for some $T>0$, and then applying this result to a
truncated sub-Gaussian feature distribution. Throughout this section, recall from \eqref{def_Lf} that, for any $p\in\cF$,
$$
L(p;y,x)=p(x)^y(1-p(x))^{1-y},
\qquad \forall ~ (y,x)\in\{0,1\}\times\cX.
$$
Then $(y,x)\mapsto L(p;y,x)$ is the joint density under regression function $p$ with respect to $\nu\otimes P_X$, where $\nu$ is the counting measure on $\{0,1\}$. Let $h_X$ denote the Hellinger distance with respect to $\nu\otimes P_X$. For any $p,q\in\cF$, one has
$$
h_X^2(L(p;\cdot,\cdot),L(q;\cdot,\cdot)) = P_X h^2(p,q).
$$

\subsection{Joint-density formulation and bracketing entropy}

For $r>0$, write
\begin{equation}\label{def_Q_r}
    \cQ_r:=\{L(p;\cdot,\cdot):p\in\cF_r\}
\end{equation}
so that for each $p\in\cF_r$,
$
\log L(p;Y_i,X_i)=Y_i\log p(X_i)+(1-Y_i)\log(1-p(X_i)).
$
Then for $\wt p_r\in\cF_r$ satisfying \eqref{eq:fixed_design_oracle_lr}, we have
\begin{equation}\label{eq:joint_density_reformulation_random}
\sum_{i=1}^n \log L(\wt p_r;Y_i,X_i)
\ge \sup_{p\in\cF_r} \sum_{i=1}^n \log L(p;Y_i,X_i) -n\tau .
\end{equation}
For $T\ge0$, write
$$
    \bar R_{r,T} := 1\vee (rT).
$$
The following lemma bounds the bracketing entropy $H_{[]}(s,\cQ_r,h_X)$, which, for any $s>0$, is the logarithm of the minimum number of brackets of $h_X$-width at most $s$ needed to cover $\cQ_r$.

\begin{lemma}
\label{lem:hellinger_bracketing_transfer_random}
Suppose
$P_X(\cB_T)=1$ for some $T\ge 0$. Then for every $s\in(0,1)$,
$$
H_{[]}(s,\cQ_r,h_X)
\le
\log\cN\left({s^2\over 2},\cF_r,\|\cdot\|_{\i,\cB_T}\right).
$$
Consequently, there exists a constant $C>0$, depending only on $d$, such that, for every
$s\in(0,e^{-1}]$,
$$
H_{[]}(s,\cQ_r,h_X)
\le
C\,\bar R_{r,T}^d
\left(\log{C\bar R_{r,T}\over s}\right)^{d+1}.
$$
\end{lemma}

\begin{proof}
Fix $s\in(0,1)$, and set $\epsilon=s^2/2$. Let
$p_1,\ldots,p_N$ be an $\epsilon$-cover of $\cF_r$ under
$\|\cdot\|_{\i,\cB_T}$. For each $j$, define bracket endpoints on
$\{0,1\}\times\cB_T$ by
$$
\ell_j(1,x):=(p_j(x)-\epsilon)\vee0,
\qquad
u_j(1,x):=(p_j(x)+\epsilon)\wedge1,
$$
and
$$
\ell_j(0,x):=1-u_j(1,x),
\qquad
u_j(0,x):=1-\ell_j(1,x).
$$
If $\|p-p_j\|_{\i,\cB_T}\le \epsilon$, then
$$
    \ell_j(y,x)\le L(p;y,x)\le u_j(y,x),
    \qquad y\in\{0,1\},\ x\in\cB_T .
$$
Moreover, since $(\sqrt a-\sqrt b)^2\le a-b$ for $0\le b\le a$,
\begin{align*}
h_X^2(\ell_j,u_j)
&=
{1\over2}
\EE\left[
\left(\sqrt{u_j(1,X)}-\sqrt{\ell_j(1,X)}\right)^2
+
\left(\sqrt{u_j(0,X)}-\sqrt{\ell_j(0,X)}\right)^2
\right]  \\
&\le
{1\over2}
\EE\left[
u_j(1,X)-\ell_j(1,X)
+
u_j(0,X)-\ell_j(0,X)
\right] \\
&\le
2\epsilon
=
s^2 .
\end{align*}
Thus each bracket has Hellinger width at most $s$, proving the first claim.
The second claim follows from \cref{thm:metric_entropy}, applied on
$\cB_T$, with $R_x$ replaced by $T$.
\end{proof}

We shall use the following likelihood-ratio inequality of
\citet[Theorem~1]{wong1995probability}. It is stated here only in the form
needed below.

\begin{proposition}
\label{prop:wong_shen_random}
Let $\cQ$ be a class of densities on a dominated measurable space, let
$q_0\in\cQ$ be the true density, and let $h$ denote the Hellinger
distance. Let $Z_1,\ldots,Z_n$ be i.i.d. observations with density $q_0$,
and write $\PP^*_{q_0}$ for outer probability under their joint law.
There exist universal constants $a\ge1$ and $c,c_1,c_2>0$ such that
the following holds. If, for some $\gamma>0$,
$$
\int_{\gamma^2/2^8}^{\sqrt{2}\gamma}
\left(H_{[]}(s/a,\cQ,h)\right)^{1/2}\,\rd s
\le
c\sqrt n\,\gamma^2,
$$
then
$$
\PP^*_{q_0}\left(
\sup_{\{q\in\cQ:\,h(q,q_0)\ge \gamma\}}
\prod_{i=1}^n {q(Z_i)\over q_0(Z_i)}
\ge
e^{-c_1n\gamma^2}
\right)
\le
4e^{-c_2n\gamma^2}.
$$
\end{proposition}

\subsection{Estimation error under bounded design}
\leavevmode\par

Recall $\epsilon_{\est}(r,T)$ from \eqref{def_epsilon_1} and $\epsilon_{\app}(r, P_X)$ from \eqref{def_chisquare_pop}.

\begin{lemma}
\label{lem:random_bounded_design}
Assume $p^*\in\cF$ and $P_X(\cB_T)=1$ for some $T>0$. Fix $r>0$ and
$\tau\ge0$. Let $\wt p_r\in\cF_r$ be any measurable estimator satisfying
\eqref{eq:fixed_design_oracle_lr}. Then there exists a constant
$C>0$, depending only on $d$, such that, for every $u\ge1$,
$$
\PP\left\{
P_X(\wt p_r-p^*)^2
>
C\left(
\epsilon_{\est}(r,T)
+
\epsilon_{\app}(r, P_X)
+
\tau
+
{u\over n}
\right)
\right\}
\le e^{-u}.
$$
\end{lemma}

\begin{proof}
 Since $p^*$ need not belong to $\cF_r$, we apply the Wong--Shen inequality
to the enlarged class
$$
    \overline{\cQ}_r:=\cQ_r\cup\{L(p^*;\cdot,\cdot)\}
$$
with $\cQ_r$ given in \eqref{def_Q_r}.
The addition of the single density $L(p^*;\cdot,\cdot)$ changes the bracketing entropy by at
most a universal additive constant. Hence, by
\cref{lem:hellinger_bracketing_transfer_random},
\begin{equation}\label{bracketing_H}
H_{[]}(s,\overline{\cQ}_r,h_X)
\le
C\,\bar R_{r,T}^d
\left(\log{C\bar R_{r,T}\over s}\right)^{d+1},
\qquad s\in(0,e^{-1}].
\end{equation}

Fix $u\ge1$, and set
$$
\gamma_{n,u,\tau}^2
:=
A\left(
\epsilon_{\est}(r,T)+\epsilon_{\app}(r, P_X)+\tau+{u\over n}
\right),
$$
where $A\ge1$ is a sufficiently large constant depending only on $d$, to be chosen below.
If
$
    \gamma_{n,u,\tau}>{e^{-1} / \sqrt2},
$
then the claim is trivial after enlarging the final constant, since
$P_X(\wt p_r-p^*)^2\le1$. We therefore assume
$$
    \gamma_{n,u,\tau}\le {e^{-1}\over\sqrt2}.
$$
By the bracketing entropy bound in \eqref{bracketing_H},
\begin{align*}
&\int_{\gamma_{n,u,\tau}^2/2^8}^{\sqrt2\gamma_{n,u,\tau}}
\left(H_{[]}(s/a,\overline{\cQ}_r,h_X)\right)^{1/2}\,\rd s  \le
C_1\gamma_{n,u,\tau}\bar R_{r,T}^{d/2}
\left(
\log{C_1\bar R_{r,T}\over \gamma_{n,u,\tau}^2}
\right)^{(d+1)/2}.
\end{align*}
Since
$$
    \gamma_{n,u,\tau}^2
    \ge
    A\,{\bar R_{r,T}^d\{\log(n\bar R_{r,T})\}^{d+1}\over n},
$$
and $\bar R_{r,T}\ge1$, one has
$$
    \log{C_1\bar R_{r,T}\over \gamma_{n,u,\tau}^2}
    \le
    C_2\log(n\bar R_{r,T}).
$$
It follows that
$$
\int_{\gamma_{n,u,\tau}^2/2^8}^{\sqrt2\gamma_{n,u,\tau}}
\left(H_{[]}(s/a,\overline{\cQ}_r,h_X)\right)^{1/2}\,\rd s
\le
{C_3\over\sqrt A}\sqrt n\,\gamma_{n,u,\tau}^2.
$$
Here $C_1,C_2,C_3>0$ depend only on $d$.
For sufficiently large $A$, the entropy integral condition in
\cref{prop:wong_shen_random} holds with
$\gamma=\gamma_{n,u,\tau}$.

For any $p^\circ\in\cF_r$,
\eqref{eq:joint_density_reformulation_random} gives
\begin{equation}\label{approx_ineq}
\prod_{i=1}^n{L(\wt p_r;Y_i,X_i)\over L(p^\circ;Y_i,X_i)}
\ge e^{-n\tau}.
\end{equation}
Write $$
A_{n,u}:=
\left\{
\sup_{\{q\in\overline{\cQ}_r:\,h_X(q,L(p^*;\cdot,\cdot))\ge\gamma_{n,u,\tau}\}}
\prod_{i=1}^n{q(Y_i,X_i)\over L(p^*;Y_i,X_i)}
\ge
e^{-c_1n\gamma_{n,u,\tau}^2}
\right\}
$$
and
$$
B_{n,u}(p^\circ):=
\left\{
\prod_{i=1}^n{L(p^*;Y_i,X_i)\over L(p^\circ;Y_i,X_i)}
\ge
\exp\{c_1n\gamma_{n,u,\tau}^2-n\tau\}
\right\}.
$$
Observe that
$$
     F_{n,u}:=\{h_X(L(\wt p_r;\cdot,\cdot),L(p^*;\cdot,\cdot))\ge \gamma_{n,u,\tau}\} \subseteq A_{n,u}\cup B_{n,u}(p^\circ)
$$
Indeed, if $F_{n,u}$ holds and $A_{n,u}$ does not hold, then
$$
\prod_{i=1}^n{L(\wt p_r;Y_i,X_i)\over L(p^*;Y_i,X_i)}
<
e^{-c_1n\gamma_{n,u,\tau}^2},
$$
while \eqref{approx_ineq} gives
$$
\prod_{i=1}^n{L(\wt p_r;Y_i,X_i)\over L(p^*;Y_i,X_i)}
\ge
e^{-n\tau}
\prod_{i=1}^n{L(p^\circ;Y_i,X_i)\over L(p^*;Y_i,X_i)},
$$
implying $B_{n,u}(p^\circ)$.

On the one hand, applying \cref{prop:wong_shen_random} with
$\cQ=\overline{\cQ}_r$, $Z_i=(Y_i,X_i)$, and $q_0=L(p^*;\cdot,\cdot)$ gives
$$
    \PP^*(A_{n,u})
    \le
    4e^{-c_2n\gamma_{n,u,\tau}^2}.
$$
On the other hand, Markov's inequality gives
\begin{align*}
\PP(B_{n,u}(p^\circ))
&\le
\exp\left\{-c_1 n\gamma_{n,u,\tau}^2+n\tau\right\}
\EE
\prod_{i=1}^n{L(p^*;Y_i,X_i)\over L(p^\circ;Y_i,X_i)}  \\
&=
\exp\left\{-c_1n\gamma_{n,u,\tau}^2+n\tau\right\}
\left[
1+
\EE{(p^*(X)-p^\circ(X))^2\over p^\circ(X)(1-p^\circ(X))}
\right]^n  \\
&\le
\exp\left\{
-nc_1\gamma_{n,u,\tau}^2+n\tau
+n\EE{(p^*(X)-p^\circ(X))^2\over p^\circ(X)(1-p^\circ(X))}
\right\}.
\end{align*}
Taking the infimum over $p^\circ\in\cF_r$, we obtain
$$
\PP(F_{n,u})
\le
4e^{-c_2n\gamma_{n,u,\tau}^2}
+
\exp\left\{
-nc_1\gamma_{n,u,\tau}^2+n\tau+n\epsilon_{\app}(r, P_X)
\right\}.
$$
By choosing $A$ sufficiently large, the right-hand side is at most $e^{-u}$.

Finally, by \cref{lem:bernoulli_hellinger_x}, integrated with respect to
$P_X$,
$$
P_X(\wt p_r-p^*)^2
\le
4h_X^2(L(\wt p_r;\cdot,\cdot),L(p^*;\cdot,\cdot)).
$$
This proves the lemma after enlarging $C$.
\end{proof}

\subsection{Proof of \cref{thm:mse_rate_random}}
\label{app:random_design_theorem_proof}

We now apply \cref{lem:random_bounded_design} after truncating the covariate distribution.
The following quadratic-form bound follows from \citet[Theorem~2.1]{Hsu2012} and will be used in the truncation argument.
\begin{lemma} \label{lem_quad}
    Let $\xi\in\RR^d$ be a sub-Gaussian random vector satisfying \cref{assum:subgaussian_x} with $m_x=0$ and $\sigma = \gamma_\xi>0$.  Then, for every symmetric positive semidefinite
    matrix $H\in\RR^{d\times d}$ and every $t\ge0$,
    $$
    \PP\left\{
    \xi^\T H\xi > \gamma_\xi^2\left(
    \sqrt{{\rm tr}(H)}+ \sqrt{2 t \|H\|_{\rm op} }
    \right)^2
    \right\} \le e^{-t}.
    $$
\end{lemma}

\begin{proof}[Proof of \cref{thm:mse_rate_random}]
Fix $u\ge1$, and recall $T:=T_u=m_x+C_{\sigma} \sqrt{d + u+\log(n)}.$
Let
$$
\cE_T:=\bigcap_{i=1}^n\{X_i\in \cB_T\}.
$$
By \cref{assum:subgaussian_x}, the vector $\xi=X_1-\EE X_1$ satisfies
the condition of \cref{lem_quad} with $\gamma_\xi=\sigma$.
Applying this lemma with $H=I_d$ and $t=u+\log(4n)$ gives
$$
    \PP\left\{
        \|X_1-\EE X_1\|_2>\sigma\left(\sqrt d+\sqrt{2t}\right)
    \right\}\le e^{-t}.
$$
Since $\|\EE X_1\|_2\le m_x$ and, for $C_\sigma\ge4\sigma$,
$$
    \sigma\left(\sqrt d+\sqrt{2t}\right)
    \le4\sigma\sqrt{d+u+\log n}
    \le T-m_x,
$$
the triangle inequality gives
$$
    P_X(\cB_T^c)
    \le \PP\left\{\|X_1-\EE X_1\|_2>T-m_x\right\}
    \le e^{-t}={e^{-u}\over4n}.
$$
Taking a union bound over $i\in[n]$, we obtain
\begin{equation}\label{bd_event_comp}
\PP(\cE_T^c)\le nP_X(\cB_T^c)\le {1\over4}e^{-u}.
\end{equation}
Furthermore, $P_X(\cB_T)\ge 1/2$.

Define
$$
    P_X^{(T)}(\cdot):=P_X(\cdot\mid\cB_T).
$$
Conditional on $\cE_T$, the observations $(X_i,Y_i)_{i=1}^n$ remain i.i.d. from the logistic-mixture model with feature marginal $P_X^{(T)}$ and the same conditional success probability $p^*(X_i)$. Indeed, conditioning on $\cE_T$ only restricts each $X_i$ to $\cB_T$, so the feature vectors remain independent with common law $P_X^{(T)}$, while the conditional law of $Y_i$ given $X_i$ is unchanged. Moreover, since neither the empirical likelihood nor the class $\cF_r$ changes on the observed sample, $\wt p_r$ remains a $\tau$-approximate sieve NPMLE over $\cF_r$. Conditional on $\cE_T$, by applying \cref{lem:random_bounded_design} to the truncated feature distribution $P_X^{(T)}$, with $u$ replaced by $u+\log 4$, and by noting that
$$
\epsilon_{\app}(r, P_X^{(T)})
=
\inf_{q\in\cF_r}
{1\over P_X(\cB_T)}P_X\left[
{(p^*-q)^2\over q(1-q)}
\mathbf 1_{\cB_T}
\right]
\le
{1\over P_X(\cB_T)}\epsilon_{\app}(r, P_X)
\le
2\epsilon_{\app}(r, P_X),
$$
we obtain
$$
\PP\left\{
P_X^{(T)}(\wt p_r-p^*)^2
>
C\left(
\epsilon_{\est}(r,T)
+
\epsilon_{\app}(r, P_X)
+
\tau
+
{u\over n}
\right)
\,\middle|\,\cE_T
\right\}
\le {e^{-u}\over4}.
$$
Since
\begin{align*}
P_X(\wt p_r-p^*)^2
&=P_X(\cB_T)P_X^{(T)}(\wt p_r-p^*)^2
+\EE\!\left[\{\wt p_r(X)-p^*(X)\}^2\mathbf 1_{\cB_T^c}\right] \\
&\le P_X^{(T)}(\wt p_r-p^*)^2+P_X(\cB_T^c),
\end{align*}
where the inequality follows from $0\le \wt p_r(x),p^*(x)\le 1$ for all $x$, in conjunction with $P_X(\cB_T^c)\le e^{-u}/(4n)\le u/n$, we obtain
$$
\PP\left\{
P_X(\wt p_r-p^*)^2
>
C\left(
\epsilon_{\est}(r,T_u)
+
\epsilon_{\app}(r, P_X)
+
\tau
+
{u\over n}
\right)
\right\}
\le e^{-u},
$$
after enlarging $C$. This completes the proof.
\end{proof}

\subsection{Proof of \cref{cor:random_tail_conditions}}
\label{app:proof_random_tail_conditions}

\begin{proof}
Throughout the proof, we apply \cref{thm:mse_rate_random} with $\tau=0$ and $u=\log n$.

For the bounded-support case, if $G^*(\cB_R)=1$ and $r\ge R$, then
\eqref{eq:pop_chisquare_tail_bound} gives
$
    \epsilon_{\app}(r, P_X)=0.
$
Applying \cref{thm:mse_rate_random} gives
$$
\PP\left\{
P_X(\wh p_{r}-p^*)^2
>
C\left(
{\bar R_{r,T_u}^d
\{\log(n\bar R_{r,T_u})\}^{d+1}\over n}
+
{u\over n}
\right)
\right\}
\le e^{-u}.
$$

For the exponential-decay case, Markov's inequality gives, for every $r>0$,
$$
    G^*(\cB_r^c)
    =
    G^*(\|\beta\|_2>r)
    \le
    \exp(-ar^q)
    \int_{\RR^d}\exp(a\|\beta\|_2^q)\,\rd G^*(\beta)
    \le
    b\exp(-ar^q).
$$
Thus $\epsilon_{\app}(r, P_X)\le b\exp(-ar^q)$. With
$
    r_n = \{ (1/a)\log n \}^{1/q},
$
we have $b\exp(-ar_n^q)=b/n\le u/n$ for sufficiently large $n$. Applying
\cref{thm:mse_rate_random} gives the stated bound.

For the finite weak $p$-th moment case, take $s=2^{1/p}\mu$, so that
$G^*(\cB_s)\ge1/2$. By \cref{assum:subgaussian_x},
\begin{align*}
    \EE e^{s\|X\|_2}
    &\le e^{sm_x}\sum_{v\in\{-1,1\}^d}
        \EE e^{sv^\T(X-\EE X)} \le 2^d\exp\left(sm_x+{\sigma^2s^2d\over2}\right).
\end{align*}
The first step uses $\|z\|_2\le\|z\|_1=\max_{v\in\{-1,1\}^d}v^\T z$;
the second uses $\|v\|_2^2=d$. Combining this bound with
\eqref{eq:pop_chisquare_tail_bound} gives
$$
    \epsilon_{\app}(r, P_X)\le C\left({\mu\over r}\right)^{2p},
    \qquad r\ge s,
$$
where $C$ depends only on $p,\mu,d,m_x,\sigma$.
Taking
$$
    r_{n}
    =
    \mu^{2p/(2p+d)}n^{1/(2p+d)}(d+\log n)^{-d/(4p+2d)}
$$
gives $r_n\ge s$ for sufficiently large $n$. Applying
\cref{thm:mse_rate_random} gives
$$
\PP\left\{
P_X(\wh p_{r_n}-p^*)^2
>
C\left(
{\bar R_{r_{n},T_u}^d
\{\log(n\bar R_{r_{n},T_u})\}^{d+1}\over n}
+
\left({\mu\over r_{n}}\right)^{2p}
+
{u\over n}
\right)
\right\}
\le e^{-u}.
$$
Since $u=\log n$, one has $T_u\lesssim\sqrt{d+\log n}$ and
$\log(n\bar R_{r_n,T_u})\lesssim\log n$ for sufficiently large $n$.
Moreover,
$$
    {r_n^d(d+\log n)^{d/2}\over n}
    =\left({\mu\over r_n}\right)^{2p}
    =\mu^{2pd/(2p+d)}(d+\log n)^{pd/(2p+d)}n^{-2p/(2p+d)}.
$$
Substituting these bounds proves the third claim.
\end{proof}

\section{Proofs of convergence for the RPG algorithm with a fixed dataset}
\label{app:particle_fixed_data_analysis}

Throughout this section, the data $\{(Y_i,x_i)\}_{i\in[n]}$ are fixed, and we have $r>0$,
$R_{x,n}>0$, and $0<\eta\le  8\rho_{n,r}^2 / R_{x,n}^2$. Recall $\bm{p}_G$,
$\bm{p}_{\delta_\beta}$, $\bm{r}(p_G)$,
$g_{n,G}(\beta)$, $\bar g_{n,G}$, and $\langle\cdot,\cdot\rangle_n$ from
\cref{sec:particle_prediction_formulation}, and $\rho_{n,r}$ and $E_t$
from \cref{sec:particle_fixed_data}. Further recall the empirical inner product $\langle u, v\rangle_n = n^{-1}\sum_{i=1}^n u_i v_i.$

The road-map  of proving \cref{thm:reweighted_particle_global} is as follows.
By \cref{lem:particle_one_block_ascent}, the likelihood gap decreases
by at least $(\eta/2)E_t$ at each outer iteration $t$. It therefore suffices to find
$b_n,c_n>0$, independent of $M$, such that
\begin{equation}
    E_t\ge c_n\{\ell_n(\wh G_r)-\ell_n(G_t)\}^2
    \quad\text{if}\quad
    0<\ell_n(\wh G_r)-\ell_n(G_t)\le b_n,
    \label{eq:app_energy_target}
\end{equation}
and to show that with high probability,
\begin{equation}\label{eq:init_G_fixed}
    \ell_n(\wh G_r)-\ell_n(G_0)\le b_n.
\end{equation}
Monotonicity in \cref{lem:particle_one_block_ascent}
then keeps every subsequent likelihood gap below $b_n$, and the recursion argument in \cref{lem:app_particle_quadratic_recursion} gives the $\cO(t^{-1})$ rate.

We start by stating a few useful results in our proof. We prove \cref{lem:particle_one_block_ascent} in \cref{app_proof_lem:particle_one_block_ascent}. The initialization requirement in \eqref{eq:init_G_fixed} is established in \cref{app:sec_initial_fixed}. Proving \eqref{eq:app_energy_target} is technically challenging and we defer its proof in \cref{app:particle_fixed_data_auxiliary}.

\subsection{Preliminaries}
\label{app:particle_fixed_data_preliminaries}

We first state some facts on the smoothness and concavity of the log-likelihood function $\ell_n(G)$ with respect to $\bm p_G$. We then prove the KKT condition \eqref{eq:exact_reweight_t_kkt} of exact weighting. Finally we review the projection
inequality used in our projected gradient-ascent step in \eqref{eq:reweighted_particle_update}.

\paragraph{Prediction and likelihood bounds}
For every $G\in\cP(\cB_r)$ with $p_G$ given in \eqref{def_p_G},  we have
\begin{equation}\label{bd_r_p_uniform}
\rho_{n,r}\le p_G(x_i)\le1-\rho_{n,r},\qquad  |r_i(p_G)|\le {1\over \rho_{n,r}},\qquad i\in [n].
\end{equation}
Since $\sigma'(z)=\sigma(z)\{1-\sigma(z)\}\le1/4$, the mean value
theorem gives, for $\beta,\wt\beta\in\cB_r$,
\begin{equation}
\begin{aligned}
    \|\bm{p}_{\delta_\beta}-\bm{p}_{\delta_{\wt\beta}}\|_n^2
    &\le{1\over16n}\sum_{i=1}^n
       \{x_i^\T(\beta-\wt\beta)\}^2 \le{R_{x,n}^2\over16}\|\beta-\wt\beta\|_2^2.
\end{aligned}
\label{eq:app_fixed_data_smoothness}
\end{equation}
For each summand in the likelihood \eqref{def_llh}, differentiation
with respect to its scalar argument gives
$$
\begin{aligned}
    \left.{d\over du}\{Y_i\log u+(1-Y_i)\log(1-u)\}\right|_{u=p_G(x_i)}
       &=r_i(p_G),\\
    -\left.{d^2\over du^2}\{Y_i\log u+(1-Y_i)\log(1-u)\}\right|_{u=p_G(x_i)}
       &={Y_i\over p_G(x_i)^2}+{1-Y_i\over\{1-p_G(x_i)\}^2}.
\end{aligned}
$$
It is easy to see that the mixed second derivatives of the log-likelihood are zero.
Since $Y_i\in\{0,1\}$, the absolute value of the second derivative of
each summand lies between $(1-\rho_{n,r})^{-2}$ and $\rho_{n,r}^{-2}$.
For $G,\wt G\in\cP(\cB_r)$ with corresponding prediction vectors $\bm{p}_G$ and $\bm{p}_{\wt G}$,
a second-order Taylor expansion along their line segment gives
\begin{align}\label{eq:app_fixed_data_ell_bounds_lower}
    \ell_n(\wt G)-\ell_n(G)
    &~ \ge~ \langle \bm{r}(p_G),\bm{p}_{\wt G}-\bm{p}_G\rangle_n
        -{\|\bm{p}_{\wt G}-\bm{p}_G\|_n^2\over2\rho_{n,r}^2},\\\label{eq:app_fixed_data_ell_bounds_upper}
    \ell_n(\wt G)-\ell_n(G)
    &~ \le~ \langle \bm{r}(p_G),\bm{p}_{\wt G}-\bm{p}_G\rangle_n
        -{\|\bm{p}_{\wt G}-\bm{p}_G\|_n^2\over2(1-\rho_{n,r})^2}.
\end{align}

Recall that $\wh G_r$ maximizes $\ell_n(G)$ over $\cP(\cB_r)$.
Viewing the likelihood as a function of the prediction vector, its
first-order optimality condition gives
$\langle \bm{r}(\bm{p}_{\wh G_r}),\bm{p}_G-\bm{p}_{\wh G_r}\rangle_n\le0$
for every $G\in\cP(\cB_r)$.
The  bound in \eqref{eq:app_fixed_data_ell_bounds_upper} therefore implies
\begin{equation}
    \ell_n(\wh G_r)-\ell_n(G)
    \ge{\|\bm{p}_G-\bm{p}_{\wh G_r}\|_n^2\over2(1-\rho_{n,r})^2}.
    \label{eq:app_empirical_gap_distance}
\end{equation}

\paragraph{Exact reweighting}
The following lemma states the optimality conditions for the exact reweighting step in \eqref{eq:exact_reweight}.

\begin{lemma}
\label{lem:app_particle_weight_kkt}
Fix $\bbeta\in\cB_r^M$ and $w\in\cM_n(\bbeta)$.
With $G=\sum_{j=1}^Mw_j\delta_{\beta_j}$, one has \begin{equation}
    g_{n,G}(\beta_j)\le\bar g_{n,G},
    \qquad
    w_j\{g_{n,G}(\beta_j)-\bar g_{n,G}\}=0,
    \qquad j\in[M].
    \label{eq:app_particle_weight_kkt}
\end{equation}
\end{lemma}
\begin{proof}

The objective $v\mapsto\ell_n(\sum_jv_j\delta_{\beta_j})$ is concave and
differentiable on $\Delta^{M-1}$, with
$$
    \left.{\partial\over\partial v_j}
       \ell_n\left(\sum_{k=1}^Mv_k\delta_{\beta_k}\right)\right|_{v=w}
    =\langle \bm{r}(p_G),\bm{p}_{\delta_{\beta_j}}\rangle_n=g_{n,G}(\beta_j).
$$
The simplex KKT conditions \citep[Section~5.5.3]{boyd2004convex}
give a scalar $\lambda$ such that $g_{n,G}(\beta_j)\le\lambda$ for
every $j$, with equality whenever $w_j>0$. Multiplying by $w_j$
and summing over $j$, we obtain
$\lambda=\sum_jw_jg_{n,G}(\beta_j)=\bar g_{n,G}$, which proves
\eqref{eq:app_particle_weight_kkt}.
\end{proof}

\paragraph{Euclidean projection}
\label{app:particle_projection_prelim}

For a nonempty closed convex set $\cC\subset\RR^d$, the Euclidean projection
$y=\Pi_{\cC}(z)$ is characterized by
\begin{equation}
    \langle z-y,v-y\rangle\le0,\qquad v\in\cC.
    \label{eq:app_projection_inequality}
\end{equation}
The normal cone at $\beta\in\cC$ is
\begin{equation}
    N_{\cC}(\beta)=\{u:\langle u,v-\beta\rangle\le0
                    \text{ for every }v\in\cC\}.
    \label{eq:app_normal_cone}
\end{equation}
Consequently, a projected ascent update leaves $\beta$ unchanged if and
only if its gradient belongs to $N_{\cC}(\beta)$. For the ball $\cB_r$,
the normal cone is $\{0\}$ in the interior and consists of the
nonnegative multiples of $\beta$ on the boundary.

\subsection{Proof  of \cref{lem:particle_one_block_ascent}: the one-step ascent lemma}\label{app_proof_lem:particle_one_block_ascent}
\begin{proof}

We first bound the likelihood gain from one location update, then sum
the gains over a block and apply exact reweighting.

Fix a block $t$ and write $w=w^t$, $\beta_j=\beta_j^{t,0}$,
$G=G_t=\sum_jw_j\delta_{\beta_j}$, and
$\bm{p}_G=\bm{p}_{G_t}=\sum_jw_j\bm{p}_{\delta_{\beta_j}}$. Let
$$
    \delta_j=\Pi_{\cB_r}\{\beta_j+\eta\nabla g_{n,G}(\beta_j)\}-\beta_j.
$$
The updated measure $G^+=\sum_jw_j\delta_{\beta_j+\delta_j}$ has prediction
vector $\bm{p}_{G^+}=\sum_{j=1}^Mw_j\bm{p}_{\delta_{\beta_j+\delta_j}}$.
Thus $\delta_j=\beta_j^{t,1}-\beta_j^{t,0}$,
$\bm{p}_{G^+}=\bm{p}_{G^{t,1}}$, and $E_t=\eta^{-2}\sum_jw_j\|\delta_j\|_2^2$.  By \eqref{eq:app_fixed_data_ell_bounds_lower}, we have
\begin{equation}\label{init_diff_ln}
     \ell_n(G^+)-\ell_n(G)
     ~ \ge~ \langle \bm{r}(p_G),\bm{p}_{G^+}-\bm{p}_G\rangle_n
       -{\|\bm{p}_{G^+}-\bm{p}_G\|_n^2\over2\rho_{n,r}^2}.
\end{equation}
By Jensen's inequality,
\begin{align} \nonumber
    \|\bm{p}_{G^+}-\bm{p}_G\|_n^2
    &=\left\|\sum_{j=1}^Mw_j\{\bm{p}_{\delta_{\beta_j+\delta_j}}-\bm{p}_{\delta_{\beta_j}}\}\right\|_n^2\\\nonumber
    &\le\sum_{j=1}^Mw_j\|\bm{p}_{\delta_{\beta_j+\delta_j}}-\bm{p}_{\delta_{\beta_j}}\|_n^2\\
    &\le{R_{x,n}^2\over16}\sum_{j=1}^Mw_j\|\delta_j\|_2^2 &&\text{by  \eqref{eq:app_fixed_data_smoothness}}.
\label{eq:app_prediction_displacement}
\end{align}
To bound from below $\langle \bm{r}(p_G),\bm{p}_{G^+}-\bm{p}_G\rangle_n$,
a second-order Taylor expansion of $g_{n,G}$ along
$[\beta_j,\beta_j+\delta_j]\subset\cB_r$  gives
$$
\begin{aligned}
    \langle \bm{r}(p_G),\bm{p}_{\delta_{\beta_j+\delta_j}}-\bm{p}_{\delta_{\beta_j}}\rangle_n
    &=g_{n,G}(\beta_j+\delta_j)-g_{n,G}(\beta_j)\\
    &=\langle\nabla g_{n,G}(\beta_j),\delta_j\rangle
      +{1\over2}\delta_j^\T\nabla^2 g_{n,G}(\bar\beta)\delta_j
\end{aligned}
$$
for some $\bar\beta$ on the segment $[\beta_j,\beta_j+\delta_j]$.
Applying \eqref{eq:app_projection_inequality} with
$z=\beta_j+\eta\nabla g_{n,G}(\beta_j)$, $y=\beta_j+\delta_j$,
and $v=\beta_j$ gives
$\langle\eta\nabla g_{n,G}(\beta_j)-\delta_j,-\delta_j\rangle\le0$ so that
\begin{equation}
    \langle\nabla g_{n,G}(\beta_j),\delta_j\rangle
    \ge{1\over\eta}\|\delta_j\|_2^2.
    \label{eq:app_projection_gain}
\end{equation}
To bound the second-order term, write $u=\sigma(z)$, so that
$\sigma''(z)=u(1-u)(1-2u)$. The derivative of this polynomial is
$1-6u+6u^2$, with roots $(3\pm\sqrt3)/6$; evaluating at these roots
and at $0$ and $1$ gives $\sup_z|\sigma''(z)|=1/(6\sqrt3)$. Hence for all $\beta\in\cB_r$,
\begin{align}\nonumber
    \|\nabla^2g_{n,G}(\beta)\|_{\mathrm{op}} &= \left\| {1\over n}\sum_{i=1}^n r_i(p_G)\sigma''(x_i^\T\beta)x_ix_i^\T\right\|_\op \\\nonumber
    &\le{1\over n}\sum_{i=1}^n
       {|\sigma''(x_i^\T\beta)|\|x_i\|_2^2\over\rho_{n,r}}\\
       &\le{R_{x,n}^2\over6\sqrt3\,\rho_{n,r}}. \label{bd_Hessian_gp}
\end{align}
Combining \eqref{eq:app_projection_gain} and \eqref{bd_Hessian_gp} gives
\begin{equation}
\begin{aligned}
    \langle \bm{r}(p_G),\bm{p}_{\delta_{\beta_j+\delta_j}}-\bm{p}_{\delta_{\beta_j}}\rangle_n
    &\ge\left({1\over\eta}-{R_{x,n}^2\over12\sqrt3\,\rho_{n,r}}\right)
       \|\delta_j\|_2^2
\end{aligned}
\label{eq:app_linear_prediction_gain}
\end{equation}
so that by multiplying by $w_j$ and summing over $j$, we obtain
\begin{equation}
    \langle \bm{r}(p_G),\bm{p}_{G^+}-\bm{p}_G\rangle_n
    \ge\left({1\over\eta}-{R_{x,n}^2\over12\sqrt3\,\rho_{n,r}}\right)
       \sum_{j=1}^Mw_j\|\delta_j\|_2^2.
    \label{eq:app_sum_linear_gain}
\end{equation}

Plugging \eqref{eq:app_sum_linear_gain} and \eqref{eq:app_prediction_displacement} into \eqref{init_diff_ln} yields
\begin{align*}
    \ell_n(G^+)-\ell_n(G)
    &\ge\left({1\over\eta}
             -{R_{x,n}^2\over12\sqrt3\,\rho_{n,r}}
             -{R_{x,n}^2\over32\rho_{n,r}^2}\right)
       \sum_{j=1}^Mw_j\|\delta_j\|_2^2\\
    &\ge{1\over2\eta}\sum_{j=1}^Mw_j\|\delta_j\|_2^2
    ={\eta\over2}E_t.
\end{align*}
The last inequality uses $0<\eta\le  8\rho_{n,r}^2 / R_{x,n}^2$ and $\rho_{n,r} \le 1/2$ so that
$\eta\{R_{x,n}^2/(6\sqrt3\,\rho_{n,r})
+R_{x,n}^2/(16\rho_{n,r}^2)\}\le1$.
The same calculation applies to every location update in the block,
so each has nonnegative likelihood gain. At the end of the block,
$w^{t+1}\in\cM_n(\bbeta^{t,H})$, whereas $G^{t,H}$ uses the feasible
weights $w^t$ at these locations. Since $G^{t,0}=G_t$, we obtain
\begin{align*}
    \ell_n(G_{t+1})-\ell_n(G_t)
    &\ge\ell_n(G^{t,H})-\ell_n(G^{t,0})\\
    &=\sum_{h=0}^{H-1}\{\ell_n(G^{t,h+1})-\ell_n(G^{t,h})\}\\
    &\ge\ell_n(G^{t,1})-\ell_n(G^{t,0})
     \ge{\eta\over2}E_t.
\end{align*}
This proves \eqref{eq:first_step_progress} and completes the proof.
\end{proof}

\subsection{Initialization and a recursive argument for convergence}\label{app:sec_initial_fixed}

Below we show that with sufficiently many uniformly sampled particles, each atom of
$\wh G_r$ has a nearby initial particle
with high probability. Then an exact reweighting step in \eqref{def_Mn} leads to the initial $G_0$ that has a small likelihood gap.

\begin{lemma}
\label{lem:app_particle_fixed_initialization}
For every $\delta\in(0,1)$, there is a finite $M_0$, depending
only on the data, $r$, and $\delta$, such that the following holds
for every $M\ge M_0$. If $\beta_1^0,\ldots,\beta_M^0$ are
independent and uniformly distributed on $\cB_r$, and
$w^0\in\cM_n(\bbeta^0)$, then with probability
at least $1-\delta$, the mixing distribution $G_0 = \sum_j w_j^{0} \delta_{\beta_j^0}$ satisfies
$\ell_n(\wh G_r)-\ell_n(G_0)\le b_n/2$ .
\end{lemma}
\begin{proof}
Recall the representation
$$
    \bm{p}_{\wh G_r}=\sum_{\ell=1}^{L^*}a_\ell^*\bm{p}_{\delta_{\beta_\ell^*}}
$$
in \eqref{eq:caratheodory_npmle_prediction} with some $L^*\le n+1$.  Set
$$
    h_n=r\wedge{2\rho_{n,r}b_n\over R_{x,n}}.
$$
Suppose that for every $\ell\in[L^*]$ an initial particle
$\wt\beta_\ell$ satisfies
$\|\wt\beta_\ell-\beta_\ell^*\|_2\le h_n$. Let
$\wt G=\sum_{\ell=1}^{L^*}a_\ell^*\delta_{\wt\beta_\ell}$, with prediction
vector $\bm{p}_{\wt G}=\sum_{\ell=1}^{L^*}a_\ell^*\bm{p}_{\delta_{\wt\beta_\ell}}$.
We first show that
$$
     \ell_n(\wh G_r)-\ell_n(\wt G) \le b_n/2.
$$
To this end, apply \eqref{eq:app_fixed_data_ell_bounds_upper} to
$\ell_n(\wh G_r)-\ell_n(\wt G)$ to obtain
\begin{align*}
    \ell_n(\wh G_r)-\ell_n(\wt G)
    &~ \le~ \langle \bm{r}(\bm{p}_{\wt G}),\bm{p}_{\wh G_r}-\bm{p}_{\wt G}\rangle_n ~ \overset{\eqref{bd_r_p_uniform}}{\le}~ {\|\bm{p}_{\wt G}-\bm{p}_{\wh G_r}\|_n\over\rho_{n,r}}.
\end{align*}
Since
\begin{align*}
    \|\bm{p}_{\wt G}-\bm{p}_{\wh G_r}\|_n
    &\le\sum_{\ell=1}^{L^*}a_\ell^*
       \|\bm{p}_{\delta_{\wt\beta_\ell}}-\bm{p}_{\delta_{\beta_\ell^*}}\|_n \overset{\eqref{eq:app_fixed_data_smoothness}}{\le}{R_{x,n}\over4}\sum_{\ell=1}^{L^*}a_\ell^*
       \|\wt\beta_\ell-\beta_\ell^*\|_2
    \le{R_{x,n}h_n\over4},
\end{align*}
we conclude
$$
    \ell_n(\wh G_r)-\ell_n(\wt G) ~ \le~ {R_{x,n}h_n\over4\rho_{n,r}}
    \le{b_n\over2}.
$$
Since $\wt G$ is supported on the initial particle locations and
$w^0\in\cM_n(\bbeta^0)$, we have $\ell_n(G_0)\ge\ell_n(\wt G)$. Therefore,
$$
     \ell_n(\wh G_r)-\ell_n(G_0) ~ \le~    \ell_n(\wh G_r)-\ell_n(\wt G) ~ \le~   {b_n\over2}.
$$
as required.

It remains to bound the probability that each atom
$\beta_\ell^*$ has an initial particle within distance $h_n$.
Since $h_n\le r$, the ball of radius $h_n/2$ centered at
$(1-h_n/(2r))\beta_\ell^*$ lies in $\cB_r$ and every point
in it is within distance $h_n$ of $\beta_\ell^*$.
Indeed, for any $\beta$ in this ball,
$$
    \|\beta\|_2 \le \left\|\beta - \left(1-{h_n \over 2r}\right)\beta_\ell^*\right\|_2 + \left(1-{h_n \over 2r}\right)\|\beta_\ell^*\|_2 \le {h_n\over 2} + r-{h_n\over 2} = r
$$
and
$$
    \|\beta-\beta_\ell^*\|_2 \le \left\|\beta - \left(1-{h_n \over 2r}\right)\beta_\ell^*\right\|_2 + \left\| \left(1-{h_n \over 2r}\right)\beta_\ell^* -  \beta_\ell^* \right\|_2  \le {h_n}.
$$
A uniform sample $\wt \beta \in \cB_r$ therefore falls within distance $h_n$
of $\beta_\ell^*$ with probability at least $(h_n/(2r))^d$.
By independence among the random draws $\wt \beta_1,\ldots, \wt \beta_M$ and taking the union bound over $\ell \in [L^*]$, one can deduce that
$$
    \PP\left\{\ell_n(\wh G_r)-\ell_n(G_0)>{b_n\over2}\right\}
    \le L^*\left\{1-\left({h_n\over2r}\right)^d\right\}^M
    \le L^*\exp\left\{-M\left({h_n\over2r}\right)^d\right\}.
$$
The claim follows by taking
$M_0=\lceil(2r/h_n)^d\log(L^*/\delta)\rceil$.
Since $b_n$ depends only on the data and $r$, so does $h_n$.
\end{proof}

We use the following elementary recursion to obtain the convergence rate.

\begin{lemma}
\label{lem:app_particle_quadratic_recursion}
Let $\{a_t\}_{t\ge0}$ be a nonnegative sequence such that, for
some $\kappa>0$, $\kappa a_0\le1/2$ and
$a_{t+1}\le a_t-\kappa a_t^2$ for every $t\ge0$. Then
\begin{equation}
    a_t\le{1\over a_0^{-1}+\kappa t},
    \qquad t\ge0,
    \label{eq:app_quadratic_recursion}
\end{equation}
where the right-hand side is interpreted as zero if $a_0=0$.
\end{lemma}
\begin{proof}

The sequence $\{a_t\}_{t\ge 0}$ is nonincreasing, so $\kappa a_t\le1/2$ for every $t$.
If $a_t=0$, all later terms are zero. Otherwise, as long as
$a_{t+1}>0$,
$$
    {1\over a_{t+1}}
    \ge{1\over a_t(1-\kappa a_t)}
    \ge{1+\kappa a_t\over a_t}
    ={1\over a_t}+\kappa.
$$
Iterating proves \eqref{eq:app_quadratic_recursion}.
\end{proof}

\subsection{Proof of \cref{thm:reweighted_particle_global}}
\label{app:particle_fixed_data_proof}

\begin{proof}

Let $M\ge M_0$, with $M_0$ given in
\cref{lem:app_particle_fixed_initialization}. Work on the event
of that lemma, which has probability at least $1-\delta$.
By \cref{lem:particle_one_block_ascent}, the likelihood is
nondecreasing, so
$\ell_n(\wh G_r)-\ell_n(G_t)\le b_n/2$ for every $t\ge0$.
Since $w^t\in\cM_n(\bbeta^t)$ for every $t$, one has
\begin{align*}
    \ell_n(\wh G_r)-\ell_n(G_{t+1})
    &\le\ell_n(\wh G_r)-\ell_n(G_t)-{\eta\over2}E_t &&\text{by \eqref{eq:first_step_progress}}\\
    &\le\ell_n(\wh G_r)-\ell_n(G_t)
       -{\eta \over2} c_n
        \{\ell_n(\wh G_r)-\ell_n(G_t)\}^2 &&\text{by \eqref{eq:app_energy_target}}.
\end{align*}
This inequality also holds when $\ell_n(\wh G_r)-\ell_n(G_t) = 0$.
By decreasing $c_n$ if necessary so that $\eta c_nb_n\le2$, applying \cref{lem:app_particle_quadratic_recursion} with
$a_t=\ell_n(\wh G_r)-\ell_n(G_t)$  and  $\kappa=\eta c_n/2$
gives
$$
    \ell_n(\wh G_r)-\ell_n(G_t)
    \le{1\over2/b_n+(\eta c_n/2)t},
    \qquad t\ge0.
$$
This proves \eqref{eq:particle_global_rate} with
$C=2/(\eta c_n)$ and $c=4/(\eta c_nb_n)$.
Both constants depend only on the fixed data, $r$, and $\eta$,
not on $M$, $H$, or the realized initialization.
Finally, as $t\to \infty$, $\ell_n(G_t)\to\ell_n(\wh G_r)
=\sup_{G\in\cP(\cB_r)}\ell_n(G)$, and
\eqref{eq:app_empirical_gap_distance} gives
$\|\bm{p}_{G_t}-\bm{p}_{\wh G_r}\|_n\to0$.
\end{proof}

\subsection{Energy lower bounds in \eqref{lb_energy}}
\label{app:particle_fixed_data_auxiliary}

We first prove the energy bound in \eqref{eq:app_energy_target} for mixing distributions with at most $n+1$
particles. \cref{lem:app_particle_stationary_values} shows that the stationary likelihood values of their stationary points form a finite set, so there is a positive  gap between the maximal likelihood and the likelihood of suboptimal  stationary points. We use this fact and construct an analytic curve in \cref{lem:app_particle_projected_curve,lem:app_particle_fixed_order_energy} to derive the quantitative bound
in \eqref{eq:app_energy_target}. To allow arbitrary particle numbers,
we express the weights as a convex combination of weight vectors with
the same prediction and at most $n+1$ positive entries; see \cref{lem:app_particle_weight_fiber}. The preceding energy bound then yields the desired result with the constants $c_n$  independent of $M$.

\paragraph{Analytic functions and curve selection}
\label{app:particle_subanalytic_prelim}

We recall some definitions from real analytic geometry. For background
on analytic functions and subanalytic sets, see
\citet[Sections 2.2 and 6.5]{krantz2002primer}.

A real function is analytic if it agrees locally with a convergent power
series. Here \emph{locally} means on some open neighborhood of each point;
the neighborhood may depend on the point. A set is semianalytic if, in
such neighborhoods, it is a finite union of sets described by finitely
many analytic equalities and inequalities. These are conditions of the
form $f(x)=0$ or $f(x)>0$, where $f$ is analytic.

A set is subanalytic if, locally, it has the form
$$
    \{z:\text{there exists }u\text{ such that }(z,u)\in B\},
$$
where $B$ is semianalytic and has compact closure within its defining
domain. Thus we keep a point $z$ if the conditions defining $B$ can
be satisfied by some $u$.

A set is semialgebraic if it is a finite union of sets described by
finitely many polynomial equalities and inequalities
\citep[Definition 1.1]{bierstone1988semianalytic}. Since polynomials are analytic,
every semialgebraic set is also semianalytic. Examples include the simplex,
the ball, and the set of pairs $(x,y)$ satisfying
$y=\Pi_{\cB_r}(x)$.\\

To describe smooth pieces of these sets, we use the notion of a
submanifold, which generalizes smooth curves and surfaces to other
dimensions. A $q$-dimensional $C^1$ submanifold of $\RR^N$ is a
subset with the following local description near each point: after
reordering coordinates if necessary, $q$ coordinates range over an open
set, and the remaining $N-q$ coordinates are continuously differentiable
functions of them. For example, a circle is a one-dimensional submanifold
of $\RR^2$, and a sphere is a two-dimensional submanifold of $\RR^3$.
A subset of $\RR^N$ is connected if it is not the union of two disjoint
nonempty sets that are open relative to the subset. For a $C^1$
submanifold, connectedness is equivalent to path-connectedness, and any
two points of a connected $C^1$ submanifold can be joined by a piecewise
$C^1$ curve in it; this is the property used in the proof of
\cref{lem:app_particle_stationary_values}.

The following two lemmas are the results from real analytic geometry
that we use. The first follows from
\citet[Theorem 4.1]{bierstone1996composite}, applied to the identity map,
and the second is the subanalytic case of
\citet[1.17]{vandenDriesMiller1996geometric}.

\begin{lemma}
\label{lem:app_subanalytic_stratification}
Let $S\subset\RR^N$ be a compact subanalytic set, and let
$A_1,\ldots,A_k$ be subanalytic subsets of $S$. Then $S$ admits a
finite partition into connected subanalytic $C^1$ submanifolds such that
each $A_i$ is a union of these pieces. In particular,
$S$ has finitely many connected components.
\end{lemma}

\begin{lemma}
\label{lem:app_subanalytic_curve_selection}
Let $S\subset\RR^N$ be subanalytic and let
$z\in\overline S\setminus S$. There exist $\epsilon>0$ and a real
analytic map $\gamma:(-\epsilon,\epsilon)\to\RR^N$ such that
$$
    \gamma(0)=z,\qquad \gamma(s)\in S,\quad 0<s<\epsilon.
$$
\end{lemma}

\Cref{lem:app_subanalytic_stratification} allows a possibly nonsmooth set
to be studied one smooth piece at a time. The partition also respects the
subsets $A_i$: each piece lies entirely inside or entirely outside each
of them. For example, if $S=[0,1]$ and $A_1=[0,1/2]$, we can use
the partition
$$
    \{0\},\quad (0,1/2),\quad \{1/2\},\quad (1/2,1),\quad \{1\}.
$$
The two open intervals are one-dimensional pieces, whereas the three
single points are zero-dimensional pieces. The first three pieces together
form $A_1$.

\Cref{lem:app_subanalytic_curve_selection} gives more than a convergent
sequence: it gives an analytic curve approaching $z$ from any point within
$S$. For example, take
$S=\{(x,x^{3/2}):0<x<1\}$ and $z=(0,0)$. The curve
$\gamma(s)=(s^2,s^3)$ lies in $S$ for $0<s<1$ and reaches
$z$ at zero. Both coordinates are polynomials in $s$, so the curve
is analytic even at its endpoint. If $F$ is analytic near $z$, then
$F(\gamma(s))=F(s^2,s^3)$ is a function of $s$ with a convergent
power series near $s=0$.

To compare how fast functions vanish along such a curve, we use the
first nonzero term of their power series. Suppose that $s\mapsto f(s)$ is
analytic on an open interval containing zero and is not identically zero.
Then, as $s\to0$,
$$
    f(s)=c s^m+\cO(s^{m+1}),\qquad c\ne0,
$$
for some integer $m\ge0$. The integer $m$ is called the
\emph{order of vanishing} of $f$ at zero
\citep[Section~1, Notation]{rainer2013perturbation}. If $m\ge1$,
differentiation gives $f'(s)=mc s^{m-1}+\cO(s^m)$, so this exponent also
determines the leading order of the derivative.

\subsubsection{Characterization of the stationary points to the likelihood}

For any $L\ge1$, consider
$(a,\bbeta)\in\Delta^{L-1}\times\cB_r^L$, where
$\bbeta=(\beta_1,\ldots,\beta_L)$,
$G=\sum_{j=1}^La_j\delta_{\beta_j}$, and
$\bm{p}_G=\sum_{j=1}^La_j\bm{p}_{\delta_{\beta_j}}$.
Define its energy by
$$
    E(a,\bbeta):={1\over\eta^2}\sum_{j=1}^La_j
       \bigl\|\Pi_{\cB_r}\{\beta_j+\eta\nabla g_{n,G}(\beta_j)\}
                     -\beta_j\bigr\|_2^2.
$$
In particular, $E(w^t,\bbeta^t)=E_t$ by
\eqref{eq:reweighted_particle_update} and \eqref{eq:particle_gap_energy}.

Denote the set of exactly reweighted configurations by
$$
    \cR_L=\{(a,\bbeta)\in\Delta^{L-1}\times\cB_r^L:
       a\in\cM_n(\bbeta)\}.
$$
By concavity and \cref{lem:app_particle_weight_kkt}, these are precisely
the configurations satisfying
\begin{equation}
    g_{n,G}(\beta_j)\le\bar g_{n,G},
    \qquad
    a_j\{g_{n,G}(\beta_j)-\bar g_{n,G}\}=0,
    \qquad j\in[L].
    \label{eq:app_reweighted_set}
\end{equation}
Let $\cS_L\subset\cR_L$ be the set of stationary configurations, defined by
\begin{equation}
   \cS_L := \left\{ (a,\bbeta) \in \cR_L:  E(a,\bbeta)=0\right\}.
    \label{eq:app_stationary_set}
\end{equation}
Since the summand defining $E(a,\bbeta)$ are nonnegative, this condition
is equivalent to
$$
    \Pi_{\cB_r}\{\beta_j+\eta\nabla g_{n,G}(\beta_j)\}=\beta_j
    \qquad\text{for every }j\text{ with }a_j>0.
$$
By \eqref{eq:app_normal_cone}, it is also equivalent to
$\nabla g_{n,G}(\beta_j)\in N_{\cB_r}(\beta_j)$ for all such $j$.
In particular, $\cS_L$ does not depend on $\eta>0$.

Both sets of $\cS_L$ and $\cR_L$ are compact: they are defined by continuous equalities and
weak inequalities on $\Delta^{L-1}\times\cB_r^L$.
They are also subanalytic. Indeed, the log-likelihood and residuals
are analytic functions of the prediction coordinates on a neighborhood of
$[\rho_{n,r},1-\rho_{n,r}]^n$, and $\sigma$ is analytic on a neighborhood
of the compact range of $x_i^\T\beta$. Thus $\bm{p}_G$, $\ell_n(G)$,
$g_{n,G}$, and $\nabla g_{n,G}$ are analytic in the configuration
variables on a neighborhood of the feasible set.
The simplex, the ball, and the graph of $\Pi_{\cB_r}$ are
semialgebraic. Including the projected locations as additional
variables expresses $\cR_L$ and $\cS_L$ as projections of compact
semianalytic sets, which proves the claim.

The following lemma shows that the likelihood takes only finitely many distinct values on $\cS_L$ whenever $L\le n+1$.
\begin{lemma}
\label{lem:app_particle_stationary_values}
For every $1\le L\le n+1$, the set
$\{\ell_n(\sum_{j=1}^La_j\delta_{\beta_j}):(a,\bbeta)\in\cS_L\}$
is finite.
\end{lemma}

\begin{proof}
We show that $\ell_n$ is constant on each of finitely many connected
pieces of $\cS_L$. Apply \cref{lem:app_subanalytic_stratification}
to the compact subanalytic set $\cS_L$ and its subanalytic subsets
$\{a_j=0\}$ and $\{\|\beta_j\|_2=r\}$, $j\in[L]$.
On each resulting piece, every weight is either positive throughout
or identically zero, and every location lies either in the interior
of $\cB_r$ throughout or on its boundary throughout.
Any two points of a connected $C^1$ submanifold can be joined by a
piecewise $C^1$ curve in it. It suffices to show that the likelihood
has derivative zero along every $C^1$ curve in a piece.

Let $s\mapsto(a(s),\bbeta(s))$ be such a curve and write
$G(s)=\sum_ja_j(s)\delta_{\beta_j(s)}$. Dots denote derivatives in $s$.
By the chain rule,
\begin{equation}
    {d\over ds}\ell_n(G(s))
    =\sum_{j=1}^Lg_{n,G}(\beta_j)\dot a_j
     +\sum_{j=1}^La_j
        \langle\nabla g_{n,G}(\beta_j),\dot\beta_j\rangle,
    \label{eq:app_stationary_curve_derivative}
\end{equation}
where $G=G(s)$, $a_j=a_j(s)$, and $\beta_j=\beta_j(s)$.
By \eqref{eq:app_reweighted_set}, $g_{n,G}(\beta_j)=\bar g_{n,G}$
when $a_j>0$. Otherwise, $a_j$ is identically zero on the piece
and $\dot a_j=0$. Hence
$$
    \sum_{j=1}^Lg_{n,G}(\beta_j)\dot a_j
    =\bar g_{n,G}\sum_{j:a_j>0}\dot a_j
    =\bar g_{n,G}\,{d\over ds}\sum_{j=1}^La_j(s)=0.
$$
For the second sum in \eqref{eq:app_stationary_curve_derivative},
fix $j$ with $a_j>0$. Stationarity gives
$\nabla g_{n,G}(\beta_j)=0$ if $\beta_j$ is in the interior.
If $\|\beta_j\|_2=r$, it gives
$\nabla g_{n,G}(\beta_j)=\lambda_j\beta_j$ for some
$\lambda_j\ge0$. The curve stays on the boundary in this case, so
$$
    2\langle\beta_j,\dot\beta_j\rangle
    ={d\over ds}\|\beta_j(s)\|_2^2=0.
$$
Thus $\langle\nabla g_{n,G}(\beta_j),\dot\beta_j\rangle=0$
in both cases. The likelihood is constant on each piece, and there
are only finitely many pieces.
\end{proof}

The optimal value is one of these stationary values. To see this,
we recall
$\bm{p}_{\wh G_r}=\sum_{\ell=1}^{L^*}a_\ell^*\bm{p}_{\delta_{\beta_\ell^*}}$
in \eqref{eq:caratheodory_npmle_prediction}, with $L^*\le n+1$.
Optimality and \eqref{eq:potential_directional_derivative} give
$g_{n,\wh G_r}(\beta)\le\bar g_{n,\wh G_r}$ for every
$\beta\in\cB_r$. Since $a^*$ maximizes the likelihood for the fixed
locations $\beta_1^*,\ldots,\beta_{L^*}^*$,
\eqref{eq:app_particle_weight_kkt} gives equality at every
$\beta_\ell^*$ with $a_\ell^*>0$. Each such atom therefore maximizes
$g_{n,\wh G_r}$ on $\cB_r$, and its gradient belongs to
$N_{\cB_r}(\beta_\ell^*)$. This representation belongs to
$\cS_{L^*}$.

Since $L$ takes only finitely many values,
\cref{lem:app_particle_stationary_values} implies that
\begin{equation}
\begin{aligned}
    \Gamma_n:=\min\bigl\{\ell_n(\wh G_r)-\ell_n(G):\
    &(a,\bbeta)\in\cS_L,\ 1\le L\le n+1,\\
    &G=\textstyle\sum_{j=1}^La_j\delta_{\beta_j},\
      \ell_n(G)<\ell_n(\wh G_r)\bigr\}
\end{aligned}
\label{eq:app_empirical_stationary_gap}
\end{equation}
is strictly positive, with $\Gamma_n=+\infty$ if the set is empty.
Define
\begin{equation}
    b_n:=1\wedge{\Gamma_n\over2}
    \label{eq:app_fixed_data_constants}
\end{equation}
which depends only on the data and $r$. We conclude from \eqref{eq:app_empirical_stationary_gap} that every configuration in $\cS_L$, $1\le L\le n+1$, either attains the maximum likelihood or has a likelihood gap of at least $b_n$.

\subsubsection{Energy bound for at most $n+1$ particles}

The next lemma compares the likelihood derivative along an analytic
curve with the projected displacement of a particle. It also applies
when the weight tends to zero and the location approaches the boundary of $\cB_r$.

\begin{lemma}
\label{lem:app_particle_projected_curve}
Let $a:(-s_0,s_0)\to\RR$ and
$\beta,v:(-s_0,s_0)\to\RR^d$ be real analytic, with
$a(s)\ge0$ and $\beta(s)\in\cB_r$ for $0\le s<s_0$. Write
$$
    d(s)=\Pi_{\cB_r}\{\beta(s)+\eta v(s)\}-\beta(s).
$$
Suppose either $\|\beta(s)+\eta v(s)\|_2\le r$ for every
$0<s<s_0$, or $\|\beta(s)+\eta v(s)\|_2\ge r$ for every
$0<s<s_0$. If $a(s)\|d(s)\|_2^2=o(s^{2m})$ as $s\downarrow0$
for some integer $m\ge1$, then
$$
    a(s)\langle v(s),\dot\beta(s)\rangle=o(s^{m-1}),
    \qquad s\downarrow0.
$$
\end{lemma}

\begin{proof}
All limits in this proof are as $s\downarrow0$. Analyticity at zero
ensures that $v$ and $\dot\beta$ are bounded for sufficiently small $s$.
We first consider $\|\beta+\eta v\|_2\le r$, in which case
$d=\eta v$. If $a$ or $v$ is identically zero, the claim is immediate.
Otherwise, let $\alpha$ be the order of vanishing of $a$ and
$\gamma$ the smallest such order among the nonzero coordinates of $v$.
Then $a(s)\asymp s^\alpha$ and $\|v(s)\|_2\asymp s^\gamma$.
The hypothesis gives $\alpha+2\gamma>2m$. Since $\alpha\ge0$,
we also have $\alpha+\gamma>m$, and therefore
$$
    |a\langle v,\dot\beta\rangle|
    \le a\|v\|_2\|\dot\beta\|_2
    =\cO(s^{\alpha+\gamma})=o(s^m).
$$

Now suppose $\|\beta+\eta v\|_2\ge r$. The projection formula gives
$$
    \beta+d={r(\beta+\eta v)\over\|\beta+\eta v\|_2},
    \qquad
    \eta v=d+\lambda(\beta+d),
    \qquad
    \lambda={\|\beta+\eta v\|_2-r\over r}.
$$
Since $\|\beta(0)+\eta v(0)\|_2\ge r>0$, these formulas extend
$d$ and $\lambda$ analytically to an open interval containing zero.
Moreover,
$$
    0\le\lambda
    \le{\|\beta\|_2+\eta\|v\|_2-r\over r}
    \le{\eta\|v\|_2\over r},
$$
so $\lambda$ is bounded. Let $s_r(s)=r-\|\beta(s)\|_2$.
Since $\|\beta+d\|_2=r$, the reverse triangle inequality gives
$0\le s_r\le\|d\|_2$.

If $s_r(0)>0$, then $\|d(s)\|_2\ge s_r(0)/2$ for sufficiently
small $s$. The hypothesis implies $a=o(s^{2m})$, and the
conclusion follows from the boundedness of $v$ and $\dot\beta$.
It remains to consider $\|\beta(0)\|_2=r$. In this case $s_r$
is analytic on an open interval containing zero.
If $a$ is identically zero, the claim is immediate.
If $d$ is identically zero, then $\|\beta(s)\|_2=r$ and
$\eta v=\lambda\beta$, which gives
$\langle v,\dot\beta\rangle
=\lambda(d/ds)\|\beta(s)\|_2^2/(2\eta)=0$.

Otherwise, let $\alpha$ be the order of $a$ and $\gamma$ the
smallest order among the nonzero coordinates of $d$. As before,
$\alpha+2\gamma>2m$ and $\alpha+\gamma>m$.
From $0\le s_r\le\|d\|_2=\cO(s^\gamma)$ and the power series
of $s_r$, we obtain $\dot s_r=\cO(s^{\gamma-1})$. Consequently,
$$
    a\langle d,\dot\beta\rangle
    =\cO(s^{\alpha+\gamma})=o(s^m),
    \qquad
    a\dot s_r=\cO(s^{\alpha+\gamma-1})=o(s^{m-1}).
$$
Since $\langle\beta,\dot\beta\rangle=-(r-s_r)\dot s_r$,
\begin{align*}
    \eta a\langle v,\dot\beta\rangle
    &=a\langle d,\dot\beta\rangle
      +\lambda a\langle\beta+d,\dot\beta\rangle\\
    &=(1+\lambda)a\langle d,\dot\beta\rangle
      -\lambda(r-s_r)a\dot s_r
    =o(s^{m-1}).
\end{align*}
This proves the claim.
\end{proof}

We now prove the energy bound for the case $L\le n+1$.

\begin{lemma}
\label{lem:app_particle_fixed_order_energy}
There exists $c_n>0$, depending only on the data, $r$, and $\eta$,
such that, for every $1\le L\le n+1$ and every
$(a,\bbeta)\in\cR_L$ with $G=\sum_{j=1}^La_j\delta_{\beta_j}$ and
$0<\ell_n(\wh G_r)-\ell_n(G)\le b_n$,
\begin{equation}
    E(a,\bbeta)
    \ge c_n\{\ell_n(\wh G_r)-\ell_n(G)\}^2.
    \label{eq:app_fixed_order_energy}
\end{equation}
\end{lemma}
\begin{proof}
We argue by contradiction. If the bound failed, compactness and
the separation of stationary values would give a sequence approaching
an optimal configuration, with energy smaller than any fixed multiple
of the squared gap. Curve selection gives an analytic curve with the
same property. Along this curve, we will obtain two incompatible
orders of vanishing for the derivative of the likelihood gap.

Since $L$ takes only finitely many values, failure of a uniform
constant gives a fixed $L$ and a sequence
$(a_k,\bbeta_k)\in\cR_L$, with $G_k=\sum_{j=1}^La_{k,j}\delta_{\beta_{k,j}}$, such that
$0<\ell_n(\wh G_r)-\ell_n(G_k)\le b_n$ and
\begin{equation}
    {E(a_k,\bbeta_k)
     \over\{\ell_n(\wh G_r)-\ell_n(G_k)\}^2}
    \le {1\over k}.
    \label{eq:app_energy_ratio_zero}
\end{equation}
By compactness, a subsequence converges to
$(a_\infty,\bbeta_\infty)\in\cR_L$.
The numerator in \eqref{eq:app_energy_ratio_zero} is at most
$b_n^2/k$ and tends to zero as $k\to \i$. Continuity therefore implies
$(a_\infty,\bbeta_\infty)\in\cS_L$ using \eqref{eq:app_stationary_set}.
Its likelihood gap lies in $[0,b_n]$, and hence is zero by
\eqref{eq:app_empirical_stationary_gap} and
\eqref{eq:app_fixed_data_constants}. In particular,
$\ell_n(G_k)\to\ell_n(\wh G_r)$.

Let $\cA$ be the set of $(a,\bbeta,u)\in\cR_L\times[0,1]$
such that, with $G=\sum_ja_j\delta_{\beta_j}$,
\begin{gather*}
    0<\ell_n(\wh G_r)-\ell_n(G)\le b_n,\\
    E(a,\bbeta)
    =u\{\ell_n(\wh G_r)-\ell_n(G)\}^2.
\end{gather*}
This set is bounded and subanalytic, by the same analytic and
semialgebraic representations used for $\cR_L$ and $\cS_L$.
Taking $u$ equal to the ratio in \eqref{eq:app_energy_ratio_zero}
shows that $(a_\infty,\bbeta_\infty,0)$ lies in
$\overline{\cA}\setminus\cA$.

Partition $\cA$ according to whether each $a_j$ is zero or positive,
and whether each
$\|\beta_j+\eta\nabla g_{n,G}(\beta_j)\|_2$ is at most $r$ or
larger than $r$. These finitely many subsets are subanalytic, and
one contains a subsequence converging to
$(a_\infty,\bbeta_\infty,0)$.
Apply \cref{lem:app_subanalytic_curve_selection} to this subset.
There is an analytic curve
$s\mapsto(a(s),\bbeta(s),u(s))$ on $(-s_0,s_0)$ that equals
$(a_\infty,\bbeta_\infty,0)$ at zero and lies in that subset for
$0<s<s_0$. Write $G(s)=\sum_ja_j(s)\delta_{\beta_j(s)}$ and
$\bm{p}_{G(s)}=\sum_ja_j(s)\bm{p}_{\delta_{\beta_j(s)}}$.
Since $\bm{p}_{G(0)}\in[\rho_{n,r},1-\rho_{n,r}]^n$, the formula
\eqref{def_llh} extends $\ell_n(G(s))$ analytically to an open interval
containing zero.
Its gap is zero at $s=0$ and positive for $s>0$. Hence, for
some $A>0$ and integer $m\ge1$,
\begin{equation}
\begin{aligned}
    \ell_n(\wh G_r)-\ell_n(G(s))
    &=As^m+\cO(s^{m+1}),\\
    E(a(s),\bbeta(s))
    &=u(s)\{\ell_n(\wh G_r)-\ell_n(G(s))\}^2
      =o(s^{2m}).
\end{aligned}
\label{eq:app_selected_analytic_curve}
\end{equation}
Let
$$
    d_j(s)=\Pi_{\cB_r}\{\beta_j(s)+\eta\nabla g_{n,G(s)}(\beta_j(s))\}
                   -\beta_j(s).
$$
Since $\eta$ is fixed, \eqref{eq:app_selected_analytic_curve} gives
$$
    0\le a_j(s)\|d_j(s)\|_2^2
       \le\sum_{k=1}^La_k(s)\|d_k(s)\|_2^2
       =\eta^2E(a(s),\bbeta(s))=o(s^{2m})
$$
for every $j$.

We next differentiate the likelihood along the curve.
Each weight is either positive throughout $(0,s_0)$ or
identically zero there. By \eqref{eq:app_reweighted_set},
$$
    \sum_{j=1}^L\dot a_j(s)g_{n,G(s)}(\beta_j(s))
    =\bar g_{n,G(s)}\sum_{j:a_j(s)>0}\dot a_j(s)
    =\bar g_{n,G(s)}\,{d\over ds}\sum_{j=1}^La_j(s)=0.
$$
The chain rule therefore gives
\begin{equation}
    {d\over ds}\ell_n(G(s))
    =\sum_{j=1}^La_j(s)
       \langle\nabla g_{n,G(s)}(\beta_j(s)),\dot\beta_j(s)\rangle.
    \label{eq:app_energy_curve_derivative}
\end{equation}
For each $j$, the functions $a_j$, $\beta_j$, and
$s\mapsto\nabla g_{n,G(s)}(\beta_j(s))$ are analytic at zero.
The partition ensures that the same case of the projection formula
applies for all $0<s<s_0$. Thus
\cref{lem:app_particle_projected_curve} applies and shows that
each summand in \eqref{eq:app_energy_curve_derivative} is
$o(s^{m-1})$. On the other hand, differentiating the first line
of \eqref{eq:app_selected_analytic_curve} gives
$$
    {d\over ds}\ell_n(G(s))=-mAs^{m-1}+\cO(s^m),
$$
with $mA>0$, a contradiction. This proves the existence of one
$c_n>0$ for all $L\le n+1$. The sets and functions used in the
argument depend only on the fixed data, $r$, and $\eta$.
\end{proof}

\subsubsection{Energy bounds for an arbitrary number of particles}

For fixed particle locations $\bbeta$, weight vectors producing the same prediction vector yield the same potential $g_{n,G}$ and hence the same updated particle locations. The particle contributions to $E(w,\bbeta)$ are therefore fixed on this set of weights, so $E(w,\bbeta)$ is linear in $w$. We use this fact to extend
\eqref{eq:app_fixed_order_energy} to any particle count.

\begin{lemma}
\label{lem:app_particle_weight_fiber}
Fix $\bbeta\in\cB_r^M$ and $w\in\cM_n(\bbeta)$.
Write $G=\sum_jw_j\delta_{\beta_j}$ and
$\bm{p}_G=\sum_jw_j\bm{p}_{\delta_{\beta_j}}$.
Then $w=\sum_s\lambda_su^{(s)}$ is a convex combination of
the extreme points of
$$
    \left\{u\in\Delta^{M-1}:\textstyle\sum_{j=1}^Mu_j\bm{p}_{\delta_{\beta_j}}=\bm{p}_G\right\}.
$$
Each $u^{(s)}$ has at most $n+1$ positive entries and belongs to
$\cM_n(\bbeta)$. Moreover,
\begin{equation}
    E(w,\bbeta)=\sum_s\lambda_s E(u^{(s)},\bbeta).
    \label{eq:app_weight_fiber_energy}
\end{equation}
\end{lemma}

\begin{proof}
The set of weights in the statement is a nonempty compact polytope,
so $w$ is a convex combination of its extreme points.
Suppose an extreme point $u$ has more than $n+1$ positive entries.
The corresponding vectors $(1,\bm{p}_{\delta_{\beta_j}}^\T)^\T\in\RR^{n+1}$ are
linearly dependent. Hence there is a nonzero vector $v$, with
$v_j=0$ whenever $u_j=0$, such that
$$
    \sum_{j=1}^Mv_j=0,
    \qquad
    \sum_{j=1}^Mv_j\bm{p}_{\delta_{\beta_j}}=0.
$$
For sufficiently small $\epsilon>0$, both $u+\epsilon v$ and
$u-\epsilon v$ have nonnegative entries and belong to the polytope.
Their average is $u$, contradicting extremality.

Every point in the polytope gives the same prediction $\bm{p}_G$ as $w$,
and hence the same likelihood $\ell_n(G)$ and potential $g_{n,G}$.
Since this likelihood is maximal
over $\Delta^{M-1}$ for the fixed locations, each $u^{(s)}$ belongs to
$\cM_n(\bbeta)$. Finally, substituting
$w_j=\sum_s\lambda_su_j^{(s)}$ and interchanging the finite sums
proves \eqref{eq:app_weight_fiber_energy}.
\end{proof}

For any $G=\sum_{j=1}^Mw_j\delta_{\beta_j}$ with
$w\in\cM_n(\bbeta)$ and $0<\ell_n(\wh G_r)-\ell_n(G)\le b_n$,
apply \cref{lem:app_particle_weight_fiber} to write
$w=\sum_s\lambda_su^{(s)}$.
Removing the zero entries of each $u^{(s)}$ and the corresponding
particle locations changes neither its prediction vector nor its energy.
Each resulting configuration remains exactly
reweighted on its support and has at most $n+1$ particles. Therefore,
\begin{align*}
    E(w,\bbeta)
    &=\sum_s\lambda_s E(u^{(s)},\bbeta) &&\text{by  \cref{lem:app_particle_weight_fiber}}\\
    &\ge\sum_s\lambda_s c_n
       \{\ell_n(\wh G_r)-\ell_n(G)\}^2 &&\text{by   \cref{lem:app_particle_fixed_order_energy}}\\
    &=c_n\{\ell_n(\wh G_r)-\ell_n(G)\}^2.
\end{align*}
Thus \eqref{eq:app_energy_target} holds with the same $b_n,c_n$
for every $M$.

\section{Identifiability and the denoising bound}\label{app_sec_ident}

\subsection{Proof of \cref{prop:identifiability} and a key lemma}
\label{app:particle_identifiability}

We first introduce notation for signed measures and then prove the
auxiliary result used in \cref{prop:identifiability}.

\paragraph{Signed measures}
For background on signed measures, see
\citet[Section~IX.3]{knapp2016basic}.

A finite signed measure may assign negative as well as positive mass to
a measurable set. For example, if $P$ and $Q$ are probability
measures on $\RR^d$, then $\mu=P-Q$ satisfies
$\mu(A)=P(A)-Q(A)$.
Every finite signed measure has a unique Jordan decomposition
$$
    \mu=\mu^+-\mu^-,\qquad \mu^+\perp\mu^-,
$$
where $\mu^+$ and $\mu^-$ are nonnegative finite measures
\citep[Theorems~9.14--9.15]{knapp2016basic}. Here $\mu^+\perp\mu^-$
means that the two measures are concentrated on disjoint measurable sets.
We use the total variation norm
$$
    \|\mu\|_{\mathrm{TV}}=\mu^+(\RR^d)+\mu^-(\RR^d).
$$
For probability measures $P$ and $Q$, it satisfies
$\|P-Q\|_{\mathrm{TV}}=2$ if and only if $P\perp Q$.
For example, $\|\delta_a-\delta_b\|_{\mathrm{TV}}=2$ whenever
$a\ne b$, however close the two points are.

The reflection of a signed measure $\mu$ is defined by
$\check\mu(A)=\mu(-A)$, where $-A:=\{-\beta:\beta\in A\}$.
In words, it moves the mass at each location $\beta$ to $-\beta$. Equivalently,
for every bounded measurable function $f$,
$$
    \int f(\beta)\,\rd \check\mu(\beta)
    =\int f(-\beta)\,\rd \mu(\beta).
$$

\paragraph{The sigmoid transform}
For a finite signed measure $\mu$ on $\cB_r$, extend the notation
for a mixture regression function by
$$
    p_\mu(x)=\int_{\cB_r}\sigma(x^\T\beta)\,\rd \mu(\beta).
$$

\begin{lemma}
\label{lem:app_sigmoid_kernel}
Let $\mu$ be a finite signed measure on $\cB_r$.
If $p_\mu$ vanishes on a nonempty open subset of $\RR^d$, then
\begin{equation}
    \mu(\cB_r)=0,
    \qquad
    \mu=\check\mu.
    \label{eq:app_sigmoid_kernel}
\end{equation}
Conversely, $p_\mu=0$ whenever $\mu=\check\mu$ and
$\mu(\cB_r)=0$.
\end{lemma}

\begin{proof}
We first show that $p_\mu$ is real analytic by constructing a convergent power series around each $x_0\in\RR^d$. Fix $x_0$ and set
$$
    \delta:=\frac{\pi}{2\sqrt d(r+1)}.
$$
For each $\beta\in\cB_r$, consider the complex-valued function
$z\mapsto\sigma(z^\T\beta)$. If
$\max_{j\in[d]}|z_j-x_{0,j}|\le\delta$, then, since $x_0$ and $\beta$ are real,
$$
    |\Im(z^\T\beta)|
    \le \|\Im z\|_2\|\beta\|_2
    \le \sqrt d\,\delta r
    <{\pi \over 2}.
$$
The poles of $\sigma(w)=(1+e^{-w})^{-1}$ are $(2k+1)\pi i$, $k\in\mathbb Z$,
so $\sigma$ is holomorphic on the strip $\{w\in\CC:|\Im w|<\pi\}$.
Consequently, $z\mapsto\sigma(z^\T\beta)$ is holomorphic on a neighborhood of
$$
    R_\delta:=\{z\in\CC^d:\max_j|z_j-x_{0,j}|\le\delta\}.
$$
We next bound its Taylor coefficients uniformly over $\beta\in\cB_r$
to justify integrating the series term by term.

For $w=a+ib$
with $|b|\le\pi/2$, since
$
    |1+e^{-w}|^2
    =1+2e^{-a}\cos b+e^{-2a}\ge1,
$
we have $|\sigma(w)|\le1$, and hence
$$
    |\sigma(z^\T\beta)|\le1
    \qquad
    \text{whenever }\beta\in\cB_r
    \text{ and } z\in R_\delta.
$$
For each $\beta$, expand around $x_0$:
$$
    \sigma(z^\T\beta)
    =\sum_{\alpha\in\mathbb N_0^d}
        c_\alpha(\beta)(z-x_0)^\alpha,
    \qquad
    c_\alpha(\beta)
    :=\frac{1}{\alpha!}
      \left.\partial_z^\alpha\sigma(z^\T\beta)\right|_{z=x_0}.
$$
Here, for a multi-index $\alpha=(\alpha_1,\ldots,\alpha_d)\in\mathbb N_0^d$,
we write $|\alpha|=\sum_j\alpha_j$, $\alpha!=\prod_j\alpha_j!$,
$z^\alpha=\prod_jz_j^{\alpha_j}$, and
$\partial_z^\alpha=\partial_{z_1}^{\alpha_1}\cdots\partial_{z_d}^{\alpha_d}$.
Cauchy's estimate on $R_\delta$ gives
$$
    |c_\alpha(\beta)|\le\delta^{-|\alpha|}
    \qquad\text{for every }\beta\in\cB_r.
$$
For any $z\in R_{\delta/2}$, we therefore
have
$$
    |c_\alpha(\beta)(z-x_0)^\alpha|
    \le2^{-|\alpha|}.
$$
Let $|\mu|$ be the total variation measure of $\mu$, and write
$\|\mu\|_{\mathrm{TV}}=|\mu|(\cB_r)$. Since this norm is finite,
\begin{align*}     \sum_{\alpha\in\mathbb N_0^d}
    \int_{\cB_r}|c_\alpha(\beta)(z-x_0)^\alpha|\,\rd |\mu|(\beta)
    &\le \|\mu\|_{\mathrm{TV}} \sum_{\alpha\in\mathbb N_0^d}2^{-|\alpha|}\\
    &= \|\mu\|_{\mathrm{TV}} \prod_{j=1}^d\left(\sum_{k=0}^\infty2^{-k}\right)\\
    &=  \|\mu\|_{\mathrm{TV}}\, 2^d<\infty.
\end{align*}
This justifies integrating the Taylor series term by term, so restricting
to real $x$ with $\max_j|x_j-x_{0,j}|<\delta/2$ gives
$$
    p_\mu(x)
    =\sum_{\alpha\in\mathbb N_0^d}
       \left(\int_{\cB_r}c_\alpha(\beta)\,\rd \mu(\beta)\right)
       (x-x_0)^\alpha.
$$
This is a convergent power series representation of $p_\mu$ around
$x_0$. Since $x_0$ was arbitrary, $p_\mu$ is real analytic on
$\RR^d$. Because $\RR^d$ is connected, the identity theorem for
real analytic functions implies that $p_\mu=0$ everywhere whenever
it vanishes on a nonempty open subset.

To prove \eqref{eq:app_sigmoid_kernel}, we first note that
$$
    p_{\mu}(0) = {1\over 2} \mu(\cB_r) = 0.
$$
It remains to prove $\mu=\check\mu$. Using
$\sigma(t)=\{1+\tanh(t/2)\}/2$, we have, for every $x\in\RR^d$,
\begin{equation}\label{eq_pmu}
\begin{aligned}
    p_\mu(x)
    &=\frac12\mu(\cB_r)
      +\frac12\int_{\cB_r}\tanh(x^\T\beta/2)\,\rd \mu(\beta) =\frac12\int_{\cB_r}\tanh(x^\T\beta/2)\,\rd \mu(\beta)=0.
\end{aligned}
\end{equation}
Let $\nu=\mu-\check\mu$. By oddness of $\tanh$,
$$
    \int_{\cB_r}\tanh(x^\T\beta/2)\,\rd \nu(\beta) = 2\int_{\cB_r}\tanh(x^\T\beta/2)\,\rd \mu(\beta) = 0.
$$
Fix $u\in\RR^d$ and set $x=tu$. Since
$\nu$ has finite total variation and compact support $\cB_r$, differentiation under
the integral gives, for $m\ge 0$,
$$
    0
    =\left.{d^{2m+1}\over dt^{2m+1}}
      \int_{\cB_r}\tanh(tu^\T\beta/2)\,\rd \nu(\beta)\right|_{t=0}
    ={\tanh^{(2m+1)}(0)\over2^{2m+1}}
      \int_{\cB_r}(u^\T\beta)^{2m+1}\,\rd \nu(\beta).
$$
To show that $\tanh^{(2m+1)}(0)\ne0$, write the Taylor series as
$\tanh t=\sum_{m\ge0}(-1)^m b_m t^{2m+1}$. The identity
$(\tanh t)'=1-(\tanh t)^2$ gives
$$
    b_0=1,\qquad
    (2m+1)b_m=\sum_{j=0}^{m-1}b_jb_{m-1-j},\quad m\ge1.
$$
Induction yields $b_m>0$. Consequently,
$\tanh^{(2m+1)}(0)=(-1)^m(2m+1)!b_m\ne0$, and
$$
    \int_{\cB_r}(u^\T\beta)^{2m+1}\,\rd \nu(\beta)=0,
    \qquad m\ge0.
$$
For any odd integer $k\ge1$, the multinomial expansion gives
$$
    0=\int_{\cB_r}(u^\T\beta)^k\,\rd \nu(\beta)
     =\sum_{|\alpha|=k}\frac{k!}{\alpha!}u^\alpha
        \int_{\cB_r}\beta^\alpha\,\rd \nu(\beta),
    \qquad u\in\RR^d.
$$
The polynomial on the right is identically zero, so each coefficient
vanishes. Thus all odd-order moments of $\nu$ vanish. For even $|\alpha|$,
including $|\alpha|=0$, the identity $\check\nu=-\nu$ gives
$$
    \int_{\cB_r}\beta^\alpha\,\rd \nu(\beta)
    =\int_{\cB_r}(-\beta)^\alpha\,\rd \nu(\beta)
    =\int_{\cB_r}\beta^\alpha\,\rd \check\nu(\beta)
    =-\int_{\cB_r}\beta^\alpha\,\rd \nu(\beta),
$$
so these moments also vanish.

It follows that $\int P\,\rd \nu=0$ for every polynomial $P$.
Since $\nu$ is a finite signed measure supported on $\cB_r$,
\cref{lem:app_compact_moment_uniqueness} gives $\nu=0$,
which proves \eqref{eq:app_sigmoid_kernel}.

Conversely, suppose that $\mu=\check\mu$ and $\mu(\cB_r)=0$. For any $x\in \RR^d$, the symmetry
of $\mu$ and the oddness of $\tanh$ give
$$
    \int_{\cB_r}\tanh(x^\T\beta/2)\,\rd \mu(\beta)=0.
$$
Using $\sigma(t)=\{1+\tanh(t/2)\}/2$ and $\mu(\cB_r)=0$
then gives $p_\mu(x)=0$. This proves the converse.
\end{proof}

\begin{proof}[Proof of \cref{prop:identifiability}]
For any $G\in \cP(\cB_r)$,  write its reflection as $\check G(A)=G(-A)$, where $-A:=\{-\beta:\beta\in A\}$.
The density lower bound in \cref{assum:particle_population}(ii) implies
that $p_G=p_{G^*}$, Lebesgue-almost everywhere on $\cB_{r_X}$.
By continuity, the two functions agree everywhere on this ball.
Applying \cref{lem:app_sigmoid_kernel}
in \cref{app:particle_identifiability} to the finite signed measure
$G-G^*$ gives
\begin{equation}
    G-\check G=G^*-\check G^*.
    \label{eq:app_reflection_identity}
\end{equation}
Let $S=\operatorname{supp}(G^*)$.
The condition $\beta_a^*+\beta_b^*\ne0$ for all $a,b\in[K^*]$ gives
$S\cap(-S)=\varnothing$. Evaluating \eqref{eq:app_reflection_identity}
on $S$ yields
$$
    G(S)-G(-S)=G^*(S)-G^*(-S)=1.
$$
Since $G(S)\le1$ and $G(-S)\ge0$, we have $G(S)=1$ and $G(-S)=0$.
For every $\beta\in S$, \eqref{eq:app_reflection_identity} then gives
$$
    G(\{\beta\})
    =G(\{\beta\})-G(\{-\beta\})
    =G^*(\{\beta\}).
$$
Both distributions are supported on the finite set $S$, so these
equalities imply $G=G^*$.
\end{proof}

\begin{remark}[Comparison of reflection conditions with the binary-choice models]
Restrictions related with our reflection condition in \cref{assum:particle_population} arise in random-coefficient binary-choice models of the form
$$
    Y=\mathbf{1}\{X^\T\beta\ge0\},
$$
where the random coefficient vector $\beta$ is independent of $X$ and is nonzero almost surely. Since positive rescaling of $\beta$ leaves the choice unchanged, these models identify the distribution of the normalized direction $\theta=\beta/\|\beta\|_2$; coefficient magnitudes are not identified.

To identify the density of $\theta$, \citet[Assumption~3.2]{gautier2013nonparametric} impose a hemisphere restriction: there exists $v\ne0$ such that $v^\T\theta>0$ almost surely. For our finite atoms, the analogous restriction would require $v^\T\beta_a^*>0$ for every $a\in[K^*]$. \citet[Assumption~A2.1]{gautier2018adaptive} use the weaker condition $f_\theta(u)f_\theta(-u)=0$ for almost every $u\in\mathbb S^{d-1}$, where $f_\theta$ is the density with respect to spherical measure. Its finite-atom analogue excludes opposite directions:
$$
    \frac{\beta_a^*}{\|\beta_a^*\|_2}
    \ne-\frac{\beta_b^*}{\|\beta_b^*\|_2},
    \qquad a,b\in[K^*].
$$

These restrictions are stronger than our reflection condition $\beta_a^*+\beta_b^*\ne0$, which excludes only exact reflected pairs and zero atoms. In particular, our condition permits opposite directions with different magnitudes, such as $\beta_b^*=-2\beta_a^*\ne0$. The comparison concerns the geometric restrictions: the binary-choice results identify a distribution of directions, whereas \cref{prop:identifiability} identifies the full mixing distribution in our logistic model, including both coefficient directions and magnitudes.
\end{remark}

\subsection{A quantitative denoising bound}
\label{app:particle_denoising}

For a candidate $G = \sum_{j=1}^M w_j \delta_{\beta_j}$, we bound its weight and particle location errors by its estimation error of $p_{G^*}$ for $G^* = \sum_{k=1}^{K^*}w_k^* \delta_{\beta_k^*}$. To accommodate several atoms of $G$ approximating a single atom of $G^*$, we define a discrepancy between $G$ and $G^*$ based on the total assigned weight and particle location errors of $G$ in a neighborhood of each true atom in $G^*$, rather than pair their atoms  individually.

\subsubsection{A discrepancy between $G$ and $G^*$}
\begin{figure}[htbp]
    \centering
    \includegraphics[width=0.88\textwidth]
        {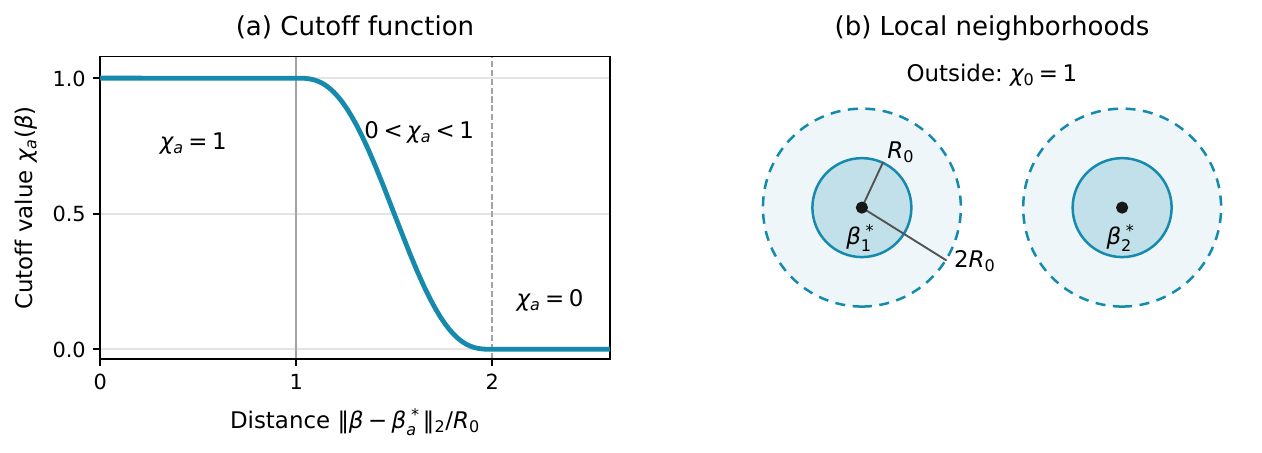}
    \caption{Cutoff functions used to define
    $\mfa(G,G^*)$. Left: the cutoff as a function of normalized
    distance from a true atom. Right: schematic neighborhoods of two true
    atoms (not to scale), with solid inner circles of radius $R_0$ and dashed outer
    circles of radius $2R_0$. The outer balls are disjoint.}
    \label{fig:population_cutoff_neighborhoods}
\end{figure}

For each $a\in[K^*]$, to define the neighborhoods around $\beta_a^*$, fix
\begin{equation}
    R_0={1\over16}\min\{s_-,s_+,\tau_0\}
    \label{eq:app_population_cutoff_radius}
\end{equation}
with $s_-$ and $s_+$ given in \cref{assum:particle_population} and $\tau_0$ given in \eqref{eq:particle_radius_choice}. We use a $C^2$ function
$\chi_a:\RR^d\to[0,1]$ that is equal to one on
$\cB(\beta_a^*,R_0)$ and is zero outside
$\cB(\beta_a^*,2R_0)$, with a smooth transition in between. We give its construction below and illustrate it in the left panel of \cref{fig:population_cutoff_neighborhoods}.  Since $4R_0<s_-$, the outer balls are disjoint. Set
$$
    \chi_0(\beta)=1-\sum_{a=1}^{K^*}\chi_a(\beta).
$$
Then $0\le\chi_0\le1$, with $\chi_0=0$ on the inner balls
and $\chi_0=1$ outside the outer balls; see
\cref{fig:population_cutoff_neighborhoods}.

For a probability measure $G$ on $\cB_r$,
$\int\chi_a\,\rd G$ is the mass assigned to the neighborhood of
$\beta_a^*$, and $\int\chi_0\,\rd G$ is the mass not assigned to the neighborhood of any true atom.
We measure the discrepancy between $G$ and $G^*$  by
\begin{equation}
\begin{aligned}
    \mfa(G,G^*)
    ={}&\int\chi_0(\beta)\,\rd G(\beta)
       +\sum_{a=1}^{K^*}\left|\int\chi_a(\beta)\,\rd G(\beta) - w_a^*\right|\\
       &+\sum_{a=1}^{K^*}\left\|\int\chi_a(\beta)(\beta-\beta_a^*)\,\rd G(\beta)\right\|_2
       +\sum_{a=1}^{K^*}\int\chi_a(\beta)\|\beta-\beta_a^*\|_2^2\,\rd G(\beta).
\end{aligned}
    \label{eq:app_population_gauge}
\end{equation}
The four terms in $\mfa(G,G^*)$ measure the mass not assigned to any true-atom neighborhood,
the weight errors within each neighborhood, the first moment particle location errors, and the second moment particle location errors. Opposite particle location errors can cancel in the third term, but not
in the fourth.
If the first and fourth terms vanish, $G$ is supported on the true
atoms. If the second term also vanishes, their weights further agree with those
of $G^*$. Consequently,
$$
    \mfa(G,G^*)=0 \quad\Longleftrightarrow\quad G=G^*.
$$

For the construction of $\chi_a$, define the nonincreasing $C^2$ function
$$
    \bar\chi(t)=
    \begin{cases}
        1, & 0\le t\le1,\\
        1-10(t-1)^3+15(t-1)^4-6(t-1)^5, & 1<t<2,\\
        0, & t\ge2,
    \end{cases}
$$
and set $\chi(u) :=\bar\chi(u/R_0)$.
For $a\in[K^*]$, take
\begin{equation}\label{def_chi_a_dist}
    \chi_a(\beta)=\chi(\|\beta-\beta_a^*\|_2).
\end{equation}
For $1<t<2$, writing $v=t-1$ gives
$$
    \bar\chi'(t)=-30v^2(1-v)^2,\qquad
    \bar\chi''(t)=-60v(1-v)(1-2v).
$$
Both derivatives vanish at the endpoints, and
$\|\bar\chi'\|_\infty\le15/8$,
$\|\bar\chi''\|_\infty\le15$, using $v(1-v)\le1/4$.
On the annulus $\{\beta:R_0<\|\beta-\beta_a^*\|_2<2R_0\}$,
differentiation gives
\begin{align*}
    \nabla\chi_a(\beta)
    &={\bar\chi'(\|\beta-\beta_a^*\|_2/R_0)\over R_0}
            {\beta-\beta_a^*\over\|\beta-\beta_a^*\|_2},\\
    \nabla^2\chi_a(\beta)
    &={\bar\chi''(\|\beta-\beta_a^*\|_2/R_0)\over R_0^2}
            {(\beta-\beta_a^*)(\beta-\beta_a^*)^\T
                      \over\|\beta-\beta_a^*\|_2^2}\\
    &\quad+{\bar\chi'(\|\beta-\beta_a^*\|_2/R_0)
                      \over R_0\|\beta-\beta_a^*\|_2}
       \left(I_d-{(\beta-\beta_a^*)(\beta-\beta_a^*)^\T
                      \over\|\beta-\beta_a^*\|_2^2}\right).
\end{align*}
Both derivatives are zero off this annulus and extend continuously to
its boundary. It is easy to see that
$$
    \|\nabla\chi_a\|_\infty\le{15\over8R_0},\qquad
    \|\nabla^2\chi_a\|_\infty\le{15\over R_0^2}.
$$

\subsubsection{Generalized functions and compactness}
\label{app:particle_uniform_preliminaries}

\paragraph{Generalized functions}
\label{app:particle_functional_prelim}

The proof of \cref{thm:app_population_inverse} uses limits of
normalized differences of mixing distributions. For background on
generalized functions, see
\citet[Sections~2.1--2.2]{hormander1990analysis}.

A sequence of scaled differences need not converge in total variation,
even when its integrals against smooth functions converge.
For example, in one dimension, let
$\mu_h=(\delta_{a+h}-\delta_a)/h$ for $h\ne0$. Although
$\|\mu_h\|_{\mathrm{TV}}=2/|h|$ diverges as $h\to0$, every
smooth function $\varphi$ satisfies
$$
    \int\varphi\,\rd \mu_h
    ={\varphi(a+h)-\varphi(a)\over h}
    \to\varphi'(a).
$$
Similarly,
$$
    \int\varphi\,\rd \left\{
       {\delta_{a+h}+\delta_{a-h}-2\delta_a\over h^2}\right\}
    ={\varphi(a+h)+\varphi(a-h)-2\varphi(a)\over h^2}
    \to\varphi''(a).
$$
In the first example, we change the location of a point mass. In the
second, we split it symmetrically about $a$. The maps
$\varphi\mapsto\varphi'(a)$ and $\varphi\mapsto\varphi''(a)$ are
linear, but neither is integration against a finite signed measure.
These maps are examples of generalized functions.

Let
$$
    U :=\cB_{r+1}^\circ := \{x\in \RR^d: \|x\|_2 < r+1\}.
$$
A smooth test
function $\varphi:U\to\RR$ is an infinitely differentiable function
that vanishes outside some compact subset of $U$. We denote the space
of these functions by $C_c^\infty(U)$. In this appendix, we need only
generalized functions of order at most two. Such a generalized function
$u$ is a linear map from $C_c^\infty(U)$ to $\RR$, whose value
at $\varphi$ is written $\langle u,\varphi\rangle$, with the
following bound: for every compact $K\subset U$, there is a finite
constant $C_K$ such that for all $\varphi \in C_c^\infty(U)$ with $\operatorname{supp}\varphi\subset K$,
\begin{equation}\label{bd_order_two}
    |\langle u,\varphi\rangle|
    \le C_K\left(\sup_{\beta\in U} |\varphi(\beta)|+
         \sup_{\beta\in U} \|\nabla_\beta\varphi(\beta)\|_2+
         \sup_{\beta\in U} \|\nabla_\beta^2\varphi(\beta)\|_\op\right).
\end{equation}
This bound depends on the first- and second-order mixed derivatives of the test
function, as in the preceding examples.

Every finite signed measure $\mu$ defines a generalized function by
$\langle\mu,\varphi\rangle=\int\varphi\,\rd \mu$. More generally,
convergence of generalized functions $u_m$ to $u$ means
$$
    \lim_{m\to\i}\langle u_m,\varphi\rangle=\langle u,\varphi\rangle
    \qquad\text{for every }\varphi\in C_c^\infty(U).
$$
We say that $u$ is supported on $\cB_r$ if
$\langle u,\varphi\rangle=0$ whenever $\varphi$ is supported in
$U\setminus\cB_r$. Any generalized function $u$ supported on $\cB_r$ can also be applied to arbitrary smooth functions on $U$. Specifically, choose a cut-off function $\zeta\in C_c^\infty(U)$ that equals one on an open neighborhood of $\cB_r$. For any $f\in C^\infty(U)$, the product $\zeta f$ belongs to $C_c^\infty(U)$, so we may define
\begin{equation}\label{def_cutoff}
    \langle u,f\rangle := \langle u,\zeta f\rangle.
\end{equation}
This value is independent of the choice of $\zeta$. To see this, let $\wt\zeta\in C_c^\infty(U)$ also equal one on an open neighborhood of $\cB_r$. Then $(\zeta-\wt\zeta)f$ is a smooth test function whose support is contained in $U\setminus\cB_r$. Since $u$ is supported on $\cB_r$, we obtain
$$
    \langle u,\zeta f\rangle-\langle u,\wt\zeta f\rangle
    =\langle u,(\zeta-\wt\zeta)f\rangle=0.
$$
In particular, if $u$ is a finite signed measure supported on $\cB_r$, then $\langle u,1\rangle=u(\cB_r)$ according to \eqref{def_cutoff}.

The proofs below also use the sigmoid transform and reflection of
generalized functions. We define these two operations as follows.
For $u$ supported on
$\cB_r$, the cutoff definition allows us to apply $u$ to the smooth
function $\beta\mapsto\sigma(x^\T\beta)$. We define its sigmoid
transform by
\begin{equation}\label{def_Tu}
    (Tu)(x):=\langle u,\sigma(x^\T\,\cdot)\rangle,
    \qquad x\in\RR^d.
\end{equation}
For a finite signed measure $\mu$, this gives
$(T\mu)(x)=\int\sigma(x^\T\beta)\,\rd \mu(\beta)$.
In particular, $TG=p_G$ for a mixing distribution $G$ on $\cB_r$.

The reflection $\check u$ is defined by applying $u$ to reflected
test functions:
$$
    \langle\check u,\varphi\rangle=\langle u,\check\varphi\rangle,
    \qquad \check\varphi(\beta)=\varphi(-\beta),
    \quad \varphi\in C_c^\infty(U).
$$
For a signed measure $\mu$, this agrees with
$\check\mu(A)=\mu(-A)$.

\paragraph{Compactness, moments, and positivity}
\label{app:particle_compactness_prelim}

We use compactness results for linear functionals,
followed by two results for generalized functions.

For generalized functions of order at most two, we use the following
compactness result. With $U=\cB^\circ_{r+1}$ and $\overline U = \cB_{r+1}$, let
$C^2(\overline U)$ consist of the real-valued functions that are
twice continuously differentiable on $U$, with the function and all
derivatives through order two extending continuously to $\overline U$.
We use the norm
\begin{equation}\label{norm_C2}
    \|\varphi\|_{C^2}
    ~ :=~  \sup_{\beta \in \overline U}|\varphi(\beta)|+\sup_{\beta \in \overline U}\|\nabla_\beta\varphi(\beta)\|_2
                        +\sup_{\beta \in \overline U}\|\nabla_\beta^2\varphi(\beta)\|_\op,
\end{equation}
Thus convergence in this norm means uniform convergence of the
functions, their gradients, and their Hessians.
To obtain the compactness result below, we use that
$C^2(\overline U)$ is a separable Banach space; see
\citet[Section~1.1.6, p.~12]{emelyanov2007homotopy}, in which an
equivalent norm is used. The next lemma follows from
\citet[Corollary~3.30]{brezis2011functional}, applied to
$C^2(\overline U)$.

\begin{lemma}
\label{lem:app_compact_functionals}
Let $u_m:C^2(\overline U)\to\RR$ be linear functionals with $m\in \NN$. Suppose
that, for some finite $C$ independent of $m$,
$$
    |\langle u_m,\varphi\rangle|\le C\|\varphi\|_{C^2},
    \qquad \varphi\in C^2(\overline U).
$$
There are a subsequence $u_{m_k}$ and a linear functional
$u:C^2(\overline U)\to\RR$ such that
$$
    \lim_{k\to\i}\langle u_{m_k},\varphi\rangle =
       \langle u,\varphi\rangle,
    \qquad |\langle u,\varphi\rangle|\le C\|\varphi\|_{C^2},
    \quad \varphi\in C^2(\overline U).
$$
The restriction of $u$ to $C_c^\infty(U)$ is a generalized
function of order at most two. In addition, if
$\varphi_k\to\varphi$ in $C^2(\overline U)$ as $k\to \i$, then
$$
    \lim_{k\to\i}\langle u_{m_k},\varphi_k\rangle =\langle u,\varphi\rangle.
$$
\end{lemma}

The same subsequence works for every test function.
The last assertion also allows the test function to vary
along the subsequence: indeed,
$$
    |\langle u_{m_k},\varphi_k\rangle-\langle u,\varphi\rangle|
    \le C\|\varphi_k-\varphi\|_{C^2}
          +|\langle u_{m_k}-u,\varphi\rangle|\to0.
$$

For a multi-index $\alpha=(\alpha_1,\ldots,\alpha_d)$ with
nonnegative integer entries, write
$$
\beta^\alpha=\prod_{j=1}^d\beta_j^{\alpha_j},
\qquad |\alpha|=\sum_{j=1}^d\alpha_j,
\qquad\text{and}\qquad
\partial^\alpha=\partial_1^{\alpha_1}\cdots\partial_d^{\alpha_d}.
$$
The polynomial moments of a compactly supported generalized function
$u$ are the values $\langle u,\beta \mapsto \beta^\alpha\rangle$, defined
using the cutoff in \eqref{def_cutoff}. The next lemma shows that if the polynomial moments of $u$ are zero, then $u$ is zero. The proof follows from Fourier inversion
and analyticity of its Fourier transform; see
\citet[Sections~7.1 and 7.3]{hormander1990analysis}.

\begin{lemma}
\label{lem:app_compact_moment_uniqueness}
Let $u$ be a generalized function of order at most two, and is supported on
$\cB_r$. If
$$
    \langle u,P\rangle=0
    \qquad\text{for every polynomial }P:\RR^d\to\RR,
$$
then $u=0$. In particular, two finite signed measures supported on
$\cB_r$ are equal if all their polynomial moments agree.
\end{lemma}

\begin{proof}
Choose any cutoff function $\zeta\in C_c^\infty(U)$ equal to one on an open neighborhood
of $\cB_r$, and let $\cK:=\operatorname{supp}(\zeta)$. Recall that, for
every smooth function $f$ on $U$, we define
$
    \langle u,f\rangle=\langle u,\zeta f\rangle
$
as in \eqref{def_cutoff}.
We also allow complex-valued functions by applying $u$ separately
to their real and imaginary parts:
$$
    \langle u,f_1+if_2\rangle
    :=\langle u,f_1\rangle+i\langle u,f_2\rangle.
$$
The order-two bound in \eqref{bd_order_two} continues to hold for complex-valued test
functions supported in $\cK$, after increasing the constant if necessary.

We first show that $\langle u,e^{-it^\T\cdot}\rangle=0$ for every $t\in \RR^d$. Fix $t\in\RR^d$ and define
$$
    P_{N,t}(\beta)
    :=\sum_{n=0}^N\frac{(-i\,t^\T\beta)^n}{n!}.
$$
Since the real and imaginary parts of $P_{N,t}$ are real polynomials, the hypothesis gives $\langle u,P_{N,t}\rangle=0$. Since the exponential series and its first two derivatives with respect to $\beta$ converge uniformly on $\cK$ as $N\to \i$, multiplication by the fixed smooth cutoff $\zeta$ therefore gives
$$
    \|\zeta(P_{N,t}-e^{-it^\T\cdot})\|_{C^2}\to0
$$
with the norm defined in \eqref{norm_C2}. Since all these functions are supported in $\cK$, the order-two bound in \eqref{bd_order_two} yields
$$
    \bigl|\langle u,P_{N,t}\rangle
          -\langle u,e^{-it^\T\cdot}\rangle\bigr|
    \le C\|\zeta(P_{N,t}-e^{-it^\T\cdot})\|_{C^2}
    \to0.
$$
Thus $\langle u,e^{-it^\T\cdot}\rangle=0$ for every $t\in \RR^d$ which implies that the Fourier transform
$$
    \wh u(t):=\langle u,e^{-it^\T\beta}\rangle
$$
vanishes for every $t\in\RR^d$. Since a compactly supported generalized function is uniquely determined by its Fourier transform, it follows that $u=0$.

Indeed, let $\varphi\in C_c^\infty(U)$ be arbitrary. Extend $\varphi$ by zero outside $U$ to obtain a smooth, compactly supported function on $\RR^d$, and define
$$
    \wh\varphi(t)
    :=\int_{\RR^d}e^{-it^\T\beta}\varphi(\beta)\,\rd \beta.
$$
Fourier inversion gives
$$
    \varphi(\beta)
    =\frac{1}{(2\pi)^d}
      \int_{\RR^d}e^{it^\T\beta}\wh\varphi(t)\,\rd t.
$$
Since $\wh\varphi$ decays faster than any inverse power of $\|t\|_2$ as $\varphi\in C_c^\infty(U)$,
$$
    \int_{\RR^d}(1+\|t\|_2^2)|\wh\varphi(t)|\,\rd t<\infty.
$$
Consequently, after multiplication by $\zeta$, the inversion integral converges in the $C^2$ norm, with all integrands supported in $\cK$. The order-two bound in \eqref{bd_order_two} therefore allows us to apply $u$ under the integral. Since $\zeta=1$ on a neighborhood of $\operatorname{supp}u$, we obtain
$$
\begin{aligned}
    \langle u,\varphi\rangle
    &=\langle u,\zeta\varphi\rangle =\frac{1}{(2\pi)^d}
      \int_{\RR^d}\wh\varphi(t)
          \langle u,\zeta e^{it^\T\beta}\rangle\,\rd t =\frac{1}{(2\pi)^d}
      \int_{\RR^d}\wh\varphi(t)\wh u(-t)\,\rd t
      =0.
\end{aligned}
$$
As $\varphi\in C_c^\infty(U)$ was arbitrary, $u=0$. The proof is complete.
\end{proof}
A generalized function $u$ is called nonnegative if
$\langle u,\varphi\rangle\ge0$ for every nonnegative
$\varphi\in C_c^\infty(U)$. The next lemma follows from the
Riesz representation theorem
\citep[Theorem~11.1]{knapp2016basic}, using the bound given below.

\begin{lemma}
\label{lem:app_positive_generalized_function}
Let $u$ be a nonnegative generalized function on $U$. There is
a unique nonnegative Borel measure $\nu$ on $U$, finite on
compact subsets, such that
$$
    \langle u,\varphi\rangle=\int_U\varphi\,\rd \nu,
    \qquad \varphi\in C_c^\infty(U).
$$
If, in addition, $u$ is supported on $\cB_r$ and
$\langle u,1\rangle=0$, then $u=0$.
\end{lemma}

\begin{proof}
Fix a compact
$K\subset U$, and take a nonnegative
$\psi\in C_c^\infty(U)$ equal to one on an open neighborhood of
$K$. For any test
function $\varphi$ supported on $K$, the inequalities
$-\|\varphi\|_\infty\psi\le\varphi\le
\|\varphi\|_\infty\psi$ give
$$
    |\langle u,\varphi\rangle|
       \le\langle u,\psi\rangle\|\varphi\|_\infty.
$$
This supremum-norm bound permits a unique extension to continuous
functions of compact support. The Riesz representation theorem gives
the nonnegative measure $\nu$. If $u$ is supported on
$\cB_r$, so is $\nu$, and
$\nu(U)=\langle u,1\rangle=0$ implies $\nu=0$.
\end{proof}

\subsubsection{The denoising bound}

We will apply the sigmoid transform to   limits of a sequence of
$\{(G_m-G^*)/\mfa(G_m,G^*)\}_{m\ge 1}$. Since such a limit may be a
generalized function, we first extend \eqref{eq:app_sigmoid_kernel} in \cref{lem:app_sigmoid_kernel}
from finite signed measures to the sigmoid transform $Tu$ in \eqref{def_Tu}.

\begin{lemma}
\label{lem:app_sigmoid_kernel_second_order}
Let $u$ be a generalized function of order at most two supported on
$\cB_r$.
If $Tu$ vanishes on a nonempty open subset of $\RR^d$, then
\begin{equation}
    \langle u,1\rangle=0,
    \qquad
    u=\check u.
    \label{eq:app_sigmoid_kernel_second_order}
\end{equation}
\end{lemma}
\begin{proof}
We follow the analyticity argument in the proof of
\cref{lem:app_sigmoid_kernel}. To apply $u$ term by term to the Taylor
expansion of $\sigma(z^\T\beta)$ around any $x_0\in\RR^d$, we establish
uniform convergence in $\beta$ of the series and its first two
derivatives with respect to $\beta$.

Fix $x_0\in\RR^d$.
Take $\delta=\pi/\{2\sqrt d(r+1)\}$, as in the proof of
\cref{lem:app_sigmoid_kernel}, and write
$Z=\|x_0\|_2+\sqrt d\,\delta$ and
$R_\delta=\{z\in\CC^d:\max_j|z_j-x_{0,j}|\le\delta\}$.
For $z\in R_\delta$ and $\beta\in\overline U$, the
facts $|\sigma(z^\T\beta)|\le1$ and
$|1-\sigma(z^\T\beta)|\le1$ give
$$
    |\sigma'(z^\T\beta)|\le1,\qquad
    |\sigma''(z^\T\beta)|\le 2,
$$
by $\sigma'=\sigma(1-\sigma)$ and
$\sigma''=\sigma(1-\sigma)(1-2\sigma)$. Therefore,
$$
    \|\sigma(z^\T\,\cdot)\|_{C^2}
    \le 1+Z+2Z^2.
$$
Here the same norm in \eqref{norm_C2} is used for complex-valued functions. Write the
local Taylor series as
$\sigma(z^\T\beta)=\sum_\alpha c_\alpha(\beta)(z-x_0)^\alpha$.
Cauchy's estimate, applied also to the first two derivatives in
$\beta$, gives
$$
    \|c_\alpha\|_{C^2}
       \le(1+Z+2Z^2)\delta^{-|\alpha|}.
$$
Consequently, when $z\in R_{\delta/2}$,
$$
    \sum_\alpha\|c_\alpha(z-x_0)^\alpha\|_{C^2}
       \le 2^d(1+Z+2Z^2).
$$
Choose a smooth cutoff $\zeta$ equal to one on an open neighborhood of
$\cB_r$, as in
\cref{app:particle_functional_prelim}. Multiplication by $\zeta$
preserves this convergence, since the product rule gives
$$
    \|\zeta f\|_{C^2}\le2\|\zeta\|_{C^2}\|f\|_{C^2}.
$$
The order-two bound on $u$ now permits applying $u$ term by term.
Thus $Tu$ is real analytic, and the identity theorem gives
$Tu=0$ on $\RR^d$.

Since $\sigma(0)=1/2$, evaluation at zero gives
$$
    0=(Tu)(0)={1\over2}\langle u,1\rangle.
$$
The identity $\sigma(t)=\{1+\tanh(t/2)\}/2$ now yields
\begin{equation}
    \langle u,\tanh(x^\T\cdot/2)\rangle=0,
    \qquad x\in\RR^d.
    \label{eq:app_second_order_tanh_zero}
\end{equation}
Set $\nu=u-\check u$. By the oddness of $\tanh$,
\eqref{eq:app_second_order_tanh_zero} implies
$$
    \langle\nu,\tanh(x^\T \cdot/2)\rangle=0,
    \qquad x\in\RR^d.
$$
Fix $v\in\RR^d$, substitute $x=t v$, and differentiate at $t=0$.
For every $m\ge0$,
$$
    0={\tanh^{(2m+1)}(0)\over2^{2m+1}}
      \langle\nu,(v^\T\,\cdot)^{2m+1}\rangle.
$$
The recurrence in the proof of \cref{lem:app_sigmoid_kernel} shows that
every odd derivative of $\tanh$ at zero is nonzero. Varying $v\in \RR^d$, as in
that proof, shows that all odd-order moments of
$\nu$ vanish. Moreover, $\check\nu=-\nu$, so all even-order moments vanish
as well.

Applying \cref{lem:app_compact_moment_uniqueness} gives $\nu=0$. Thus
$u=\check u$, as required.
\end{proof}

The following theorem is the denoising result that bounds $\mfa(G,G^*)$ from above by the
estimation error of $p_G$ for any $G^*$ satisfying \cref{assum:particle_population}.

\begin{theorem}
\label{thm:app_population_inverse}
Grant \cref{assum:particle_population} and \eqref{eq:particle_radius_choice}.
There is a constant $c_{\mathrm{inv}}>0$, depending only on $d,r,\tau_0$ and
the constants in \cref{assum:particle_population}, such that,
for every  probability measure $G$ on $\cB_r$,
\begin{equation}
    \left\{P_X(p_G-p^*)^2\right\}^{1/2}
    \ge c_{\mathrm{inv}}\,\mfa(G,G^*).
    \label{eq:app_population_inverse}
\end{equation}
The constant $c_{\mathrm{inv}}$ does not depend on the support size of $G$.
\end{theorem}
\begin{proof}
We argue by contradiction. We normalize $G-G^*$ by
$\mfa(G,G^*)$ and show that a subsequence converges to a nonzero
generalized function. Failure of the bound forces its sigmoid transform
to vanish, so \cref{lem:app_sigmoid_kernel_second_order} makes the limit
symmetric with zero total mass. Away from the limiting true atoms, the
limit is nonnegative. Reflection separation and symmetry then imply
nonnegativity everywhere, contradicting its nonzero value.

Keep $d,r,\tau_0$ and all constants in
\cref{assum:particle_population} fixed. If the claimed uniform bound fails,
there are admissible sequences $G_m^*$, $P_{X,m}$, and $G_m$ such
that, with $\mfa_m=\mfa(G_m,G_m^*)>0$,
\begin{equation}
    {1\over \mfa_m}\{P_{X,m}(p_{G_m}-p_{G_m^*})^2\}^{1/2}
    \to0.
    \label{eq:app_inverse_contradiction_sequence}
\end{equation}
Throughout the proof, we pass to subsequences without changing the index
$m$.

Since the number of true atoms is bounded by $K_0$, we may first pass
to a subsequence with a fixed number of true atoms. After labeling the
atoms, compactness of $\cB_{r-\tau_0}$
allows a further subsequence such that
$$
    G_m^*=\sum_{a=1}^{K^*}
        w_{a,m}^*\delta_{\beta_{a,m}^*},
    \qquad
    \lim_{m\to\i}\beta_{a,m}^* =\beta_{a,\infty}^*,
    \quad a\in[K^*].
$$
Both separation conditions are preserved in the limit. The cutoff radius
$R_0$ is fixed, and we write
$$
    \chi_{a,m}(\beta)=\chi(\|\beta-\beta_{a,m}^*\|_2),
    \qquad
    \chi_{0,m}=1-\sum_{a=1}^{K^*}\chi_{a,m}.
$$
Write the four components in
\eqref{eq:app_population_gauge} as
$$
\begin{aligned}
    F_{0,m}&=\int\chi_{0,m}(\beta)\,\rd G_m(\beta),&
    A_{a,m}&=\int\chi_{a,m}(\beta)\,\rd G_m(\beta)-w_{a,m}^*,\\
    B_{a,m}&=\int\chi_{a,m}(\beta)(\beta-\beta_{a,m}^*)\,\rd G_m(\beta),&
    C_{a,m}&=\int\chi_{a,m}(\beta)\|\beta-\beta_{a,m}^*\|_2^2\,\rd G_m(\beta).
\end{aligned}
$$
In particular,
$$
    \mfa_m=F_{0,m}+\sum_{a=1}^{K^*}|A_{a,m}|+
        \sum_{a=1}^{K^*}\|B_{a,m}\|_2+\sum_{a=1}^{K^*} C_{a,m}.
$$

To obtain a limit without restricting the support size of $G_m$, set
$$
    \mu_m={G_m-G_m^*\over \mfa_m}
$$
and use $U=\cB^\circ_{r+1}$ and the $C^2$ norm from
\cref{app:particle_compactness_prelim}. For
$\varphi\in C^2(\overline U)$, expand $\varphi$ to first order around
each $\beta_{a,m}^*$. Multiplying the expansions by the cutoffs
and integrating gives
\begin{align}
    \langle G_m-G_m^*,\varphi\rangle
    ={}&\int\chi_{0,m}\varphi\,\rd G_m
       +\sum_{a=1}^{K^*} A_{a,m}\varphi(\beta_{a,m}^*)\notag\\
       &+\sum_{a=1}^{K^*} B_{a,m}^\T\nabla_\beta\varphi(\beta_{a,m}^*)
       +R_m(\varphi),
    \label{eq:app_inverse_taylor_decomposition}
\end{align}
where, writing $h=\beta-\beta_{a,m}^*$ inside the $a$th integral,
$$
    R_m(\varphi)
    =\sum_{a=1}^{K^*}\int\chi_{a,m}(\beta)
       \int_0^1(1-t)h^\T
       \nabla^2\varphi(\beta_{a,m}^*+th)h\,\rd t\,\rd G_m(\beta).
$$
Since $\int_0^1(1-t)\,\rd t=1/2$, this gives
$$
    |R_m(\varphi)|
    \le {1\over2}\|\nabla_\beta^2\varphi\|_\infty\sum_{a=1}^{K^*} C_{a,m}.
$$
Dividing by $\mfa_m$ therefore gives
\begin{align*}
    |\langle\mu_m,\varphi\rangle|
    \le{}&{F_{0,m}+\sum_{a=1}^{K^*}|A_{a,m}|\over \mfa_m}\|\varphi\|_\infty
       +{\sum_{a=1}^{K^*}\|B_{a,m}\|_2\over \mfa_m}
          \|\nabla_\beta\varphi\|_\infty\\
       &+{\sum_{a=1}^{K^*} C_{a,m}\over2\mfa_m}
          \|\nabla_\beta^2\varphi\|_\infty
    \le \|\varphi\|_{C^2}.
\end{align*}
Applying \cref{lem:app_compact_functionals} gives a bounded linear
functional $\mu: C^2(\overline U) \to \RR$ such that
$$
    \langle\mu_m,\varphi\rangle\to\langle\mu,\varphi\rangle,
    \qquad \varphi\in C^2(\overline U),
    \qquad |\langle\mu,\varphi\rangle|\le\|\varphi\|_{C^2}.
$$
Restricting each $\varphi\in C_c^\infty(\RR^d)$ to
$\overline U$ defines a generalized function, also denoted by $\mu$.
The bound $|\langle\mu,\varphi\rangle|\le\|\varphi\|_{C^2}$
shows that its order is at most two. If a smooth test function is supported
away from $\cB_r$, its pairing with every $\mu_m$ is zero, so its
pairing with $\mu$ is also zero. Thus $\mu$ is supported on $\cB_r$.
Its action on any smooth function agrees with the cutoff convention in
\cref{app:particle_functional_prelim}: multiplication by a smooth cutoff
equal to one on an open neighborhood of $\cB_r$ leaves every
$\mu_m$-pairing unchanged, and
this identity passes to the limit.

We next find a function $\varphi_\infty\in C^2(\overline U)$ such that
$\langle\mu,\varphi_\infty\rangle>0$. The nonnegative numbers
$$
    {F_{0,m}\over \mfa_m},\qquad
    {|A_{a,m}|\over \mfa_m},\qquad
    {\|B_{a,m}\|_2\over \mfa_m},\qquad
    {C_{a,m}\over \mfa_m},\qquad a\in[K^*],
$$
sum to one. Passing to a subsequence makes every coordinate converge.
There are $1+3K^*$ coordinates, so the largest limit satisfies
$\gamma\ge(1+3K^*)^{-1}\ge(1+3K_0)^{-1}>0$.
To construct the required functions, recall that
$$
    \mfa_m\langle\mu_m,\varphi\rangle
    =\int\varphi\,\rd G_m
       -\sum_{b=1}^{K^*}w_{b,m}^*\varphi(\beta_{b,m}^*).
$$
The cutoff functions satisfy
$\chi_{0,m}(\beta_{b,m}^*)=0$ and
$\chi_{a,m}(\beta_{b,m}^*)=\mathbf1\{a=b\}$. Hence
\begin{align*}
    \mfa_m\langle\mu_m,\chi_{0,m}\rangle
       &=\int\chi_{0,m}\,\rd G_m=F_{0,m},\\
    \mfa_m\langle\mu_m,\chi_{a,m}\rangle
       &=\int\chi_{a,m}\,\rd G_m-w_{a,m}^*=A_{a,m}.
\end{align*}
For any $v\in\RR^d$, both
$\chi_{a,m}(\beta)v^\T(\beta-\beta_{a,m}^*)$ and
$\chi_{a,m}(\beta)\|\beta-\beta_{a,m}^*\|_2^2$ vanish at every
true atom: the displacement is zero when $b=a$, and the cutoff is
zero when $b\ne a$. Their integrals under $G_m^*$ are therefore
zero, giving
\begin{align*}
    \mfa_m\big\langle\mu_m,
       \chi_{a,m}(\cdot)v^\T(\cdot-\beta_{a,m}^*)\big\rangle
       &=v^\T\int\chi_{a,m}(\beta)(\beta-\beta_{a,m}^*)\,\rd G_m(\beta)
         =v^\T B_{a,m},\\
    \mfa_m\big\langle\mu_m,
       \chi_{a,m}(\cdot)\|\cdot-\beta_{a,m}^*\|_2^2\big\rangle
       &=\int\chi_{a,m}(\beta)\|\beta-\beta_{a,m}^*\|_2^2\,\rd G_m(\beta)
         =C_{a,m}.
\end{align*}
Multiplying $\chi_{a,m}$ by $\operatorname{sign}(A_{a,m})$
gives $|A_{a,m}|$ in place of $A_{a,m}$. When $B_{a,m}\ne0$,
choosing $v=B_{a,m}/\|B_{a,m}\|_2$ gives
$v^\T B_{a,m}=\|B_{a,m}\|_2$. Dividing by $\mfa_m$ yields
the four identities below:
$$
\begin{array}{c|c}
    \varphi_m(\beta) & \langle\mu_m,\varphi_m\rangle\\[2pt]
    \hline
    \chi_{0,m}(\beta) & F_{0,m}/\mfa_m\\[3pt]
    \operatorname{sign}(A_{a,m})\chi_{a,m}(\beta)
        & |A_{a,m}|/\mfa_m\\[3pt]
    \displaystyle\chi_{a,m}(\beta)
       {B_{a,m}^\T(\beta-\beta_{a,m}^*)\over\|B_{a,m}\|_2}
        & \|B_{a,m}\|_2/\mfa_m\\[7pt]
    \chi_{a,m}(\beta)\|\beta-\beta_{a,m}^*\|_2^2
        & C_{a,m}/\mfa_m
\end{array}
$$
Choose the row corresponding to the positive limit $\gamma$.
In the second or third row, the selected coefficient is nonzero for
all sufficiently large $m$. Pass to a further subsequence so that its
sign is constant in the second row, or its normalized vector
$B_{a,m}/\|B_{a,m}\|_2$ converges in the third row, using compactness of
the unit sphere. Since the cutoffs are translates of a fixed $C^2$
function and their centers converge, the
selected functions converge in $C^2(\overline U)$ to a function
$\varphi_\infty$. The last assertion of
\cref{lem:app_compact_functionals} gives
$$
    \langle\mu,\varphi_\infty\rangle
    =\lim_{m\to\infty}\langle\mu_m,\varphi_m\rangle
    =\gamma>0.
$$

Since $\varphi_\infty$ need only be $C^2$, we next show that $\mu$ is nonzero on some smooth, compactly supported test function. To construct such function based on $\varphi_\i$, first choose a cutoff $\xi\in C_c^\infty(U)$ equal to one on an open neighborhood of $\cB_r$. Since $\mu_m$ is supported on $\cB_r$,
$
    \langle\mu_m,\xi\varphi_\infty\rangle
    =\langle\mu_m,\varphi_\infty\rangle.
$
Passing to the limit gives $\langle\mu,\xi\varphi_\infty\rangle=\gamma$. Extend $\xi\varphi_\infty$ by zero outside $U$, obtaining a function in $C_c^2(\RR^d)$. Then choose a nonnegative $\omega\in C_c^\infty(\cB_1^\circ)$ with $\int_{\RR^d}\omega(z)\,\rd z=1$, and set
$$
    \varphi_\epsilon(\beta) :=\int_{\RR^d}\omega(z)
          (\xi\varphi_\infty)(\beta-\epsilon z)\,\rd z ={1\over \epsilon^d}\int_{\RR^d}
          \omega\!\left(\frac{\beta-y}{\epsilon}\right)
          (\xi\varphi_\infty)(y)\,\rd y.
$$
The second expression shows that $\varphi_\epsilon$ is smooth, since derivatives of every order can be taken on the kernel $\omega$ under the integral. Its support lies within distance $\epsilon$ of $\operatorname{supp}(\xi\varphi_\infty)\subset U$, so it remains compactly supported in $U$ for sufficiently small $\epsilon>0$. Uniform continuity of the derivatives of $\xi\varphi_\infty$ through order two gives
$$
    \|\varphi_\epsilon-\xi\varphi_\infty\|_{C^2}\to0.
$$
The order-two bound therefore implies $\langle\mu,\varphi_\epsilon\rangle\to\gamma>0$. Thus $\mu$ is nonzero on a smooth, compactly supported test function, proving that $\mu\ne0$ as a generalized function.

For every fixed
$x$, the function $\beta\mapsto\sigma(x^\T\beta)$ belongs to
$C^2(\overline U)$, so
$$
    (T\mu_m)(x)
    ={p_{G_m}(x)-p_{G_m^*}(x)\over \mfa_m}
    \to(T\mu)(x).
$$
Fatou's lemma and the common design density lower bound give
\begin{align*}
    \int_{\cB_{r_X}}|(T\mu)(x)|^2\,\rd x
    &\le\liminf_{m\to\infty}
        \int_{\cB_{r_X}}|(T\mu_m)(x)|^2\,\rd x\\
    &\le {1\over m_X}\lim_{m\to\infty}
        {P_{X,m}(p_{G_m}-p_{G_m^*})^2\over \mfa_m^2}
      =0,
\end{align*}
where the last equality uses
\eqref{eq:app_inverse_contradiction_sequence}. The map
$x\mapsto\sigma(x^\T\,\cdot)$ is continuous in
$C^2(\overline U)$, so $T\mu$ is continuous and vanishes everywhere
on this ball. Applying \cref{lem:app_sigmoid_kernel_second_order} yields
\begin{equation}
    \mu=\check\mu,
    \qquad
    \langle\mu,1\rangle=0.
    \label{eq:app_inverse_limit_kernel}
\end{equation}

Finally, we show that $\mu$ is nonnegative. Since its total mass is
zero, this will contradict $\mu\ne0$.
Set
$$
    S=\{\beta_{1,\infty}^*,\ldots,
        \beta_{K^*,\infty}^*\}.
$$
Reflection separation gives
$\operatorname{dist}(S,-S)\ge s_+$. If
$\varphi\in C_c^\infty(\RR^d)$ is nonnegative and supported away from
$S$, then all true atoms lie outside its support for sufficiently large
$m$. Hence
$$
    \langle\mu,\varphi\rangle
    =\lim_{m\to\infty}{1\over \mfa_m}
       \left(\int\varphi\,\rd G_m-\int\varphi\,\rd G_m^*\right)
    =\lim_{m\to\infty}{1\over \mfa_m}\int\varphi\,\rd G_m
    \ge0.
$$
For a nonnegative test function $\varphi$ supported in
$B(\beta_{a,\infty}^*,s_+/4)$, its reflection
$\check\varphi(\beta)=\varphi(-\beta)$ is supported in
$B(-\beta_{a,\infty}^*,s_+/4)$, which is separated
from $S$. Thus $\langle\mu,\check\varphi\rangle\ge0$, and
\eqref{eq:app_inverse_limit_kernel} gives
$$
    \langle\mu,\varphi\rangle
    =\langle\mu,\check\varphi\rangle\ge0.
$$

To apply these inequalities to an arbitrary nonnegative test function,
choose smooth functions $0\le\zeta_a\le1$ equal to one on
$B(\beta_{a,\infty}^*,R_0)$, with compact support in
$B(\beta_{a,\infty}^*,2R_0)$, for $a\in[K^*]$.
These outer balls are disjoint and $2R_0<s_+/4$.
Thus $\zeta_0=1-\sum_{a=1}^{K^*}\zeta_a$ is nonnegative and vanishes on an
open neighborhood of $S$.
For any nonnegative $\varphi\in C_c^\infty(\RR^d)$, the test
$\zeta_0\varphi$ is supported away from $S$, and each
$\zeta_a\varphi$, $a\ge1$, is supported in the corresponding ball.
Since these tests are nonnegative and sum to $\varphi$,
$$
    \langle\mu,\varphi\rangle
    =\sum_{a=0}^{K^*}\langle\mu,\zeta_a\varphi\rangle\ge0.
$$
Thus $\mu$ is a nonnegative generalized function with compact support and,
by \eqref{eq:app_inverse_limit_kernel}, zero total mass. Applying
\cref{lem:app_positive_generalized_function} gives $\mu=0$, contrary
to $\mu\ne0$. This proves
\eqref{eq:app_population_inverse}. All subsequence arguments used only
the fixed radius, separation, and design bounds, and $K^*\le K_0$.
Thus the same positive $c_{\mathrm{inv}}$ applies to every admissible
$G^*$, $P_X$, and $G$, without a restriction on the support size
of $G$.
\end{proof}

\subsection{Dependence of the quantity $c_{\mathrm{inv}}$ on  $s_+$}
The quantity $c_{\mathrm{inv}}$ in \cref{thm:app_population_inverse} necessarily deteriorates as
$s_+$ decreases. To see this, take $d=1$, $r=2$,
$P_X=\mathrm{Unif}[-1,1]$, and
$$
    G_\varepsilon^*
    =\frac{1}{2}\delta_{1+\varepsilon}
     +\frac{1}{2}\delta_{-1+\varepsilon},
    \qquad G=\delta_0,
    \qquad 0<\varepsilon<\frac12.
$$
Here $s_+=2\varepsilon$,
$\mfa(G,G_\varepsilon^*)=2$, and $p_G(x)=1/2$. Since
$\{\sigma(x)+\sigma(-x)\}/2=1/2$, the fundamental theorem of
calculus gives
$$
    p_{G_\varepsilon^*}(x)-\frac12
    ={\varepsilon x\over2}\int_0^1
       \{\sigma'((1+t\varepsilon)x)
          +\sigma'((-1+t\varepsilon)x)\}\,\rd t.
$$
The bound $\sup_t\sigma'(t)=1/4$ and
$\EE X^2=\frac12\int_{-1}^1x^2\,\rd x=1/3$ imply
$$
    \{P_X(p_G-p_{G_\varepsilon^*})^2\}^{1/2}
    \le {\varepsilon\over4}(\EE X^2)^{1/2}
    =\frac{\varepsilon}{4\sqrt3}.
$$
Consequently, every inverse constant valid for this family must satisfy
\begin{equation}
    c_{\mathrm{inv}}
    \le \frac{s_+}{16\sqrt3}.
    \label{eq:app_population_inverse_upper_example}
\end{equation}

Unfortunately, the linear dependence on $s_+$ is not the worst when the support size of
$G$ is unrestricted. To see this, take $d=2$ and let $P_X$ be
uniform on the unit disk. Write $e_1,e_2$ for the canonical basis of $\RR^2$,
fix an integer $q\ge0$ and a constant $0<a<r-\tau_0$,
and consider the two-atom mixing distribution
$$
    G_\varepsilon^*
    ={1\over2}\delta_{ae_1+\varepsilon e_2}
     +{1\over2}\delta_{-ae_1+\varepsilon e_2},
    \qquad \varepsilon>0.
$$
For all sufficiently small $\varepsilon$, the true atoms lie in
$\cB_{r-\tau_0}$, their distance is $2a$, and their minimum
reflection distance is $s_+=2\varepsilon$. Thus all other model
constants can be held fixed as $\varepsilon$ tends to zero. Define
\begin{equation}\label{def_Q_delta}
    Q_\delta(x)
    ={1\over2}\{\sigma(ax_1+\delta x_2)
                  +\sigma(-ax_1+\delta x_2)\}
\end{equation}
such that
$$
    p_{G_\varepsilon^*}(x) =  Q_\varepsilon(x).
$$
Since $Q_{-\delta}(x)=1-Q_\delta(x)$, the function
$\delta\mapsto Q_\delta(x)-1/2$ is odd. Taylor's formula at zero
through order $2q+2$ therefore gives
\begin{equation}\label{taylor_Q_delta}
    Q_\delta(x)-\frac12
    =\sum_{j=0}^q
       \frac{\partial_\delta^{2j+1}Q_0(x)}{(2j+1)!}
       \delta^{2j+1}+R_q(\delta,x),
\end{equation}
where, uniformly over $\|x\|_2\le1$,
\begin{equation}\label{bd_remainder_Q_delta}
    |R_q(\delta,x)|
    ~ \le ~ \frac{\|\sigma^{(2q+3)}\|_\infty}{(2q+3)!}
       |\delta|^{2q+3}.
\end{equation}
Here the sup-norm of $\sigma^{(2q+3)}$ is taken over $\RR$.
Indeed,
$$
    \partial_\delta^kQ_\delta(x)
    =\frac{x_2^k}{2}
       \{\sigma^{(k)}(ax_1+\delta x_2)
          +\sigma^{(k)}(-ax_1+\delta x_2)\},
$$
so $|\partial_\delta^kQ_\delta(x)|\le\|\sigma^{(k)}\|_\infty$
when $\|x\|_2\le1$, and all positive even-order derivatives at $\delta=0$
vanish by oddness of $Q_\delta(x) - 1/2$.

We next construct a
positive measure $\nu_q$, with total mass $1/2$, such that
\begin{equation}\label{nu_odd_moment}
    \int \lambda^{2j+1}\,\rd \nu_q(\lambda)=1,
    \qquad j=0,\ldots,q.
\end{equation}
For $R>1$ to be specified below, take $t_\ell=R^{2(\ell+1)}$ and the Lagrange
interpolation weights
$$
    d_\ell=\prod_{k\ne\ell}{1-t_k\over t_\ell-t_k},
    \qquad \ell=0,\ldots,q.
$$
They satisfy $\sum_{\ell=0}^q d_\ell t_\ell^j=1$ for
$j=0,\ldots,q$, by evaluating the interpolation formula for the
polynomial $t\mapsto t^j$ at $t=1$. Set
$z_\ell=\sqrt{t_\ell}$ and place mass
$\lvert d_\ell\rvert/z_\ell$ at
$\operatorname{sign}(d_\ell)z_\ell$. If $R\ge\sqrt2$, each factor
in $d_\ell$ satisfies
$$
    \left|{1-t_k\over t_\ell-t_k}\right|
    \le
    \begin{cases}
      \displaystyle {t_k\over t_\ell-t_k}
        ={1\over R^{2(\ell-k)}-1}\le1, ~~ & k<\ell,\\[8pt]
      \displaystyle {t_k\over t_k-t_\ell}
        ={1\over1-R^{-2(k-\ell)}}\le2, ~~& k>\ell.
    \end{cases}
$$
Thus $|d_\ell|\le2^q$, and
$$
    \sum_{\ell=0}^q{|d_\ell|\over z_\ell}
    \le{(q+1)2^q\over R}.
$$
Choose
$$
    R = 2^{q+2}(q+1) .
$$
Then
$\sum_{\ell=0}^q|d_\ell|/z_\ell\le 1/4$.
Place mass $1/2-\sum_{\ell=0}^q|d_\ell|/z_\ell$ at zero so that
$\nu_q(\RR)=1/2$. Its odd moments satisfy
$$
    \int\lambda^{2j+1}\,\rd \nu_q(\lambda)
    =\sum_{\ell=0}^q {|d_\ell|\over z_\ell}
       \{\operatorname{sign}(d_\ell)z_\ell\}^{2j+1}
    =\sum_{\ell=0}^q d_\ell t_\ell^j=1,
    \qquad 0\le j\le q,
$$
as desired.

For this $\nu_q$, define
$$
    G_{\varepsilon,q}
    =\int\left\{
       {1\over2}\delta_{ae_1+\lambda\varepsilon e_2}
       +{1\over2}\delta_{-ae_1+\lambda\varepsilon e_2}
      \right\}\,\rd\nu_q(\lambda)+{1\over 2}\delta_0.
$$
Its total mass is $\nu_q(\RR)+1/2=1$.
For $0<\varepsilon\le\sqrt{r^2-a^2}/R^{q+1}$, its nonzero atoms have norm at most
$\sqrt{a^2+R^{2q+2}\varepsilon^2}\le r$, hence $G_{\varepsilon,q}$ is supported on $\cB_r$. Its induced regression function is
\begin{align*}
    p_{G_{\varepsilon,q}}(x) &=  \int \left\{
     {1\over2}\sigma(ax_1+\lambda\varepsilon x_2)
       +{1\over2}\sigma(-ax_1+\lambda\varepsilon x_2)
    \right\}\,\rd\nu_q(\lambda) + {1\over 4}\\
    &= \int Q_{\lambda\varepsilon}(x)\,\rd\nu_q(\lambda) + {1\over 4} &&\text{by \eqref{def_Q_delta}}.
\end{align*}
By using \eqref{taylor_Q_delta} twice for both $p_{G_\varepsilon^*} = Q_\varepsilon$ and $p_{G_{\varepsilon,q}}$, we have
\begin{align*}
    p_{G_{\varepsilon,q}}(x)-p_{G_\varepsilon^*}(x) &= \int \left\{{1\over 2} + \sum_{j=0}^q
       \frac{\partial_\delta^{2j+1}Q_0(x)}{(2j+1)!}
       (\lambda \varepsilon)^{2j+1}+R_q(\lambda \varepsilon,x)\right\}\,\rd\nu_q(\lambda) + {1\over 4}\\
       &\qquad -  {1\over 2 } - \sum_{j=0}^q
       \frac{\partial_\delta^{2j+1}Q_0(x)}{(2j+1)!}
       \varepsilon^{2j+1}-R_q(\varepsilon,x)\\
       &=  \sum_{j=0}^q
       \frac{\partial_\delta^{2j+1}Q_0(x)}{(2j+1)!} \varepsilon^{2j+1}\left(
       \int     \lambda^{2j+1}\,\rd\nu_q(\lambda)   - 1\right)\\
       &\qquad + \int R_q(\lambda \varepsilon,x)\,\rd\nu_q(\lambda)-R_q(\varepsilon,x).
\end{align*}
The second equality uses $\nu_q(\RR)=1/2$, so that the constant terms
sum to $\nu_q(\RR)/2+1/4-1/2=0$.
Invoking \eqref{nu_odd_moment} and \eqref{bd_remainder_Q_delta} therefore gives
$$
\begin{aligned}
    \sup_{\|x\|_2\le1}\left|p_{G_{\varepsilon,q}}(x)-p_{G_\varepsilon^*}(x)\right| &\le {\|\sigma^{(2q+3)}\|_\infty\over(2q+3)!}
       \left(1+\int|\lambda|^{2q+3}\,\rd \nu_q(\lambda)\right)
       \varepsilon^{2q+3}\\
    &\le \left({2\over\pi}\right)^{2q+3}
          \left\{1+{1\over 2}R^{(q+1)(2q+3)}\right\} \varepsilon^{2q+3}.
\end{aligned}
$$
Since $\|X\|_2\le1$ almost surely, this also bounds
$\{P_X(p_{G_{\varepsilon,q}}-p_{G_\varepsilon^*})^2\}^{1/2}$.
The complex sigmoid bound in the proof of \cref{lem:app_sigmoid_kernel} gives
$|\sigma(z)|\le1$ on every circle $|z-t|=\pi/2$, $t\in\RR$.
Cauchy's derivative estimate therefore yields
$\|\sigma^{(k)}\|_\infty\le k!(2/\pi)^k$, while
$\int|\lambda|^k\,\rd \nu_q(\lambda)\le R^{k(q+1)}/2$.
Moreover, $\chi_0(0)=1$ for all sufficiently small $\varepsilon$, so
the mass $(1/2)\delta_0$ gives
$\mfa(G_{\varepsilon,q},G_\varepsilon^*)\ge1/2$.

Consequently, \eqref{eq:app_population_inverse} gives
$$
    c_{\mathrm{inv}}
    \le 2\left({2\over\pi}\right)^{2q+3}
       \left\{1+{1\over2}R^{(q+1)(2q+3)}\right\}
       \left({s_+\over2}\right)^{2q+3},
$$
where we used $\varepsilon=s_+/2$.
Substituting $R=2^{q+2}(q+1)$ gives, for every $q\ge0$,
\begin{equation}
    c_{\mathrm{inv}}
     ~ \le ~  \frac{2+\{2^{q+2}(q+1)\}^{(q+1)(2q+3)}}{\pi^{2q+3}}
       s_+^{2q+3}.
    \label{eq:app_population_inverse_arbitrary_order}
\end{equation}
For any fixed $\alpha<\infty$, choose $q$ with $2q+3>\alpha$.
Then \eqref{eq:app_population_inverse_arbitrary_order} gives
$c_{\mathrm{inv}}/s_+^\alpha\to0$ as $s_+\to0$.
Hence no positive constant $c_\alpha$, independent of $s_+$, can
give $c_{\mathrm{inv}}\ge c_\alpha s_+^\alpha$ uniformly over
probability measures $G$ of unrestricted support size. Notice that the true mixing distribution in this
construction remains discrete with two atoms whereas the support size of
$G_{\varepsilon,q}$  increases with
$q$.

\section{Proof of Theorem~\ref{thm:particle_population}: the population-level  convergence of the finite particle approximation}
\label{app:particle_uniform_analysis}

Throughout this section, we work under \cref{assum:particle_population}
and \eqref{eq:particle_radius_choice}, with $d,r,\tau_0$ and the
assumption constants fixed. In particular, every true atom belongs to
$\cB_{r-\tau_0}$.

Fix any $t\ge 0$. Recall that the mixing distribution at the beginning of the $t$-th  outer iteration is
$$
    G_t = \sum_{j=1}^M w_j^t\delta_{\beta_j^t}
$$
with  $w^t\in\cM(\bbeta^t)$. To enter the $H$ inner steps of updating locations, let $\bbeta^{t,0} := \bbeta^t$ and $G^{t,0} := G_t$. The $H$ location updates yield the particles $\bbeta^{t,1},\ldots, \bbeta^{t,H}$ and we write $\bbeta^{t+1} := \bbeta^{t,H}$, yielding
\begin{equation}\label{def_G_tH}
    G^{t,H}  = \sum_{j=1}^M w_j^t\delta_{\beta_j^{t,H}} = \sum_{j=1}^M w_j^t\delta_{\beta_j^{t+1}}.
\end{equation}
The exact
reweighting then updates the mixing weights $w^t$ to   $w^{t+1}\in\cM(\bbeta^{t+1})$ so that  the final mixing distribution by the end of the $t$-th iteration is
$$
    G_{t+1} = \sum_{j=1}^M w_j^{t+1}\delta_{\beta_j^{t+1}}.
$$

\paragraph{Proof strategy for \cref{thm:particle_population}}

To prove the convergence of $\cL(G_{t+1}) - \cL(G_t)$, since the definition of exact reweighting yields
$$
    \cL(G_{t+1}) - \cL(G_{t}) ~ \ge ~   \cL(G^{t,H}) -  \cL(G_{t}) ~ = ~ \cL(G^{t,H}) -  \cL(G^{t,0}),
$$
it is critical to bound from below the increment $\cL(G^{t,H}) -  \cL(G^{t,0})$ due to location updates. The following lemma establishes such result. Since the proof is rather complicated, we present and explain it separately in \cref{app:proof:prop:app_population_block_gain}.

\begin{proposition}
\label{prop:app_population_block_gain}
There are some constants $\kappa\in(0,1)$ and $\epsilon_0>0$, depending only on $d,r,\tau_0$ and the
constants in \cref{assum:particle_population}, such that the following
holds for any $t\ge 0$.
For any integer
$H\ge1$, let $G^{t,H}$ be defined in \eqref{def_G_tH} using the location updates \eqref{eq:reweighted_particle_update} with the step size satisfying $0<\eta\le27\sigma(-rR_X)/(4R_X^2)$.
If $\cL(G^*)-\cL(G_t)\le 2\epsilon_0^2$ holds, then
\begin{equation}
    \left\{\cL(G^{t,H})-\cL(G^{t,0})\right\}
    ~ \ge~ \kappa \left\{\cL(G^*)-\cL(G^{t,0})\right\}
       \min\left\{\eta H\sqrt{\cL(G^*)-\cL(G^{t,0})},1\right\}.
    \label{eq:app_population_block_gain}
\end{equation}
\end{proposition}
For choosing $M\ge M_0$ large,
\cref{lem:app_population_initialization} in \cref{app:sec:lemmas_init_recursive_population} ensures that  with high probability, the gap $\cL(G^*) - \cL(G_0)$ can be made small, which together with the monotonicity  of $\cL(G_t)$ ensures the requirement in \cref{prop:app_population_block_gain}.
Combining these results with the following recursion bound, proved in
\cref{app:sec:lemmas_init_recursive_population}, gives the claimed convergence
rate in \cref{thm:particle_population}.

\begin{lemma}
\label{lem:app_particle_saturated_recursion}
Let $\{a_t\}_{t\ge0}$ be nonnegative and nonincreasing, with
$a_0\le\Gamma$ for some $\Gamma>0$. Suppose that, for constants
$b\ge0$, $s>0$, and $\kappa\in(0,1)$,
$$
    a_t>b
    \qquad\Longrightarrow\qquad
    a_{t+1}\le a_t-\kappa a_t\min\{s\sqrt{a_t},1\}.
$$
Then there are constants $C,c>0$, depending only on $\Gamma$ and
$\kappa$, such that
\begin{equation}
    a_t\le C\{(1+st)^{-2}+e^{-ct}\}+b,
    \qquad t\ge0.
    \label{eq:app_saturated_recursion_rate}
\end{equation}
\end{lemma}

\subsection{Proof of \cref{thm:particle_population}}
\label{app:particle_population_proof}

\begin{proof}
Take $\epsilon_0$ and $\kappa$ from \cref{prop:app_population_block_gain} and
work on the event
$$
    \cE_0 := \left\{\cL(G^*)-\cL(G_0)\le \epsilon_0^2\right\}
$$ which, according to  \cref{lem:app_population_initialization}, holds with probability at least $1-\delta$ for every $M\ge M_0$. Fix any $t\ge 0$.
Since the exact reweighting ensures $\cL(G_{t+1}) \ge \cL(G^{t,H})$,  we obtain
\begin{align}\nonumber
    \cL(G^*)-\cL(G_{t+1})~ &\le ~  \cL(G^*)-\cL(G^{t,H})\\\nonumber
    ~&\le ~ \cL(G^*)-\cL(G^{t,0}) &&\text{by \cref{prop:app_population_block_gain}}\\
    ~&=~ \cL(G^*)-\cL(G_t).\label{monotone_G_t}
\end{align}
Consequently, $\cE_0$ ensures
$$
    \cL(G^*)-\cL(G_t) ~ \le ~ \epsilon_0^2,\qquad t\ge 0.
$$
One can then apply
\cref{prop:app_population_block_gain} to bound from below $\cL(G^{t,H}) - \cL(G^{t,0})$:   for any $t\ge 0$, similar arguments as in \eqref{monotone_G_t} yield
\begin{align*}
    \cL(G^*)-\cL(G_{t+1})
    &\le\cL(G^*)-\cL(G^{t,H})\\
    &\le\{\cL(G^*)-\cL(G_t)\}-\kappa\{\cL(G^*)-\cL(G_t)\}\\
    &\qquad\qquad\qquad\qquad\times
       \min\{\eta H\sqrt{\cL(G^*)-\cL(G_t)},1\}.
\end{align*}
Now apply
\cref{lem:app_particle_saturated_recursion} with
$a_t=\{\cL(G^*)-\cL(G_t)\}$, $s=\eta H$, $b=0$, and
$\Gamma=\Gamma_0/2$ to obtain
$$
    \cL(G^*)-\cL(G_t)
    ~ \le~ C\{(1+\eta Ht)^{-2}+e^{-ct}\},\qquad t\ge0.
$$
Finally, the bound for $P_X(p_{G_t} - p^*)^2$  follows from \eqref{eq:app_population_kl_l2}.

The constants $C,c$ depend only on $\Gamma_0,\kappa$, hence only on
$d,r,\tau_0$ and the constants in \cref{assum:particle_population}.
The initialization threshold $M_0$ may additionally depend on
$\delta$; none of these quantities depends on $M,H$, or $\eta$.
\end{proof}

\subsection{Proof of \cref{rem:particle_converged_locations}}
\label{app:particle_reweighting_remark_proof}

\begin{proof}
Take $C,c$ from \cref{thm:particle_population} and fix $\eta>0$.
For $Ht\ge2$, choosing $t\ge(2/c)\log(Ht)$ gives
$$
    e^{-ct}\le(Ht)^{-2},
    \qquad
    (1+\eta Ht)^{-2}\le\eta^{-2}(Ht)^{-2}.
$$
On the event in \cref{lem:app_population_initialization}, the theorem
therefore yields
$$
    \cL(G^*)-\cL(G_t)
       \le C(1+\eta^{-2})(Ht)^{-2}.
$$
Such choices can be made with $t$ of order $\log(Ht)$. For example,
take $t=\lceil(4/c)\log H\rceil$ for all sufficiently large integers
$H$. Then $t\le H$, so
$(2/c)\log(Ht)\le(4/c)\log H\le t$.
Since $G_t$ uses $Ht$ location updates and $t+1$ exact
reweightings, this proves the first assertion. The initialization
event and its particle threshold do not depend on $H$ or $t$.

For the assertion with location updates run to convergence, work on the
event in \cref{lem:app_population_initialization}.
With the weights held fixed, the gap along the location iterates is
nonincreasing in $h$ and nonnegative, so it has a limit. If the starting gap is
positive, the minimum in \eqref{eq:app_population_block_gain} equals one
for all sufficiently large $H$. Hence
$$
    \lim_{H\to\infty}\{\cL(G^*)-\cL(G_{\bbeta^{(h)}})\}
    \le(1-\kappa)\{\cL(G^*)-\cL(G_{\bbeta^{(0)}})\}.
$$
The feasible location set is compact. By continuity of $\cL(G_{\bbeta^{(h)}})$ in $\bbeta^{(h)}$, every accumulation
point of the iterates has this limiting gap. Since
reweighting does not  increase the gap, this proves the claim for $\cL(G^*) > \cL(G_{\bbeta^{(0)}})$. When $\cL(G^*) = \cL(G_{\bbeta^{(0)}})$, the result follows
from monotonicity, and repeating the same argument.
\end{proof}

\subsection{Initialization and  proof of the recursive lemma}\label{app:sec:lemmas_init_recursive_population}

We next show that the likelihood gap of $G_0$ is small with high
probability.

\begin{lemma}
\label{lem:app_population_initialization}
Let $\epsilon_0$ be specified in \eqref{def_epsilon_0}.
For every $\delta\in(0,1)$, there is a finite $M_0$, depending only on
$d,r,\tau_0$, the constants in \cref{assum:particle_population}, and $\log(1/\delta)$, such
that the following holds whenever $M\ge M_0$. If
$\beta_1^0,\ldots,\beta_M^0$ are independent and uniformly distributed
on $\cB_r$, and $w^0\in\cM(\bbeta^0)$, then  for $G_0 = \sum_{j=1}^M w_j^0 \delta_{\beta_j^0}$, the following holds
with probability at least $1-\delta$,
$$
    \cL(G^*)-\cL(G_0)\le \epsilon_0^2
$$
\end{lemma}

\begin{proof}
It suffices that each true atom has an initial particle close to it:
assigning the true weights to those particles gives a small likelihood
gap, and exact reweighting can only decrease that gap. We first determine
the required distance, then bound the probability of this event.

By the integral form of the mean value theorem,
$$
    |p_{\delta_\beta}(x)-p_{\delta_{\wt\beta}}(x)|
    =\left|\int_0^1\sigma'(x^\T(\wt\beta+t(\beta-\wt\beta)))
                   x^\T(\beta-\wt\beta)\,\rd t\right|
    \overset{\eqref{eq:app_population_location_derivatives}}{\le} \frac{R_X}{4}\|\beta-\wt\beta\|_2.
$$
If the initial configuration contains a
particle $\wt\beta_a$ within distance $h$ of each true atom, then
the feasible candidate measure
$$
    \wt G
    =\sum_{a=1}^{K^*}w_a^*\delta_{\wt\beta_a}
$$
has prediction $\wt p=p_{\wt G}$ satisfying
\begin{equation}\label{lip_p_beta}
    \{P_X(\wt p-p^*)^2\}^{1/2}
    \le\|\wt p-p^*\|_\infty
    \le\sum_{a=1}^{K^*} w_a^* \frac{R_X}{4}\|\wt\beta_a-\beta_a^*\|_2
    \le \frac{R_X}{4}h.
\end{equation}
Since $f_{p^*(x)}'(p^*(x))=0$, Taylor's integral formula gives
\begin{align*}
    \cL(G^*)-\cL(\wt G)
    &=\EE_X\int_0^1(1-s)
       f_{p^*(X)}''\{p^*(X)+s(\wt p(X)-p^*(X))\}
                       (\wt p(X)-p^*(X))^2\,\rd s\\
    &\le {1\over2\rho^2}P_X(\wt p-p^*)^2
     \le {R_X^2\over32\rho^2}h^2.
\end{align*}
The inequalities follow from \eqref{eq:app_population_scalar_loss_bounds}
and \eqref{lip_p_beta}, respectively. Choose
$$
    h_0={\tau_0\over2}\wedge{s_-\over4}
        \wedge\frac{4\sqrt{2}\rho\, \epsilon_0}{R_X}.
$$
Each ball $\cB(\beta_a^*,h_0)$ is contained in $\cB_r$, and one
uniform particle lies in it with probability $(h_0/r)^d$, by the
ratio of the two ball volumes. Independence and a union bound give
\begin{align*}
 \PP\{\exists\ a\in[K^*]\text{ such that } \cB(\beta_a^*,h_0)\text{ contains no initial particle}\}
    &     \le K^*\{1-(h_0/r)^d\}^M\\
    &  \le K_0\exp\{-M(h_0/r)^d\}.
\end{align*}
The last bound is at most $\delta$ when
$$
    M\ge M_0
    :=K_0\vee
      \left\lceil
      \left({r\over h_0}\right)^d\log{K_0\over\delta}
      \right\rceil.
$$
Thus, with probability at least $1-\delta$, every ball contains
an initial particle. The balls are disjoint because $2h_0<s_-$,
so these particles are distinct. Assigning weight $w_a^*$ to the
particle in the $a$th ball and zero weight to all others makes
$\wt G$ feasible on the initial locations. Since
$w^0\in\cM(\bbeta^0)$, we have
$$
    \cL(G^*)-\cL(G_0)
    \le\cL(G^*)-\cL(\wt G)
    \le\frac{R_X^2h_0^2}{32\rho^2}
    \le \epsilon_0^2,
$$
which completes the proof.
\end{proof}

We now prove the recursive \cref{lem:app_particle_saturated_recursion}.

\begin{proof}
We apply the recursion to $a_j^{-1/2}$. This quantity increases
geometrically while it is below $s$, and at least linearly thereafter,
giving the two terms in \eqref{eq:app_saturated_recursion_rate}.

Once the sequence reaches $b$, it remains at most $b$ by
monotonicity. It therefore suffices to consider an index $t$ for which
$a_t>b$. Then $a_j>0$ for $0\le j\le t$. Write
$x_j=a_j^{-1/2}$. The inequality
$(1-z)^{-1/2}\ge1+z/2$, for $0\le z<1$, gives
$$
    x_{j+1}
    \ge x_j\{1-\kappa\min(s/x_j,1)\}^{-1/2}
    \ge x_j+{\kappa\over2}\min\{s,x_j\},
    \qquad 0\le j<t.
$$
Thus $x_j$ grows geometrically while it is below $s$, and by at
least $\kappa s/2$ per iteration thereafter.

For $t\ge2$, set $m=\lfloor t/2\rfloor$. If $x_m\ge s$, the
last $t-m$ iterations give
$$
    x_t\ge\Gamma^{-1/2}+{\kappa s\over2}(t-m)
       \ge\Gamma^{-1/2}+{\kappa st\over4} \ge \min\{\Gamma^{-1/2},\kappa/4\}(1+st).
$$
Hence
$$
    a_t\le \max\{\Gamma,16/\kappa^2\}(1+st)^{-2}.
$$
If $x_m<s$, then $x_j<s$ for
every $j\le m$, so
$$
    x_t\ge x_m\ge\Gamma^{-1/2}(1+\kappa/2)^m
$$
and, by $2m\ge t-1$,
$$
    a_t\le\Gamma(1+\kappa/2)^{-2m} \le  \Gamma(1+\kappa/2)^{1-t} = \Gamma(1+\kappa/2)e^{-ct}
$$
with
$c=\log(1+\kappa/2)$. Taking
$C=\max\{16/\kappa^2,\Gamma(1+\kappa/2)\}$ also covers
$t=0,1$, since $Ce^{-c}\ge\Gamma$.
This proves \eqref{eq:app_saturated_recursion_rate}, with $C,c$
independent of $s,b$.
\end{proof}

\subsection{Proofs and intermediate results of \cref{prop:app_population_block_gain}: block location updates improvement and ascent directions}\label{app:proof:prop:app_population_block_gain}

\subsubsection{Preliminaries on population-level likelihood function and its derivatives}

We use the notation of \cref{assum:particle_population} and write
$\rho=\sigma(-rR_X)$. For every $G\in\cP(\cB_r)$ and
$x\in\cB_{R_X}$, we have $\rho\le p_G(x)\le1-\rho$.
Since $p_{G^*}=p^*$, the
definition of $\cL$ in \eqref{eq:population_likelihood_gap} yields
\begin{equation}
    \cL(G^*)-\cL(G)
    =\EE_X\kl{\Bern(p^*(X))}{\Bern(p_G(X))},
    \label{eq:app_population_kl_identity}
\end{equation}
where $\kl{P}{Q}$ denotes the Kullback--Leibler divergence from $P$ to
$Q$. For $p,q\in[\rho,1-\rho]$, the inequality
$\log t\le t-1$ gives
$$
\begin{aligned}
    \kl{\Bern(p)}{\Bern(q)}
    &\le p\left({p\over q}-1\right)
       +(1-p)\left({1-p\over1-q}-1\right) ={ (p-q)^2\over q(1-q)}
     \le\rho^{-2}(p-q)^2.
\end{aligned}
$$
Combining this bound with Pinsker's inequality
\citep[Lemma 2.5]{tsybakov2009introduction} and taking expectation in
\eqref{eq:app_population_kl_identity} yields
\begin{equation}
    2P_X(p_G-p^*)^2
    \le\cL(G^*)-\cL(G)
    \le\rho^{-2}P_X(p_G-p^*)^2.
    \label{eq:app_population_kl_l2}
\end{equation}
The population counterpart of $g_{n,G}$ is obtained by replacing the sample
average by expectation and using $\EE(Y\mid X)=p^*(X)$. Define
\begin{equation}
    r_G(x)={p^*(x)-p_G(x)\over p_G(x)\{1-p_G(x)\}},
    \qquad
    g_G(\beta)=\EE_X\{r_G(X)\sigma(X^\T\beta)\}.
    \label{eq:app_population_residual_potential}
\end{equation}
Differentiating $g_G$ with respect to $\beta$, with $G$ fixed,
gives
$$
    \nabla g_G(\beta)
    =\EE_X\{r_G(X)\sigma'(X^\T\beta)X\}.
$$
If $G=\sum_{j=1}^Mw_j\delta_{\beta_j}$ with $w\in\cM(\bbeta)$, the
population version of \eqref{eq:app_particle_weight_kkt} gives
\begin{equation}
    g_G(\beta_j)\le\sum_{k=1}^Mw_kg_G(\beta_k),
    \qquad\text{with equality whenever }w_j>0.
    \label{eq:app_population_weight_kkt}
\end{equation}
Indeed, $g_G(\beta_j)$ is the partial derivative of $\cL(G)$
with respect to $w_j$, and $\cL(G)$ is concave in the weights.
\paragraph{Derivative bounds}

Write
$z=\sigma(t)$ and $v=z(1-z)$. Differentiation gives
$$
    \sigma'(t)=z(1-z),\qquad
    \sigma''(t)=z(1-z)(1-2z),\qquad
    \sigma'''(t)=v(1-6v).
$$
Since $0\le v\le1/4$, the first derivative is at most $1/4$.
The derivative of $z(1-z)(1-2z)$ in $z$ is
$1-6z+6z^2$, whose zeros are $(3\pm\sqrt3)/6$.
Evaluating there and at the endpoints gives the second bound below.
Finally, $v(1-6v)$ has maximum $1/24$ at $v=1/12$
and minimum $-1/8$ on $[0,1/4]$. Thus
\begin{equation}
    \|\sigma'\|_\infty={1\over4},\qquad
    \|\sigma''\|_\infty={1\over6\sqrt3},\qquad
    \|\sigma'''\|_\infty={1\over8}.
    \label{eq:app_population_scalar_derivatives}
\end{equation}
For a single atom, $p_{\delta_\beta}(x)=\sigma(x^\T\beta)$.
For fixed $x$, change the location from $\beta$ to $\beta+tu$,
where $u\in\RR^d$ is fixed and $t$ is a scalar. We write
$$
    D_\beta p_{\delta_\beta}[u](x)
    :=\left.{d\over dt}p_{\delta_{\beta+tu}}(x)\right|_{t=0}
    =\sigma'(x^\T\beta)x^\T u.
$$
The vector in brackets specifies the direction of the location change.
More generally, for $u_1,\ldots,u_k\in\RR^d$, we write
$D_\beta^k p_{\delta_\beta}[u_1,\ldots,u_k]$ for successive
differentiation with respect to $\beta$ in these $k$ directions,
with $x$ fixed. The chain rule
and $|x^\T u_i|\le R_X\|u_i\|_2$ give, for $k=1,2,3$,
\begin{equation}
    |D_\beta^k p_{\delta_\beta}[u_1,\ldots,u_k](x)|
    =|\sigma^{(k)}(x^\T\beta)|
       \prod_{i=1}^k|x^\T u_i|
    \le R_X^k\|\sigma^{(k)}\|_\infty
       \prod_{i=1}^k\|u_i\|_2.
    \label{eq:app_population_location_derivatives}
\end{equation}
These bounds are uniform over $\|x\|_2\le R_X$ and $\beta\in\RR^d$.

To control the Taylor expansions of the population likelihood, we use
the following derivative bounds. For $a\in[0,1]$ and $p\in(0,1)$, let
$f_a(p)=-a\log p-(1-a)\log(1-p)$.
Setting $a=p^*(x)$ and $p=p_G(x)$, we have
$$
    -\cL(G)=\EE_X\left[f_{p^*(X)}\{p_G(X)\}\right].
$$
For $\rho\le p\le1-\rho$, direct differentiation gives
\begin{align}
\label{eq:app_population_scalar_loss_bounds} \nonumber
    f_a'(p)&={p-a\over p(1-p)},
       &&|f_a'(p)|\le\min\left\{{1\over \rho},{|p-a| \over \rho^2}\right\},\\
    f_a''(p)&={a\over p^2}+{1-a\over(1-p)^2},
       &&0<f_a''(p)\le {1\over \rho^2},\\
    f_a'''(p)&=-{2a\over p^3}+{2(1-a)\over(1-p)^3},
       &&|f_a'''(p)|\le {2\over \rho^3}. \nonumber
\end{align}
The first bound uses the fact that
$f_a'(p)=(1-a)/(1-p)-a/p$ is a convex combination of
$1/(1-p)$ and $-1/p$. The other two bounds follow by bounding
each denominator below by the corresponding power of $\rho$.

\subsubsection{Likelihood gains from inner-loop location updates}

We first construct a configuration with unchanged weights and a smaller
likelihood gap. We then compare each location update with this
configuration. For $M$ particles with fixed weights
$w=(w_1,\ldots,w_M)^\T\in\Delta^{M-1}$ and locations
$\bbeta=(\beta_1,\ldots,\beta_M)\in\cB_r^M$, write
\begin{equation}\label{def_G_beta}
    G_{\bbeta}=\sum_{j=1}^Mw_j\delta_{\beta_j}.
\end{equation}
For location directions $\boldsymbol u=(u_1,\ldots,u_M)$ and
$\boldsymbol v=(v_1,\ldots,v_M)$, with $u_j,v_j\in\RR^d$, set
$$
    \langle\boldsymbol u,\boldsymbol v\rangle_w
       =\sum_{j=1}^Mw_j u_j^\T v_j,
    \qquad
    \|\boldsymbol u\|_w^2=\sum_{j=1}^Mw_j\|u_j\|_2^2.
$$
These define an inner product and its norm on the coordinates with
$w_j>0$. Coordinates with zero weight do not affect $\cL(G_{\bbeta})$ and may be
omitted from the fixed-weight calculations. The directional derivative of the likelihood is
\begin{equation}
    \left.{d\over ds}\cL(G_{\bbeta+s\boldsymbol u})\right|_{s=0}
       =\sum_{j=1}^M w_j\langle\nabla g_{G_{\bbeta}}(\beta_j),u_j\rangle.
    \label{eq:app_fixed_weight_directional}
\end{equation}
The population location update uses
$\nabla g_{G_{\bbeta}}(\beta_j)$, without a factor $w_j$.
Coordinatewise projection onto $\cB_r$ minimizes the weighted
squared distance. The
algorithm updates zero-weight particles as well, although their
locations do not affect the fixed-weight likelihood.

In the geometry of $\langle\cdot,\cdot\rangle_w$, the map
$\bbeta\mapsto\cL(G_{\bbeta})$ therefore has gradient
$\{\nabla g_{G_{\bbeta}}(\beta_j)\}_{j\in[M]}$, and each location
update is a projected gradient-ascent step in this geometry.
Recall from \cref{lem:particle_one_block_ascent} that, for fixed data and a
step size satisfying $0<\eta\le  8\rho_{n,r}^2 / R_{x,n}^2$, one location update
increases $\ell_n$ by at least $(\eta/2)E_t$, where $\eta^2E_t$ is the
weighted squared displacement of the particles. The next lemma, part $(ii)$ in particular, is the
population counterpart of this result and has the same form. Its proof is
organized as that of \cref{lem:particle_one_block_ascent}. The only
difference is the bound on the second-order term, which is obtained in one
calculation and yields a step-size condition of order $\rho$ rather than
$\rho^2$. Part (ii) follows from part (i). We state part (i) separately,
since the proof of \cref{prop:app_population_block_gain} uses it directly
to compare the iterates with a reference configuration. Recall that
$\rho=\sigma(-rR_X)$.

\begin{lemma}
\label{lem:app_population_one_step_ascent}
Fix $w\in\Delta^{M-1}$. Let $G_\bbeta$ be defined in \eqref{def_G_beta} for any $\bbeta\in \cB_r^M$.
\begin{enumerate}[label=(\roman*),noitemsep,leftmargin=7mm]
\item For any $\bbeta,\wt\bbeta\in\cB_r^M$,
\begin{equation}
    \cL(G_{\wt\bbeta})-\cL(G_{\bbeta})
    ~ \ge~ \sum_{j=1}^Mw_j\langle\nabla g_{G_{\bbeta}}(\beta_j),
       \wt\beta_j-\beta_j\rangle
      - {2R_X^2 \over 27\rho}\|\wt\bbeta-\bbeta\|_w^2
    \label{eq:app_fixed_weight_upper_model}
\end{equation}
and
\begin{equation}
    \cL(G_{\wt\bbeta})-\cL(G_{\bbeta})
    ~ \le~ \sum_{j=1}^Mw_j
       \langle\nabla g_{G_{\bbeta}}(\beta_j),\wt\beta_j-\beta_j\rangle
      +\frac{R_X^2  \{\cL(G^*)-\cL(G_{\bbeta})\}^{1/2}}{12\sqrt6\,\rho^2}
       \|\wt\bbeta-\bbeta\|_w^2.
    \label{eq:app_fixed_weight_lower_model}
\end{equation}
\item Let $0<\eta\le {27\rho / 4R_X^2}$ and
$\beta_j^+=\Pi_{\cB_r}\{\beta_j+\eta\nabla g_{G_{\bbeta}}(\beta_j)\}$
for $j\in[M]$. Then
\begin{equation}
    \cL(G_{\bbeta^+})-\cL(G_{\bbeta})
    ~ \ge ~ {1\over2\eta}\|\bbeta^+-\bbeta\|_w^2.
    \label{eq:app_fixed_weight_descent}
\end{equation}
\end{enumerate}
If $\max_{i\in[n]}\|X_i\|_2\le R_X$, then \eqref{eq:app_fixed_weight_upper_model} and \eqref{eq:app_fixed_weight_descent} also hold with
$\cL$ and $g_{G_{\bbeta}}$ replaced by $\ell_n$ and
$g_{n,G_{\bbeta}}$, respectively.
\end{lemma}

\begin{proof}
We follow the proof of \cref{lem:particle_one_block_ascent}.
Write $\delta_j=\wt\beta_j-\beta_j$ and
$\bdelta=(\delta_1,\ldots,\delta_M)$. Let
$$
    L_X := {4R_X^2 \over 27\rho}.
$$
Since $\cB_r$ is convex,
$\bbeta+s\bdelta\in\cB_r^M$ for all $0\le s\le1$. Let
$$
F(s):=\cL(G_{\bbeta+s\bdelta}),\qquad \text{ for } 0\le s\le 1.
$$
For every $x\in\cB_{R_X}$,
$p_{G_{\bbeta+s\bdelta}}(x)\in[\rho,1-\rho]$, and by
\eqref{eq:app_population_location_derivatives} its first two derivatives in
$s$ are bounded uniformly in $x$ and $s$. By
\eqref{eq:app_population_scalar_loss_bounds}, $F$ is therefore twice
continuously differentiable on $[0,1]$ and may be differentiated under the
expectation, and \eqref{eq:app_fixed_weight_directional} gives
$F'(0)=\sum_jw_j\langle\nabla g_{G_{\bbeta}}(\beta_j),\delta_j\rangle$.
Taylor's formula with integral remainder gives
\begin{equation}\label{init_diff_popu}
    \cL(G_{\wt\bbeta})-\cL(G_{\bbeta})
    =F(1)-F(0)
    =\sum_{j=1}^Mw_j\langle\nabla g_{G_{\bbeta}}(\beta_j),\delta_j\rangle
     +\int_0^1(1-s)F''(s)\,\rd s.
\end{equation}
Since $\int_0^1(1-s)\,\rd s=1/2$, to prove \eqref{eq:app_fixed_weight_upper_model},  it suffices to show
\begin{equation}\label{bd_second_popu}
    F''(s)\ge-L_X\|\bdelta\|_w^2,\qquad 0\le s\le1.
\end{equation}
To this end, fix $x\in\cB_{R_X}$, $y\in\{0,1\}$, and
$s\in[0,1]$, and consider the log-likelihood of a single observation along
the path,
$$
    \varphi(s)=y\log q(s)+(1-y)\log\{1-q(s)\},\qquad
    q(s)=\sum_{j=1}^Mw_j\sigma\{x^\T(\beta_j+s\delta_j)\}.
$$
Since $1-\sigma(t)=\sigma(-t)$, writing
$q_j(s)=\sigma\{(2y-1)x^\T(\beta_j+s\delta_j)\}$ and
$\bar q(s)=\sum_jw_jq_j(s)$ gives $\varphi(s)=\log\bar q(s)$ for both
values of $y$. With $t_j=(2y-1)x^\T\delta_j$, differentiation gives
$$
    \bar q'=\sum_{j=1}^Mw_jq_j(1-q_j)t_j,
    \qquad
    \bar q''=\sum_{j=1}^Mw_jq_j(1-q_j)(1-2q_j)t_j^2,
$$
where the dependence on $s$ is suppressed. Since the numbers
$w_jq_j/\bar q$, $j\in[M]$, are nonnegative and sum to one, the
Cauchy--Schwarz inequality gives
$$
    \left({\bar q'\over\bar q}\right)^2
    =\left\{\sum_{j=1}^M{w_jq_j\over\bar q}(1-q_j)t_j\right\}^2
    \le{1\over\bar q}\sum_{j=1}^Mw_jq_j(1-q_j)^2t_j^2.
$$
Consequently,
\begin{align}\nonumber
    \varphi''
    ={\bar q''\over\bar q}-\left({\bar q'\over\bar q}\right)^2
    &\ge{1\over\bar q}\sum_{j=1}^Mw_jq_j
       \{(1-q_j)(1-2q_j)-(1-q_j)^2\}t_j^2\\\nonumber
    &=-{1\over\bar q}\sum_{j=1}^Mw_jq_j^2(1-q_j)t_j^2\\
    &\ge- {4R_X^2\over27\rho}\sum_{j=1}^Mw_j\|\delta_j\|_2^2
    =-L_X\|\bdelta\|_w^2.
    \label{bd_second_pointwise}
\end{align}
The last inequality uses $\bar q\ge\rho$, $|t_j|\le R_X\|\delta_j\|_2$,
and $\max_{0\le z\le1}z^2(1-z)=4/27$; the maximum is attained at
$z=2/3$, since the derivative of $z^2(1-z)$ is $z(2-3z)$.
The two terms in $\varphi''$ are the curvature of the mixture $\bar q$
along the path and the curvature of the Bernoulli log-likelihood in its
success probability. In the proof of \cref{lem:particle_one_block_ascent},
these were bounded separately, in \eqref{bd_Hessian_gp} and
\eqref{eq:app_prediction_displacement}; bounding their difference in one
calculation keeps the factor $q_j$ inside the sum and yields a constant
of order $1/\rho$ rather than $1/\rho^2$.
Taking the conditional expectation of \eqref{bd_second_pointwise} given
$X=x$, which replaces $y$ by $p^*(x)$, and then the expectation over
$X$ proves \eqref{bd_second_popu}, and hence (i).

We next prove the upper bound in
\eqref{eq:app_fixed_weight_lower_model}. Recall the scalar loss
$f_a(p)=-a\log p-(1-a)\log(1-p)$ from
\eqref{eq:app_population_scalar_loss_bounds}, so that
$-\cL(G)=\EE_X[f_{p^*(X)}\{p_G(X)\}]$ and
$f_a'(p)=(p-a)/\{p(1-p)\}$. Since $f_a$ is convex, the tangent
inequality $f_a(p')\ge f_a(p)+f_a'(p)(p'-p)$, applied with
$a=p^*(x)$, $p=p_{G_{\bbeta}}(x)$, and
$p'=p_{G_{\wt\bbeta}}(x)$, and integrated over $X$, gives
\begin{equation}\label{eq:block_tangent}
    \cL(G_{\wt\bbeta})-\cL(G_{\bbeta})
    \le\EE_X\left[r_{G_{\bbeta}}(X)
       \{p_{G_{\wt\bbeta}}(X)-p_{G_{\bbeta}}(X)\}\right],
\end{equation}
with the residual $r_G$ defined in
\eqref{eq:app_population_residual_potential}.
Using $\delta_j=\wt\beta_j-\beta_j$, Taylor's formula with integral
remainder in the locations gives, for every $x$,
$$
    p_{G_{\wt\bbeta}}(x)-p_{G_{\bbeta}}(x)
    =\sum_{j=1}^Mw_jD_\beta p_{\delta_{\beta_j}}[\delta_j](x)+R(x),
$$
where
$$
    R(x)=\sum_{j=1}^Mw_j\int_0^1(1-t)
       \sigma''\{x^\T(\beta_j+t\delta_j)\}(x^\T\delta_j)^2\,\rd t.
$$
By the definition of $\nabla g_G$ in
\eqref{eq:app_population_residual_potential}, the first-order term
contributes
$$
    \EE_X\Big[r_{G_{\bbeta}}(X)
       \sum_{j=1}^Mw_jD_\beta p_{\delta_{\beta_j}}[\delta_j](X)\Big]
    =\sum_{j=1}^Mw_j
       \langle\nabla g_{G_{\bbeta}}(\beta_j),\delta_j\rangle
$$
to \eqref{eq:block_tangent}. For the remainder,
\eqref{eq:app_population_location_derivatives} and
$\int_0^1(1-t)\,\rd t=1/2$ give
$$
\|R\|_\infty\le\{R_X^2/(12\sqrt3)\}\|\bdelta\|_w^2,
$$
whereas the Cauchy--Schwarz inequality, $p_G(1-p_G)\ge\rho^2$, and
\eqref{eq:app_population_kl_l2} give
$$
    \EE_X|r_{G_{\bbeta}}(X)|
    \le\rho^{-2}\{P_X(p_{G_{\bbeta}}-p^*)^2\}^{1/2}
    \le\frac{1}{\sqrt2\,\rho^2}
       \{\cL(G^*)-\cL(G_{\bbeta})\}^{1/2}.
$$
Plugging these bounds into \eqref{eq:block_tangent} proves
\eqref{eq:app_fixed_weight_lower_model}.

To prove (ii), apply \eqref{eq:app_fixed_weight_upper_model} in part (i) with $\wt\bbeta=\bbeta^+$, so that
$\delta_j=\beta_j^+-\beta_j$. To bound from below the first-order term in
\eqref{eq:app_fixed_weight_upper_model}, apply
\eqref{eq:app_projection_inequality} with
$z=\beta_j+\eta\nabla g_{G_{\bbeta}}(\beta_j)$, $y=\beta_j^+$, and
$v=\beta_j$, which gives
$\langle\eta\nabla g_{G_{\bbeta}}(\beta_j)-\delta_j,-\delta_j\rangle\le0$,
so that
\begin{equation}\label{eq:app_projection_gain_popu}
    \langle\nabla g_{G_{\bbeta}}(\beta_j),\delta_j\rangle
    \ge{1\over\eta}\|\delta_j\|_2^2,\qquad j\in[M].
\end{equation}
Multiplying \eqref{eq:app_projection_gain_popu} by $w_j$, summing over
$j$, and plugging the result into \eqref{eq:app_fixed_weight_upper_model}
yields
$$
    \cL(G_{\bbeta^+})-\cL(G_{\bbeta})
    \ge\left({1\over\eta}-{L_X\over2}\right)\|\bbeta^+-\bbeta\|_w^2
    \ge{1\over2\eta}\|\bbeta^+-\bbeta\|_w^2,
$$
where the last inequality uses $\eta\le L_X^{-1}$. This proves (ii).

Finally, \eqref{bd_second_pointwise} holds for every $(x,y)$ with
$\|x\|_2\le R_X$. Averaging it over the observations $(X_i,Y_i)$,
$i\in[n]$, instead of taking expectations, gives \eqref{bd_second_popu}
with $F(s)=\ell_n(G_{\bbeta+s\bdelta})$, whose derivative at zero is
$\sum_jw_j\langle\nabla g_{n,G_{\bbeta}}(\beta_j),\delta_j\rangle$.
The remaining steps are unchanged.
\end{proof}

\subsubsection{Feasible weight and location directions, and construction of the reference configuration}

We use the  denoising result of \cref{thm:app_population_inverse} in \cref{app:particle_denoising} to construct a joint direction for the
weights and locations in \cref{lem:app_population_direction}.
For an exactly reweighted iterate with sufficiently small likelihood
gap, \cref{lem:app_population_reference} shows that moving only the
locations along this direction reduces the gap by a fixed fraction,
without changing the weights.
We use this result to bound the likelihood increase over each block in
\cref{app:prop:app_population_block_gain} by linking the gradient-ascent step to the aforementioned direction.

Fix a finite mixing distribution
$G=\sum_{j=1}^Mw_j\delta_{\beta_j}$ on $\cB_r$. We start by defining the feasible directions (relative to this $G$).
\begin{definition}[Feasible directions]\label{def_feasible}
For a weight direction
$\dot w=(\dot w_1,\ldots,\dot w_M)^\T\in \RR^M$ and location directions
$\xi=(\xi_1,\ldots,\xi_M)$, we call $(\dot w,\xi)$ {\em feasible} if, for all sufficiently small
$s\ge0$, the weights $w_j + s\dot w_j$, $j\in [M]$, are nonnegative and sum to one, and all
locations $\beta_j + s\xi_j \in \cB_r$.
\end{definition}
\noindent For any feasible $(\dot w,\xi)$, consider the perturbation
\begin{equation}\label{def_G_s}
    G_s=\sum_{j=1}^M(w_j+s\dot w_j)
        \delta_{\beta_j+s\xi_j} \in \cP(\cB_r).
\end{equation}
We use $D$ for the derivative of the prediction function $p_{G_s}$ in $s$ evaluated at
$s=0$. The limit may be taken from the right if the path is feasible only for $s\ge0$.
With $x$ held fixed, the prediction derivative is
\begin{equation}
    Dp_G[\dot w,\xi](x)
    :=\left.{d\over ds}p_{G_s}(x)\right|_{s=0}
    =\sum_{j=1}^M\dot w_j\sigma(x^\T\beta_j)
     +\sum_{j=1}^Mw_j\sigma'(x^\T\beta_j)x^\T\xi_j.
    \label{eq:app_population_prediction_derivative}
\end{equation}
For sufficiently small $s$, the finite sum defining
$p_{G_s}$ and \eqref{eq:app_population_location_derivatives} bound
its first two derivatives in $s$ uniformly over $x\in\cB_{R_X}$.
Since $p_{G_s}(x)\in[\rho,1-\rho]$,
\eqref{eq:app_population_scalar_loss_bounds} allows
differentiation under the expectation through order two.
The population likelihood derivative along the same path is
\begin{equation}  \label{eq:app_population_lh_derivative}
    D\cL(G)[\dot w,\xi]
    :=\left.{d\over ds}\cL(G_s)\right|_{s=0} = \EE \left\{
        {p^*(X)-p_G(X)\over p_G(X)\{1-p_G(X)\}}  Dp_G[\dot w,\xi](X)
    \right\}.
\end{equation}
For the location derivatives, we use
\eqref{eq:app_population_location_derivatives}. In particular, for
$u,v\in\RR^d$,
$$
    D_\beta p_{\delta_\beta}[u](x)=\sigma'(x^\T\beta)x^\T u,
    \qquad
    D_\beta^2 p_{\delta_\beta}[u,v](x)
       =\sigma''(x^\T\beta)(x^\T u)(x^\T v).
$$
Recall that $\rho=\sigma(-rR_X)$.

\begin{lemma}
\label{lem:app_population_direction}
Let $G=\sum_{j=1}^Mw_j\delta_{\beta_j}$ be any finite
probability measure on $\cB_r$, with positive weights and arbitrary
$M$. There are constants $\epsilon_0>0$ and $C_0<\infty$, depending only on
$d,r,\tau_0$ and the constants in \cref{assum:particle_population}, such that the
following holds.  If
\begin{equation}\label{init_p_G}
    0< P_X\left(p_G-p^*\right)^2 \le\epsilon_0^2,
\end{equation}
then there is a feasible direction $(\dot w,\xi)$ satisfying
\begin{align}
    \sum_{j=1}^M\dot w_j&=0,
    \label{eq:app_direction_mass}\\
    \sum_{j=1}^Mw_j\|\xi_j\|_2^2
    &\le\mfa(G,G^*).
    \label{eq:app_direction_norm}
\end{align}
Moreover, the induced first-order change in the prediction function
approximates  $p^*-p_G$:
\begin{equation}
    P_X\left(Dp_G[\dot w,\xi]-(p^*-p_G)\right)^2
    ~ \le~ \rho^4  P_X\left(p^*-p_G\right)^2 .
    \label{eq:app_direction_prediction}
\end{equation}
For every $j\in[M]$ and $0\le s\le1$, the location
$\beta_j+s\xi_j$ remains in $\cB_r$. If $\xi_j\ne0$, it also
lies at least $\tau_0/2$ from the boundary throughout this interval.

For the prediction function along this path,
\begin{equation}\label{def_ps}
    p_s=\sum_{j=1}^M w_jp_{\delta_{\beta_j+s\xi_j}}
\end{equation}
with $p_s'=\partial p_s/\partial s$ and
$p_s''=\partial^2p_s/\partial s^2$,
we also have
\begin{equation}
    \sup_{0\le s\le1}
    \{\|p_s'\|_\infty+\|p_s''\|_\infty\}
    \le C_0 \mfa(G,G^*),
    \label{eq:app_direction_path_derivatives}
\end{equation}
where the sup-norm is taken over $\cB_{R_X}$.
\end{lemma}

\begin{proof}
We start by grouping the particles in $G$ according to their closeness to each of the true atoms in $G^*$. With $R_0$ given in
\eqref{eq:app_population_cutoff_radius}, let $h_* \in (0, R_0\wedge1)$ be fixed below, and define
\begin{equation}\label{def_Ia}
    I_a=\left\{j\in [M]:\|\beta_j-\beta_a^*\|_2\le h_*\right\},\quad a\in [K^*],
    \qquad
    I_0=[M]\setminus\bigcup_a I_a.
\end{equation}
It is easy to see that
\begin{enumerate}[noitemsep]
    \item[(a)] the groups $\{I_0,I_1,\ldots,I_{K^*}\}$ are disjoint because $2h_*<s_-$;
    \item[(b)] for any $a\in[K^*]$, every particle in
$I_a$ lies at least
$\tau_0-h_*\ge\tau_0/2$ from the boundary, using
$h_*<R_0\le\tau_0/16$.
\end{enumerate}
For each $a\in[K^*]$, write
$$
    \delta_j=\beta_j-\beta_a^*,\quad j\in I_a,\qquad
    \bar\delta_a=\sum_{j\in I_a}w_j\delta_j,\qquad \bar w_a=\sum_{j\in I_a}w_j.
$$
Recalling the function $\chi_a$ in \eqref{def_chi_a_dist}, we have
$$
    \chi_a(\beta_j)=  \chi(\|\beta_j-\beta_a^*\|_2) = 1,\qquad j\in I_a,
$$
and
$$
    \|\beta_j -\beta_a^*\|_2 > h_*,\qquad j\notin I_a.
$$
It then follows  that
$$
    \int\chi_a(\beta) \,\rd G(\beta) = \sum_{j\in I_a} w_j \chi_a(\beta_j) + \sum_{j \notin I_a} w_j  \chi_a(\beta_j) =    \bar w_a +\int_{\beta: \|\beta-\beta_a^*\|_2> h_*} \chi_a(\beta) \,\rd G(\beta),
$$
so that by the definition of $\mfa(G,G^*)$ in \eqref{eq:app_population_gauge}
\begin{align*}
    |\bar w_a-w_a^*|
    &\le  \left|\int\chi_a(\beta) \,\rd G(\beta) -  w_a^*\right| + \left| \int_{\beta: \|\beta-\beta_a^*\|_2> h_*} \chi_a(\beta) \,\rd G(\beta)\right|\\
    &\le \left|\int\chi_a(\beta) \,\rd G(\beta) -  w_a^*\right|
       +{1\over h_*^2}\int\chi_a(\beta)\|\beta-\beta_a^*\|_2^2\,\rd G(\beta) \\
       &\le(1+h_*^{-2})\mfa(G,G^*).
\end{align*}
Similarly, since
$$
    \int \chi_a(\beta)(\beta-\beta_a^*) \,\rd G(\beta) =  \sum_{j\in I_a} w_j (\beta_j - \beta_a^*)  + \int_{\beta: \|\beta-\beta_a^*\|_2> h_*} \chi_a(\beta)(\beta-\beta_a^*) \,\rd G(\beta),
$$
we have
\begin{align*}
     \|\bar\delta_a\|_2 &~ =~  \left\| \sum_{j\in I_a}w_j (\beta_j - \beta_a^*) \right\|_2\\  &~ \le~    \left\|\int\chi_a(\beta)(\beta-\beta_a^*)\,\rd G(\beta)\right\|_2
       + {1\over h_*}\int\chi_a(\beta)\|\beta-\beta_a^*\|_2^2\,\rd G(\beta)\\
    &~ \le~ (1+h_*^{-1})\,\mfa(G,G^*).
\end{align*}
It is also easy to see that
\begin{equation*}
    \sum_{j\in I_a}w_j\|\delta_j\|_2^2 ~ \le\int\chi_a(\beta)\|\beta-\beta_a^*\|_2^2\,\rd G(\beta)
     \le\mfa(G,G^*).
\end{equation*}
Summing over $a\in [K^*]$ concludes
$$
\begin{aligned}
    \sum_{a=1}^{K^*}|\bar w_a-w_a^*|&\le(1+h_*^{-2})\mfa(G,G^*),\\
    \sum_{a=1}^{K^*}\|\bar\delta_a\|_2&\le(1+h_*^{-1})\mfa(G,G^*),\\
      \sum_{a=1}^{K^*}\sum_{j\in I_a}w_j\|\delta_j\|_2^2 &\le \sum_{a=1}^{K^*}\int\chi_a(\beta)\|\beta-\beta_a^*\|_2^2\,\rd G(\beta)
       \le\mfa(G,G^*).
\end{aligned}
$$
By \cref{thm:app_population_inverse},
\begin{equation}
    \mfa(G,G^*)
    ~ \le ~  c_{\mathrm{inv}}^{-1}\{P_X(p_G-p^*)^2\}^{1/2},
    \label{eq:app_gauge_by_prediction}
\end{equation}
so that  choosing
$\epsilon_0\le w_0c_{\mathrm{inv}}/(2(1+h_*^{-2}))$ ensures
\begin{equation*}
\bar w_a\ge w_a^*-(1+h_*^{-2})\epsilon_0/c_{\mathrm{inv}}\ge w_0/2.
\end{equation*}

We first construct the weight direction $\dot w$. We decrease the
weights in $I_0$ and adjust the total weight $\bar w_a$ of each
$I_a$ towards $w_a^*$, keeping the relative weights within each
group unchanged:
\begin{equation}\label{def_dot_w}
    \dot w_j=
    \begin{cases}
       -w_j,&~ j\in I_0,\\
       (w_a^*/\bar w_a-1)w_j,&~j\in I_a.
    \end{cases}
\end{equation}
The fact that $\sum_{a=1}^{K^*} w_a^*=1$ ensures
$$
    \sum_{j=1}^M \dot w_j
       =-\sum_{j\in I_0}w_j+\sum_{a=1}^{K^*}(w_a^*-\bar w_a)=0.
$$
For $0\le s\le1$, the new weights are
$$
    w_j+s\dot w_j =
    \begin{cases}
       (1-s)w_j,&~ j\in I_0,\\
     (1-s+sw_a^*/\bar w_a)w_j,&~j\in I_a, ~ a\in [K^*]
    \end{cases}
$$
and hence are nonnegative. This proves that $\dot w$ is a feasible weight direction for any $0\le s\le 1$.

For the location direction $\xi = (\xi_1,\ldots,\xi_M)$, take
$$
    \xi_j= \begin{cases}
       0,&~ j\in I_0,\\
       -{1\over 2}(\delta_j + \bar\delta_a/\bar w_a),&~ j\in I_a, ~ a\in [K^*].
    \end{cases}
$$
By the definitions of $\bar w_a$ and $\bar\delta_a$,
\begin{equation}
    \sum_{j\in I_a}w_j\xi_j
       =-\frac12\sum_{j\in I_a}w_j\delta_j
        -\frac{\bar\delta_a}{2\bar w_a}\sum_{j\in I_a}w_j
       =-\bar\delta_a.
    \label{eq:app_direction_first_moment}
\end{equation}
Since
$\|\bar\delta_a/\bar w_a\|_2\le h_*$, for $0\le s\le1$ we also have
$$
    \|\beta_j+s\xi_j-\beta_a^*\|_2
    =\left\|(1-s/2)\delta_j
                -{s \bar\delta_a\over2\bar w_a}\right\|_2
    \le h_*.
$$
Thus the location direction $\xi$  is  also feasible on the whole interval $0\le s\le 1$. Furthermore, since
the weighted Cauchy--Schwarz gives
$\|\bar\delta_a\|_2^2\le \bar w_a\sum_{j\in I_a}w_j\|\delta_j\|_2^2$, one has
\begin{equation}
    \begin{aligned}
    \sum_{j\in I_a}w_j\|\xi_j\|_2^2
    &={1\over4}\sum_{j\in I_a}w_j\|\delta_j\|_2^2
       +{3\over4\bar w_a}\|\bar\delta_a\|_2^2\\
    &\le\sum_{j\in I_a}w_j\|\delta_j\|_2^2
     \le \int\chi_a(\beta)\|\beta-\beta_a^*\|_2^2\,\rd G(\beta).
    \end{aligned}
    \label{eq:app_direction_cluster_norm}
\end{equation}
Summing over the groups thus yields
$$
    \sum_{j=1}^M w_j\|\xi_j\|_2^2
    \le\sum_{a=1}^{K^*} \int\chi_a(\beta)\|\beta-\beta_a^*\|_2^2\,\rd G(\beta)\le\mfa(G,G^*),
$$
which proves \eqref{eq:app_direction_norm}.

To prove \eqref{eq:app_direction_prediction}, first note that
$\dot w_j=-w_j$ and $\xi_j=0$ for $j\in I_0$. Thus
$$
    \sum_{j\in I_0}
       \{(\dot w_j+w_j)p_{\delta_{\beta_j}}
          +w_jD_\beta p_{\delta_{\beta_j}}[\xi_j]\}=0.
$$
Using this equality and substituting the weights on each $I_a$ into
\eqref{eq:app_population_prediction_derivative} gives
\begin{align*}
    &Dp_G[\dot w,\xi]-(p^*-p_G)\\
    &\quad=\sum_{a=1}^{K^*}\left\{
       {w_a^*\over \bar w_a}\sum_{j\in I_a}w_jp_{\delta_{\beta_j}}
       -w_a^* p_{\delta_{\beta_a^*}}
       +\sum_{j\in I_a}w_jD_\beta p_{\delta_{\beta_j}}[\xi_j]\right\}\\
    &\quad=\sum_{a=1}^{K^*}\left[
       {w_a^*-\bar w_a\over \bar w_a}
          \sum_{j\in I_a}w_j(p_{\delta_{\beta_j}}-p_{\delta_{\beta_a^*}})
       +\sum_{j\in I_a}w_j
          \{p_{\delta_{\beta_j}}-p_{\delta_{\beta_a^*}}+D_\beta p_{\delta_{\beta_j}}[\xi_j]\}
       \right].
\end{align*}
We bound the two terms for each group separately. For the first term,
Taylor's theorem at $\beta_a^*$ and
$\sum_{j\in I_a}w_j\delta_j=\bar\delta_a$ give
$$
\begin{aligned}
    \left\{P_X\Bigl[\sum_{j\in I_a}w_j
        (p_{\delta_{\beta_j}}-p_{\delta_{\beta_a^*}})
       -D_\beta p_{\delta_{\beta_a^*}}[\bar\delta_a]\Bigr]^2\right\}^{1/2}
    &\le {R_X^2\over12\sqrt3}\sum_{j\in I_a}w_j\|\delta_j\|_2^2\\
     & \le {R_X^2\over12\sqrt3}
          \int\chi_a(\beta)\|\beta-\beta_a^*\|_2^2\,\rd G(\beta).
\end{aligned}
$$
Indeed, Taylor's integral remainder for a displacement $u$ is
$\int_0^1(1-t)D_\beta^2p_{\delta_{\beta+tu}}[u,u]\,\rd t$, whose
sup-norm is at most $ R_X^2\|u\|_2^2/(12\sqrt3)$ by
\eqref{eq:app_population_location_derivatives}. The triangle inequality
therefore yields
\begin{equation}
    \begin{aligned}
    &\left\{P_X\Bigl[\sum_{j\in I_a}w_j
       (p_{\delta_{\beta_j}}-p_{\delta_{\beta_a^*}})\Bigr]^2\right\}^{1/2}
    \\&\quad\le\{P_X(D_\beta p_{\delta_{\beta_a^*}}[\bar\delta_a])^2\}^{1/2}
       +{R_X^2\over12\sqrt3}\int\chi_a(\beta)\|\beta-\beta_a^*\|_2^2\,\rd G(\beta)\\
    &\quad\le \frac{R_X}{4}\|\bar\delta_a\|_2+{R_X^2\over12\sqrt3}\int\chi_a(\beta)\|\beta-\beta_a^*\|_2^2\,\rd G(\beta).
    \end{aligned}
    \label{eq:app_direction_weight_remainder}
\end{equation}

For the second term, expand $p_{\delta_{\beta_j}}$ to second order and
$D_\beta p_{\delta_{\beta_j}}[\xi_j]$ to first order about the true location.
By \eqref{eq:app_population_location_derivatives}, the two integral
remainders have factors $\int_0^1(1-t)^2/2\,\rd t=1/6$ and
$\int_0^1(1-t)\,\rd t=1/2$, respectively. Thus, for every
$j\in I_a$,
\begin{align*}
    &\left\{P_X\left[p_{\delta_{\beta_j}}-p_{\delta_{\beta_a^*}}
       -D_\beta p_{\delta_{\beta_a^*}}[\delta_j]
       -{1\over2}D_\beta^2 p_{\delta_{\beta_a^*}}
          [\delta_j,\delta_j]
       \right]^2\right\}^{1/2} \le~  {R_X^3\over48}\|\delta_j\|_2^3,\\
    &\left\{P_X\left[D_\beta p_{\delta_{\beta_j}}[\xi_j]-D_\beta p_{\delta_{\beta_a^*}}[\xi_j]
       -D_\beta^2 p_{\delta_{\beta_a^*}}[\delta_j,\xi_j]
       \right]^2\right\}^{1/2} ~~  \le~  {R_X^3\over16}\|\delta_j\|_2^2\|\xi_j\|_2.
\end{align*}
Recall that $h_*$ is the radius defining $I_a$, so
$\|\delta_j\|_2\le h_*$. The definition of $\xi_j$
and $\|\bar\delta_a/\bar w_a\|_2\le h_*$ also give $\|\xi_j\|_2\le h_*$.
Consequently, the sum of the above two remainders is bounded by
\begin{align}\nonumber
     \sum_{j\in I_a}w_j
       \left\{{R_X^3\over48} \|\delta_j\|_2^3
              + {R_X^3\over16}\|\delta_j\|_2^2\|\xi_j\|_2 \right\}
    &\le {R_X^3h_*\over12}\sum_{j\in I_a}w_j\|\delta_j\|_2^2\\
    &\le {R_X^3h_*\over12}
          \int\chi_a(\beta)\|\beta-\beta_a^*\|_2^2\,\rd G(\beta).\label{bd_remainder}
\end{align}
By \eqref{eq:app_direction_first_moment}, the terms containing the
first derivative sum to zero:
$$
    \sum_{j\in I_a}w_j
       D_\beta p_{\delta_{\beta_a^*}}[\delta_j+\xi_j]
    =D_\beta p_{\delta_{\beta_a^*}}\left[\bar\delta_a+\sum_{j\in I_a}w_j\xi_j\right]
    =0.
$$
For the quadratic terms, substituting the definition of $\xi_j$ yields
\begin{align*}
     {1\over2}\sum_{j\in I_a}w_jD_\beta^2 p_{\delta_{\beta_a^*}}
       [\delta_j,\delta_j]
       +\sum_{j\in I_a}w_jD_\beta^2 p_{\delta_{\beta_a^*}}
          [\delta_j,\xi_j]
    & =\sum_{j\in I_a}w_jD_\beta^2 p_{\delta_{\beta_a^*}}
       \left[\delta_j,
          {\delta_j\over2}+\xi_j
       \right]\\
    & =-{1\over2\bar w_a}D_\beta^2 p_{\delta_{\beta_a^*}}
       \left[\sum_{j\in I_a}w_jd_j,\bar\delta_a\right]\\
     &  =-{1\over2\bar w_a}D_\beta^2 p_{\delta_{\beta_a^*}}[\bar\delta_a,\bar\delta_a].
\end{align*}
Thus the terms involving
$\sum_{j\in I_a}w_jd_jd_j^\T$ cancel, and the remaining second
derivative is applied only to $\bar\delta_a$.

Combining these cancellations with the bound in \eqref{bd_remainder} and
$\bar w_a\ge w_0/2$ gives
$$
\begin{aligned}
    &\left\{P_X\left[\sum_{j\in I_a}w_j
       \{p_{\delta_{\beta_j}}-p_{\delta_{\beta_a^*}}
          +D_\beta p_{\delta_{\beta_j}}[\xi_j]\}\right]^2\right\}^{1/2}\\
    &\qquad\le {R_X^2\over6\sqrt3\,w_0}\|\bar\delta_a\|_2^2
       +{R_X^3h_*\over12}
          \int\chi_a(\beta)\|\beta-\beta_a^*\|_2^2\,\rd G(\beta).
\end{aligned}
$$
Combining the two groupwise bounds gives
\begin{align*}
    &\{P_X(Dp_G[\dot w,\xi]-(p^*-p_G))^2\}^{1/2}\\
    &\quad\le {2\over w_0}\sum_{a=1}^{K^*}|w_a^*-\bar w_a|
       \left\{{R_X\over4}\|\bar\delta_a\|_2
          +{R_X^2\over12\sqrt3}
             \int\chi_a(\beta)\|\beta-\beta_a^*\|_2^2\,\rd G(\beta)\right\}\\
    &\qquad+{R_X^2\over6\sqrt3\,w_0}\sum_{a=1}^{K^*}\|\bar\delta_a\|_2^2
       +{R_X^3h_*\over12}
          \sum_{a=1}^{K^*}\int\chi_a(\beta)\|\beta-\beta_a^*\|_2^2\,\rd G(\beta)\\
    &\quad\le {R_X^3h_*\over12}\mfa(G,G^*)
          +C_{h_*}\mfa(G,G^*)^2,
\end{align*}
where we used the summed moment bounds and may take
$$
    C_{h_*}={2(1+h_*^{-2})\over w_0}\left(\frac{R_X}{4}(1+h_*^{-1})+{R_X^2\over12\sqrt3}\right)
                +{ R_X^2 (1+h_*^{-1})^2\over 6\sqrt3\,w_0}.
$$
In particular, these constants do not depend on the support size.
Choose
\begin{equation}\label{def_epsilon_0}
    h_*={1\over2}\min\left\{R_0,1,
                   {6\rho^2c_{\mathrm{inv}}\over R_X^3}\right\},
    \qquad
    \epsilon_0 = \min\left\{c_{\mathrm{inv}},\,
       {w_0c_{\mathrm{inv}}\over2(1+h_*^{-2})},\,
       {\rho^2c_{\mathrm{inv}}^2\over2C_{h_*}}\right\}.
\end{equation}
With these choices, \eqref{eq:app_gauge_by_prediction} bounds the
approximation error by
$$
    \left({R_X^3h_*\over12c_{\mathrm{inv}}}
       +{C_{h_*}\epsilon_0\over c_{\mathrm{inv}}^2}\right)
       \{P_X(p^*-p_G)^2\}^{1/2}
    \le\rho^2\{P_X(p^*-p_G)^2\}^{1/2}.
$$
This proves \eqref{eq:app_direction_prediction}, with both choices
depending only on $d,r,\tau_0$ and the constants in
\cref{assum:particle_population}.

To prove \eqref{eq:app_direction_path_derivatives}, first use
\eqref{eq:app_direction_first_moment} to write
$$
\begin{aligned}
    p_0'
    &=-\sum_{a=1}^{K^*}D_\beta p_{\delta_{\beta_a^*}}[\bar\delta_a] +\sum_{a=1}^{K^*}\sum_{j\in I_a}w_j
       \int_0^1D_\beta^2p_{\delta_{\beta_a^*+td_j}}
                    [\delta_j,\xi_j]\,\rd t.
\end{aligned}
$$
By \eqref{eq:app_population_location_derivatives}, the integral terms
for group $I_a$ have supremum norm at most
$$
\begin{aligned}
    \frac{R_X^2}{6\sqrt3}\sum_{j\in I_a}w_j\|\delta_j\|_2\|\xi_j\|_2
    &\le \frac{R_X^2}{6\sqrt3}
             \left(\sum_{j\in I_a}w_j\|\delta_j\|_2^2\right)^{1/2}
             \left(\sum_{j\in I_a}w_j\|\xi_j\|_2^2\right)^{1/2}\\
    &\le \frac{R_X^2}{6\sqrt3}
             \int\chi_a(\beta)\|\beta-\beta_a^*\|_2^2\,\rd G(\beta).
\end{aligned}
$$
The last step uses \eqref{eq:app_direction_cluster_norm}.
Combining this bound with the summed moment bounds gives
$$
\begin{aligned}
    \|p_0'\|_\infty
    &\le \frac{R_X}{4}\sum_{a=1}^{K^*}\|\bar\delta_a\|_2
       +\frac{R_X^2}{6\sqrt3}\sum_{a=1}^{K^*}
             \int\chi_a(\beta)\|\beta-\beta_a^*\|_2^2\,\rd G(\beta)\\
    &\le\left\{\frac{R_X}{4}(1+h_*^{-1})
                  +\frac{R_X^2}{6\sqrt3}\right\}\mfa(G,G^*).
\end{aligned}
$$
Differentiating the fixed-weight path twice also gives
$$
    p_s''(x)=\sum_{j=1}^M w_j\sigma''(x^\T(\beta_j+s\xi_j))(x^\T\xi_j)^2,
    \qquad
    \sup_{0\le s\le1}\|p_s''\|_\infty\le \frac{R_X^2}{6\sqrt3}\mfa(G,G^*).
$$
Since $p_s'=p_0'+\int_0^s p_t''\,\rd t$, we obtain
$$
    \sup_{0\le s\le1}\{\|p_s'\|_\infty+\|p_s''\|_\infty\}
    \le\left\{\frac{R_X}{4}(1+h_*^{-1})+\frac{R_X^2}{2\sqrt3}\right\}\mfa(G,G^*),
$$
proving \eqref{eq:app_direction_path_derivatives} with
$C_0=R_X(1+h_*^{-1})/4+R_X^2/(2\sqrt3)$.
Our choice $\epsilon_0\le c_{\mathrm{inv}}$ also gives
$\mfa(G,G^*)\le\epsilon_0/c_{\mathrm{inv}}\le1$.
The constants $h_*,\epsilon_0,C_0$ depend only on the fixed model
parameters, and not on $M$ or the particle weights.
\end{proof}

The next lemma constructs a feasible reference with unchanged weights
and a smaller gap. This comparison point is not an algorithmic update.

\begin{lemma}
\label{lem:app_population_reference}
Let $\epsilon_0$ be given in \eqref{def_epsilon_0} of
\cref{lem:app_population_direction}, which depends only on $d,r,\tau_0$
and the constants in \cref{assum:particle_population}. For every $K>0$, there are constants $\alpha\in(0,1)$
and $A>0$, depending only on the same parameters and $K$, such that
$KA\le\alpha/4$ and, whenever
$w\in\cM(\bbeta)$ and
$0<\cL(G^*)-\cL(G_{\bbeta})\le 2\epsilon_0^2$, there exists
$\bbeta^\dagger\in\cB_r^M$ satisfying
\begin{align}
    \cL(G^*)-\cL(G_{\bbeta^\dagger})
       &~ \le~ (1-\alpha)\{\cL(G^*)-\cL(G_{\bbeta})\},
       \label{eq:app_reference_gap}\\
    \|\bbeta^\dagger-\bbeta\|_w^2
       &~ \le~  A\sqrt{\cL(G^*)-\cL(G_{\bbeta})}.
       \label{eq:app_reference_distance}
\end{align}
Both $\alpha$ and $A$ are independent of $M$, $w$, and
$\bbeta$.
\end{lemma}

\begin{proof}
Since by \eqref{eq:app_population_kl_l2}, $0<\cL(G^*)-\cL(G_{\bbeta})\le 2\epsilon_0^2$ implies
$$
    0 < P_X(p_G-p^*)^2  \le \epsilon_0^2,
$$
the direction $(\dot w,\xi)$ in \cref{lem:app_population_direction} is applicable.
Below we use such direction to construct new particle locations $\bbeta^\dagger$, with unchanged weights, that meet both \eqref{eq:app_reference_gap} and \eqref{eq:app_reference_distance}.

We first omit the particles with $w_j=0$, which does not change
$G_{\bbeta}$. All remaining weights are positive, so we apply
\cref{lem:app_population_direction} to obtain $(\dot w_j,\xi_j)$
for the indices with $w_j>0$. For every omitted particle, we define
$\dot w_j=0$ and $\xi_j=0$.

By \eqref{eq:app_population_lh_derivative},
$$
    D\cL(G_{\bbeta})[\dot w,\xi]
       =\sum_{j=1}^M\dot w_jg_{G_{\bbeta}}(\beta_j)
          +D\cL(G_{\bbeta})[0,\xi],
$$
where the two terms on the right are the contributions from weight and
location changes, respectively.

Since the weights $w\in\cM(\bbeta)$
maximize the likelihood for fixed $\bbeta$,
\eqref{eq:app_population_weight_kkt} gives the same value of
$g_{G_{\bbeta}}(\beta_j)$ for all $j$ with $w_j>0$.
Using $\dot w_j=0$ whenever $w_j=0$ and $\sum_j\dot w_j=0$, we obtain
$$
    \sum_{j=1}^M\dot w_jg_{G_{\bbeta}}(\beta_j)
       =\left\{\sum_{k=1}^M w_kg_{G_{\bbeta}}(\beta_k)\right\}
          \sum_{j=1}^M\dot w_j=0.
$$
Thus moving only the locations along $\xi$, with all weights fixed,
gives the same likelihood derivative at $s=0$ as the joint direction
$(\dot w,\xi)$. We use the feasible path
$$
    G_{\bbeta+s\xi}=\sum_{j=1}^M w_j\delta_{\beta_j+s\xi_j},
    \qquad 0\le s\le1.
$$

To find a suitable step size on the  direction in $\xi$, by Taylor's integral formula, for any $0<\lambda\le1$, the increment of likelihood follows
$$
    \cL(G_{\bbeta+\lambda\xi})-\cL(G_{\bbeta})
    ~ = ~ \lambda\left.\frac{d}{ds}
                \cL(G_{\bbeta+s\xi})\right|_{s=0} +\int_0^\lambda(\lambda-s)\frac{d^2}{ds^2}
                \cL(G_{\bbeta+s\xi})\,\rd s.
$$
We bound the first term from below and control the absolute value of
the integral remainder.

First, using \eqref{eq:app_population_lh_derivative} and the definition of $r_G$ in  \eqref{eq:app_population_residual_potential} gives
\begin{align}\nonumber
    \left.{d\over ds}\cL(G_{\bbeta+s\xi})\right|_{s=0} &=  P_X\left\{r_{G_{\bbeta}}Dp_{G_{\bbeta}}[\dot w,\xi]\right\}\\
    &= P_X\frac{(p^*-p_{G_{\bbeta}})^2}
                    {p_{G_{\bbeta}}(1-p_{G_{\bbeta}})}
       +P_X\left\{r_{G_{\bbeta}}
          \left(Dp_{G_{\bbeta}}[\dot w,\xi]-(p^*-p_{G_{\bbeta}})\right)\right\}.\label{eq:app_population_directional_gain}
\end{align}
The fact $p_{G_{\bbeta}}(X)(1-p_{G_{\bbeta}}(X)) \le1/4$ ensures
$$
    P_X\frac{(p^*-p_{G_{\bbeta}})^2}
                    {p_{G_{\bbeta}}(1-p_{G_{\bbeta}})} \ge 4P_X(p_{G_{\bbeta}}-p^*)^2.
$$
Since $p_{G_{\bbeta}}(X)\{1-p_{G_{\bbeta}}(X)\}\ge\rho^2$,
we have $\lvert r_{G_{\bbeta}}(X)\rvert
\le\rho^{-2}\lvert p^*(X)-p_{G_{\bbeta}}(X)\rvert$.
Cauchy--Schwarz then gives
\begin{align*}
    & \left|P_X\left\{r_{G_{\bbeta}}
          \left(Dp_{G_{\bbeta}}[\dot w,\xi]-(p^*-p_{G_{\bbeta}})\right)\right\}\right|\\
        &\le \rho^{-2}\{P_X(p_{G_{\bbeta}}-p^*)^2\}^{1/2}
          \{P_X(Dp_{G_{\bbeta}}[\dot w,\xi]-(p^*-p_{G_{\bbeta}}))^2\}^{1/2} \\
        &\le P_X(p_{G_{\bbeta}}-p^*)^2 &&\text{by \eqref{eq:app_direction_prediction}}.
\end{align*}
Combining these bounds with \eqref{eq:app_population_directional_gain} yields
\begin{equation}
    \left.{d\over ds}\cL(G_{\bbeta+s\xi})\right|_{s=0}
       \ge 3 P_X(p_{G_{\bbeta}}-p^*)^2 \overset{\eqref{eq:app_population_kl_l2}}{\ge} 3\rho^2 \{\cL(G^*)-\cL(G_{\bbeta})\}.
    \label{eq:app_reference_initial_derivative}
\end{equation}

Second, we bound the absolute value of the second derivative of
$\cL(G_{\bbeta+s\xi})$ in $s$. With $p_s$ defined in
\eqref{def_ps}, differentiation gives
\begin{align*}
    - {d^2\over ds^2} \cL(G_{\bbeta+s\xi})
    =\EE_X\left[
       \left\{{p^*\over p_s^2}
              +{1-p^*\over(1-p_s)^2}\right\}(p_s')^2
       +{p_s-p^*\over p_s(1-p_s)}p_s''
       \right].
\end{align*}
Note that
\begin{align*}
    \sup_{0\le s\le1}
       \{\|p_s'\|_\infty+\|p_s''\|_\infty\} &\le C_0\, \mfa(G_{\bbeta},G^*) &&\text{by \eqref{eq:app_direction_path_derivatives}}\\
       &\le {C_0\over c_{\mathrm{inv}}} \{P_X(p_{G_\bbeta}-p_{G^*})^2\}^{1/2} &&\text{by \eqref{eq:app_gauge_by_prediction}}\\
       &\le {C_0\over\sqrt2\,c_{\mathrm{inv}}}\sqrt{\cL(G^*)-\cL(G_{\bbeta})}  &&\text{by \eqref{eq:app_population_kl_l2}}.
\end{align*}
Also, $p_s-p^*=p_0-p^*+\int_0^s p_t'\,\rd t$ and \eqref{eq:app_population_kl_l2} imply
$$
\begin{aligned}
    \{P_X(p_s-p^*)^2\}^{1/2}
    &\le\{P_X(p_0-p^*)^2\}^{1/2}
          +\int_0^s\|p_t'\|_\infty\,\rd t\\
          &\le\left( {1\over \sqrt2}+{C_0\over\sqrt2\,c_{\mathrm{inv}}}\right)
       \{\cL(G^*)-\cL(G_{\bbeta})\}^{1/2}.
\end{aligned}
$$
By  \eqref{eq:app_population_scalar_loss_bounds} and Cauchy--Schwarz,
\begin{align*}
    \sup_{0\le s\le1}
       \left|{d^2\over ds^2} \cL(G_{\bbeta+s\xi}) \right|
    & ~ \le {1\over \rho^2}\sup_{0\le s\le 1}\left\{
       P_X(p_s')^2
       +\{P_X(p_s-p^*)^2\}^{1/2}\{P_X(p_s'')^2\}^{1/2}\right\}\\
    &~ \le C_F \{\cL(G^*)-\cL(G_{\bbeta})\}
\end{align*}
where
$$
C_F=  {2C_0^2+C_0c_{\mathrm{inv}}\over 2c_{\mathrm{inv}}^2\rho^2}  .
$$

Combining this  second-derivative bound with
\eqref{eq:app_reference_initial_derivative} yields, for
$0<\lambda\le\min\{1,3\rho^2/C_F\}$,
\begin{equation}\label{lh_bd_ref}
    \cL(G_{\bbeta+\lambda\xi})-\cL(G_{\bbeta})
     \ge\left(3\rho^2\lambda-{C_F\lambda^2\over2}\right)
       \{\cL(G^*)-\cL(G_{\bbeta})\} \ge {3\rho^2\lambda\over2}
       \{\cL(G^*)-\cL(G_{\bbeta})\}.
\end{equation}
Moreover,
\begin{equation}
    \lambda^2 \|\xi\|_w^2
       \overset{\eqref{eq:app_direction_norm}}{\le}\lambda^2\mfa(G_{\bbeta},G^*)
       \le {\lambda^2\over\sqrt2\,c_{\mathrm{inv}}}\sqrt{\cL(G^*)-\cL(G_{\bbeta})}.\label{norm_bd_ref}
\end{equation}
Choose
$$
    \lambda=\min\left\{1,\frac{6c_{\mathrm{inv}}^2\rho^4}{2C_0^2+C_0c_{\mathrm{inv}}} ,
                 \frac{3\sqrt2\,\rho^2c_{\mathrm{inv}}}{8K}\right\},
$$
and set
$$
\bbeta^\dagger=\bbeta+\lambda\xi,\qquad \alpha={3\rho^2\lambda\over 2},\qquad  A={\lambda^2 \over \sqrt2\,c_{\mathrm{inv}}}.
$$
Since the location path is feasible, the above $\bbeta^\dagger\in\cB_r^M.$
The bounds in \eqref{lh_bd_ref} and \eqref{norm_bd_ref} prove
\eqref{eq:app_reference_gap} and \eqref{eq:app_reference_distance}, respectively.
Moreover,
$$
    {KA\over\alpha}
    ={2K\lambda\over3\sqrt2\,\rho^2c_{\mathrm{inv}}}\le{1\over4}.
$$
Since $\rho\le1/2$ and $\lambda\le1$, we also have
$0<\alpha\le3/8<1$. \end{proof}

\subsubsection{Proof of \cref{prop:app_population_block_gain}}\label{app:prop:app_population_block_gain}

\begin{proof}
Fixing any $t\ge 0$, we simplify the notation by writing
$$
    \bbeta^{(h)} := \bbeta^{t,h},\qquad G_{\bbeta^{(h)}} := G^{t,h},\qquad 0\le h\le H.
$$
The weights $w := w^t \in \cM(\bbeta^{(0)})$ are therefore fixed as those of $G_{\bbeta^{(0)}}$. Write
$$
    \Delta_0:=\cL(G^*)-\cL(G_{\bbeta^{(0)}}),
    \qquad
    \bdelta_h:=\bbeta^{(h+1)}-\bbeta^{(h)},\quad 0\le h<H.
$$
The step size satisfies $0<\eta\le27\rho/(4R_X^2)=:1/L_X$, so
\eqref{eq:app_fixed_weight_descent} of \cref{lem:app_population_one_step_ascent}(ii) gives
\begin{equation}\label{eq:block_one_step}
    \cL(G_{\bbeta^{(h+1)}})-\cL(G_{\bbeta^{(h)}})
    \ge{1\over2\eta}\|\bdelta_h\|_w^2,\qquad 0\le h<H.
\end{equation}
In particular, $\cL(G_{\bbeta^{(h)}})\ge\cL(G_{\bbeta^{(0)}})$, which implies,
\begin{equation}\label{eq:block_monotone}
    \cL(G^*)-\cL(G_{\bbeta^{(h)}}) \le  \cL(G^*)-\cL(G_{\bbeta^{(0)}})\le\Delta_0,
    \qquad 0\le h\le H.
\end{equation}
Note that if $\Delta_0 = 0$, then the above ensures $\cL(G^*) - \cL(G_{\bbeta^{(h)}}) = 0$ for all $h\le H$ so that the claim in \eqref{eq:app_population_block_gain} holds.

It thus remains to focus on $0<\Delta_0 \le \Gamma_0 := 2\epsilon_0^2$ with $\epsilon_0$ specified in \eqref{def_epsilon_0} of
\cref{lem:app_population_direction}. In this case, applying \cref{lem:app_population_reference} with $\bbeta = \bbeta^{(0)}$ and
\begin{equation}\label{def_K}
    K:={R_X^2 \over  12\sqrt6\,\rho^2}
\end{equation}
gives
$\bbeta^\dagger\in\cB_r^M$ and constants $\alpha\in(0,1)$, $A>0$
with $KA\le\alpha/4$ such that
\begin{equation}\label{eq:block_reference}
    \cL(G_{\bbeta^\dagger})-\cL(G_{\bbeta^{(0)}})
    \ge\alpha\Delta_0,
    \qquad
    \|\bbeta^\dagger-\bbeta^{(0)}\|_w^2\le A\sqrt{\Delta_0}.
\end{equation}
Set
\begin{equation}\label{def_kappa}
\kappa=\min\left\{{\alpha\over 2},\, {\alpha^2 \over 32A}\right\}.
\end{equation}
If
$\cL(G_{\bbeta^{(H)}})-\cL(G_{\bbeta^{(0)}})
\ge(\alpha/2)\Delta_0$, then \eqref{eq:app_population_block_gain} holds trivially as  $\kappa\le\alpha/2$.

Now suppose
\begin{equation}\label{contra_Fw}
    \cL(G_{\bbeta^{(H)}})-\cL(G_{\bbeta^{(0)}})
    <{\alpha\over2}\Delta_0
\end{equation}
which is also bounded from below by $\Delta_0>0$.
We show that, in this case,
\begin{equation}\label{eq:block_large_step}
    \|\bdelta_h\|_w^2 ~ \ge ~ {\eta^2\alpha^2\over16A}\Delta_0^{3/2},
    \qquad 0\le h<H.
\end{equation}
If this is the case, then summing \eqref{eq:block_one_step} over $h$ and using
\eqref{eq:block_large_step} gives
$$
    \cL(G_{\bbeta^{(H)}})-\cL(G_{\bbeta^{(0)}})
    \ge{1\over2\eta}\sum_{h=0}^{H-1}\|\bdelta_h\|_w^2
    \ge{\alpha^2\over32A}\,\eta H\,\Delta_0^{3/2}
    \ge\kappa\Delta_0\min\{\eta H\sqrt{\Delta_0},1\},
$$
which is the claim in \eqref{eq:app_population_block_gain}.

To prove \eqref{eq:block_large_step},  we compare the iterates with the reference $\bbeta^\dagger$.
By \eqref{eq:block_one_step},
$\cL(G_{\bbeta^{(h+1)}})\le\cL(G_{\bbeta^{(H)}})$, so
 that for all $ 0\le h<H$,
\begin{align}\nonumber
    \cL(G_{\bbeta^\dagger})-\cL(G_{\bbeta^{(h+1)}})
    &\ge\cL(G_{\bbeta^\dagger})-\cL(G_{\bbeta^{(H)}})\\\nonumber
    &= \cL(G_{\bbeta^\dagger})-\cL(G_{\bbeta^{(0)}}) + \cL(G_{\bbeta^{(0)}}) -\cL(G_{\bbeta^{(H)}})\\
    &> {\alpha\over2}\Delta_0
    \label{eq:app_reference_remaining_gap}
\end{align}
where the last step uses \eqref{eq:block_reference} and \eqref{contra_Fw}.
Next, fix $0\le h<H$ and compare the $h$-th update $\bbeta^{(h+1)}$ with $\bbeta^\dagger$.
Both \eqref{eq:app_fixed_weight_upper_model}, applied with
$(\bbeta,\wt\bbeta)=(\bbeta^{(h)},\bbeta^{(h+1)})$, and
\eqref{eq:app_fixed_weight_lower_model}, applied with
$(\bbeta,\wt\bbeta)=(\bbeta^{(h)},\bbeta^\dagger)$, are
expansions at the same point $\bbeta^{(h)}$. Subtracting the former from
the latter, and using \eqref{eq:block_monotone} in the coefficient of the
quadratic term, gives
\begin{align}
    &\cL(G_{\bbeta^\dagger})-\cL(G_{\bbeta^{(h+1)}})\notag\\
    &\quad\le\sum_{j=1}^Mw_j
       \langle\nabla g_{G_{\bbeta^{(h)}}}(\beta_j^{(h)}),
               \beta_j^\dagger-\beta_j^{(h+1)}\rangle
       +{L_X\over2}\|\bdelta_h\|_w^2
       +K\sqrt{\Delta_0}\,\|\bbeta^\dagger-\bbeta^{(h)}\|_w^2\notag\\
    &\quad\le{1\over\eta}\langle\bdelta_h,\bbeta^\dagger-\bbeta^{(h+1)}\rangle_w
       +{L_X\over2}\|\bdelta_h\|_w^2
       +K\sqrt{\Delta_0}\,\|\bbeta^\dagger-\bbeta^{(h)}\|_w^2\notag\\
    &\quad={1\over\eta}\langle\bdelta_h,\bbeta^\dagger-\bbeta^{(h)}\rangle_w
       -\left({1\over\eta}-{L_X\over2}\right)\|\bdelta_h\|_w^2
       +K\sqrt{\Delta_0}\,\|\bbeta^\dagger-\bbeta^{(h)}\|_w^2\notag\\
    &\quad\le{1\over\eta}\langle\bdelta_h,\bbeta^\dagger-\bbeta^{(h)}\rangle_w
       -{1\over2\eta}\|\bdelta_h\|_w^2
       +K\sqrt{\Delta_0}\,\|\bbeta^\dagger-\bbeta^{(h)}\|_w^2.
    \label{eq:app_reference_step_comparison}
\end{align}
The second inequality applies \eqref{eq:app_projection_inequality} with
$z=\beta_j^{(h)}+\eta\nabla g_{G_{\bbeta^{(h)}}}(\beta_j^{(h)})$,
$y=\beta_j^{(h+1)}$, and $v=\beta_j^\dagger$, which gives
$\langle\eta\nabla g_{G_{\bbeta^{(h)}}}(\beta_j^{(h)})-\delta_{h,j},
\beta_j^\dagger-\beta_j^{(h+1)}\rangle\le0$ for each $j$; the equality
uses $\bbeta^\dagger-\bbeta^{(h+1)}=\bbeta^\dagger-\bbeta^{(h)}-\bdelta_h$;
and the last inequality uses $\eta\le L_X^{-1}$.

From \eqref{eq:app_reference_step_comparison}, we argue by induction that
\begin{equation}\label{eq:block_distance}
    \|\bbeta^{(h)}-\bbeta^\dagger\|_w^2\le A\sqrt{\Delta_0},
    \qquad 0\le h\le H.
\end{equation}
Indeed, for $h=0$, this is \eqref{eq:block_reference}. Suppose
\eqref{eq:block_distance} holds for some $h<H$. Since
$\bbeta^{(h+1)}-\bbeta^\dagger=(\bbeta^{(h)}-\bbeta^\dagger)+\bdelta_h$,
expanding the squared norm gives the identity
$$
    {2\over\eta}\langle\bdelta_h,\bbeta^\dagger-\bbeta^{(h)}\rangle_w
    ={1\over\eta}\left\{\|\bbeta^{(h)}-\bbeta^\dagger\|_w^2
       -\|\bbeta^{(h+1)}-\bbeta^\dagger\|_w^2
       +\|\bdelta_h\|_w^2\right\}.
$$
Substituting this identity into \eqref{eq:app_reference_step_comparison}
cancels the terms in $\|\bdelta_h\|_w^2$ and yields
$$
    \|\bbeta^{(h+1)}-\bbeta^\dagger\|_w^2
    \le\|\bbeta^{(h)}-\bbeta^\dagger\|_w^2
       +2\eta K\sqrt{\Delta_0}\,\|\bbeta^{(h)}-\bbeta^\dagger\|_w^2
       -2\eta\{\cL(G_{\bbeta^\dagger})
               -\cL(G_{\bbeta^{(h+1)}})\}.
$$
By the induction hypothesis and $KA\le\alpha/4$, the second term on the
right is at most $2\eta KA\Delta_0\le(\eta\alpha/2)\Delta_0$, whereas by
\eqref{eq:app_reference_remaining_gap} the last term is less than
$-\eta\alpha\Delta_0$. Hence
$$
    \|\bbeta^{(h+1)}-\bbeta^\dagger\|_w^2
    \le\|\bbeta^{(h)}-\bbeta^\dagger\|_w^2-{\eta\alpha\over2}\Delta_0
    \le A\sqrt{\Delta_0},
$$
which completes the induction and proves \eqref{eq:block_distance}.

Finally, to prove \eqref{eq:block_large_step}, invoking \eqref{eq:app_reference_step_comparison} again and using Cauchy-Schwarz inequality give
\begin{align*}
\cL(G_{\bbeta^\dagger})-\cL(G_{\bbeta^{(h+1)}}) &\le {1\over \eta} \|\bdelta_h\|_w \|\bbeta^{(h+1)}-\bbeta^\dagger\|_w + K\sqrt{\Delta_0} \|\bbeta^{(h+1)}-\bbeta^\dagger\|_w^2\\
  &  \le{1\over\eta}\|\bdelta_h\|_w\,A^{1/2}\Delta_0^{1/4}
       +KA\Delta_0\\
 &   \le{1\over\eta}\|\bdelta_h\|_w\,A^{1/2}\Delta_0^{1/4}
       +{\alpha\over4}\Delta_0.
\end{align*}
Combining this bound with \eqref{eq:app_reference_remaining_gap} yields
\eqref{eq:block_large_step}.  This completes the proof.
\end{proof}

\section{Proofs of the sample-level convergence of the finite particle approximation}\label{app:proof_thm:particle_uniform_empirical}

\subsection{Empirical comparisons and block improvement}
\label{app:particle_empirical_auxiliary}

We use the direction and fixed-weight comparison from the proof of
\cref{thm:particle_population}. We first bound the differences between
the empirical and population likelihoods and their derivatives,
uniformly over the candidate mixtures. When the population likelihood
gap exceeds a constant multiple of $\varepsilon_n(\delta)$, these
bounds preserve the improvement obtained in the population proof.

The empirical likelihood does not decrease along the iterates, but
the population likelihood need not increase. We therefore use the uniform
likelihood comparison to show that all iterates remain in the
neighborhood where the direction can be constructed. We then obtain
a recursion for $\ell_n(\wh G_r)-\ell_n(G_t)$ until it reaches
a constant multiple of $\varepsilon_n(\delta)$. Monotonicity keeps it
at most that level thereafter.

\subsection{Empirical process bounds on  likelihoods and their directional derivatives}

Let $P$ denote the joint law of $(X,Y)$, and let $P_n$ be its
empirical measure. For a function $f$ of $(X,Y)$, we write
$Pf=\EE f(X,Y)$ and $P_nf=n^{-1}\sum_{i=1}^n f(X_i,Y_i)$.
Recall $\rho=\sigma(-rR_X)$. We use $\varepsilon_n(\delta)$ as
defined in \eqref{eq:particle_statistical_floor}. Every $p\in\cF_r$ takes values
in $[\rho,1-\rho]$ on $\cB_{R_X}$. Supremum norms for functions of
$(x,y)$ are taken over $\cB_{R_X}\times\{0,1\}$.
The derivative of the Bernoulli log-likelihood with respect to its
success probability is
$$
    r_y(t)={y-t\over t(1-t)}={y\over t}-{1-y\over1-t},
    \qquad y\in\{0,1\},\quad t\in[\rho,1-\rho].
$$
Recall the likelihood $L(p;y,x)$ from \eqref{def_Lf} and the derivative of $\log L(p;y,x)$ with respect to $p(x)$ is $r_y(p(x))$. Moreover,
$$
    P_n\log L(p_G;Y,X)=\ell_n(G),
    \qquad
    P\log L(p_G;Y,X)=\cL(G).
$$

For the directional derivative of $\cL(G_s)$ in
\eqref{eq:app_population_lh_derivative} with $G_s$ given in \eqref{def_G_s}, its empirical
counterpart  is
\begin{equation}
    D\ell_n(G)[\dot w,\xi]
    :={\left.{d\over ds}\ell_n(G_s)\right|}_{s=0}
    =P_n\{r_Y(p_G)Dp_G[\dot w,\xi]\}.
    \label{eq:app_empirical_parameter_derivative}
\end{equation}
For the second-order location remainder, define
\begin{equation}\label{def_cJ}
    \mathcal J
    :=\left\{
      (x,y)\mapsto r_y\{p(x)\}\sigma''(x^\T\beta)(x^\T v)^2:
      p\in\cF_r,\ \beta\in\cB_r,\ \|v\|_2=1
    \right\}.
\end{equation}

The next lemma bounds the differences between the sample and population versions of the likelihood, the derivative along the constructed direction, and the second-order location terms. All three bounds hold on the same event.

\begin{lemma}
\label{lem:app_particle_empirical_comparison}
There is a constant $C>0$, depending only on $d,r,\tau_0$ and the constants in \cref{assum:particle_population}, such that, for every $\delta\in(0,1)$, the following holds on an event of probability at least $1-\delta/2$.
\begin{enumerate}
\item[(i)] For every $G\in\cP(\cB_r)$,
\begin{equation}
    \left|
      \ell_n(G)-\ell_n(G^*)+\{\cL(G^*)-\cL(G)\}
    \right|
    \le C\left[
       \{\{\cL(G^*)-\cL(G)\}\varepsilon_n(\delta)\}^{1/2}
       +\varepsilon_n(\delta)
    \right].
    \label{eq:app_empirical_likelihood_comparison}
\end{equation}
\item[(ii)] For every $G=\sum_{j=1}^Mw_j\delta_{\beta_j}$ satisfy \eqref{init_p_G} of \cref{lem:app_population_direction}, with the direction
$(\dot w,\xi)$ constructed in its proof,
\begin{equation}
    \left|
      D\ell_n(G)[\dot w,\xi]
      -D\cL(G)[\dot w,\xi]
    \right|
    \le C\left[
       \{\{\cL(G^*)-\cL(G)\}\varepsilon_n(\delta)\}^{1/2}
       +\varepsilon_n(\delta)
    \right].
    \label{eq:app_empirical_direction_comparison}
\end{equation}

\item[(iii)] With the set $\cJ$ given in \eqref{def_cJ}, one has
\begin{equation}
    \sup_{J\in\mathcal J}|(P_n-P)J|
    \le C\sqrt{\varepsilon_n(\delta)}.
    \label{eq:app_empirical_curvature_comparison}
\end{equation}
\end{enumerate}
\end{lemma}

\begin{proof}
We apply Bernstein's inequality to finite covers of the three
function classes separately. The covers approximate functions, rather than
particle configurations, and their sizes do not depend on $M$.
For the first two bounds, we also show that the second moment of each
function is at most a constant times the corresponding population
gap. This gives the factor
$\{\{\cL(G^*)-\cL(G)\}\varepsilon_n(\delta)\}^{1/2}$. We start by constructing finite covers and  deriving envelope functions and bounds for the corresponding variance terms.\\

\noindent{\bf The likelihood difference.}
For $t\in[\rho,1-\rho]$, note that
\begin{equation}\label{bd_ry}
    r_y'(t)=-\frac{y}{t^2}-\frac{1-y}{(1-t)^2},
    \qquad |r_y(t)|\le\rho^{-1},
    \qquad |r_y'(t)|\le\rho^{-2}.
\end{equation}
Consequently, for any $p,\wt p\in\cF_r$, any $x\in \cB_{R_X}$ and $y\in \{0,1\}$,
\begin{equation}\label{lip_ry}
    |r_y\{p(x)\}-r_y\{\wt p(x)\}|
       \le\rho^{-2}|p(x)-\wt p(x)|.
\end{equation}
Furthermore, integrating the derivative of
$y\log p+(1-y)\log(1-p)$ with respect to $p$ gives
\begin{equation}\label{lip_log_L}
    \left|\log L(p_G;y,x)-\log L(p_{\wt G};y,x)\right|
       \le\rho^{-1}|p_G(x)-p_{\wt G}(x)|
\end{equation}
as well as
\begin{equation}
    P\{\log L(p_G;Y,X)-\log L(p^*;Y,X)\}^2
    \le{1\over \rho^2}P_X(p_G-p^*)^2
    \overset{\eqref{eq:app_population_kl_l2}}{\le}{1\over2\rho^2}\{\cL(G^*)-\cL(G)\}.
    \label{eq:app_likelihood_contrast_variance}
\end{equation}
In view of \eqref{lip_log_L}, a cover of $\cF_r$ of radius
$\rho u$ therefore gives a cover of $\{\log L(p;y,x): p\in \cF_r\}$ of radius $u$.
Applying \cref{thm:metric_entropy}, with its constants depending only
on $d,r,R_X$, yields
\begin{align}\nonumber
    &\log\mathcal N
    \left(u,\{\log L(p_G;\cdot,\cdot)-\log L(p^*;\cdot,\cdot):G\in\cP(\cB_r)\},
      \|\cdot\|_{\infty}\right)\\
    &\qquad\le \log\mathcal N(\rho u,\cF_r,\|\cdot\|_{\infty,\cB_{R_X}})\notag\\
    &\qquad \le \{C\log(C/u)\}^{d+1},
    && 0<u\le1.
    \label{eq:app_likelihood_contrast_entropy}
\end{align}
Finally,  it is easy to see that
\begin{equation}\label{sup-rho}
    |\log L(p_G;y,x)-\log L(p^*;y,x)| \le
 \log\{(1-\rho)/\rho\},\qquad y\in \{0,1\}, x\in \cB_{R_X}.
\end{equation}

\noindent{\bf The directional derivative difference.}
For the directional derivative comparison, recall from \eqref{eq:app_empirical_parameter_derivative} and \eqref{eq:app_population_lh_derivative} that
$$
     D\ell_n(G)[\dot w,\xi] = P_n\left\{r_Y(p_G)Dp_G[\dot w,\xi]\right\}, \qquad D\cL(G)[\dot w,\xi] = P\left\{r_Y(p_G) Dp_G[\dot w,\xi]\right\}.
$$
We  first quantify the directional derivatives of prediction function $Dp_G[\dot w,\xi]$.
For notational simplicity, write
\begin{equation}
    h(x) :=  Dp_G[\dot w,\xi](x) =\int\sigma(x^\T\beta)\,\rd \nu(\beta)
       +\sum_{\ell=1}^d x_\ell
          \int\sigma'(x^\T\beta)\,\rd \nu_\ell(\beta).
    \label{eq:app_direction_signed_measure_form}
\end{equation}
where
$$
\nu=\sum_{j=1}^M \dot w_j\delta_{\beta_j},\qquad \nu_\ell=\sum_{j=1}^M  w_j \xi_{j,\ell}\delta_{\beta_j}$$
are finite signed
measures on $\cB_r$.

Recall from \eqref{def_Ia} that the groups $I_0, I_1,\ldots, I_{K^*}$ form a partition of $[M]$ and  $\bar w_a = \sum_{j\in I_a} w_j$ for $a\in [K^*]$. Since
$\sum_{a=1}^{K^*} w_a^*=1$ and
$\sum_{a=1}^{K^*} \bar w_a+\sum_{j\in I_0}w_j=1$,
$$
    \|\nu\|_{\mathrm{TV}}
    \le \sum_{j=1}^M |\dot w_j|
      =\sum_{j\in I_0}w_j+\sum_{a=1}^{K^*}|w_a^*-\bar w_a|
    \le\sum_{j\in I_0}w_j+\sum_{a=1}^{K^*}(w_a^*+\bar w_a)=2.
$$
The choice $\epsilon_0\le c_{\mathrm{inv}}$ in
\cref{lem:app_population_direction} and
\eqref{eq:app_gauge_by_prediction} give $\mfa(G,G^*)\le1$.
By Cauchy--Schwarz and \eqref{eq:app_direction_norm}, we also have
$$
    \begin{aligned}
    \sum_{\ell=1}^d\|\nu_\ell\|_{\mathrm{TV}}
    &\le\sum_{j=1}^M w_j\|\xi_j\|_1
     \le\sqrt d\sum_{j=1}^M  w_j\|\xi_j\|_2
    \\
    &\le\Bigl(d\sum_{j=1}^M  w_j\|\xi_j\|_2^2\Bigr)^{1/2}
     \le\{d\mfa(G,G^*)\}^{1/2}\le\sqrt d.
    \end{aligned}
$$
Define the function class
$$
    \cV := \left\{h \text{ in } \eqref{eq:app_direction_signed_measure_form}: \|\nu\|_{\mathrm{TV}}+\sum_{\ell=1}^d \|\nu_\ell\|_{\mathrm{TV}}\le 2+\sqrt d \right\}.
$$
Its envelope satisfies
\begin{equation}
    \|h\|_{\infty,\cB_{R_X}}
    \le\|\nu\|_{\mathrm{TV}}
       +{R_X\over4}\sum_{\ell=1}^d \|\nu_\ell\|_{\mathrm{TV}}
    \le (2+\sqrt d)(1+R_X/4) =: B_V.
    \label{eq:app_empirical_direction_envelope}
\end{equation}

To obtain a finite cover of $\cV$, we use the polynomial approximation argument in
\cref{prop:finite_dim_approx_class}, now for both $\sigma$ and
$\sigma'$. Put $R=rR_X$, and take the ellipse parameter
$\rho_R$ from \cref{lem:logistic_chebyshev}. That lemma verifies
holomorphy and the bound $|\sigma|\le1$ on
$R\overline E_{\rho_R}$, where $|\Im z|\le\pi/2$.
The derivative is also holomorphic there and, for $z=t+iv$,
$$
    |\sigma'(z)|
       ={e^{-t}\over1+e^{-2t}+2e^{-t}\cos v}
       \le{e^{-t}\over1+e^{-2t}}\le{1\over2}.
$$
By \cref{lem:chebyshev_projection_analytic}, there are polynomials
$\pi_{0,k},\pi_{1,k}$ of degree at most $k$ such that
$$
    \max_{\ell\in\{0,1\}}
    \sup_{|t|\le R}|\sigma^{(\ell)}(t)-\pi_{\ell,k}(t)|
       \le {2\over\rho_R-1}\rho_R^{-k}.
$$
For $0<u\le1$, choose
$$
    k=1\vee\left\lceil
       {\log\{8(2+\sqrt d)(1+R_X)/((\rho_R-1)u)\}\over\log\rho_R}
       \right\rceil.
$$
Replacing $\sigma,\sigma'$ in
\eqref{eq:app_direction_signed_measure_form} by these polynomials gives
a polynomial $h_k$ in $x$ of degree at most $k+1$, with
\begin{align*}
    \|h-h_k\|_{\infty,\cB_{R_X}}
    &\le {2\rho_R^{-k}\over\rho_R-1}
          \left(\|\nu\|_{\mathrm{TV}}
             +R_X\sum_{\ell=1}^d \|\nu_\ell\|_{\mathrm{TV}}\right) \le {2(2+\sqrt d)(1+R_X)\over\rho_R-1}\rho_R^{-k}
       \le {u\over 4},
\end{align*}
and, by \eqref{eq:app_empirical_direction_envelope},
\begin{equation}\label{sup-hk}
    \|h_k\|_{\infty,\cB_{R_X}}\le B_V+1/2.
\end{equation}
The polynomial space has dimension $\binom{d+k+1}{d}\le(k+2)^d$.
Apply the finite-dimensional covering argument in
\cref{prop:finite_dim_approx_class} with polynomial degree $k+1$,
norm bound $B_V+1/2$, and covering radius $u/4$.
Together with the approximation error $u/4$, this gives a
$u/2$-cover of $\cV$. Choosing a function $h\in\cV$ from each ball
that intersects $\cV$ doubles the covering radius, giving a
$u$-cover with centers in $\cV$. Hence
\begin{equation}
    \log\mathcal N(u,\cV,\|\cdot\|_{\infty,\cB_{R_X}})
    \le \binom{d+k+1}{d}
          \log\{1+16(B_V+1/2)/u\}
    \le\{C'\log(C'/u)\}^{d+1}.
    \label{eq:app_direction_class_entropy}
\end{equation}
Here $C'$ depends only on $d,r,R_X$. The same cover applies
regardless of the support sizes of the signed measures.

For the likelihood directional derivatives, consider
$r_y\{p(x)\}h(x)$, with $p\in\cF_r$ and
$h\in\cV$. By \eqref{bd_ry} and \eqref{lip_ry}, one has
\begin{equation}\label{sup-sh}
    \|r_Y(p)h\|_\infty\le{ B_V\over \rho},
    \qquad
    \|r_Y(p)h-r_Y(\wt p)\wt h\|_\infty
    \le {1\over\rho}\|h-\wt h\|_\infty
          +{B_V\over\rho^2}\|p-\wt p\|_\infty.
\end{equation}
Thus the product of a cover of $\cF_r$ of radius
$\rho^2u/(2B_V)$ and a cover of the directional function $h$ of radius $\rho u/2$
has radius $u$ for this class. The sum of their log-cardinalities
is at most $\{C\log(C/u)\}^{d+1}$, by
\cref{thm:metric_entropy} and
\eqref{eq:app_direction_class_entropy}. For the constructed
$h=Dp_G[\dot w,\xi]$, \eqref{eq:app_direction_prediction} also gives
$$
    \{P_Xh^2\}^{1/2}
       \le(1+\rho^2)\{P_X(p_G-p^*)^2\}^{1/2}
$$
such that
\begin{equation}
    P\{r_Y(p_G)h\}^2
    \le\rho^{-2}P_Xh^2
    \le{(1+\rho^2)^2\over2\rho^2}\{\cL(G^*)-\cL(G)\}.
    \label{eq:app_direction_score_variance}
\end{equation}
Moreover,
$D\ell_n(G)[\dot w,\xi]-D\cL(G)[\dot w,\xi] = (P_n-P)\{r_Y(p_G)h\}$.\\

\noindent{\bf The second-order location remainder difference.}
For the function in $\cJ$ in \eqref{def_cJ},  write
$J_{p,\beta,v}(x,y)=r_y(p(x))\sigma''(x^\T\beta)(x^\T v)^2$. Since
$$
    |(x^\T v)^2-(x^\T\wt v)^2|
    \le |x^\T(v-\wt v)|
         \{|x^\T v|+|x^\T\wt v|\}
    \le2R_X^2\|v-\wt v\|_2,
$$
by \eqref{eq:app_population_location_derivatives},
one has
\begin{align}\label{sup-J}
    \|J_{p,\beta,v}\|_\infty&\le \frac{R_X^2}{6\sqrt3\,\rho},\\\nonumber
    \|J_{p,\beta,v}-J_{\wt p,\wt\beta,\wt v}\|_\infty
    &\le {R_X^2\over6\sqrt3\,\rho^2}\|p-\wt p\|_\infty
          +{R_X^3\over8\rho}\|\beta-\wt\beta\|_2
          +{R_X^2\over3\sqrt3\,\rho}\|v-\wt v\|_2.
\end{align}
The second term uses
$|\sigma''(x^\T\beta)-\sigma''(x^\T\wt\beta)|
 \le(R_X/8)\|\beta-\wt\beta\|_2$.
Using covers of respective radii
$2\sqrt3\,\rho^2u/R_X^2$, $8\rho u/(3R_X^3)$, and
$\sqrt3\,\rho u/R_X^2$, we obtain
\begin{align*}
    &\log\mathcal N(u,\mathcal J,\|\cdot\|_\infty)\\
    &\le
       \log\mathcal N\left({2\sqrt3\,\rho^2u\over R_X^2},
                    \cF_r,\|\cdot\|_{\infty,\cB_{R_X}}\right) +d\log\left(1+{3rR_X^3\over4\rho u}\right)
          +d\log\left(1+{2R_X^2\over\sqrt3\,\rho u}\right)\\
    &\le\{C\log(C/u)\}^{d+1}.
\end{align*}
If the covering radius for $\cF_r$ exceeds one, we use radius one
instead, which only increases this upper bound. Also
\begin{equation}\label{bd_PJsq}
    PJ^2\le \frac{R_X^4}{108\rho^2}.
\end{equation}
\medskip

\noindent {\bf Application of Bernstein's inequalities.}
We now apply Bernstein's inequality to the three classes. In view of the sup-norm bounds in \eqref{sup-rho}, \eqref{sup-sh} and \eqref{sup-J}, set
$$
    B_*=\max\{1,\log\{(1-\rho)/\rho\},B_V/\rho,
                     R_X^2/(6\sqrt3\,\rho)\}.
$$
Choose $A,B$ in \eqref{eq:particle_statistical_floor} sufficiently
large that each class has an $n^{-2}$-net with log-cardinality at most
$A\{\log(Bn)\}^{d+1}$, and set
$$
    t_n=A\{\log(Bn)\}^{d+1}
             +\log(12/\delta).
$$
For each net member $f^\circ$, Bernstein's inequality, with
$\|f^\circ - Pf^\circ\|_\i \le 2B_*$ and variance at most $P(f^\circ)^2$, gives
$$
    \PP\left\{ |(P_n-P)f^\circ|
       \le\sqrt{2P(f^\circ)^2t_n\over n}
              +{4B_*t_n\over3n}\right\} \ge  1 - 2e^{-t_n}.
$$
Taking the union over the three nets gives failure probability at most
$$
    6\exp\left[A\{\log(Bn)\}^{d+1}-t_n\right]=\delta/2.
$$

When $\varepsilon_n(\delta)<1$, we have
$t_n/n\le\varepsilon_n(\delta)\le1$ and
$n^{-2}\le\varepsilon_n(\delta)$.
For an arbitrary function $f$ in one of the three classes, choose
$f^\circ$ with $\|f-f^\circ\|_\infty\le n^{-2}$.
Then
$$
    |(P_n-P)(f-f^\circ)|\le {2 \over n^2},
    \qquad
    P(f^\circ)^2\le2Pf^2+{2\over n^4}.
$$
Consequently, on the same event,
\begin{align*}
    |(P_n-P)f| &\le |(P_n-P)(f-f^\circ)| + |(P_n-P)f^\circ| \\
    &\le {2\over n^2}
          +2\sqrt{\Bigl(Pf^2+{1\over n^4}\Bigr)\varepsilon_n(\delta)}
          +{4B_*\over3}\varepsilon_n(\delta)\\
    &\le2\sqrt{Pf^2\,\varepsilon_n(\delta)}
          +\left(4+{4B_*\over3}\right)\varepsilon_n(\delta).
\end{align*}
This inequality holds for every function in each class on the same
event, with its own second moment $Pf^2$ and $B^*$. Substituting
\eqref{eq:app_likelihood_contrast_variance},
\eqref{eq:app_direction_score_variance}, and
\eqref{bd_PJsq} proves the three claims. For example, a common
comparison constant is
\begin{equation}\label{def_C_rel}
    C_{\mathrm{cmp}}
    :=\max\left\{
       2B_*,\,{\sqrt2(1+\rho^2)\over\rho},\,
       4+{4B_*\over3}+{R_X^2\over3\sqrt3\,\rho}
       \right\}.
\end{equation}

When $\varepsilon_n(\delta)=1$, the deterministic bound
$|(P_n-P)f|\le2B_*$ proves all three claims with the same constant.
The covers, envelopes, and variance bounds depend only on the stated
model parameters. Thus the constant and failure probability are
uniform over the model class and independent of $M$.
We use $C_{\mathrm{cmp}}$ for the constant in
\eqref{eq:app_empirical_likelihood_comparison}--\eqref{eq:app_empirical_curvature_comparison} below.
\end{proof}

\subsection{Reference configuration and  likelihood improvement}

We first compare the empirical and population likelihood gaps.
These bounds will keep the empirical iterates in the required
population neighborhood and convert the block improvement into a
recursion for the empirical gap. By $ab\le a^2/8+2b^2$,
\eqref{eq:app_empirical_likelihood_comparison} implies that for every $G\in\cP(\cB_r)$,
\begin{equation}
    \left|
      \ell_n(G)-\ell_n(G^*)+\{\cL(G^*)-\cL(G)\}
    \right|
    \le {1\over8}\{\cL(G^*)-\cL(G)\}+C_{\mathrm{rel}}\varepsilon_n(\delta)
    \label{eq:app_empirical_relative_likelihood}
\end{equation}
where
$C_{\mathrm{rel}}=2C_{\mathrm{cmp}}^2+C_{\mathrm{cmp}}$ with $C_{\mathrm{cmp}}$ given in \eqref{def_C_rel}.
Since $G^*\in\cP(\cB_r)$, taking the supremum over $G\in \cP(\cB_r)$ gives
$$
    0\le\ell_n(\wh G_r)-\ell_n(G^*)
       \le C_{\mathrm{rel}}\varepsilon_n(\delta).
$$
Combining this with \eqref{eq:app_empirical_relative_likelihood} yields
\begin{equation}
    {7\over8}\{\cL(G^*)-\cL(G)\}-C_{\mathrm{rel}}\varepsilon_n(\delta)
    \le
    \ell_n(\wh G_r)-\ell_n(G)
    \le
    {9\over8}\{\cL(G^*)-\cL(G)\}+2C_{\mathrm{rel}}\varepsilon_n(\delta).
    \label{eq:app_empirical_gap_comparison}
\end{equation}

We next construct new particle locations that increase the empirical
likelihood without changing the weights. This is the empirical version
of \cref{lem:app_population_reference}, using the direction from
\cref{lem:app_population_direction}. The additional requirement is that
the population likelihood gap be large enough relative to
$\varepsilon_n(\delta)$ to control the sampling errors.
For fixed weights $w$, we write
$G_{\bbeta}=\sum_jw_j\delta_{\beta_j}$ and use the weighted norm
$\|\cdot\|_w$ from the population proof.

\begin{lemma}
\label{lem:app_particle_empirical_reference}
Work on the event in \cref{lem:app_particle_empirical_comparison},
and let $\epsilon_0$ be as in \cref{lem:app_population_reference}.
There is a constant $C_E\ge1$ such that the following holds for every
$K>0$. There are constants
$\alpha\in(0,1)$ and $A>0$, with $KA\le\alpha/4$, such that,
whenever $w\in\cM_n(\bbeta)$ and
\begin{equation}\label{log_L_gap_emp}
    C_E\varepsilon_n(\delta)\le
       \cL(G^*)-\cL(G_{\bbeta})\le\Gamma_0,
\end{equation}
there exists $\bbeta^\dagger\in\cB_r^M$ satisfying
\begin{align}
    \ell_n(\wh G_r)-\ell_n(G_{\bbeta^\dagger})
       &\le \{\ell_n(\wh G_r)-\ell_n(G_{\bbeta})\}-\alpha \{\cL(G^*)-\cL(G_{\bbeta})\},
       \label{eq:app_empirical_reference_gain}\\
    \|\bbeta^\dagger-\bbeta\|_w^2
       &\le A\sqrt{\cL(G^*)-\cL(G_{\bbeta})}.
       \label{eq:app_empirical_reference_distance}
\end{align}
The constants $C_E,\Gamma_0$ depend only on $d,r,\tau_0$ and the constants
in \cref{assum:particle_population}; $\alpha,A$ may additionally
depend on $K$. None of these depends on $n,M$, or the realized sample.
\end{lemma}

\begin{proof}
We use the direction in \cref{lem:app_population_direction} to construct
$\bbeta^\dagger$. As in \cref{lem:app_population_reference}, we first
bound the likelihood derivative from below and then bound the second
derivative to choose a step size. Here we also use
\cref{lem:app_particle_empirical_comparison} to bound the sampling
errors in both derivatives.

Set $\Gamma_0=2\epsilon_0^2$ as in the proof of
\cref{lem:app_population_reference}. With $C_{\mathrm{cmp}}$ in \eqref{def_C_rel}, choose
$$
    C_E=\max\left\{
       1,\left({4C_{\mathrm{cmp}}\over3\rho^2}\right)^2,
       {4C_{\mathrm{cmp}}\over3\rho^2}\right\}.
$$
By \eqref{eq:app_population_kl_l2}, \eqref{log_L_gap_emp} implies
$$
    0<P_X(p_{G_{\bbeta}}-p^*)^2\le\epsilon_0^2.
$$
We omit the particles with $w_j=0$, which does not change
$G_{\bbeta}$, and apply \cref{lem:app_population_direction}
to the remaining particles. This gives $(\dot w_j,\xi_j)$ for
each index with $w_j>0$. For every omitted index, define
$\dot w_j=0$ and $\xi_j=0$.

By \eqref{eq:app_empirical_parameter_derivative} and the definition
of $g_{n,G}$ in \eqref{eq:particle_potential},
$$
    D\ell_n(G_{\bbeta})[\dot w,\xi]
       =\sum_{j=1}^M\dot w_jg_{n,G_{\bbeta}}(\beta_j)
          +D\ell_n(G_{\bbeta})[0,\xi].
$$
The two terms on the right are the contributions from weight and
location changes, respectively. The weights $w\in\cM_n(\bbeta)$
maximize $\ell_n$ for fixed $\bbeta$, so the first term is zero:
By \eqref{eq:app_particle_weight_kkt}, the partial derivatives with
respect to the positive weights are all equal. Since $\dot w_j=0$
on zero-weight coordinates and $\sum_j\dot w_j=0$,
$$
    \sum_j\dot w_jg_{n,G_{\bbeta}}(\beta_j)
       =\left\{\sum_kw_kg_{n,G_{\bbeta}}(\beta_k)\right\}
          \sum_j\dot w_j=0.
$$
Thus we keep the weights fixed and use only the location direction
$\xi$. The feasible path
$$
    G_{\bbeta+s\xi}=\sum_{j=1}^Mw_j\delta_{\beta_j+s\xi_j},
    \qquad p_s=p_{G_{\bbeta+s\xi}},\qquad 0\le s\le1,
$$
has the same likelihood derivative at $s=0$ as the joint direction
$(\dot w,\xi)$.

Taylor's integral formula gives, for $0<\lambda\le1$,
$$
\begin{aligned}
    \ell_n(G_{\bbeta+\lambda\xi})-\ell_n(G_{\bbeta})
    ={}&\lambda\left.\frac{d}{ds}\ell_n(G_{\bbeta+s\xi})\right|_{s=0}
        +\int_0^\lambda(\lambda-s)\frac{d^2}{ds^2}\ell_n(G_{\bbeta+s\xi})\,\rd s.
\end{aligned}
$$
We bound the first term from below and control the absolute value of
the integral remainder.

First, the Cauchy--Schwarz calculation following
\eqref{eq:app_population_directional_gain} uses
\eqref{eq:app_direction_prediction} to give
$$
\begin{aligned}
    D\cL(G_{\bbeta})[\dot w,\xi]
    &=P_X\{r_{G_{\bbeta}}Dp_{G_{\bbeta}}[\dot w,\xi]\}\\
    &\ge3P_X(p_{G_{\bbeta}}-p^*)^2
     \ge3\rho^2\{\cL(G^*)-\cL(G_{\bbeta})\},
\end{aligned}
$$
where the last inequality follows from \eqref{eq:app_population_kl_l2}.
This bound on the joint derivative does not require population weight
optimality. Applying \eqref{eq:app_empirical_direction_comparison} and
then the lower bound in \eqref{log_L_gap_emp} gives
\begin{align*}
    D\ell_n(G_{\bbeta})[\dot w,\xi]
    &\ge3\rho^2\{\cL(G^*)-\cL(G_{\bbeta})\}
          -C_{\mathrm{cmp}}
             \{\sqrt{\{\cL(G^*)-\cL(G_{\bbeta})\}\varepsilon_n(\delta)}+\varepsilon_n(\delta)\}\\
    &\ge\left\{3\rho^2
          -{C_{\mathrm{cmp}}\over\sqrt{C_E}}
          -{C_{\mathrm{cmp}}\over C_E}\right\}\{\cL(G^*)-\cL(G_{\bbeta})\}
     \ge {3\rho^2\over2}\{\cL(G^*)-\cL(G_{\bbeta})\}.
\end{align*}
The choice of $C_E$ bounds each of
$C_{\mathrm{cmp}}/\sqrt{C_E}$ and $C_{\mathrm{cmp}}/C_E$
by $3\rho^2/4$. Since the fixed-weight path has the same initial
derivative, we obtain
\begin{equation}
    \left.{d\over ds} \ell_n(G_{\bbeta+s\xi})\right|_{s=0} = D\ell_n(G_{\bbeta})[\dot w,\xi]  \ge {3\rho^2\over2}\{\cL(G^*)-\cL(G_{\bbeta})\}.
    \label{eq:app_empirical_reference_initial_derivative}
\end{equation}

Second, we bound the absolute value of the second derivative of
$\ell_n(G_{\bbeta+s\xi})$. By
\eqref{eq:app_direction_path_derivatives},
\eqref{eq:app_gauge_by_prediction}, and
\eqref{eq:app_population_kl_l2},
$$
    \sup_{0\le s\le1}
       \{\|p_s'\|_\infty+\|p_s''\|_\infty\}
    \le C_0\mfa(G_{\bbeta},G^*)
    \le {C_0\over\sqrt2\,c_{\mathrm{inv}}}
       \sqrt{\cL(G^*)-\cL(G_{\bbeta})}.
$$
Also, $p_s-p^*=p_0-p^*+\int_0^s p_t'\,\rd t$ implies, for
$0\le s\le1$,
\begin{align*}
    \{P_X(p_s-p^*)^2\}^{1/2}
    &\le\{P_X(p_0-p^*)^2\}^{1/2}
          +\int_0^s\|p_t'\|_\infty\,\rd t\\
    &\le\left({1\over\sqrt2}+{C_0\over\sqrt2\,c_{\mathrm{inv}}}\right)
       \sqrt{\cL(G^*)-\cL(G_{\bbeta})}.
\end{align*}
Differentiating twice gives
$$
   - {d^2\over ds^2}\ell_n(G_{\bbeta+s\xi})
    =P_n\left[
       \left\{{Y\over p_s^2}+{1-Y\over(1-p_s)^2}\right\}(p_s')^2
       +{p_s-Y\over p_s(1-p_s)}p_s''
       \right].
$$
The first term inside the empirical average is nonnegative and at most
$$
\{C_0^2/(2\rho^2c_{\mathrm{inv}}^2)\}
\{\cL(G^*)-\cL(G_{\bbeta})\}.
$$
For the second term, taking the conditional expectation of $Y$ and
applying Cauchy--Schwarz gives
$$
\begin{aligned}
    \left|P\left\{\frac{p_s-Y}{p_s(1-p_s)}p_s''\right\}\right|
    &=\left|P_X\left\{\frac{p_s-p^*}{p_s(1-p_s)}p_s''\right\}\right|\\
    &\le\rho^{-2}\{P_X(p_s-p^*)^2\}^{1/2}
                    \{P_X(p_s'')^2\}^{1/2}\\
    &\le{C_0(c_{\mathrm{inv}}+C_0)\over2\rho^2c_{\mathrm{inv}}^2}
                  \{\cL(G^*)-\cL(G_{\bbeta})\}.
\end{aligned}
$$
To bound its empirical deviation, write
$$
    r_Y(p_s)p_s''
       =\sum_{j:\xi_j\ne0}w_j\|\xi_j\|_2^2\,
          r_Y(p_s)\sigma''\{X^\T(\beta_j+s\xi_j)\}
             \left(X^\T\frac{\xi_j}{\|\xi_j\|_2}\right)^2.
$$
Each function multiplying $w_j\|\xi_j\|_2^2$ belongs to
$\mathcal J$ in \eqref{def_cJ}: $p_s\in\cF_r$, the locations
stay in $\cB_r$, and the direction is a unit vector. Moreover,
\eqref{eq:app_direction_norm}, \eqref{eq:app_gauge_by_prediction},
and \eqref{eq:app_population_kl_l2} give
$$
    \sum_{j=1}^M w_j\|\xi_j\|_2^2
       \le\mfa(G_{\bbeta},G^*)
       \le{\sqrt{\cL(G^*)-\cL(G_{\bbeta})}\over\sqrt2\,c_{\mathrm{inv}}}.
$$
Since
$(p_s-Y)/\{p_s(1-p_s)\}=-r_Y(p_s)$,
\eqref{eq:app_empirical_curvature_comparison} yields, uniformly in $s$,
\begin{align*}
    \left|(P_n-P)\left[
       {p_s-Y\over p_s(1-p_s)}p_s''\right]\right|
    &\le C_{\mathrm{cmp}}\sqrt{\varepsilon_n(\delta)}
             \sum_jw_j\|\xi_j\|_2^2\\
    &\le{C_{\mathrm{cmp}}\over\sqrt2\,c_{\mathrm{inv}}}
             \sqrt{\{\cL(G^*)-\cL(G_{\bbeta})\}\varepsilon_n(\delta)}\\
    & \le {C_{\mathrm{cmp}}\over\sqrt2\,c_{\mathrm{inv}}}\{\cL(G^*)-\cL(G_{\bbeta})\},
\end{align*}
where the last inequality follows from \eqref{log_L_gap_emp} and
$C_E\ge1$, which imply
$\varepsilon_n(\delta)\le\cL(G^*)-\cL(G_{\bbeta})$.
Combining the bounds on the two terms gives, with
$C_{\mathrm{ref}}=C_F+C_{\mathrm{cmp}}/(\sqrt2\,c_{\mathrm{inv}})$
and the same $C_F=(2C_0^2+C_0c_{\mathrm{inv}})/(2c_{\mathrm{inv}}^2\rho^2)$
as in \cref{lem:app_population_reference},
$$
    \sup_{0\le s\le1}
       \left|{d^2\over ds^2}\ell_n(G_{\bbeta+s\xi})\right|
       \le C_{\mathrm{ref}}\{\cL(G^*)-\cL(G_{\bbeta})\}.
$$

Combining this second-derivative bound with
\eqref{eq:app_empirical_reference_initial_derivative} in Taylor's
integral formula gives, for $0<\lambda\le1$,
$$
    \ell_n(G_{\bbeta+\lambda\xi})-\ell_n(G_{\bbeta})
    \ge\left(\frac{3\rho^2\lambda}{2}
                  -\frac{C_{\mathrm{ref}}\lambda^2}{2}\right)
                \{\cL(G^*)-\cL(G_{\bbeta})\}.
$$
Choose
$$
    \lambda=\min\left\{
       {1\over2},\,{3\rho^2\over2C_{\mathrm{ref}}},
       {3\rho^2\sqrt2\,c_{\mathrm{inv}}\over16K}
       \right\}.
$$
Then $3\rho^2\lambda/2-C_{\mathrm{ref}}\lambda^2/2
\ge3\rho^2\lambda/4$. Together with the bound on
$\sum_jw_j\|\xi_j\|_2^2$, this proves the claims with
$$
    \bbeta^\dagger=\bbeta+\lambda\xi,\qquad
    \alpha={3\rho^2\lambda\over4},\qquad
    A={\lambda^2\over\sqrt2\,c_{\mathrm{inv}}}.
$$
The first restriction on $\lambda$ ensures $0<\alpha<1$, and the
last gives $KA/\alpha=4K\lambda/(3\sqrt2\,\rho^2c_{\mathrm{inv}})
\le1/4$. Only $\lambda,\alpha,A$ depend on $K$;
$C_E,\Gamma_0$ were fixed before $K$ was chosen.
All constants are independent of $n,M$, the weights, and the sample.
\end{proof}

We next compare the fixed-weight iterates with this reference and show that the likelihood improves over the inner steps of location updates. This is analogous to \cref{prop:app_population_block_gain}.

\begin{lemma}
\label{lem:app_particle_empirical_block_gain}
On the event in \cref{lem:app_particle_empirical_comparison}, there
are constants $C_E\ge1$, $\Gamma_0>0$, and $\kappa\in(0,1)$,
depending only on $d,r,\tau_0$ and the constants in
\cref{assum:particle_population}, such that the following holds.
Suppose $w\in\cM_n(\bbeta^{(0)})$ and
$$
    C_E\varepsilon_n(\delta)
       \le \cL(G^*)-\cL(G_{\bbeta^{(0)}})\le\Gamma_0.
$$
For any integer $H\ge1$ and $0<\eta\le27\rho/(4R_X^2)$, let
$\bbeta^{(h)}$, $0\le h\le H$, be the subsequent
fixed-weight empirical location iterates. Then
\begin{equation}
    \{\ell_n(G_{\bbeta^{(H)}})-\ell_n(G_{\bbeta^{(0)}})\}
       \ge\kappa \{\cL(G^*)-\cL(G_{\bbeta^{(0)}})\}\min\{\eta H\sqrt{\cL(G^*)-\cL(G_{\bbeta^{(0)}})},1\}.
    \label{eq:app_empirical_block_gain}
\end{equation}
The constants are independent of $n,M,H,\eta$, and the sample.
\end{lemma}

\begin{proof}
This is the empirical version of \cref{prop:app_population_block_gain}.
We follow its proof, but must also control the population likelihood
gaps along the empirical iterates and the sampling error in the
likelihood comparison.

Take $C_E,\Gamma_0$ from
\cref{lem:app_particle_empirical_reference}. The weights
$w\in\cM_n(\bbeta^{(0)})$ remain fixed throughout the $H$
location updates. As in the population proof, write
$$
    \Delta_0:=\cL(G^*)-\cL(G_{\bbeta^{(0)}}),
    \qquad
    \bdelta_h:=\bbeta^{(h+1)}-\bbeta^{(h)},\quad 0\le h<H.
$$
Since $\|X_i\|_2\le R_X$ and
$\eta\le L_X^{-1}$, with $L_X=4R_X^2/(27\rho)$, the empirical
version of \eqref{eq:app_fixed_weight_descent} in
\cref{lem:app_population_one_step_ascent}(ii) gives
\begin{equation}
    \ell_n(G_{\bbeta^{(h+1)}})-\ell_n(G_{\bbeta^{(h)}})
    \ge{1\over2\eta}\|\bdelta_h\|_w^2,
    \qquad 0\le h<H.
    \label{eq:app_empirical_location_gain}
\end{equation}
Thus the empirical likelihood does not decrease. The population
likelihood need not increase, so we use
\eqref{eq:app_empirical_relative_likelihood} to obtain, for
$0\le h\le H$,
\begin{align*}
    {7\over8}\{\cL(G^*)-\cL(G_{\bbeta^{(h)}})\}
    &\le\ell_n(G^*)-\ell_n(G_{\bbeta^{(h)}})
          +C_{\mathrm{rel}}\varepsilon_n(\delta)\\
    &\le\ell_n(G^*)-\ell_n(G_{\bbeta^{(0)}})
          +C_{\mathrm{rel}}\varepsilon_n(\delta)\\
    &\le{9\over8}\Delta_0+2C_{\mathrm{rel}}\varepsilon_n(\delta).
\end{align*}
The middle inequality uses \eqref{eq:app_empirical_location_gain}.
Since $\Delta_0\ge C_E\varepsilon_n(\delta)\ge\varepsilon_n(\delta)$,
it follows that
$$
    \cL(G^*)-\cL(G_{\bbeta^{(h)}})
       \le{9+16C_{\mathrm{rel}}\over7}\Delta_0,
    \qquad 0\le h\le H.
$$

Next, we bound the empirical likelihood difference between two
configurations with the same weights. For
$\boldsymbol z,\boldsymbol v\in\cB_r^M$, put $d_j=v_j-z_j$.
Taylor's integral formula gives
$$
\begin{aligned}
    R_{\boldsymbol z,\boldsymbol v}(x)
    &:=p_{G_{\boldsymbol v}}(x)-p_{G_{\boldsymbol z}}(x)
       -\sum_jw_j\sigma'(x^\T z_j)x^\T d_j\\
    &=\sum_jw_j\int_0^1(1-s)
       \sigma''\{x^\T(z_j+sd_j)\}(x^\T d_j)^2\,\rd s.
\end{aligned}
$$
As in the proof of \eqref{eq:app_fixed_weight_lower_model},
\eqref{eq:app_population_location_derivatives} gives
$\|R_{\boldsymbol z,\boldsymbol v}\|_\infty
\le\{R_X^2/(12\sqrt3)\}\|\boldsymbol v-\boldsymbol z\|_w^2$.
Taking the conditional expectation of $Y$ given $X$, then using
Cauchy--Schwarz and \eqref{eq:app_population_kl_l2}, yields
\begin{align*}
    \left|P\{r_Y(p_{G_{\boldsymbol z}})
                       R_{\boldsymbol z,\boldsymbol v}\}\right|
    &=\left|P_X\{r_{G_{\boldsymbol z}}
                        R_{\boldsymbol z,\boldsymbol v}\}\right|\\
    &\le\rho^{-2}\{P_X(p_{G_{\boldsymbol z}}-p^*)^2\}^{1/2}
                  \|R_{\boldsymbol z,\boldsymbol v}\|_\infty\\
    &\le {R_X^2\over12\sqrt6\,\rho^2}
           \sqrt{\cL(G^*)-\cL(G_{\boldsymbol z})}
           \|\boldsymbol v-\boldsymbol z\|_w^2.
\end{align*}
For the empirical deviation, $p_{G_{\boldsymbol z}}\in\cF_r$,
the segment $z_j+sd_j$ stays in $\cB_r$, and
$d_j/\|d_j\|_2$ is a unit vector whenever $d_j\ne0$.
The corresponding functions therefore belong to $\mathcal J$ in
\eqref{def_cJ}. Applying \eqref{eq:app_empirical_curvature_comparison}
gives
\begin{align*}
    \left|(P_n-P)\{r_Y(p_{G_{\boldsymbol z}})
                          R_{\boldsymbol z,\boldsymbol v}\}\right|
    &\le\sum_{j:d_j\ne0}w_j\|d_j\|_2^2
         \int_0^1(1-s)C_{\mathrm{cmp}}
                  \sqrt{\varepsilon_n(\delta)}\,\rd s\\
    &= {C_{\mathrm{cmp}}\over2}\sqrt{\varepsilon_n(\delta)}
          \|\boldsymbol v-\boldsymbol z\|_w^2.
\end{align*}
For each observation, the function
$t\mapsto Y\log t+(1-Y)\log(1-t)$ is concave. Its tangent
inequality gives
$$
    \ell_n(G_{\boldsymbol v})-\ell_n(G_{\boldsymbol z})
       \le P_n\{r_Y(p_{G_{\boldsymbol z}})
                 (p_{G_{\boldsymbol v}}-p_{G_{\boldsymbol z}})\}.
$$
Substituting the Taylor expansion and the two remainder bounds yields
\begin{align}
    \ell_n(G_{\boldsymbol v})-\ell_n(G_{\boldsymbol z})
    &\le\sum_jw_j
       \langle\nabla g_{n,G_{\boldsymbol z}}(z_j),v_j-z_j\rangle
       \notag\\
    &\quad+\left\{
       {R_X^2\over12\sqrt6\,\rho^2}
          \{\cL(G^*)-\cL(G_{\boldsymbol z})\}^{1/2}
       +{C_{\mathrm{cmp}}\over2}\sqrt{\varepsilon_n(\delta)}
       \right\}\|\boldsymbol v-\boldsymbol z\|_w^2.
\label{eq:app_empirical_fixed_weight_lower_model}
\end{align}

For $\boldsymbol z=\bbeta^{(h)}$, the bound on
$\cL(G^*)-\cL(G_{\bbeta^{(h)}})$ and
$\varepsilon_n(\delta)\le\Delta_0$ show that the coefficient of
the squared distance in \eqref{eq:app_empirical_fixed_weight_lower_model}
is at most $K\sqrt{\Delta_0}$, where
$$
    K={R_X^2\over12\sqrt6\,\rho^2}
          \sqrt{{9+16C_{\mathrm{rel}}\over7}}
          +{C_{\mathrm{cmp}}\over2}.
$$
Since $w\in\cM_n(\bbeta^{(0)})$ and
$C_E\varepsilon_n(\delta)\le\Delta_0\le\Gamma_0$, applying
\cref{lem:app_particle_empirical_reference} at $\bbeta^{(0)}$
with this $K$ gives $\bbeta^\dagger\in\cB_r^M$ and
$\alpha,A$ with $KA\le\alpha/4$. By
\eqref{eq:app_empirical_reference_gain} and
\eqref{eq:app_empirical_reference_distance},
$$
    \ell_n(G_{\bbeta^\dagger})-\ell_n(G_{\bbeta^{(0)}})
       \ge\alpha\Delta_0,
    \qquad
    \|\bbeta^\dagger-\bbeta^{(0)}\|_w^2\le A\sqrt{\Delta_0}.
$$
We keep $\bbeta^\dagger,\alpha,A$ fixed throughout the block and set
$$
    \kappa=\min\{\alpha/2,\alpha^2/(32A)\}.
$$

If $\ell_n(G_{\bbeta^{(H)}})-\ell_n(G_{\bbeta^{(0)}})
\ge(\alpha/2)\Delta_0$, then \eqref{eq:app_empirical_block_gain}
follows from $\kappa\le\alpha/2$. Now suppose
$$
    \ell_n(G_{\bbeta^{(H)}})-\ell_n(G_{\bbeta^{(0)}})
       <{\alpha\over2}\Delta_0.
$$
We show that, in this case, every update satisfies
$$
    \|\bdelta_h\|_w^2
       \ge{\eta^2\alpha^2\over16A}\Delta_0^{3/2},
    \qquad 0\le h<H.
$$
Indeed, this bound and \eqref{eq:app_empirical_location_gain} imply
$$
\begin{aligned}
    \ell_n(G_{\bbeta^{(H)}})-\ell_n(G_{\bbeta^{(0)}})
    &\ge{1\over2\eta}\sum_{h=0}^{H-1}\|\bdelta_h\|_w^2\\
    &\ge{\alpha^2\over32A}\eta H\Delta_0^{3/2}
     \ge\kappa\Delta_0\min\{\eta H\sqrt{\Delta_0},1\},
\end{aligned}
$$
which is \eqref{eq:app_empirical_block_gain}.

To prove the lower bound on $\|\bdelta_h\|_w^2$, we compare each
update with $\bbeta^\dagger$. Empirical likelihood monotonicity
and \eqref{eq:app_empirical_reference_gain} give, for $0\le h<H$,
\begin{align*}
    \ell_n(G_{\bbeta^\dagger})-\ell_n(G_{\bbeta^{(h+1)}})
    &\ge\ell_n(G_{\bbeta^\dagger})-\ell_n(G_{\bbeta^{(H)}})\\
    &=\{\ell_n(G_{\bbeta^\dagger})-\ell_n(G_{\bbeta^{(0)}})\}
       -\{\ell_n(G_{\bbeta^{(H)}})-\ell_n(G_{\bbeta^{(0)}})\}\\
    &>{\alpha\over2}\Delta_0.
\end{align*}
For the $h$-th update, \eqref{eq:app_projection_inequality} gives
$$
    \left\langle
       \eta\nabla g_{n,G_{\bbeta^{(h)}}}(\beta_j^{(h)})-\delta_{h,j},
       \beta_j^\dagger-\beta_j^{(h+1)}\right\rangle\le0,
    \qquad j\in[M].
$$
Apply the empirical version of \eqref{eq:app_fixed_weight_upper_model}
from \cref{lem:app_population_one_step_ascent}(i) to
$(\bbeta^{(h)},\bbeta^{(h+1)})$, and apply
\eqref{eq:app_empirical_fixed_weight_lower_model} to
$(\boldsymbol z,\boldsymbol v)=(\bbeta^{(h)},\bbeta^\dagger)$.
The former bounds the actual likelihood increase from below, while
the latter bounds the increase to $\bbeta^\dagger$ from above.
Subtracting the two bounds and using the projection inequality yields
\begin{equation}
\begin{aligned}
    &\ell_n(G_{\bbeta^\dagger})-\ell_n(G_{\bbeta^{(h+1)}})\\
    &\quad\le\sum_jw_j
       \langle\nabla g_{n,G_{\bbeta^{(h)}}}(\beta_j^{(h)}),
               \beta_j^\dagger-\beta_j^{(h+1)}\rangle
       +{L_X\over2}\|\bdelta_h\|_w^2
       +K\sqrt{\Delta_0}\|\bbeta^\dagger-\bbeta^{(h)}\|_w^2\\
    &\quad\le{1\over\eta}\langle\bdelta_h,\bbeta^\dagger-\bbeta^{(h+1)}\rangle_w
       +{L_X\over2}\|\bdelta_h\|_w^2
       +K\sqrt{\Delta_0}\|\bbeta^\dagger-\bbeta^{(h)}\|_w^2\\
    &\quad={1\over\eta}\langle\bdelta_h,\bbeta^\dagger-\bbeta^{(h)}\rangle_w
       -\left({1\over\eta}-{L_X\over2}\right)\|\bdelta_h\|_w^2
       +K\sqrt{\Delta_0}\|\bbeta^\dagger-\bbeta^{(h)}\|_w^2\\
    &\quad\le{1\over\eta}\langle\bdelta_h,\bbeta^\dagger-\bbeta^{(h)}\rangle_w
       -{1\over2\eta}\|\bdelta_h\|_w^2
       +K\sqrt{\Delta_0}\|\bbeta^\dagger-\bbeta^{(h)}\|_w^2.
\end{aligned}
\label{eq:app_empirical_step_comparison}
\end{equation}
The equality uses
$\bbeta^\dagger-\bbeta^{(h+1)}
=\bbeta^\dagger-\bbeta^{(h)}-\bdelta_h$, and the last inequality
uses $\eta\le L_X^{-1}$.

We next prove by induction that
$$
    \|\bbeta^{(h)}-\bbeta^\dagger\|_w^2\le A\sqrt{\Delta_0},
    \qquad 0\le h\le H.
$$
For $h=0$, this is \eqref{eq:app_empirical_reference_distance}.
Suppose the bound holds for some $h<H$. Since
$\bbeta^{(h+1)}=\bbeta^{(h)}+\bdelta_h$, expanding the squared
distance gives
$$
    {2\over\eta}\langle\bdelta_h,\bbeta^\dagger-\bbeta^{(h)}\rangle_w
    ={1\over\eta}\left\{\|\bbeta^{(h)}-\bbeta^\dagger\|_w^2
       -\|\bbeta^{(h+1)}-\bbeta^\dagger\|_w^2
       +\|\bdelta_h\|_w^2\right\}.
$$
Substituting this identity into \eqref{eq:app_empirical_step_comparison}
cancels the terms in $\|\bdelta_h\|_w^2$ and yields
$$
\begin{aligned}
    \|\bbeta^{(h+1)}-\bbeta^\dagger\|_w^2
    &\le\|\bbeta^{(h)}-\bbeta^\dagger\|_w^2\\
    &\quad+2\eta K\sqrt{\Delta_0}\,
       \|\bbeta^{(h)}-\bbeta^\dagger\|_w^2
       -2\eta\{\ell_n(G_{\bbeta^\dagger})
          -\ell_n(G_{\bbeta^{(h+1)}})\}.
\end{aligned}
$$
By the induction hypothesis, the second term on the right is at most
$2\eta KA\Delta_0$. The last term is less than
$-\eta\alpha\Delta_0$, since the likelihood at $\bbeta^\dagger$
exceeds that at $\bbeta^{(h+1)}$ by more than
$\alpha\Delta_0/2$. Thus $KA\le\alpha/4$ gives
$$
\begin{aligned}
    \|\bbeta^{(h+1)}-\bbeta^\dagger\|_w^2
    &\le\|\bbeta^{(h)}-\bbeta^\dagger\|_w^2
          +(2KA-\alpha)\eta\Delta_0\\
    &\le\|\bbeta^{(h)}-\bbeta^\dagger\|_w^2
          -{\alpha\eta\over2}\Delta_0
     \le A\sqrt{\Delta_0},
\end{aligned}
$$
which completes the induction.

Finally, applying Cauchy--Schwarz to
\eqref{eq:app_empirical_step_comparison} and using
$-\|\bdelta_h\|_w^2/(2\eta)\le0$ gives
\begin{align*}
    {\alpha\over2}\Delta_0
    &<\ell_n(G_{\bbeta^\dagger})-\ell_n(G_{\bbeta^{(h+1)}})\\
    &\le{1\over\eta}\|\bdelta_h\|_w
          \|\bbeta^{(h)}-\bbeta^\dagger\|_w
       +K\sqrt{\Delta_0}\|\bbeta^{(h)}-\bbeta^\dagger\|_w^2\\
    &\le{\sqrt A\over\eta}\Delta_0^{1/4}\|\bdelta_h\|_w
       +KA\Delta_0\\
    &\le{\sqrt A\over\eta}\Delta_0^{1/4}\|\bdelta_h\|_w
       +{\alpha\over4}\Delta_0.
\end{align*}
Rearranging and squaring yields
$\|\bdelta_h\|_w^2\ge\eta^2\alpha^2\Delta_0^{3/2}/(16A)$
for every $0\le h<H$, as required. Combining the two cases proves
\eqref{eq:app_empirical_block_gain}. All constants depend
only on the stated population class.
\end{proof}

\subsection{Proof of \cref{thm:particle_uniform_empirical}}
\label{app:particle_uniform_empirical_proof}

\begin{proof}
Take $\Gamma_0,C_E,\kappa$ from
\cref{lem:app_particle_empirical_block_gain}. Since
$\varepsilon_n(\delta)\to0$ for fixed $\delta$, there is a finite
$n_0$, depending only on the stated class and $\delta$, such that
$$
    C_{\mathrm{rel}}\varepsilon_n(\delta)\le\Gamma_0/16,
    \qquad n\ge n_0.
$$
We first show that random initialization and likelihood
monotonicity keep the population likelihood gaps below $\Gamma_0$.
The block bound then applies whenever the empirical gap exceeds a
constant multiple of $\varepsilon_n(\delta)$. Below that level,
likelihood monotonicity alone suffices.

Use the construction in \cref{lem:app_population_initialization} with
$$
    h_{\mathrm{emp}}
       ={\tau_0\over2}\wedge{s_-\over4}
            \wedge\sqrt{{2\Gamma_0\rho^2\over R_X^2}},
    \qquad
    M_0=K_0\vee
       \left\lceil
         (r/h_{\mathrm{emp}})^d\log(2K_0/\delta)
       \right\rceil.
$$
The restrictions $h_{\mathrm{emp}}\le\tau_0/2$ and
$h_{\mathrm{emp}}\le s_-/4$ ensure that the balls around the true
atoms are contained in $\cB_r$ and disjoint. The probability
calculation in that proof shows that all these balls contain initial
particles with probability at least $1-\delta/2$ when $M\ge M_0$.
By the union bound, this event and the event in
\cref{lem:app_particle_empirical_comparison} hold jointly with
probability at least $1-\delta$. We work on their intersection
for the rest of the proof.

For each $a\in[K^*]$, choose an index $j_a$ such that
$\beta_{j_a}^0\in\cB(\beta_a^*,h_{\mathrm{emp}})$.
The balls are disjoint, so the indices $j_1,\ldots,j_{K^*}$ are
distinct. Assign weight $w_a^*$ to $\beta_{j_a}^0$ and zero
weight to every other initial particle, and define
$$
    \wt G=\sum_{a=1}^{K^*}w_a^*\delta_{\beta_{j_a}^0}.
$$
The weights are nonnegative and sum to one, so $\wt G$ is feasible
on the initial locations. We use $\wt G$ only for comparison;
the algorithm does not need the true atoms or weights.

Write $\wt p=p_{\wt G}$. For $x\in\cB_{R_X}$, the mean value
theorem and \eqref{eq:app_population_scalar_derivatives} give
\begin{align*}
    |\wt p(x)-p^*(x)|
    &\le\sum_{a=1}^{K^*}w_a^*
       |\sigma(x^\T\beta_{j_a}^0)-\sigma(x^\T\beta_a^*)|\\
    &\le{R_X\over4}\sum_{a=1}^{K^*}w_a^*
       \|\beta_{j_a}^0-\beta_a^*\|_2
     \le{R_X\over4}h_{\mathrm{emp}}.
\end{align*}
Consequently,
$P_X(\wt p-p^*)^2\le R_X^2h_{\mathrm{emp}}^2/16$.
For the scalar loss $f_a$ in
\eqref{eq:app_population_scalar_loss_bounds},
$f_{p^*(x)}'(p^*(x))=0$. Both $p^*(x)$ and $\wt p(x)$
belong to $[\rho,1-\rho]$, so the segment between them remains
in this interval. Taylor's integral formula and
$f_a''\le\rho^{-2}$ on this interval give
\begin{align*}
    \cL(G^*)-\cL(\wt G)
    &=\EE_X\int_0^1(1-s)
       f_{p^*(X)}''\{p^*(X)+s(\wt p(X)-p^*(X))\}
       (\wt p(X)-p^*(X))^2\,\rd s\\
    &\le{1\over2\rho^2}P_X(\wt p-p^*)^2.
\end{align*}

Assigning each true weight to a particle in the corresponding ball
gives the same feasible candidate $\wt G $, now satisfying
\begin{equation}
    \cL(G^*)-\cL(\wt G)
       \le{R_X^2\over32\rho^2}h_{\mathrm{emp}}^2
       \le{\Gamma_0\over16}.
    \label{eq:app_empirical_initial_candidate}
\end{equation}
The last inequality follows from
$h_{\mathrm{emp}}^2\le2\Gamma_0\rho^2/R_X^2$.
Since $\wt G$ uses only the initial locations and
$w^0\in\cM_n(\bbeta^0)$ maximizes $\ell_n$ over their weights,
we have
$\ell_n(G_0)\ge\ell_n(\wt G)$.

As in the proof of \cref{thm:particle_population}, we first verify
the conditions for applying the block improvement bound along the
iterates. Here the empirical likelihood is nondecreasing, but the
population likelihood need not be. We therefore use the likelihood
comparison in \eqref{eq:app_empirical_relative_likelihood}.
Since $\max_i\|X_i\|_2\le R_X$ and
$0<\eta\le27\rho/(4R_X^2)$, the empirical version of
\cref{lem:app_population_one_step_ascent}(ii) gives, for every
$t\ge0$ and $0\le h<H$,
$$
    \ell_n(G^{t,h+1})-\ell_n(G^{t,h})
    \ge{1\over2\eta}
       \|\bbeta^{t,h+1}-\bbeta^{t,h}\|_{w^t}^2\ge0.
$$
This bound does not require a small population likelihood gap.
Exact reweighting then gives
$$
    \ell_n(G_{t+1})\ge\ell_n(G^{t,H})
       \ge\ell_n(G^{t,0})=\ell_n(G_t).
$$
Consequently, $\ell_n(G_t)\ge\ell_n(G_0)\ge\ell_n(\wt G)$
for all $t\ge0$.

Write the population and empirical likelihood gaps as
$$
    \Delta_t:=\cL(G^*)-\cL(G_t),
    \qquad
    a_t:=\ell_n(\wh G_r)-\ell_n(G_t).
$$
Applying \eqref{eq:app_empirical_relative_likelihood} to
$G_t$ and $\wt G$ gives
$$
\begin{aligned}
    {7\over8}\Delta_t
    &\le\ell_n(G^*)-\ell_n(G_t)
          +C_{\mathrm{rel}}\varepsilon_n(\delta)\\
    &\le\ell_n(G^*)-\ell_n(\wt G)
          +C_{\mathrm{rel}}\varepsilon_n(\delta)\\
    &\le{9\over8}\{\cL(G^*)-\cL(\wt G)\}
          +2C_{\mathrm{rel}}\varepsilon_n(\delta)\\
    &\le {9\Gamma_0\over128}+{\Gamma_0\over8}
     ={25\Gamma_0\over128}.
\end{aligned}
$$
The last step uses \eqref{eq:app_empirical_initial_candidate} and
the choice of $n_0$. Thus
\begin{equation}
    \Delta_t\le{25\over112}\Gamma_0\le\Gamma_0,
    \qquad t\ge0.
    \label{eq:app_empirical_iterates_local}
\end{equation}

We next derive a recursion for $a_t$ using
\cref{lem:app_particle_empirical_block_gain}, the empirical
counterpart of \cref{prop:app_population_block_gain}.
Set $C_1=\max\{4C_{\mathrm{rel}},9C_E/4\}$, and suppose that
$$
    a_t>C_1\varepsilon_n(\delta).
$$
The upper bound in \eqref{eq:app_empirical_gap_comparison} and
$C_1\ge4C_{\mathrm{rel}}$ give
\begin{align*}
    \Delta_t
    &\ge{8\over9}\{a_t-2C_{\mathrm{rel}}\varepsilon_n(\delta)\}\\
    &\ge{4\over9}a_t
     >{4C_1\over9}\varepsilon_n(\delta)
     \ge C_E\varepsilon_n(\delta).
\end{align*}
Together with \eqref{eq:app_empirical_iterates_local}, this gives
$C_E\varepsilon_n(\delta)<\Delta_t\le\Gamma_0$.
Moreover, $w^t\in\cM_n(\bbeta^t)$ at the start of block $t$,
and these weights remain fixed during its $H$ location updates.
Hence all conditions of \cref{lem:app_particle_empirical_block_gain}
hold. Applying that lemma and then using exact reweighting at the
end of the block, we obtain
\begin{align*}
    a_t-a_{t+1}
    &=\ell_n(G_{t+1})-\ell_n(G_t)\\
    &\ge\ell_n(G^{t,H})-\ell_n(G^{t,0})\\
    &\ge\kappa\Delta_t\min\{\eta H\sqrt{\Delta_t},1\}\\
    &\ge{4\kappa\over9}a_t
         \min\{(2/3)\eta H\sqrt{a_t},1\}\\
    &\ge{8\kappa\over27}a_t\min\{\eta H\sqrt{a_t},1\}.
\end{align*}
The last two steps use $\Delta_t\ge(4/9)a_t$ and
$\min\{(2/3)x,1\}\ge(2/3)\min\{x,1\}$ for $x\ge0$.
Thus
$$
    a_t>C_1\varepsilon_n(\delta)
    \quad\Longrightarrow\quad
    a_{t+1}\le a_t-{8\kappa\over27}a_t
                       \min\{\eta H\sqrt{a_t},1\}.
$$

Since $\wh G_r$ maximizes $\ell_n$ over $\cP(\cB_r)$,
the sequence $\{a_t\}$ is nonnegative. It is also nonincreasing
by the likelihood monotonicity proved above. Moreover,
$\ell_n(G)\in[\log\rho,0]$ for every $G\in\cP(\cB_r)$, so
$a_0\le-\log\rho$. If $a_t\le C_1\varepsilon_n(\delta)$,
then $a_k\le C_1\varepsilon_n(\delta)$ for every $k\ge t$.
We can therefore apply \cref{lem:app_particle_saturated_recursion}
as in the population proof, now with
$$
    s=\eta H,\qquad b=C_1\varepsilon_n(\delta),\qquad
    \Gamma=-\log\rho,
$$
and decrement constant $8\kappa/27\in(0,1)$.
This gives
$$
    a_t\le C\{(1+\eta Ht)^{-2}+e^{-ct}\}
                  +C_1\varepsilon_n(\delta),\qquad t\ge0.
$$
Finally, \eqref{eq:app_population_kl_l2} and the lower bound in
\eqref{eq:app_empirical_gap_comparison} yield a bound on the
estimation error of the regression function:
\begin{align*}
    2P_X(p_{G_t}-p^*)^2
    &\le\Delta_t
     \le{8\over7}\{a_t+C_{\mathrm{rel}}\varepsilon_n(\delta)\}\\
    &\le{8C\over7}\{(1+\eta Ht)^{-2}+e^{-ct}\}
          +{8(C_1+C_{\mathrm{rel}})\over7}
                      \varepsilon_n(\delta).
\end{align*}
Since $8C/7\ge C$ and
$8(C_1+C_{\mathrm{rel}})/7\ge C_1$, both bounds in the theorem
hold with
$$
    C_{\mathrm{opt}}={8C\over7},\qquad
    C_{\mathrm{stat}}={8(C_1+C_{\mathrm{rel}})\over7}.
$$
The constants $C,c$ in the recursion depend only on
$-\log\rho$ and $\kappa$. Hence
$C_{\mathrm{opt}},C_{\mathrm{stat}},c$ depend only on the model
parameters, whereas $n_0,M_0$ may also depend on $\delta$.
None of these constants or thresholds depends on $n,M,H,\eta$, or
the realized sample.
The event in \cref{lem:app_particle_empirical_comparison} gives the
comparison bounds for every candidate measure, hence for all iterates
computed from the same data. Its intersection with the initialization
event has probability at least $1-\delta$. Therefore both bounds
hold simultaneously for all $t\ge0$, without an additional union
bound over blocks or location updates. This proves the theorem.
\end{proof}

\section{Proof of the lower bound of algorithmic convergence rate}
\label{app:particle_lower_bound_construction}

We construct a fixed one-dimensional model with one true atom.
\Cref{lem:app_sharpness_bracket} stated below identifies the two particles that
receive positive optimal weights and bounds their location gradients.
With probability bounded below independently of $M$, these particles
start at distance of order $M^{-1}$ from the true atom and remain at
that distance for order $M^2$ blocks. The resulting gap is of order
$M^{-4}$, which gives the inverse-square lower bound.

\paragraph{A fixed model and exact reweighting}

Take $d=1$, $G^*=\delta_1$, and
$X\sim\operatorname{Unif}[-1,1]$. The radius $r>1$ will be fixed
below, independently of $M$. The sup-norm of functions of $x$
below are taken over $[-1,1]$. Write
\begin{gather*}
    p^*(x)=\sigma(x),\qquad
    a_1(x)=x\sigma'(x),\qquad a_2(x)=x^2\sigma''(x).
\end{gather*}
Thus $a_1$ and $a_2$ are the first two location derivatives of
$p_{\delta_\beta}$, evaluated at the true location $\beta=1$.
The quadratic likelihood approximation uses the
following weighted inner product:
\begin{gather*}
    \langle f,g\rangle_*=\EE\{f(X)g(X)/\sigma'(X)\},\qquad
    I=\langle a_1,a_1\rangle_*,\\
    c_1={\langle a_1,a_2\rangle_*\over I},\qquad
    \psi=a_2-c_1a_1,\qquad
    D=\langle\psi,\psi\rangle_*.
\end{gather*}
The function $\psi$ is the part of the second-order prediction change
orthogonal to the first-order location change. Since
$\sigma'(1)\le \sigma'(x)\le1/4$ for $|x|\le1$, we have
$$
    {\sigma'(1)\over3}
    \le I=\EE\{X^2\sigma'(X)\}\le {1\over12}.
$$
Moreover, $a_2(x)/a_1(x)=-x\tanh(x/2)$ for $x\ne0$.
This ratio is not constant on $(0,1)$. Thus $a_2$ is not a
constant multiple of $a_1$ in the weighted inner product, and
$$
    D=\langle a_2,a_2\rangle_*
       -{\langle a_1,a_2\rangle_*^2\over I}>0.
$$
Consider the function
\begin{equation}
    K(\beta)=\langle\psi,p_{\delta_\beta}-p^*\rangle_*.
    \label{eq:app_sharpness_certificate}
\end{equation}
Orthogonality gives $K(1)=K'(1)=0$ and $K''(1)=D$.
Primes on $K$ and on a potential $g_G$ denote differentiation in
the scalar location $\beta$, with $G$ held fixed in the latter case.
The sign of $K$ will exclude all particles outside the active pair
from receiving positive optimal weight.

We first verify the sign of $K$ on a fixed interval. Let $\nu$ be
the probability measure on $[0,1]$ defined by
$$
    \int f(x)\,\rd \nu(x)
    ={1\over I}\EE\left\{{a_1(X)^2\over \sigma'(X)}f(|X|)\right\}.
$$
For $f(x)=-x\tanh(x/2)$ and
$h_\beta(x)=\sigma'(\beta x)/\sigma'(x)$, we have
$\psi(x)=a_1(x)\{f(|x|)-c_1\}$. Consequently,
\begin{gather*}
    c_1=\int f\,\rd \nu,\qquad
    D=I\operatorname{Var}_\nu(f),\\
    K'(\beta)
    =\EE\left\{{a_1(X)^2\over \sigma'(X)}
            \{f(|X|)-c_1\}h_\beta(|X|)\right\}
    =I\operatorname{Cov}_\nu(f,h_\beta),\\
    {d\over dx}\log h_\beta(x)
    =\tanh(x/2)-|\beta|\tanh(|\beta|x/2).
\end{gather*}
For $x>0$,
$$
    f'(x)=-\tanh(x/2)-{x\over2\cosh^2(x/2)}<0.
$$
Also, $a\mapsto a\tanh(ax/2)$ is strictly increasing on
$[0,\infty)$, because its derivative is
$$
\tanh(ax/2)+ax/\{2\cosh^2(ax/2)\}>0
$$
for $a>0$.
It follows that $h_\beta$ is strictly increasing when
$|\beta|<1$ and strictly decreasing when $|\beta|>1$.
The measure $\nu$ has density $x^2\sigma'(x)/I>0$ on $(0,1)$.
The identity
$$
    2\operatorname{Cov}_\nu(f,h)
    =\iint\{f(x)-f(y)\}\{h(x)-h(y)\}\,\rd \nu(x)\,\rd \nu(y)
$$
therefore shows that $K'(\beta)<0$ for $|\beta|<1$ and
$K'(\beta)>0$ for $|\beta|>1$. In particular,
$$
    K(-1)=K(1)-\int_{-1}^1K'(b)\,\rd b
          =-\int_{-1}^1K'(b)\,\rd b>0.
$$
By continuity, choose a fixed $\tau\in(0,1/4]$ such that $K>0$
on $[-1-\tau,-1]$, and set
\begin{equation}
    r=1+\tau.
    \label{eq:app_sharpness_radius}
\end{equation}
On the rest of $[-r,r]\setminus\{1\}$, positivity follows from the
derivative signs above. Thus
\begin{equation}
    K(\beta)>0\qquad
    \text{for all }\beta\in[-r,r]\setminus\{1\}.
    \label{eq:app_sharpness_certificate_positive}
\end{equation}
This fixes the model. It satisfies \cref{assum:particle_population} with
$K_0=w_0=s_-=s_+=R_X=r_X=1$ and $m_X=1/2$; the ordinary separation
condition is vacuous for a single atom. The radius in
\eqref{eq:app_sharpness_radius} satisfies \eqref{eq:particle_radius_choice}
with $\tau_0=\tau$.

\begin{lemma}
\label{lem:app_sharpness_bracket}
For the model above, there is $\epsilon_0\in(0,\tau)$ with the following
property. Suppose $u,v>0$, $u+v\le\epsilon_0$, and
$1/6\le u/v\le6$. Let $w=w(u,v)$ maximize the population
likelihood over mixtures supported at $1-u$ and $1+v$, with
$w$ denoting the weight at $1-u$. This maximizer is unique and
belongs to $(0,1)$. Write
$$
    G_{u,v}=w\delta_{1-u}+(1-w)\delta_{1+v},
$$
where $g_G$ is the population potential in
\eqref{eq:app_population_residual_potential}. Uniformly over these
$u,v$, as $u+v\to0$,
\begin{align}
    w(u,v)
    &={v\over u+v}+{c_1uv\over2(u+v)}+O\{(u+v)^2\},
    \notag\\
    p_{G_{u,v}}-p^*
    &={uv\over2}\psi+O\{uv(u+v)\},
    \label{eq:app_sharpness_pair_expansion}\\
    \cL(G^*)-\cL(G_{u,v})
    &={D\over8}u^2v^2\{1+\cO(u+v)\},\notag\\
    g'_{G_{u,v}}(1-u)
    &={D\over4}uv(u+v)\{1+\cO(u+v)\},\notag\\
    g'_{G_{u,v}}(1+v)
    &=-{D\over4}uv(u+v)\{1+\cO(u+v)\}.
    \notag
\end{align}
The function-valued remainder is uniform on $[-1,1]$.
Moreover, if a finite particle configuration in $[-r,r]$ contains
exactly one particle at each endpoint and all its other particles lie outside
$[1-u,1+v]$, its unique population-optimal weight vector assigns
weights $w(u,v)$ and $1-w(u,v)$ to these endpoints and zero to
every other particle. After decreasing $\epsilon_0$, we have
\begin{equation}
\begin{gathered}
    \cL(G^*)-\cL(G_{u,v})\ge {D\over16}u^2v^2,\\
    \max\{|g'_{G_{u,v}}(1-u)|,|g'_{G_{u,v}}(1+v)|\}
    \le {D\over2}uv(u+v).
\end{gathered}
    \label{eq:app_sharpness_pair_bounds}
\end{equation}
\end{lemma}

\begin{proof}
Optimal weighting of two particles on opposite sides of the true atom
cancels the first-order term in the regression function error. The remaining error is of
second order in their distances, so its likelihood contribution is of
fourth order. The weight optimality equation also cancels the leading
term in the location gradients, which are of third order. We derive
these expansions and then verify that every other particle receives
zero optimal weight.

Set $s=u+v$ and $\theta=v/s$, so that
$\theta\in[1/7,6/7]$. Parameterize the left weight by
$w=\theta+sz$, and let $p$ denote the corresponding prediction
function. To extend these functions to a neighborhood of $s=0$,
write $u=(1-\theta)s$ and $v=\theta s$ also for negative $s$.
In this proof, write $\rho=\sigma(-r)$ and
$Z=1+|c_1|/8$, and restrict $|z|\le Z$ and
$|s|\le\min\{1,\tau/2,(14Z)^{-1}\}$. These restrictions give
$w\in[1/14,13/14]$ and $p(x)\in[\rho,1-\rho]$.
For $|t|\le\tau/2$, Taylor's formula gives
$$
    p_{\delta_{1+t}}=p^*+ta_1+{t^2\over2}a_2+R_t,\qquad
    R_t={t^3\over2}\int_0^1(1-a)^2
                       \partial_\beta^3p_{\delta_{1+at}}\,\rd a,
    \qquad \|R_t\|_\infty\le {|t|^3\over48}.
$$
Here we used
$\partial_\beta^3p_{\delta_\beta}(x)=x^3\sigma'''(\beta x)$ and
$\sup|\sigma'''|\le1/8$. Since
$$
    -wu+(1-w)v=-s^2z,\qquad
    wu^2+(1-w)v^2=s^2\theta(1-\theta)+s^3z(1-2\theta),
$$
substitution yields
\begin{align*}
    p-p^*
      &=s^2\left\{-za_1+{\theta(1-\theta)\over2}a_2\right\}
           +R_p,\\
    R_p&={s^3z(1-2\theta)\over2}a_2+wR_{-u}+(1-w)R_v,\\
    \|R_p\|_\infty
      &\le A_Z|s|^3,\qquad
        A_Z={Z\over12\sqrt3}+{1\over48}.
\end{align*}
For the secant derivative, the same expansion gives
\begin{gather*}
    a_{\mathrm{sec}}:={p_{\delta_{1+v}}-p_{\delta_{1-u}}\over s}
       =a_1+{v-u\over2}a_2+R_{\mathrm{sec}},\qquad
    \|R_{\mathrm{sec}}\|_\infty
       \le {|u|^3+|v|^3\over48|s|}
       \le {s^2\over48}.
\end{gather*}
Alternatively,
$a_{\mathrm{sec}}=\int_0^1D_\beta p_{\delta_{1-u+as}}\,\rd a$. In particular,
$$
    \|a_{\mathrm{sec}}\|_\infty\le {1\over4},
    \qquad
    \|a_{\mathrm{sec}}-a_1\|_\infty
       \le {|s|\over6\sqrt3}.
$$
An interior optimum satisfies
\begin{equation}
    {1\over s^2}\EE\left\{
      {p(X)-p^*(X)\over p(X)\{1-p(X)\}}
      a_{\mathrm{sec}}(X)\right\}=0.
    \label{eq:app_sharpness_weight_equation}
\end{equation}
Let $\Phi(s,\theta,z)$ denote the left side of
\eqref{eq:app_sharpness_weight_equation}. The integral remainders
above show that $(p-p^*)/s^2$ and $a_{\mathrm{sec}}$ extend
smoothly to $s=0$. The denominator is bounded below by
$\rho^2$, so differentiation under the expectation is valid
on this parameter set. At $s=0$, we obtain
$$
    \Phi(0,\theta,z)
       =-Iz+{c_1I\over2}\theta(1-\theta)
       =-I\{z-c_1\theta(1-\theta)/2\}.
$$
To bound the remainder uniformly, put
$$
    B_Z={Z\over4}+{1\over48\sqrt3},\qquad E_Z=B_Z+A_Z.
$$
The leading coefficient of $s^2$ in $p-p^*$ has supremum
norm at most $B_Z$, and $\|p-p^*\|_\infty\le E_Zs^2$.
For $\rho\le t\le1-\rho$,
$$
    {1\over t(1-t)}\le\rho^{-2},\qquad
    \left|{d\over dt}{1\over t(1-t)}\right|
       ={|2t-1|\over t^2(1-t)^2}\le\rho^{-4}.
$$
By adding and subtracting
the terms with $p=p^*$ and $a_{\mathrm{sec}}=a_1$, we find
\begin{align*}
    |\Phi(s,\theta,z)-\Phi(0,\theta,z)|
    &\le {A_Z\over4\rho^2}|s|
          +{B_ZE_Z\over4\rho^4}s^2
          +{B_Z\over6\sqrt3\,\rho^2}|s| \le C_\Phi|s|
\end{align*}
where
$$
    C_\Phi :={A_Z\over4\rho^2}
          +{B_ZE_Z\over4\rho^4}
          +{B_Z\over6\sqrt3\,\rho^2}.
$$
Further, $\partial_zp=-s^2a_{\mathrm{sec}}$, and hence
$$
    \partial_z\Phi(s,\theta,z)
       =-\EE\left[
          \left\{{p^*(X)\over p(X)^2}
                 +{1-p^*(X)\over\{1-p(X)\}^2}\right\}
                    a_{\mathrm{sec}}(X)^2\right]<0.
$$
At $s=0$ this derivative equals $-I$. Since
$|c_1\theta(1-\theta)/2|\le |c_1|/8$, both $c_1\theta(1-\theta)/2-1/2$
and $c_1\theta(1-\theta)/2+1/2$ belong to $[-Z,Z]$. If
$C_\Phi|s|\le I/4$, the values of $\Phi$ at these two
points have opposite signs. Thus there is a unique root between
them. The implicit function theorem gives its smooth dependence
on $s,\theta$, and the preceding bound gives
$$
    |z-c_1\theta(1-\theta)/2|\le {C_\Phi\over I}|s|.
$$
For $s>0$, its weight belongs to $(0,1)$. Moreover,
$$
    \begin{aligned}
    &{d^2\over dw^2}\cL
        \{w\delta_{1-u}+(1-w)\delta_{1+v}\}
    \\&\quad=-\EE\left[
       \left\{{p^*(X)\over p(X)^2}
              +{1-p^*(X)\over\{1-p(X)\}^2}\right\}
           \{p_{\delta_{1-u}}(X)-p_{\delta_{1+v}}(X)\}^2\right]<0.
    \end{aligned}
$$
The root is therefore the unique maximizer.
Substituting it into the prediction expansion gives
$$
    \begin{gathered}
    p_{G_{u,v}}-p^*={uv\over2}\psi+R_e,\qquad
    \|R_e\|_\infty
       \le\left(A_Z+{C_\Phi\over4I}\right)s^3
       \le C_euvs,
    \\
    C_e={49\over6}\left(A_Z+{C_\Phi\over4I}\right),
    \end{gathered}
$$
where $6s^2/49\le uv=s^2\theta(1-\theta)\le s^2/4$.
This proves the first two expansions in
\eqref{eq:app_sharpness_pair_expansion}, with constants independent
of $\theta$.

We use the scalar loss in \eqref{eq:app_population_scalar_loss_bounds}
with $a=p^*(x)$. Holding $x$ fixed, we have
\begin{gather*}
    f_{p^*}'(p^*)=0,\qquad
    f_{p^*}''(p^*)={1\over p^*(1-p^*)}={1\over\sigma'(x)},\\
    f_{p^*}'''(t)=-{2p^*\over t^3}+{2(1-p^*)\over(1-t)^3},
    \qquad
    |f_{p^*}'''(t)|\le {2\over\rho^3}
          \quad (\rho\le t\le1-\rho).
\end{gather*}
Taylor's formula and expectation therefore give
$$
    \left|\{\cL(G^*)-\cL(G_{u,v})\}
             -{1\over2}\langle p_{G_{u,v}}-p^*,p_{G_{u,v}}-p^*\rangle_*\right|
       \le {\|p_{G_{u,v}}-p^*\|_\infty^3\over3\rho^3}.
$$
Put $B_e=\|\psi\|_\infty/2+C_e$, so that
$\|p_{G_{u,v}}-p^*\|_\infty\le B_euv$ for $s\le1$. By expanding the
quadratic term and using $\|R_e\|_\infty\le C_euvs$, we obtain
\begin{align*}
    \left|\{\cL(G^*)-\cL(G_{u,v})\}-{D\over8}u^2v^2\right|
    &\le {uv\over2}|\langle\psi,R_e\rangle_*|
          +{1\over2}\langle R_e,R_e\rangle_*
          +{B_e^3(uv)^3\over3\rho^3} \le C_\Delta(uv)^2s,
\end{align*}
where
$$
  C_\Delta ={\|\psi\|_\infty C_e\over2\rho^2} +{C_e^2\over2\rho^2}+{B_e^3\over12\rho^3}.
$$
In the last step, we used $uv\le s^2/4$ and $s\le1$.
This proves the likelihood expansion.

For the location gradients, the weight equation gives
$\EE\{r_{G_{u,v}}(X)a_{\mathrm{sec}}(X)\}=0$. Consequently, at either
endpoint,
$$
    g'_{G_{u,v}}(\beta)
    =\EE\left[r_{G_{u,v}}(X)
        \{\partial_\beta p_{\delta_\beta}(X)-a_{\mathrm{sec}}(X)\}\right].
$$
To bound each factor, Taylor's formula for the location derivative gives
$$
    D_\beta p_{\delta_{1+t}}=a_1+ta_2+\wt R_t,\qquad
    \|\wt R_t\|_\infty\le {t^2\over16}.
$$
Combining this with the expansion of $a_{\mathrm{sec}}$ yields
\begin{gather*}
    \left.\partial_\beta p_{\delta_\beta}\right|_{\beta=1-u}-a_{\mathrm{sec}}=-{s\over2}a_2+R_-,\qquad
    \left.\partial_\beta p_{\delta_\beta}\right|_{\beta=1+v}-a_{\mathrm{sec}}={s\over2}a_2+R_+,\\
    \max\{\|R_-\|_\infty,\|R_+\|_\infty\}
       \le {s^2\over16}+{s^2\over48}={s^2\over12}.
\end{gather*}
Since $r_{G_{u,v}}=-(p_{G_{u,v}}-p^*)/\{p_{G_{u,v}}(1-p_{G_{u,v}})\}$ and
$\sigma'=p^*(1-p^*)$,
\begin{align*}
    r_{G_{u,v}}&=-{uv\over2}{\psi\over\sigma'}+R_r,\\
    \|R_r\|_\infty
    &\le {\|R_e\|_\infty\over\rho^2}
          +{uv\|\psi\|_\infty\over2\rho^4}\|p_{G_{u,v}}-p^*\|_\infty
      \le C_ruvs &&\text{by }uv\le s^2/4\le s/4,
\end{align*}
with
$$
    C_r={C_e\over\rho^2}
              +{\|\psi\|_\infty B_e\over8\rho^4}.
$$
As $\langle\psi,a_2\rangle_*
=\langle\psi,\psi+c_1a_1\rangle_*=D$, multiplication gives
\begin{gather*}
    \left|g'_{G_{u,v}}(1-u)-{D\over4}uvs\right|
       \le C_guvs^2,\qquad
    \left|g'_{G_{u,v}}(1+v)+{D\over4}uvs\right|
       \le C_guvs^2,
\end{gather*}
with
$$
    C_g ~ =~ {\|\psi\|_\infty\over24\rho^2} +C_r\left({1\over12\sqrt3}+{1\over12}\right).
$$
Indeed, the product of $-uv\psi/(2\sigma')$ with $R_\pm$ is
bounded by $uvs^2\|\psi\|_\infty/(24\rho^2)$, while the
product involving $R_r$ is bounded by
$C_ruvs\{s/(12\sqrt3)+s^2/12\}$.
This proves the last two expansions in
\eqref{eq:app_sharpness_pair_expansion}.

The residual expansion implies
\begin{equation}
    {g_{G_{u,v}}(\beta)-g_{G_{u,v}}(1)\over uv}
    \to -{1\over2}K(\beta)
    \label{eq:app_sharpness_potential_limit}
\end{equation}
uniformly on $[-r,r]$, together with its first two derivatives in
$\beta$. More explicitly, the difference between the two sides is
$\EE\{R_r(X)(p_{\delta_\beta}(X)-p^*(X))\}/(uv)$. Since
$$
    \sup_\beta\|p_{\delta_\beta}-p^*\|_\infty\le1,\qquad
    \sup_\beta\|\partial_\beta p_{\delta_\beta}\|_\infty\le{1\over4},\qquad
    \sup_\beta\|\partial_\beta^2p_{\delta_\beta}\|_\infty
          \le{1\over6\sqrt3},
$$
the supremum norm of this difference, and of each of its first two
location derivatives, is at most $C_rs$.
Differentiation gives
$|K'''(\beta)|\le\|\psi\|_\infty/(8\rho^2)$, using
$\|\sigma'''\|_\infty\le1/8$. Set
$$
    h=\min\left\{{\tau\over2},
                    {4D\rho^2\over\|\psi\|_\infty}\right\},\qquad
    U=(1-h,1+h).
$$
Since $K''(1)=D$, the mean value theorem gives $K''\ge D/2$
on $U$, whose closure lies in $(-r,r)$.
For $s\le h/2$, both endpoints belong to $U$.
If also $C_rs\le D/8$, then
$$
    {g''_{G_{u,v}}(\beta)\over uv}
       \le-{1\over2}K''(\beta)+C_rs
       \le-{D\over8},\qquad \beta\in U.
$$
Thus $g_{G_{u,v}}$ is strictly concave on $U$. Weight
optimality gives
$g_{G_{u,v}}(1-u)=g_{G_{u,v}}(1+v)$. Strict concavity therefore gives
$$
    g_{G_{u,v}}(\beta)<g_{G_{u,v}}(1-u)
    \quad\text{for }\beta\in U\setminus[1-u,1+v].
$$
On the compact complement $[-r,r]\setminus U$, put
$$
    k_U=\min_{\beta\in[-r,r]\setminus U}K(\beta)>0,\qquad
    B_K={\|\psi\|_\infty\over6\sqrt3\,\rho^2}.
$$
The positivity follows from
\eqref{eq:app_sharpness_certificate_positive}, and
$\sup_\beta|K''(\beta)|\le B_K$ follows by differentiation
under the expectation. Since $K(1)=K'(1)=0$,
$$
    \left|{g_{G_{u,v}}(1-u)-g_{G_{u,v}}(1)\over uv}\right|
       \le {1\over2}|K(1-u)|+C_rs
       \le {B_Ks^2\over4}+C_rs.
$$
Taking $s\le\sqrt{k_U/(4B_K)}$ and $s\le k_U/(16C_r)$
makes the last expression at most $k_U/8$.
For $\beta\in[-r,r]\setminus U$, we then have
$$
    {g_{G_{u,v}}(\beta)-g_{G_{u,v}}(1)\over uv}
       \le-{k_U\over2}+C_rs
       \le-{3k_U\over8}
       <-{k_U\over8}
       \le {g_{G_{u,v}}(1-u)-g_{G_{u,v}}(1)\over uv}.
$$
This proves the strict inequality on the complement as well.
All these choices are independent of $\theta$.

The two-endpoint mixture consequently satisfies the full weight KKT
conditions, with strict inequality at every other particle. By concavity,
any competing weight vector giving positive weight to another particle
has strictly smaller likelihood. To see this, let $\wt w$
be such a vector and $c=g_{G_{u,v}}(1-u)=g_{G_{u,v}}(1+v)$.
The directional derivative from the two-endpoint optimum is
$$
    \sum_j\wt w_j\{g_{G_{u,v}}(\beta_j)-c\}<0,
$$
and concavity bounds the difference of the two likelihoods from
above by this derivative. A competing vector supported
on the endpoints alone has smaller likelihood unless it equals the
unique two-point optimum. This proves uniqueness of the full weight
vector, even if some inactive particles coincide with one another.
Finally, take $\epsilon_0$ small enough to satisfy all the
restrictions above, as well as
$$
    C_\Delta\epsilon_0\le {D\over16},
    \qquad C_g\epsilon_0\le {D\over4}.
$$
The likelihood and gradient error bounds then give
\eqref{eq:app_sharpness_pair_bounds}.
\end{proof}

\subsection{Proof of \cref{prop:particle_inverse_square_sharp}: lower bound of the convergence rate of RPG}
\label{app:particle_inverse_square_sharp_proof}

\begin{proof}
We choose an initialization event on which the two particles adjacent
to the true atom are both at distance of order $M^{-1}$.
Their gradients are then of order $M^{-3}$, so they remain at this
distance for order $M^2$ blocks. Order preservation prevents other
particles from replacing this pair. The resulting likelihood difference
of order $M^{-4}$ gives the claimed lower bound.

Fix an integer $H\ge1$, and set
\begin{gather*}
    \rho=\sigma(-r),\qquad
    L={1\over6\sqrt3\,\rho}+{1\over16\rho^2},\\
    0<\eta\le\min\{L^{-1},27\rho/4\},\qquad
    B_H=\sum_{h=0}^{H-1}(1+\eta L)^h.
\end{gather*}
Every mixture on $[-r,r]$
has prediction values in $[\rho,1-\rho]$. Its population residual
satisfies
$$
    r_G(x)={p^*(x)\over p_G(x)}
                 -{1-p^*(x)\over1-p_G(x)},\qquad
    -{1\over1-p_G(x)}\le r_G(x)\le {1\over p_G(x)}.
$$
Thus $\|r_G\|_\infty\le\rho^{-1}$, and differentiation gives
$$
    |g_G''(\beta)|
       =|\EE\{r_G(X)X^2\sigma''(\beta X)\}|
       \le {1\over6\sqrt3\,\rho}<L.
$$
Consequently,
$$
    {d\over d\beta}\{\beta+\eta g_G'(\beta)\}
       \ge1-{\eta\over6\sqrt3\,\rho}
       \ge1-{1\over6\sqrt3\,\rho L}>0.
$$
The map is strictly increasing, and its projection onto
$[-r,r]$ is nondecreasing. If the
updated location of a given particle is interior, every particle strictly
to its left or right remains strictly on that side after the update.
Indeed, strict order holds before projection, and clipping an outside
location to a boundary cannot make it equal to an interior location.
This applies to every particle, including those with zero weight.

We next bound displacement within a block without assuming that the
weights remain optimal between reweightings. Fix the two weights
$w,1-w$, and let $F(b_-,b_+)$ be the vector of their two
unweighted population location gradients. For
$p=wp_{\delta_{b_-}}+(1-w)p_{\delta_{b_+}}$, the derivative of
$p^*(x)/p-(1-p^*(x))/(1-p)$ with respect to $p$ is
$$
    -A_p(x),\qquad
    A_p(x)={p^*(x)\over p(x)^2}
                +{1-p^*(x)\over\{1-p(x)\}^2}
           \le \rho^{-2}.
$$
Write $w_-=w$ and $w_+=1-w$.
Since $F_i=\EE\{r_G(X)D_\beta p_{\delta_{b_i}}(X)\}$, we have
$$
    \partial_{b_j}F_i
       =-w_j\EE\{A_p(X)D_\beta p_{\delta_{b_j}}(X)D_\beta p_{\delta_{b_i}}(X)\}
          +\mathbf1_{\{i=j\}}\EE\{r_G(X)D_\beta^2 p_{\delta_{b_i}}(X)\}.
$$
Using $\|D_\beta p_{\delta_{b}}\|_\infty\le1/4$,
$\|D_\beta^2 p_{\delta_{b}}\|_\infty\le1/(6\sqrt3)$, and
$\|r_G\|_\infty\le\rho^{-1}$, we obtain
$$
    \sum_{j\in\{-,+\}}|\partial_{b_j}F_i(b_-,b_+)|
    \le {1\over6\sqrt3\,\rho}
       +{w+(1-w)\over16\rho^2}=L,
    \qquad i\in\{-,+\}.
$$
For $b,c\in[-r,r]^2$, integration along their line segment now gives
$$
    |F_i(b)-F_i(c)|
       \le\int_0^1\sum_j
          |\partial_{b_j}F_i\{c+t(b-c)\}|\,|b_j-c_j|\,\rd t
       \le L\|b-c\|_\infty.
$$
Thus $F$ is $L$-Lipschitz in the maximum norm. If
$b^h=(b_-^h,b_+^h)$ are the projected fixed-weight iterates in a
block, nonexpansiveness of projection gives
$$
    \|b^{h+1}-b^h\|_\infty\le\eta\|F(b^h)\|_\infty,
    \qquad
    \|F(b^{h+1})\|_\infty
    \le(1+\eta L)\|F(b^h)\|_\infty.
$$
By induction,
$\|F(b^h)\|_\infty\le(1+\eta L)^h\|F(b^0)\|_\infty$.
Summing the displacement bound over $h$ therefore gives
\begin{equation}
    \sum_{h=0}^{H-1}\|b^{h+1}-b^h\|_\infty
    \le\eta B_H\|F(b^0)\|_\infty.
    \label{eq:app_sharpness_block_displacement}
\end{equation}
For sufficiently large $M$, let $\mathcal E_M$ be the event that
exactly one initial particle belongs to each of
$$
    [1-2r/M,1-r/M],\qquad [1+r/M,1+2r/M],
$$
and all other particles lie outside $[1-2r/M,1+2r/M]$. These
intervals lie in $[-r,r]$ for all sufficiently large $M$.
Each selected interval has length $r/M$, and hence probability
$1/(2M)$ under the uniform initialization on $[-r,r]$.
The interval $[1-2r/M,1+2r/M]$ has length $4r/M$; its complement
has probability $1-2/M$. There are $M(M-1)$ choices of the
ordered pair of particles in the two selected intervals.
Independence therefore gives
\begin{equation}
    \begin{aligned}
    \PP(\mathcal E_M)
    &=M(M-1)\left({1\over2M}\right)^2
                  \left(1-{2\over M}\right)^{M-2}
       \\
    &={M-1\over4M}\left(1-{2\over M}\right)^{M-2}
    \to {e^{-2}\over4}.
    \end{aligned}
    \label{eq:app_sharpness_initial_event}
\end{equation}
The limit follows from
$(M-2)\log(1-2/M)=-2+\cO(M^{-1})$.
On this event, denote the two selected locations at block starts by
$1-u_t$ and $1+v_t$. Initially
$r/M\le u_0,v_0\le2r/M$, and these particles are adjacent in
the ordered configuration. Define the deterministic horizon
\begin{equation}
    a={1\over128\eta B_HD r^2},\qquad T_M=\lfloor aM^2\rfloor.
    \label{eq:app_sharpness_horizon}
\end{equation}

We prove by induction over blocks that these same particles are the only
active ones and that their distances from $1$ remain in
$[r/(2M),3r/M]$, including all intermediate location updates up to
block $T_M$. Take $M$ large enough that
$6r/M\le\epsilon_0$ and
$[1-3r/M,1+3r/M]\subset(-r,r)$. Whenever a block starts with
the claimed configuration, \cref{lem:app_sharpness_bracket} applies:
the distance ratio belongs to $[1/6,6]$, and exact reweighting
uses only these two particles. Its initial maximum gradient is at most
$$
    {D\over2}u_tv_t(u_t+v_t)\le {27Dr^3\over M^3}.
$$
By \eqref{eq:app_sharpness_block_displacement}, the total displacement
of either selected particle over the first $t+1\le T_M$ blocks,
assuming the induction assertions at their starts, is at most
\begin{equation}
    (t+1)\eta B_H{27Dr^3\over M^3}
    \le {27r\over128M}<{r\over4M}.
    \label{eq:app_sharpness_total_displacement}
\end{equation}
The same bound holds for every partial sum of the inner-iteration
displacements. In particular, the distance of each selected particle
from its initial location is less than $r/(4M)$ throughout these
blocks. Since its initial distance from $1$ is between $r/M$
and $2r/M$, its distance from $1$ remains between
$$
    {r\over M}-{r\over4M}={3r\over4M}
    \quad\text{and}\quad
    {2r\over M}+{r\over4M}={9r\over4M}.
$$
These bounds imply the claimed interval $[r/(2M),3r/M]$;
they also show that neither selected particle crosses $1$.
Both selected particles remain in the interior of $[-r,r]$.
The order argument prevents any other
particle from entering their interval or coinciding with an endpoint.
By \cref{lem:app_sharpness_bracket}, reweighting at the end of the block
again gives the same active pair. This proves the induction, starting with the
configuration on $\mathcal E_M$; it also includes the final
reweighting defining $G_{T_M}$.

Applying \eqref{eq:app_sharpness_pair_bounds} at the final reweighting gives
$$
    \cL(G^*)-\cL(G_{T_M})\ge {D\over16}u_{T_M}^2v_{T_M}^2
        \ge {D\over16}\left({r\over2M}\right)^4
        ={Dr^4\over256M^4}.
$$
The estimation error bound follows from
$P_X(p_G-p^*)^2\ge\rho^2\{\cL(G^*)-\cL(G)\}$, by
\eqref{eq:app_population_kl_l2}. Thus, on $\mathcal E_M$,
\begin{equation}
\begin{gathered}
    \cL(G^*)-\cL(G_{T_M})\ge {b_\Delta\over M^4},\qquad
    P_X(p_{G_{T_M}}-p^*)^2\ge {b_R\over M^4},\\
    b_\Delta={Dr^4\over256},\quad b_R=\rho^2b_\Delta.
\end{gathered}
    \label{eq:app_sharpness_population_lower}
\end{equation}
These constants are positive and independent of $M$. Choose
\begin{equation}
    \delta_0={e^{-2}\over16},\qquad
    c={a^2\over2}\min\{b_\Delta,b_R\}.
    \label{eq:app_sharpness_common_constants}
\end{equation}
For all sufficiently large $M$,
$\PP(\mathcal E_M)\ge2\delta_0$, and
$1+T_M\ge aM^2$. Thus
\eqref{eq:app_sharpness_population_lower} implies the claimed
population lower bound, as well as the same lower bound for squared
estimation error. Also $T_M\to\infty$.
\end{proof}

\end{document}